%% file: main.tex
\documentclass[preprint]{imsart}

\RequirePackage[OT1]{fontenc}
\RequirePackage{amsthm,amsmath,amssymb}
\RequirePackage{natbib}
\RequirePackage[colorlinks,citecolor=blue,urlcolor=blue]{hyperref}
\usepackage{graphicx}
\usepackage{amsthm}
\usepackage{pgf,tikz}
\usepackage{url}
\usepackage{thmtools}
\usetikzlibrary{arrows}
\usepackage{pdfcomment}
\usepackage[final]{pdfpages} 

\usetikzlibrary{arrows}
\usetikzlibrary{matrix}
\usepackage{diagbox}
\usepackage{multirow}
\usepackage{booktabs}
\usepackage{adjustbox}
\usepackage{hhline}
\usepackage{mathtools}

\numberwithin{equation}{section}
\theoremstyle{plain}
\declaretheorem[name=Theorem,numberwithin=section]{Th}
\declaretheorem[name=Lemma,numberwithin=section]{Lem}

\declaretheorem[name=Corollary,numberwithin=section]{Co}

\theoremstyle{definition}
\newtheorem{Def}{Definition}[section]
\newtheorem{Not}{Notation}[section]

\newtheorem{Ass}{Assumption}[section]

\theoremstyle{remark}

 \newtheorem{Remth}{Remark}[Th]

 \newtheorem{Rem}{Remark}[section]

\input{macro}

\begin{document}

\begin{frontmatter}
\title{Statistical applications of Random matrix theory:\\ comparison of two populations II\thanksref{T1}}
\runtitle{Comparison of two populations}
\thankstext{T1}{This paper is constructed from the Thesis of R\'emy Mari\'etan that will be divided in three parts.}

\begin{aug}
\author{\fnms{R\'emy} \snm{Mari\'etan}\thanksref{t1}\ead[label=e1]{remy.marietan@alumni.epfl.ch}}
\and
\author{\fnms{Stephan} \snm{Morgenthaler}\thanksref{t2}\ead[label=e2]{stephan.morgenthaler@epfl.ch}}

\thankstext{t1}{PhD Student at EPFL in mathematics department}
\thankstext{t2}{Professor at EPFL in mathematics department}
\runauthor{R\'emy Mari\'etan and Stephan Morgenthaler}

\affiliation{\'Ecole Polytechnique F\'ed\'eral de Lausanne, EPFL}

\address{Department of Mathematics\\
\'Ecole Polytechnique F\'ed\'eral de Lausanne\\1015 Lausanne\\
\printead{e1}, 
\printead*{e2}}
\end{aug}

\begin{abstract}
This paper investigates a statistical procedure for testing the equality of two independent estimated covariance matrices when the number of potentially dependent data vectors is large and proportional to the size of the vectors, that is, the number of variables.
Inspired by the spike models used in random matrix theory, we concentrate on the largest eigenvalues of the matrices in order to determine significance. To avoid false rejections we must guard against residual spikes and need a sufficiently precise description of the behaviour of the largest eigenvalues under the null hypothesis. 

In this paper we propose an ``invariance" theorems that allows us to extend the test of \cite{mainarticle} for a perturbation of order $1$ to a general tests for order $k$.
The statistics introduced in this paper allow the user to test the equality of two populations based on high-dimensional multivariate data. Furthermore, simulations show that these tests have more power of detection than standard multivariate methods.
\end{abstract}



\end{frontmatter}

\section{Introduction} 

Random matrix theory (RMT) can be used to describe the asymptotic spectral properties of estimators of high-dimensional covariance matrices. The theory has been applied to multi-antenna channels in wireless communication engineering and to financial mathematics models. In other data-rich and high-dimensional areas where statistics is used, such as brain imaging or genetic research, it has not found widespread use. The main barrier to the adoption of RMT may be the lack of concrete statistical results from the probability side. Simply using classical multivariate theory in the high dimension setting can sometimes lead to success, but such  procedures are valid only under strict assumptions about the data such as normality or independence. Even minor differences between the model assumptions and the actual data distribution typically lead to catastrophic results and such procedures do also often have little to no power.

This paper proposes a statistical procedure for testing the equality of two covariance matrices $\Sigma_X$ and $\Sigma_Y$ when the number of potentially dependent data vectors $n$ and the number of variables $m$ are large. RMT tells us what happens to the eigenvalues and eigenvectors of estimators of covariance matrices $\hat{\Sigma}$ when both $n$ and $m$ tend to infinity in such a way that $\lim \frac{m}{n}=c>0$. The classical case, when $m$ is finite and $n$ tends to infinity, is presented in the books of \cite{multi3}, \cite{multi} and \cite{multi2} (or its original version \cite{multi22}).
In the RMT case, the behaviour is more complex, but by now, results of interest are known. \cite{Alice}, \cite{Tao} and more recently \cite{bookrecent} contain comprehensive introductions to RMT and \cite{Appliedbook} covers the case of empirical (estimated) covariance matrices.

Although the existing theory builds a good intuition of the behaviour of these matrices, it does not provide enough of a basis to construct a statistical test with good power.
Inspired by the existing theory, we extend the residual spikes introduced in \cite{mainarticle} and provide a description of the behaviour of diverse types of statistics under a null hypothesis when the perturbation is of order $k$. These results enable the user to test the equality of two populations as well as other null hypotheses such as the independence of two sets of variables.
The remainder of the paper is organized as follows. First, we review the main theorem of \cite{mainarticle} and then indicate how to generalize the test (see Section \ref{sec:test}). We next look at case studies and a compare the new test with alternatives. Finally, in Section \ref{sec:Theorems}, we present the main theorems. The proofs themselves are technical and presented in the supplementary material appendix \ref{appendixproof}.

\section{Statistical test}\label{sec:test}
\subsection{Introduction}
\subsubsection{Hypotheses}
We compare the spectral properties of two covariance estimators $\hat{\Sigma}_X$ and $\hat{\Sigma}_Y$ of dimension $m\times m$ which can be represented as
\begin{Ass}\label{Ass=matrice}
\begin{eqnarray*}
\hat{\Sigma}_X=P_X^{1/2} W_X P_X^{1/2} \text{ and } \hat{\Sigma}_Y=P_Y^{1/2} W_Y P_Y^{1/2}.
\end{eqnarray*}
In this equation, $W_X$ and $W_Y$ are of the form
\begin{eqnarray*}
W_X=O_X \Lambda_X O_X \text{ and } W_Y=O_Y \Lambda_Y O_Y,
\end{eqnarray*}
with $O_X$ and $O_Y$ being independent unit orthonormal random matrices whose distributions are invariant under rotations, while  $\Lambda_X$ and $\Lambda_Y$ are independent positive random diagonal matrices, independent of  $O_X, O_Y$ with trace equal to m and a bound on the diagonal elements. Note that the usual RMT assumption, $\frac{m}{n}=c$ is replaced by this bound! The \textit{(multiplicative) spike model of order} $k$ determines the form of the perturbation $P_X$ (and $P_Y$), which satisfies
$$P_X= \I_m + \sum_{s=1}^k (\theta_{X,s}-1) u_{X,s} u_{X,s}^t\,,$$ 
where $\theta_{X,1}>\theta_{X,2}>...>\theta_{X,k}$ and the scalar product $\left\langle u_{X,s} ,u_{X,r}  \right\rangle=\delta_{s,r}$. $P_Y$ is of the same form. 
\end{Ass}
Some results require large value for $\theta$ and others not. To be precise, we will make use of the following types of hypotheses:
\begin{Ass}\label{Ass=theta}

\begin{itemize} 
\item[(A1)] $\frac{\theta}{\sqrt{m} }\rightarrow \infty.$
\item[(A2)] $\theta \rightarrow \infty.$
\item[(A3)] $\theta_i=p_i \theta$, where $p_i$ is fixed different from $1$.
\item[(A4)] For  $i=1,...,k_\infty, \ \theta_i=p_i \theta$, $\theta \rightarrow \infty$  according to  (A1) or (A2),\\
For $i=k_\infty+1,...,k, \ \theta_i=p_i \theta_0$.\\
For all $i \neq j$, $p_i \not= p_j$.
\end{itemize}
\end{Ass}
The result of this paper will apply to finite eigenvalues $\theta_s$. However, they must be \textbf{detectable}.
\begin{Def}\ \label{Def:detectable} 
\begin{enumerate}
\item We assume that a perturbation $P=\I_m+ (\theta-1) u u^t$ is \textbf{detectable} in $\hat{\Sigma}=P^{1/2} W P^{1/2}$ if the perturbation creates a largest isolated eigenvalue, $\hat{\theta}$.
\item We say that a finite perturbation of order $k$ is \textbf{detectable} if it creates $k$  large eigenvalues separated from the spectrum of $W$.
\end{enumerate}
\end{Def}
Finally, we generalize the filtered estimator of the covariance matrix introduced in \cite{mainarticle}.
\begin{Def}\ \label{Def=unbiased} \\
Suppose $\hat{\Sigma}$ is of the form given in Assumption \ref{Ass=matrice}.\\
The \textbf{unbiased estimator of }$\theta_s$ ($s=1,\ldots,k$) is defined as 
$$ \hat{\hat{\theta}}_s=1+\frac{1}{\frac{1}{m-k} \sum_{i=k+1}^{m} \frac{\hat{\lambda}_{\hat{\Sigma},i}}{\hat{\theta}_s-\hat{\lambda}_{\hat{\Sigma},i}}},$$
where $\hat{\lambda}_{\hat{\Sigma},i}$ is the $i^{\text{th}}$ eigenvalue counting from largest to smallest of $\hat{\Sigma}$.\\
Suppose that $\hat{u}_i$ denotes the eigenvector of $\hat{\Sigma}$ corresponding to the $i^{\text{th}}$ eigenvalue, the \textbf{filtered estimated covariance matrix } is then defined as 
$$\hat{\hat{\Sigma}}= \I_m+\sum_{i=1}^k (\hat{\hat{\theta}}_i-1) \hat{u}_i \hat{u}_i^t.$$
Under Assumption \ref{Ass=matrice}, this estimator is asymptotically equivalent to the theoretical estimator using 
$$ \hat{\hat{\theta}}_s=1+\frac{1}{\frac{1}{m-k} \sum_{i=k+1}^m \frac{\hat{\lambda}_{W,i}}{\hat{\theta}_{s}-\hat{\lambda}_{W,i}}},$$
where $\hat{\lambda}_{W,i}$ is the $i^{\text{th}}$ eigenvalue of $W$.
\end{Def}

Our results will apply to any two centered data matrices $\X \in \mathbb{R}^{m\times n_X}$ and $\Y \in \mathbb{R}^{m\times n_Y}$ which are such that 
\begin{eqnarray*}
\hat{\Sigma}_X= \frac{1}{n_X}\X \X^t \text{ and } \hat{\Sigma}_Y= \frac{1}{n_Y}\Y \Y^t
\end{eqnarray*}
can be decomposed in the manner indicated. This is the basic assumption concerning the covariance matrices. 

We will assume throughout the paper that $n_X\geq n_Y$. 

Because $O_X$ and $O_Y$ are independent and invariant by rotation we can assume without loss of generality that for $s=1,2,...,k$, $u_{X,s}=e_s$ as in \cite{deformedRMT}. Under the null hypothesis, $P_X=P_Y$, we use the simplified notation $P_k$ for both matrices, where for $s=1,2,...,k$, $\theta_{X,s}=\theta_{Y,s}=\theta_s$ and $u_{X,s}=u_{Y,s}(=e_s)$. 

\subsection{The case of $k=1$}
This paper generalises \cite{mainarticle}, in which the following key result was established.

\begin{Th} \label{jointdistribution}\ 
Suppose $W_X$ and $W_Y$ satisfy \ref{Ass=matrice} with $P=P_X=P_Y$, a detectable perturbation of order $k=1$. Moreover, we assume as known the spectra \linebreak $S_{W_X}=\left\lbrace \hat{\lambda}_{W_X,1},\hat{\lambda}_{W_X,2},...,\hat{\lambda}_{W_X,m} \right\rbrace$ and $S_{W_Y}=\left\lbrace \hat{\lambda}_{W_Y,1},\hat{\lambda}_{W_Y,2},...,\hat{\lambda}_{W_Y,m} \right\rbrace$. 
If $\left(\hat{\theta}_X,\hat{\theta}_Y\right)$ converges to $\left(\rho_X,\rho_Y\right)$ in $O_p\left(\theta/\sqrt{m}\right)$ and $\E\left[ \hat{\theta}_X \right]=\rho_X+o\left( \frac{\theta}{\sqrt{m}} \right)$ and $\E\left[ \hat{\theta}_Y \right]=\rho_X+o\left( \frac{\theta}{\sqrt{m}} \right)$, then we have
\begin{equation*}
\resizebox{.95\hsize}{!}{$ \left. \begin{pmatrix}
\hat{\hat{\theta}}_X  \\ 
\hat{\hat{\theta}}_Y  \\ 
\left\langle \hat{u}_X,\hat{u}_Y \right\rangle^2
\end{pmatrix} \right| S_{W_X}, S_{W_Y} \sim \Normal \left(
\begin{pmatrix}
\theta \\
\theta  \\ 
\alpha_{X,Y}^2
\end{pmatrix}
,\frac{1}{m}\begin{pmatrix}
\sigma_{\theta,X}^2 & 0 & \sigma_{\theta,\alpha^2,X} \\ 
0 & \sigma_{\theta,Y}^2  &  \sigma_{\theta,\alpha^2,Y} \\
\sigma_{\theta,\alpha^2,X} & \sigma_{\theta,\alpha^2,Y} & \sigma_{\alpha^2,X,Y}^2
\end{pmatrix}
 \right)+\begin{pmatrix}
o_p\left(\frac{\theta}{\sqrt{m}}\right) \\ 
o_p\left(\frac{\theta}{\sqrt{m}}\right) \\ 
o_p\left( \frac{1}{\theta \sqrt{m}}\right)
\end{pmatrix},$}
\end{equation*}

\noindent where all the parameters in the limit law depend on 
\begin{eqnarray*}
 \resizebox{.95\hsize}{!}{$ M_{s,r,X}(\rho_X)= \frac{1}{m} \sum_{i=1}^m \frac{\hat{\lambda}_{W_X,i}^s}{(\rho_X-\hat{\lambda}_{W_X,i})^r} \text{ and } M_{s,r,Y}(\rho_Y)=\frac{1}{m} \sum_{i=1}^m \frac{\hat{\lambda}_{W_X,i}^s}{(\rho_Y-\hat{\lambda}_{W_X,i})^r}  .$}
\end{eqnarray*}
\end{Th}

\subsection{Generalization}

Suppose $\hat{\Sigma}_X$ and $\hat{\Sigma}_Y$ are two random matrices that verify Assumption \ref{Ass=matrice}. We want to test
\begin{eqnarray*}
{\rm H}_0 : P_X=P_Y,\text{ against }{\rm H}_1 : P_X \neq P_Y.
\end{eqnarray*}
When $P_X=P_Y=P_1$ are perturbation of order $1$, we can 
use Theorem \ref{jointdistribution} to study any test statistic which is a function of the three statistics
$\hat{\hat{\theta}}_X,\hat{\hat{\theta}}_Y,\left\langle \hat{u}_{X},\hat{u}_{Y} \right\rangle^2$,
where $\hat{\hat{\theta}}_X$ and $\hat{\hat{\theta}}_Y$ are asymptotic unbiased estimator of $\theta_X$ and $\theta_Y$ defined in \ref{Def=unbiased} and $\left\langle\hat{u}_{X},\hat{u}_{Y}\right\rangle$ is the scalar product between the two largest eigenvectors of $\hat{\Sigma}_X$ and $\hat{\Sigma}_Y$.\\
In this paper we want to generalise such test statistics to perturbations of order $k$ by considering functions of 
\begin{equation}
\hat{\hat{\theta}}_{X,1},...,\hat{\hat{\theta}}_{X,k},\hat{\hat{\theta}}_{Y,1},...,\hat{\hat{\theta}}_{Y,k}, \sum_{i=1}^k \left\langle \hat{u}_{X,1},\hat{u}_{Y,i} \right\rangle^2,...,\sum_{i=1}^k \left\langle \hat{u}_{X,k},\hat{u}_{Y,i} \right\rangle^2\,.
\label{eq:tig}
\end{equation}
Some possible tests are:
\begin{itemize}
\item $T_1=m \sum_{i=1}^k \left( \frac{\hat{\hat{\theta}}_{X,i}-\hat{\hat{\theta}}_{Y,i}} {\sigma_{\theta_i}} \right)^2 $, where $\sigma_{\theta_i}^2$ is the asymptotic variance of $\hat{\hat{\theta}}_{X,i}-\hat{\hat{\theta}}_{Y,i}$.
\item \begin{displaymath}
\resizebox{.85\hsize}{!}{$T_2=\sum_{i=1}^m\begin{pmatrix}
\hat{\hat{\theta}}_{X,i}-\hat{\hat{\theta}}_{Y,i}\\
\sum_{j=1}^k \left\langle \hat{u}_{X,i},\hat{u}_{Y,j} \right\rangle^2- \hat{\alpha}_{X,Y,i}^2
\end{pmatrix}^t \Sigma_{T_2}^{-1}\begin{pmatrix}
\hat{\hat{\theta}}_{X,i}-\hat{\hat{\theta}}_{Y,i}\\
\sum_{j=1}^k \left\langle \hat{u}_{X,i},\hat{u}_{Y,j} \right\rangle^2- \alpha_{X,Y,i}^2,
\end{pmatrix}$}
\end{displaymath}
where $\Sigma_{T_2}$ is the asymptotic variance of \linebreak $\left( \hat{\hat{\theta}}_{X,i}-\hat{\hat{\theta}}_{Y,i}, \ 
\sum_{j=1}^k \left\langle \hat{u}_{X,i},\hat{u}_{Y,j} \right\rangle^2- \hat{\alpha}_{X,Y,i}^2 \right)$.
\item 
\begin{displaymath}
\resizebox{.85\hsize}{!}{$T_3^{\pm}(s)=\lambda^{\pm}\left( \hat{\hat{\Sigma}}_X^{-1/2} \left( \left( \hat{\hat{\theta}}_{Y,s}-1\right) \hat{u}_{Y,s} \hat{u}_{Y,s}^t \right)  \hat{\hat{\Sigma}}_X^{-1/2} + \I_m \left( \frac{1}{\hat{\hat{\theta}}_{X,s}}-1 \right) \hat{u}_{X,s}\hat{u}_{X,s}^t \right)$}
\end{displaymath}
are also statistics of this form, where $\lambda^{\pm}()$ gives the extreme eigenvalues and $\hat{\hat{\Sigma}}_X$ is the filtered estimator defined in \ref{Def=unbiased}.
\item $\sum_{i=1}^m \lambda_i\left( \hat{\hat{\Sigma}}_X^{-1/2} \hat{\hat{\Sigma}}_Y \hat{\hat{\Sigma}}_X^{-1/2} \right)$
\item $ \sum_{i=1}^k \lambda_i\left( \hat{\hat{\Sigma}}_X^{-1/2} \hat{\hat{\Sigma}}_Y \hat{\hat{\Sigma}}_X^{-1/2} \right) $

\end{itemize}
In order to understand such statistics, we need to understand the joint properties of all the components in (\ref{eq:tig}).\\
The results of this paper show that the distributions of $\hat{\hat{\theta}}_X$ ,$\hat{\hat{\theta}}_Y$ and $\left\langle \hat{u}_{X},\hat{u}_{Y} \right\rangle^2$ we found for perturbation of order $1$ describe also the general case.

\subsection{Test statistic $T_1$}
Based on Theorem \ref{jointdistribution}, Theorem \ref{Thinvarianteigenvalue} and the fact that all the terms are uncorrelated by Theorem {\color{red}3.1} of \cite{mainarticle}, we can show that 
\begin{eqnarray*}
T_1 \sim \chi_k^2 + o_p(1),
\end{eqnarray*}
where 
\begin{equation*}
\resizebox{0.95\hsize}{!}{$\sigma_{\theta_i}^2=\sigma_{\theta_i,X}^2+\sigma_{\theta_i,Y}^2 =\frac{2 \left(M_{2,2,X}(\rho_{X,i})-M_{1,1,X}(\rho_{X,i})^2\right)}{M_{1,1,X}(\rho_{X,i})^4} +\frac{2 \left(M_{2,2,Y}(\rho_{Y,i})-M_{1,1,Y}(\rho_{Y,i})^2\right)}{M_{1,1,Y}(\rho_{Y,i})^4}.$}
\end{equation*}
Finally we can estimate $\sigma_{\theta_i}$ with $\hat{\sigma}_{\theta_i}$ by replacing $(\rho_{X,i},\rho_{Y,i})$ by \linebreak $(\hat{\theta}_{X,i},\hat{\theta}_{Y,i})$.
\subsection{Test statistic $T_2$}
We can show that 
\begin{equation*}
\resizebox{.95\hsize}{!}{$\sum_{i=1}^m\begin{pmatrix}
\hat{\hat{\theta}}_{X,i}-\hat{\hat{\theta}}_{Y,i}\\
\sum_{j=1}^k \left\langle \hat{u}_{X,i},\hat{u}_{Y,j} \right\rangle^2- \hat{\alpha}_{X,Y,i}^2
\end{pmatrix}^t \Sigma_{T_2}^{-1}\begin{pmatrix}
\hat{\hat{\theta}}_{X,i}-\hat{\hat{\theta}}_{Y,i}\\
\sum_{j=1}^k \left\langle \hat{u}_{X,i},\hat{u}_{Y,j} \right\rangle^2- \alpha_{X,Y,i}^2
\end{pmatrix} \sim \chi^2_{2k}+o(1),$}\,,
\end{equation*}
where 
\begin{equation*}
\resizebox{0.95\hsize}{!}{$M_{s_1,s_2,X}(\rho_X)=\frac{1}{m-k} \sum_{i=k+1}^m \frac{\hat{\lambda}_{\hat{\Sigma}_X,i}^{s_1}}{\left( \rho_X -\hat{\lambda}_{\hat{\Sigma}_X,i} \right)^{s_2}}  \text{ and } \ M_{s_1,s_2,Y}(\rho_Y)=\frac{1}{m-k} \sum_{i=k+1}^m \frac{\hat{\lambda}_{\hat{\Sigma}_Y,i}^{s_1}}{\left( \rho_Y -\hat{\lambda}_{\hat{\Sigma}_Y,i} \right)^{s_2}},$}
\end{equation*}
\begin{eqnarray*}
\hat{\alpha}_{X,Y,i}&=&\frac{\hat{\hat{\theta}}_{X,i}\hat{\hat{\theta}}_{Y,i}}{(\hat{\hat{\theta}}_{X,i}-1)^2(\hat{\hat{\theta}}_{Y,i}-1)^2} \frac{1}{\hat{\theta}_{X,i} \hat{\theta}_{Y,i} M_{1,2,X}\left(\hat{\theta}_{X,i}\right) M_{1,2,Y}\left(\hat{\theta}_{Y,i}\right)}, \\
\hat{\hat{\theta}}_{X,i}&=&\frac{1}{M_{1,1,X}\left(\hat{\theta}_{X,i}\right)}+1,\\
\\
\hat{\hat{\theta}}_{Y,i}&=&\frac{1}{M_{1,1,Y}\left(\hat{\theta}_{Y,i}\right)}+1.
\end{eqnarray*}
Moreover
\begin{eqnarray*}
\Sigma_{T_2}=\nabla\left(G\right)^t \Sigma \nabla\left(G\right),
\end{eqnarray*}
where $G: \mathbb{R}^3\rightarrow \mathbb{R}^2$ is such that 
\begin{eqnarray*}
\begin{pmatrix}
\hat{\hat{\theta}}_{X,i}-\hat{\hat{\theta}}_{Y,i}\\
\sum_{j=1}^k \left\langle \hat{u}_{X,i},\hat{u}_{Y,j} \right\rangle^2- \hat{\alpha}_{X,Y,i}^2
\end{pmatrix}= G\begin{pmatrix}
\hat{\theta}_{X,i}\\
\hat{\theta}_{Y,i}\\
\sum_{j=1}^k \left\langle \hat{u}_{X,i},\hat{u}_{Y,j} \right\rangle^2
\end{pmatrix}.
\end{eqnarray*}
and  
\begin{eqnarray*}
\begin{pmatrix}
\hat{\theta}_{X,i}\\
\hat{\theta}_{Y,i}\\
\sum_{j=1}^k \left\langle \hat{u}_{X,i},\hat{u}_{Y,j} \right\rangle^2
\end{pmatrix} \sim \Normal\left( 
 \begin{pmatrix}
\rho_{X,i} \\
\rho_{Y,i}\\
\frac{\theta_{X,i}\theta_{Y,i}}{(\theta_{X,i}-1)^2(\theta_{Y,i}-1)^2} \frac{1}{\rho_X \rho_Y M_{1,2,X}M_{1,2,Y}}
\end{pmatrix}, \Sigma \right),
\end{eqnarray*}
Using similar argument as in the proof of Theorem \ref{jointdistribution}, we can show that \\
\scalebox{0.4}{
\begin{minipage}{1\textwidth}
\begin{eqnarray*}
\Sigma_{1,1}&=&-\frac{2  (M_{1,1,X}(\rho_X)+M_{1,1,X}(\rho_X)^2-M_{1,2,X}(\rho_X)  \rho_X)}{M_{1,1,X}(\rho_X)^4}\\
\Sigma_{1,2}&=&0\\
\Sigma_{1,3}&=&\frac{2  (M_{1,1,X}(\rho_X)  (1+M_{1,1,X}(\rho_X))  M_{1,2,X}(\rho_X)+M_{1,1,X}(\rho_X)  (M_{1,2,X}(\rho_X)^2-2  (1+M_{1,1,X}(\rho_X))  M_{1,3,X}(\rho_X))  \rho_X+M_{1,2,X}(\rho_X)  M_{1,3,X}(\rho_X)  \rho_X^2)  \theta_{X,i}  \theta_{Y,i}}{M_{1,1,X}(\rho_X)^2  M_{1,2,X}(\rho_X)^3  M_{1,2,Y}(\rho_Y)  \rho_X^2  \rho_Y  (-1+\theta_{X,i})^2  (-1+\theta_{Y,i})^2}\\
\Sigma_{2,2}&=&-\frac{2  (M_{1,1,Y}(\rho_Y)+M_{1,1,Y}(\rho_Y)^2-M_{1,2,Y}(\rho_Y)  \rho_Y)}{M_{1,1,Y}(\rho_Y)^4}\\
\Sigma_{2,3}&=&\frac{2  (M_{1,1,Y}(\rho_Y)  (1+M_{1,1,Y}(\rho_Y))  M_{1,2,Y}(\rho_Y)+M_{1,1,Y}(\rho_Y)  (M_{1,2,Y}(\rho_Y)^2-2  (1+M_{1,1,Y}(\rho_Y))  M_{1,3,Y}(\rho_Y))  \rho_Y+M_{1,2,Y}(\rho_Y)  M_{1,3,Y}(\rho_Y)  \rho_Y^2)  \theta_{X,i}  \theta_{Y,i}}{M_{1,2,X}(\rho_X)  M_{1,1,Y}(\rho_Y)^2  M_{1,2,Y}(\rho_Y)^3  \rho_X  \rho_Y^2  (-1+\theta_{X,i})^2  (-1+\theta_{Y,i})^2}\\
\Sigma_{3,3}&=&\Bigg( \left.2  \theta_{X,i}  \theta_{Y,i}  (2  M_{1,2,X}(\rho_X)^5  M_{1,2,Y}(\rho_Y)^4  \rho_X^3  \rho_Y^2  (-1 + 
      \theta_{X,i})^2  (M_{1,2,Y}(\rho_Y)  \rho_Y  (-1 + \theta_{Y,i})^2 - \theta_{Y,i}) - (1 + 
      2  M_{1,1,X}(\rho_X))  M_{1,2,X}(\rho_X)^3  M_{1,2,Y}(\rho_Y)^4  \rho_X  \rho_Y^2  \theta_{X,i}  \theta_{Y,i} + \right.\\
   && \hspace{0.5cm} \left.
   4  M_{1,1,X}(\rho_X)  (1 + M_{1,1,X}(\rho_X))  M_{1,2,X}(\rho_X)  M_{1,3,X}(\rho_X)  M_{1,2,Y}(\rho_Y)^4  \rho_X  \rho_Y^2  \theta_{X,i}  \theta_{Y,i} 
   -4  M_{1,1,X}(\rho_X)  (1 + M_{1,1,X}(\rho_X))  M_{1,3,X}(\rho_X)^2  M_{1,2,Y}(\rho_Y)^4  \rho_X^2  \rho_Y^2  \theta_{X,i}  \theta_{Y,i} + 
   \right.\\
   && \hspace{0.5cm} \left.
   M_{1,2,X}(\rho_X)^2  M_{1,2,Y}(\rho_Y)^4  (-M_{1,1,X}(\rho_X)  (1 + M_{1,1,X}(\rho_X)) + (1 + 4  M_{1,1,X}(\rho_X))  M_{1,3,X}(\rho_X)  \rho_X^2 + 
      M_{1,4,X}(\rho_X)  \rho_X^3)  \rho_Y^2  \theta_{X,i}  \theta_{Y,i} - \right.\\
   && \hspace{0.5cm} \left.
   M_{1,2,X}(\rho_X)^4  \rho_X^2  \theta_{X,i}  (2  M_{1,2,Y}(\rho_Y)^5  \rho_Y^3  (-1 + \theta_{Y,i})^2 + (1 + 
         2  M_{1,1,Y}(\rho_Y))  M_{1,2,Y}(\rho_Y)^3  \rho_Y  \theta_{Y,i} - 
      4  M_{1,1,Y}(\rho_Y)  (1 + M_{1,1,Y}(\rho_Y))  M_{1,2,Y}(\rho_Y)  M_{1,3,Y}(\rho_Y)  \rho_Y  \theta_{Y,i} +\right.\\
   && \hspace{0.5cm} \left. 
      4  M_{1,1,Y}(\rho_Y)  (1 + M_{1,1,Y}(\rho_Y))  M_{1,3,Y}(\rho_Y)^2  \rho_Y^2  \theta_{Y,i} + 
      M_{1,2,Y}(\rho_Y)^2  (M_{1,1,Y}(\rho_Y)  (1 + M_{1,1,Y}(\rho_Y)) - (2  M_{1,3,Y}(\rho_Y) + 4  M_{1,1,Y}(\rho_Y)  M_{1,3,Y}(\rho_Y) + 
            M_{2,4,Y}(\rho_Y))  \rho_Y^2)  \theta_{Y,i}))\right. \Bigg) \\
             &&\hspace{2cm} \Bigg/ \left(M_{1,2,X}(\rho_X)^6  M_{1,2,Y}(\rho_Y)^6  \rho_X^4  \rho_Y^4  (-1 + 
   \theta_{X,i})^4  (-1 + \theta_{Y,i})^4\right)
\end{eqnarray*}.
\end{minipage}}\\
Finally we can estimate $\Sigma$ with $\hat{\Sigma}$ by replacing $(\rho_{X,i},\theta_{X,i},\rho_{Y,i},\theta_{X,i})$ by \linebreak $(\hat{\theta}_{X,i},\hat{\hat{\theta}}_{X,i},\hat{\theta}_{Y,i},\hat{\hat{\theta}}_{Y,i})$ and $\Sigma_{T_2}$ with
\begin{eqnarray*}
\hat{\Sigma}_{T_2}=\nabla\left(G\right)^t \hat{\Sigma} \nabla\left(G\right).
\end{eqnarray*}

\subsection{Test statistic $T_3$}
Elementary linear algebra in conjunction with the theorems of \cite{mainarticle} and this paper show that
\begin{eqnarray*}
\resizebox{1\hsize}{!}{$T_3^{\pm}(s)=\frac{1}{2} \left(\hat{\hat{\theta}}_{Y,s}+\sum_{i=1}^k \left\langle \hat{u}_{Y,s},\hat{u}_{X,i} \right\rangle^2-\hat{\hat{\theta}}_{Y,s} \left( \sum_{i=1}^k \left\langle \hat{u}_{Y,s},\hat{u}_{X,i} \right\rangle^2 \right)+\frac{1+(\hat{\hat{\theta}}_{Y,s}-1)\left( \sum_{i=1}^k \left\langle \hat{u}_{Y,s},\hat{u}_{X,i} \right\rangle^2 \right) \pm \sqrt{-4 \hat{\hat{\theta}}_{Y,s} \hat{\hat{\theta}}_{X,s}+ \left( 1+\hat{\hat{\theta}}_{Y,s} \hat{\hat{\theta}}_{X,s} - (\hat{\hat{\theta}}_{X,s}-1) (\hat{\hat{\theta}}_{X,s}-1) \left(\sum_{i=1}^k \left\langle \hat{u}_{Y,s},\hat{u}_{X,i} \right\rangle^2 \right) \right)^2}}{\hat{\hat{\theta}}_{X,s}} \right)+O\left( \frac{1}{m} \right).$}
\end{eqnarray*}
This result can be obtain by looking at the trace and the square of the matrix. This statistic is the residual spike defined in \cite{mainarticle}. Therefore $T_3$ is bounded by
\begin{eqnarray*}
\Normal \left( \lambda^+  , \frac{{\sigma^{+}}^2}{m} \right)+o\left( \frac{1}{\sqrt{m}} \right) \text{ and } \Normal \left( \lambda^-  , \frac{{\sigma^{-}}^2}{m} \right)+o\left( \frac{1}{\sqrt{m}} \right)\,,
\end{eqnarray*}
with the parameters as defined in Theorem 2.1 of \cite{mainarticle}.

\subsection{Simulation}
Assume $\mathbf{X} \in \mathbb{R}^{m\times n_X}$ and $\mathbf{Y} \in \mathbb{R}^{m\times n_Y}$ with $\mathbf{X}=\left(X_1,X_2,...,X_{n_X}\right)$ and $\mathbf{Y}=\left(Y_1,Y_2,...,Y_{n_Y}\right)$. The components of the random vectors are independent and the covariance between the vectors is as follows:\\
\scalebox{0.65}{
\begin{minipage}{1\textwidth}
\begin{eqnarray*}
&&X_i \sim \Normal_m\left(\vec{0},\sigma^2 \I_m \right) \text{ with } X_1 = \epsilon_{X,1}\text{ and } 
X_{i+1}=\rho X_i+\sqrt{1-\rho^2} \ \epsilon_{X,i+1}, \text{ where }  \epsilon_{X,i} \overset{i.i.d}{\sim} \Normal_m\left(\vec{0},\sigma^2 \I_m\right),\\
&&Y_i \sim \Normal_m\left(\vec{0},\sigma^2 \I_m \right) \text{ with } Y_1 = \epsilon_{Y,1}\text{ and } Y_{i+1}=\rho Y_i+\sqrt{1-\rho^2} \ \epsilon_{Y,i+1}, \text{ where }  \epsilon_{Y,i} \overset{i.i.d}{\sim} \Normal_m\left(\vec{0},\sigma^2 \I_m\right)
\end{eqnarray*}
\end{minipage}} \vspace{0.2cm}

Let $P_X= \I_m + \sum_{i=1}^k (\theta_{X,i}-1) u_{X,i} u_{X,i}^t$ and $P_Y= \I_m+ \sum_{i=1}^k (\theta_{Y,i}-1) u_{Y,i} u_{Y,i}^t$ be two perturbations in  $\mathbb{R}^{m\times m}$ and put
\begin{eqnarray*}
\mathbf{X}_P=P_X^{1/2} \mathbf{X} \text{ and  } \mathbf{Y}_P=P_Y^{1/2} \mathbf{Y},\\
\hat{\Sigma}_X=\frac{\mathbf{X}_P^t \mathbf{X}_P}{n_X} \text{ and } \hat{\Sigma}_Y=\frac{\mathbf{Y}_P^t \mathbf{Y}_P}{n_Y}.
\end{eqnarray*}


\subsubsection{Comparison with existing tests}
In the classical multivariate theory, the trace or the determinant of $\hat{\Sigma}_X^{-1/2} \hat{\Sigma}_Y \hat{\Sigma}_X^{-1/2}$ are used to test the equality of two covariance matrices (see, for example, \cite{multi22}).\\
Suppose
\begin{eqnarray*}
X_1,X_2,...,X_{n_X} \overset{i.i.d.}\sim \Normal_m(0,\Sigma_X),\\
Y_1,Y_2,...,Y_{n_Y} \overset{i.i.d.}\sim \Normal_m(0,\Sigma_Y).
\end{eqnarray*}
We want to test 
\begin{eqnarray*}
{\rm H}_0 : \Sigma_X=\Sigma_Y\text{ against }{\rm H}_1 : \Sigma_X \neq \Sigma_Y,
\end{eqnarray*} 
In this section we show that any test statistic using either the log-determinant \linebreak $T_4=\log\left| \hat{\Sigma}_X^{-1/2} \hat{\Sigma}_Y \hat{\Sigma}_X^{-1/2} \right|$ or $T_5=\Tr \left( \hat{\Sigma}_X^{-1/2} \hat{\Sigma}_Y \hat{\Sigma}_X^{-1/2} \right)$ have difficulties to detect differences between the finite perturbations $P_X$ and $P_Y$.
To explore this problem, we compare the performance of these tests with $T_1$, $T_2$ and $T_3$ by simulation. Table \ref{Tabletest} shows the power of these tests to detect under a variety of alternatives and sample sizes. For $T_1$ and $T_2$ the critical values are based on the asymptotic chi-squared distributions, for $T_3$ the following two-sided power is used 
\begin{equation*}
\resizebox{.9\hsize}{!}{$P_{{\rm H}_1}\left( \underset{s=1,2,...,k}{\max}\left(\sqrt{m}\frac{T_3^+(s) -\lambda^+}{\sigma^+}\right) < q_{\Normal(0,1)} (1-0.025/k) \text{ or }  \underset{s=1,2,...,k}{\min}\left(\sqrt{m}\frac{T_3^-(s) -\lambda^-}{\sigma^-}\right) < q_{\Normal(0,1)} (0.025/k) \right),$}
\end{equation*}
with the parameters of Theorem 2.1 of \cite{mainarticle}.
For the tests $T_4$ and $T_5$ the critical values are determined by simulation.
In order to apply these tests to degenerated matrices, the determinant is defined as the product of the non-null eigenvalues of the matrix and the inverse is the generalised inverse.

\begin{table}[hbt] \scalebox{0.6}{
\begin{tabular}{ccccccc}
\begin{minipage}{1cm}
\begin{eqnarray*}
&&m=500,\\
&&n_X=n_Y=250
\end{eqnarray*}
\end{minipage}
 & \\
&\begin{minipage}{3cm}
\begin{eqnarray*}
&&\theta_X=7, u_X=e_1,\\
&&\theta_Y=7, u_Y=e_1,
\end{eqnarray*}
\end{minipage}&\begin{minipage}{3cm}
\begin{eqnarray*}
&&\theta_X=7, u_X=e_1,\\
&&\theta_Y=7, u_Y=e_2,
\end{eqnarray*}
\end{minipage} &
\begin{minipage}{3cm}
\begin{eqnarray*}
&&\theta_X=50, u_X=e_1,\\
&&\theta_Y=50, u_Y=e_2,
\end{eqnarray*}
\end{minipage}  &
\begin{minipage}{3cm}
\begin{eqnarray*}
&&\theta_X =7, u_X=e_1,\\
&&\theta_Y =17, u_Y=e_1,
\end{eqnarray*}
\end{minipage} &
\begin{minipage}{3cm}
\begin{eqnarray*}
&&\theta_X=300, u_X=e_1,\\
&&\theta_Y=600, u_Y=e_1,
\end{eqnarray*} 
\end{minipage} &\\
\hline
$T_1$  & 0.05 & 0.04 & 0.05 & 1    & 0.91\\
$T_2$ & 0.06 & 1     & 1    & 1    & 0.99\\
$T_3$ & 0 & 0.37     & 1    & 0.85 & 0.995\\
$T_4$ & 0.04 & 0.06  & 0.11 & 0.06 & 0.06\\
$T_5$ & 0.035 & 0.04 & 1    & 0.12 & 0.07\\
\hline 
\hline
\begin{minipage}{1cm}
\begin{eqnarray*}
&&m=500,\\
&&n_X=1000,\\
&&n_Y=250
\end{eqnarray*}
\end{minipage}
 & \\
&\begin{minipage}{3cm}
\begin{eqnarray*}
&&\theta_X=7, u_X=e_1,\\
&&\theta_Y=7, u_Y=e_1,
\end{eqnarray*}
\end{minipage}&\begin{minipage}{3cm}
\begin{eqnarray*}
&&\theta_X=7, u_X=e_1,\\
&&\theta_Y=7, u_Y=e_2,
\end{eqnarray*}
\end{minipage} &
\begin{minipage}{3cm}
\begin{eqnarray*}
&&\theta_X=50, u_X=e_1,\\
&&\theta_Y=50, u_Y=e_2,
\end{eqnarray*}
\end{minipage}  &
\begin{minipage}{3cm}
\begin{eqnarray*}
&&\theta_X =7, u_X=e_1,\\
&&\theta_Y =17, u_Y=e_1,
\end{eqnarray*}
\end{minipage} &
\begin{minipage}{3cm}
\begin{eqnarray*}
&&\theta_X=300, u_X=e_1,\\
&&\theta_Y=600, u_Y=e_1,
\end{eqnarray*} 
\end{minipage} &\\
\hline
$T_1$ & 0.06 & 0.045 & 0.06 & 1    & 1 \\
$T_2$ & 0.06 & 1     & 1    & 1    & 0.99 \\
$T_3$ & 0.01 & 1     & 1    & 0.96 & 1 \\
$T_4$ & 0.075 & 0.12 & 0.55 & 0.08 & 0.05\\
$T_5$ & 0.05 & 0.25  & 1    & 0.06 & 0.09\\
\hline 
\hline
\begin{minipage}{1cm}
\begin{eqnarray*}
&&m=500,\\
&&n_X=1000,\\
&&n_Y=1000
\end{eqnarray*}
\end{minipage}
 & \\
&\begin{minipage}{3cm}
\begin{eqnarray*}
&&\theta_X=7, u_X=e_1,\\
&&\theta_Y=7, u_Y=e_1,
\end{eqnarray*}
\end{minipage}&\begin{minipage}{3cm}
\begin{eqnarray*}
&&\theta_X=7, u_X=e_1,\\
&&\theta_Y=7, u_Y=e_2,
\end{eqnarray*}
\end{minipage} &
\begin{minipage}{3cm}
\begin{eqnarray*}
&&\theta_X=50, u_X=e_1,\\
&&\theta_Y=50, u_Y=e_2,
\end{eqnarray*}
\end{minipage}  &
\begin{minipage}{3cm}
\begin{eqnarray*}
&&\theta_X =7, u_X=e_1,\\
&&\theta_Y =17, u_Y=e_1,
\end{eqnarray*}
\end{minipage} &
\begin{minipage}{3cm}
\begin{eqnarray*}
&&\theta_X=300, u_X=e_1,\\
&&\theta_Y=600, u_Y=e_1,
\end{eqnarray*} 
\end{minipage} &\\
\hline
$T_1$ & 0.03 & 0.05 & 0.06  & 1   & 1\\
$T_2$ & 0.03 & 1    & 1     & 1   & 1\\
$T_3$ & 0    & 1    & 1     & 1   & 1\\
$T_4$ & 0.045& 0.07 & 0.01  & 0.1  & 0.04\\
$T_5$ & 0.04 & 0.56 & 1     & 0.09 & 0.04
\end{tabular}}
\caption{Probability to detect the alternative with a test at level $0.05$ when $P_X=\I_m+ 500 e_3 e_3^2 + 150 e_4 e_4^t + (\theta_X-1) u_X u_X^t$ and $P_Y=\I_m+ 500 e_3 e_3^2 + 150 e_4 e_4^t +(\theta_Y-1) u_Y u_Y^t$ for the different tests. The distribution of $T_4$ and $T_5$ is computed empirically by assuming the same perturbation $P_X$ for the two groups.} \label{Tabletest}
\end{table}
In the simulated cases, the trace and the determinant have difficulties to catch the alternatives. On the other hand, our procedures easily detect even small effects. These classical statistics $T_4$ and $T_5$ would presumably do well with global perturbations such as a multiplicative change of the covariance matrix.

\begin{Rem}
\begin{enumerate}

\item Under the assumption that $\hat{\Sigma}_X=P^{1/2}_X W_X P^{1/2}_X$ and $\hat{\Sigma}_Y=P^{1/2}_Y W_Y P^{1/2}_Y$ satisfy Assumption \ref{Ass=matrice}, the procedures $T_1$, $T_2$ and $T_3$ required the estimation of $M_{s,r,X}=\frac{1}{m}\sum_{i=1}^m \frac{\hat{\lambda}_{W_X,i}^s}{\left(\rho-\hat{\lambda}_{W_X,i}\right)^r}$ and $M_{s,r,Y}=\frac{1}{m}\sum_{i=1}^m \frac{\hat{\lambda}_{W_Y,i}^s}{\left(\rho-\hat{\lambda}_{W_Y,i}\right)^r}$ for $s,r=1,2,3,4$. By Cauchy's interlacing law and the upper bound on the eigenvalues of $\hat{\lambda}_{W_X,i}$ and $\hat{\lambda}_{W_Y,i}$, we can use the following estimator
\begin{eqnarray*}
\hat{M}_{s,r,X}=\frac{1}{m-k}\sum_{i=k+1}^m \frac{\hat{\lambda}_{\hat{\Sigma}_X,i}^s}{\left(\rho-\hat{\lambda}_{\hat{\Sigma}_X,i}\right)^r}= M_{s,r,X} +O\left(\frac{1}{m}\right).
\end{eqnarray*}
\item The theorems of this paper always assume perturbations with distinct eigenvalues. When $\theta_1=\theta_2=...\theta_k$, the results fail and most of the procedures are not conservative.  
\end{enumerate}
\end{Rem}

\section{Theorems} \label{sec:Theorems}
\subsection{Notation and definition}

\begin{Not}\ \label{Not=Theorem}\\
We use a precise notation to enunciate the theorems, the proofs, however, often use a simpler notation when no confusion is possible. This difference is always specified at the beginning of the proofs. 
\begin{itemize}
\item For any symmetric random matrix $A$ we denote by $\left(\hat{\lambda}_{A,i}, \hat{u}_{A,i} \right)$ its $i^{\rm th}$ eigenvalue and eigenvector.
\item A finite perturbation of order $k$ is denoted by $P_k= \I_m  +\sum_{i=1}^k (\theta_i-1) u_i u_i^t \in \mathbb{R}^{m \times m}$ with $u_1,u_2,...,u_k \in \mathbb{R}^{m \times m}$ orthonormal vectors.
\item $W \in  \mathbb{R}^{m\times m}$ denotes a random matrix as defined in Assumption \ref{Ass=matrice} which is invariant under rotation. Moreover, the estimated covariance matrix is $\hat{\Sigma}=P_k^{1/2} W P_k^{1/2}$. \\
When comparing two groups, we use $W_X$, $W_Y$ and $\hat{\Sigma}_X$, $\hat{\Sigma}_Y$.
\item When we consider only one group, $\hat{\Sigma}_{P_r}=P_r^{1/2} W P_r^{1/2}$ is the perturbation of order $r$ of the matrix $W$ and: 
\begin{itemize}
\item $\hat{u}_{P_r,i}$ is its $i^{\text{th}}$ eigenvector. When $r=k$ we just use the simpler notation $\hat{u}_{i}=\hat{u}_{P_k,i}$ after an explicit statement.
\item  $\hat{u}_{P_r,i,j}$ is the $j^{\rm th}$ component of the $i^{\rm th}$ eigenvector.
\item $\hat{\lambda}_{P_r,i}$ is the $i^{\text{th}}$ eigenvalue. If $\theta_1>\theta_2>...>\theta_r$, then for $i=1,2,...,r$ we use also the notation $\hat{\theta}_{P_r,i}=\hat{\lambda}_{P_r,i}$. We call these eigenvalues the spikes. When $r=k$, we just use the simpler notation $\hat{\theta}_i=\hat{\theta}_{P_k,i}$ after an explicit statement.
\item $\hat{\alpha}_{P_r,i}^2=\sum_{j=1}^r \left\langle \hat{u}_{P_r,i},u_j \right\rangle^2$ is called the \textbf{general angle}.
\end{itemize}  
With this notation, we have $\hat{\Sigma}=\hat{\Sigma}_{P_k}=P_k^{1/2} W P_k^{1/2}$.
\item When we consider two groups $X$ and $Y$, we use a notation similar to the above. The perturbation of order $r$ of the matrices $W_X$ and $W_Y$ are $\hat{\Sigma}_{X,P_r}=P_r^{1/2} W_X P_r^{1/2}$ and $\hat{\Sigma}_{Y,P_r}=P_r^{1/2} W_Y P_r^{1/2}$ respectively. Then, we define for the group $\hat{\Sigma}_{X,P_r}$ (and similarly for $\hat{\Sigma}_{Y,P_r}$):
\begin{itemize}
\item $\hat{u}_{\hat{\Sigma}_{X,P_r},i}$ is its $i^{\text{th}}$ eigenvector. When $r=k$ we use the simpler notation $\hat{u}_{X,i}=\hat{u}_{\hat{\Sigma}_{X,P_k},i}$.
\item  $\hat{u}_{\hat{\Sigma}_{X,P_r},i,j}$ is the $j^{\rm th}$ component of the $i^{\rm}$ eigenvector.
\item $\hat{\lambda}_{\hat{\Sigma}_{X,P_r},i}$ is its $i^{\text{th}}$ eigenvalue. If $\theta_1>\theta_2>...>\theta_r$, then for $i=1,2,...,r$ we use the notation $\hat{\theta}_{\hat{\Sigma}_{X,P_r},i}=\hat{\lambda}_{\hat{\Sigma}_{X,P_r},i}$. When $r=k$, we use the simpler notation $\hat{\theta}_{X,i}=\hat{\theta}_{\hat{\Sigma}_{X,P_k},i}$.
\item $\hat{\alpha}_{\hat{\Sigma}_{X,P_r},i}^2=\sum_{j=1}^r \left\langle \hat{u}_{\hat{\Sigma}_{X,P_r},i},u_j \right\rangle^2$. 
\item $\hat{\alpha}_{X,Y,P_r,i}^2=\sum_{j=1}^r \left\langle \hat{u}_{\hat{\Sigma}_{X,P_r},i},\hat{u}_{\hat{\Sigma}_{Y,P_r},j} \right\rangle^2$ is the \textbf{double angle} and, when no confusion is possible, we use the simpler notation $\hat{\alpha}_{P_r,i}^2$.
\end{itemize} 
\item Some theorems assume the sign convention 
\begin{eqnarray*}
\hat{u}_{P_s,i,i}>0, \text{ for $s=1,2,...,k$ and $i=1,2,...,s$,}
\end{eqnarray*}
 as in Theorem \ref{ThInvariantdot} or \ref{Thcomponentdistribution}. Others assume the convention 
\begin{eqnarray*}
\hat{u}_{P_s,i,s}>0, \text{ for $s=1,2,...,k$ and $i=1,2,...,s$,}
\end{eqnarray*}
as in Theorem \ref{Theoremcaraceigenstructure}.\\
Theorems that are not affected by this convention do not specify it precisely. Nevertheless, the convention will be mentioned in the proofs when confusion is possible.
\item We define the function $M_{s_1,s_2,X}(\rho_X)$, $M_{s_1,s_2,Y}(\rho_Y)$ and $M_{s_1,s_2}(\rho_X,\rho_Y)$ as 
\begin{eqnarray*}
M_{s_1,s_2,X}(\rho_X)&=&\frac{1}{m} \sum_{i=1}^m \frac{\hat{\lambda}_{W_X,i}^{s_1}}{\left( \rho_X -\hat{\lambda}_{W_X,i} \right)^{s_2}},\\
M_{s_1,s_2,Y}(\rho_Y)&=&\frac{1}{m} \sum_{i=1}^m \frac{\hat{\lambda}_{W_Y,i}^{s_1}}{\left( \rho_Y -\hat{\lambda}_{W_Y,i} \right)^{s_2}},\\
M_{s_1,s_2}(\rho_X,\rho_Y)&=& \frac{M_{s_1,s_2,X}(\rho_X)+M_{s_1,s_2,Y}(\rho_Y)}{2}.
\end{eqnarray*}
In particular, when $s_2=0$, we use $M_{s_1,X}=M_{s_1,0,X}$. When we only study one group, we use the simpler notation $M_{s_1,s_2}(\rho)$ when no confusion is possible.
\item We use two transforms inspired by the T-transform:
\begin{itemize}
\item $T_{W,u}(z)= \sum_{i=1}^m \frac{\hat{\lambda}_{W,i}}{z-\hat{\lambda}_{W,i}} \left\langle \hat{u}_{W,i},u \right\rangle^2$ is the T-transform in direction $u$ using the random matrix $W$.
\item $\hat{T}_{\hat{\Sigma}_X}(z)= \frac{1}{m}\sum_{i=k+1}^m \frac{\hat{\lambda}_{\hat{\Sigma}_X,i}}{z-\hat{\lambda}_{\hat{\Sigma}_X,i}}$,
 and  $\hat{T}_{W_X}(z)= \frac{1}{m} \sum_{i=1}^m \frac{\hat{\lambda}_{W_X,i}}{z-\hat{\lambda}_{W_X,i}}$, the estimated T-transforms using $\hat{\Sigma}_X$ and $W$ respectively.
\end{itemize}
\item In some theorems we use the notation $\overset{\scalebox{0.5}{order}}{\sim}$ to describe the order size in probability of a positive random variable. For example, $X_m \overset{\scalebox{0.5}{order}}{\sim} 1/m$ if
$\frac{X_m}{1/m}$ tends to a random variable $X$ independent of $m$, with $P\left\lbrace X>\epsilon_j \right\rbrace \overset{j\rightarrow \infty}{\longrightarrow} 1$ for any sequences $\epsilon_j$ tending to $0$.
\end{itemize}
\end{Not}

This paper extends previous results to perturbations of order $k>1$ for some \textbf{invariant} statistics.
\begin{Def}\ \label{Def:Invariant}\\
Suppose $W$ is a random matrix. Moreover, define $P_1=\I_m+(\theta_1-1) u_1 u_1^t$ and $P_k=\I_m+\sum_{i=1}^k (\theta_i-1)u_i u_i^t$ some perturbations of order $1$ and $k>1$, respectively. We say that a statistic $T\left( W_m,P_1\right)$ is \textbf{invariant} with respect to $k$, if $T\left( W_m,P_k\right)$ is such that 
\begin{eqnarray*}
\resizebox{.9\hsize}{!}{$T\left( W_m,P_k\right) = T\left( W_m,P_1\right) + \epsilon_m, \text{ where }
\max\left(\frac{\epsilon_m}{\E\left[ T\left( W,P_1\right)\right]},\frac{\epsilon_m^2}{\var\left( T\left( W,P_1\right) \right)}\right) \rightarrow 0.$}
\end{eqnarray*}
\end{Def}

\subsection{Invariant Eigenvalue Theorem}
Theorem \ref{jointdistribution} provides distributions of statistics for perturbations of order $1$. This estimated eigenvalue is an invariant statistics as defined in \ref{Def:Invariant}.

\begin{Th} \label{Thinvarianteigenvalue} Suppose that $W$ satisfies Assumption \ref{Ass=matrice} and
\begin{eqnarray*}
&&\tilde{P}_s=\I_m+(\theta_s-1) e_s e_s^t, \text{ for } s=1,2,...,k,\\
&&P_k=\I_m+\sum_{i=1}^k (\theta_i-1) e_i e_i^t \text{ satisfies \ref{Ass=theta} (A4),}
\end{eqnarray*} 
where $\theta_1>\theta_2>...>\theta_k$.
We define
\begin{eqnarray*}
&&\hat{\Sigma}_{\tilde{P}_s}=\tilde{P}_s^{1/2} W \tilde{P}_s^{1/2},\\
&&\hat{\Sigma}_{P_k}= P_k^{1/2} W P_k^{1/2}.
\end{eqnarray*}
Moreover, for $s=1,2,...,k$, we define
\begin{eqnarray*}
\hat{u}_{\tilde{P}_s,1}, \hat{\theta}_{\tilde{P}_s,1} &\text{ s.t. }& \hat{\Sigma}_{\tilde{P}_1} \hat{u}_{\tilde{P}_s,1} =\hat{\theta}_{\tilde{P}_s,1} \hat{u}_{\tilde{P}_s,1},  \\
\hat{u}_{P_k,s},\hat{\theta}_{P_k,s} &\text{ s.t. }& \hat{\Sigma}_{P_k} \hat{u}_{P_k,s}= \hat{\theta}_{P_k,s}\hat{u}_{P_k,s},  
\end{eqnarray*}
where $\hat{\theta}_{\tilde{P}_s,1}= \hat{\lambda}_{\hat{\Sigma}_{\tilde{P}_s,1}}$ and
 $\hat{\theta}_{P_k,s}=\hat{\lambda}_{\hat{\Sigma}_{P_k},s}$. 
\begin{enumerate}
\item  Then, for $s>1$,
$$ \boxed{\hat{\theta}_{P_{k},s} - \hat{\theta}_{\tilde{P}_{s},1} \overset{\scalebox{0.5}{order}}{\sim} \frac{\theta_s}{m}} $$
and
$$ \boxed{\hat{\theta}_{P_{k},1} - \hat{\theta}_{\tilde{P}_{1},1} \overset{\scalebox{0.5}{order}}{\sim} \frac{\theta_2}{m},}.$$
The distribution of $ \hat{\theta}_{P_{k},s} $ is therefore asymptotically the same as the distribution of $\hat{\theta}_{\tilde{P}_{s},1}$ studied in Theorem \ref{jointdistribution}.
\item 
More precisely we define for $r,s \in \left\lbrace 1,2,...,k \right\rbrace$ with $r \neq s $ ,
$$ P_{-r}=\I_m + \sum_{\underset{i\neq r}{i=1}}^k \left( \theta_i-1 \right) e_i e_i^t. $$
\begin{itemize}
\item If $\theta_s> \theta_r$, then 
\begin{equation*}
 \resizebox{.8\hsize}{!}{$\hat{\theta}_{P_{k},s}-\hat{\theta}_{P_{-r},s}
= 
-\frac{\hat{\theta}_{P_{-r},s} \hat{\theta}_{P_{k},s} (\theta_r-1)}{\theta_r-1 -\hat{\theta}_{P_{k},s}} \hat{u}_{P_{-r},s,r}^2 + O_p \left(\frac{1}{m} \right)+ O_p \left(\frac{\theta_r}{m^{3/2}} \right) .$}
\end{equation*}
\item If $\theta_s < \theta_r$, then 
\begin{equation*}
\resizebox{.8\hsize}{!}{$ \hat{\theta}_{P_{k},s}-\hat{\theta}_{P_{-r},s-1}
= 
-\frac{\hat{\theta}_{P_{-r},s-1} \hat{\theta}_{P_{k},s} (\theta_r-1)}{\theta_r-1 -\hat{\theta}_{P_{k},s}} \hat{u}_{P_{-r},s-1,r}^2 + O_p \left(\frac{1}{m} \right)+ O_p \left(\frac{\theta_s}{m^{3/2}} \right).$}
\end{equation*}
\end{itemize}

\end{enumerate} 
 
\begin{Rem}\ \\
In this manuscript, we are interested in the unbiased estimation of $\hat{\hat{\theta}}_{P_{k},1}$. The invariance of $\hat{\hat{\theta}}_{P_{k},1}$ is a direct consequence of the theorem. Moreover, 
Theorem \ref{jointdistribution} provides the distribution of $\hat{\hat{\theta}}_{P_{1},1}$.
\end{Rem}

\end{Th}
(Proof in appendix \ref{appendixproof}.)

\subsection{Invariant Angle Theorem}\label{Invariantanglesection}

The cosine of the angle between two vectors is linked to $\left\langle u,v \right\rangle$. We need the more general notion of the angle between a vector and a subspace of dimension $k$ associated with $\sum_{i=1}^k \left\langle u,v_i \right\rangle^2$, where $(v_1,...,v_k)$ is a orthonormal basis of the subspace. This generalization of the angle used with the correct subspace leads to an invariance in the sense of Definition \ref{Def:Invariant}.

\begin{Th}\label{Invariantth}\ \\
Using the same notation as Theorem \ref{Thinvarianteigenvalue},
\begin{enumerate}
\item The general angle is invariant in the sense of Definition \ref{Def:Invariant},
\begin{eqnarray*}
&& \boxed{\sum_{i=1}^k \hat{u}_{P_k,s,i}^2=\hat{u}_{\tilde{P}_s,1,s}^2 + O_p\left( \frac{1}{\theta_s m} \right).}
\end{eqnarray*}
Therefore,  the distribution of $\sum_{i=1}^k \hat{u}_{P_k,s,i}^2 $ is asymptotically the same as the distribution of $\hat{u}_{\tilde{P}_s,1,s}^2$ studied in Theorem \ref{jointdistribution}.
\item Moreover, 
$$ \hat{u}_{P_k,s,s}^2 = \hat{u}_{\tilde{P}_s,1,s}^2 + O_p\left( \frac{1}{m} \right). $$
\end{enumerate}

\begin{Remth}\ 
\begin{enumerate}
\item If
\begin{eqnarray*}
\hat{u}_{P_1,1,1}^2 \sim \Normal \left(\alpha^2, \frac{\sigma^2_{\alpha^2}}{\theta_1^2 m} \right)+o_p\left( \frac{1}{\theta_1 \sqrt{m}}\right),
\end{eqnarray*}
then
\begin{eqnarray*}
\sum_{i=1}^k \hat{u}_{P_k,1,i}^2 \sim \Normal \left(\alpha^2, \frac{\sigma^2_{\alpha^2}}{\theta_1^2 m} \right)+o_p\left( \frac{1}{\theta_1 \sqrt{m}}\right),
\end{eqnarray*}
where the parameter can be computed as in Theorem \ref{jointdistribution} in \cite{mainarticle}.

\item Assuming that $c=m/n$ and that $W$ is a Wishart random matrix of dimension $m$ with $n$ degree of freedom, $\alpha^2=\frac{1-\frac{c}{(\theta_1-1)^2}}{1+\frac{c}{\theta_1-1}}$ and $\sigma^2_{\alpha^2}=2 c^2(c+1)+o_{\theta}(1)$. \\
In particular if $\frac{\theta_1}{\sqrt{m}}$ is large, then  $\alpha^2 \approx 1-c/\theta_1$,
\item In the general case, if $\frac{\theta_1}{\sqrt{m}}$ is large,$$\alpha^2 \approx 1+\frac{1-M_{2,X}}{\theta_1}  \text{ and } \sigma^2_{\alpha^2} \approx 2\left( 4 M_{2,X}^3- M_{2,X}^2-4 M_{2,X} M_{3,X}+ M_{4,X} \right).$$
\end{enumerate}
\end{Remth}
\end{Th}
(Proof in appendix \ref{appendixproof}.)

\subsection{Asymptotic distribution of the dot product}
In this section, we compute the distribution of a dot product used in this paper to prove Theorem \ref{ThInvariantdouble} and in a future work to compute the distributions of the residual spikes defined in \cite{mainarticle} for perturbation of order $k$. 

\begin{Th}
\label{ThDotproduct} 
Suppose that $W$ satisfies Assumption \ref{Ass=matrice} and $P_2=\I_m+\sum_{i=1}^2 (\theta_i-1) e_i e_i^t$ with $\theta_1>\theta_2$.
We define
\begin{eqnarray*}
&&\hat{\Sigma}_{P_2}= P_2^{1/2} W P_2^{1/2} \text{ and } \hat{\Sigma}_{P_1}= P_1^{1/2} W P_1^{1/2}.
\end{eqnarray*}
Moreover, for $s,k=1,2$ and $s \leqslant k$, we define
\begin{eqnarray*}
\hat{u}_{P_k,s},\hat{\theta}_{P_k,s} &\text{ s.t. }& \hat{\Sigma}_{P_k} \hat{u}_{P_k,s}= \hat{\theta}_{P_k,s}\hat{u}_{P_k,s},  
\end{eqnarray*}
where $\hat{\theta}_{P_k,s}=\hat{\lambda}_{\hat{\Sigma}_{P_k},s}$. Finally the present theorem uses the convention:
\begin{eqnarray*}
\text{For $s=1,2,...,k$ and $i=1,2,...,s$, } \hat{u}_{P_s,i,i}>0.
\end{eqnarray*}
\begin{enumerate}
\item Assuming that the conditions \ref{Ass=theta} (A2) and (A3) $(\theta_i=p_i\theta \rightarrow \infty)$ hold, we have \\
\scalebox{0.73}{
\begin{minipage}{1\textwidth}
\begin{eqnarray*}
\sum_{s=3}^m  \hat{u}_{P_2,1,s}  \hat{u}_{P_2,2,s} &=&\hat{u}_{P_2,1,2} \left(  \frac{ 1}{\theta_1} - \frac{1}{\theta_2}\right) -
 \frac{1}{\theta_2^{1/2}} \sum_{j>1 }^m \hat{\lambda}_{P_{1},j} \hat{u}_{P_{1},j,1} \hat{u}_{P_{1},j,2} \\
 &&\hspace{1cm}+ O_p\left(\frac{1}{\theta_1^{1/2}\theta_2^{1/2} m} \right)+ O_p\left(\frac{1}{\theta_1^{1/2} \theta_2^{3/2} m^{1/2}} \right)\\
  &=&\frac{-\left(1+M_2 \right) W_{1,2} +\left(W^2\right)_{1,2}}{\sqrt{\theta_1 \theta_2}}+ O_p\left(\frac{1}{\theta_1^{1/2}\theta_2^{1/2} m} \right)+ O_p\left(\frac{1}{\theta_1^{1/2} \theta_2^{3/2} m^{1/2}} \right).
\end{eqnarray*}
\end{minipage}}\\

Thus,  we can estimate the distribution conditional on the spectrum of $W$, \\
\scalebox{0.75}{
\begin{minipage}{1\textwidth}
\begin{eqnarray*}
 \sum_{s=3}^m  \hat{u}_{P_2,1,s}  \hat{u}_{P_k,2,s} &\sim & \Normal \left(0, \frac{\left(1+M_2 \right)^2(M_2-1)+\left(M_4-(M_2)^2\right)-2 \left(1+M_2\right)\left(M_3-M_2\right)}{\theta_1 \theta_2 m} \right)\\
&&\hspace{1cm} + O_p\left(\frac{1}{\theta_1^{1/2}\theta_2^{1/2} m} \right)+ O_p\left(\frac{1}{\theta_1^{1/2} \theta_2^{3/2} m^{1/2}} \right). 
\end{eqnarray*}
\end{minipage}} \vspace{0.2cm}

\item If $\theta_2$ is finite, then 
\begin{eqnarray*}
\sum_{s=3}^m  \hat{u}_{P_2,1,s}  \hat{u}_{P_2,2,s} &=& O_p\left( \frac{1}{\sqrt{\theta_1 m}}\right).
\end{eqnarray*}
\end{enumerate}
\begin{Rem}\ 
\begin{enumerate}
\item 
We can easily show \\
\scalebox{0.9}{
\begin{minipage}{1\textwidth}
\begin{eqnarray*}
&&\hspace{-1cm} \hat{u}_{P_2,1,2}  \left(  \frac{ 1}{\theta_1} - \frac{1}{\theta_2}\right)\delta +\sum_{s=3}^m  \hat{u}_{P_2,1,s}  \hat{u}_{P_2,2,s}\\
 &&= \frac{-\left(\delta +M_2 \right) W_{1,2} +\left(W^2\right)_{1,2}}{\sqrt{\theta_1 \theta_2}}+O_p\left(\frac{1}{\theta m} \right)+O_p\left(\frac{1}{\theta^2 m^{1/2}} \right)\\
&&\sim  \Normal \left(0, \frac{\left(\delta+M_2 \right)^2(M_2-1)+\left(M_4-(M_2)^2\right)-2 \left(\delta +M_2 \right)\left(M_3-M_2\right)}{\theta_1 \theta_2 m} \right)\\
&&\hspace{1cm} +O_p\left(\frac{1}{\theta m} \right)+O_p\left(\frac{1}{\theta^2 m^{1/2}} \right).
\end{eqnarray*}
\end{minipage}} \vspace{0.2cm}

\item If $W$ is a standard Wishart random matrix and Assumptions \ref{Ass=theta} (A2) and (A3) is verified, then
$$ \sum_{s=3}^m  \hat{u}_{P_2,1,s}  \hat{u}_{P_2,2,s} \sim \Normal \left(0, \frac{(1-\alpha_1^2)(1-\alpha_2^2)}{m} \right)+ o_p\left(\frac{1}{\theta \sqrt{m}} \right), $$
\end{enumerate}
where $\alpha_s^2=\underset{m \rightarrow \infty}{\lim} \sum_{i=1}^2\left\langle \hat{u}_{P_2,s},u_i \right\rangle^2$

\end{Rem}
\end{Th}
(Proof in appendix \ref{appendixproof}.)

\subsection{Invariant Dot Product Theorem}

\begin{Th}\label{ThInvariantdot} 
Suppose that $W$ satisfies Assumption \ref{Ass=matrice} and
\begin{eqnarray*}
&&P_{s,r}=\I_m+\sum_{i=s,r}^2 (\theta_i-1) e_i e_i^t\\
&&P_k=\I_m+\sum_{i=1}^k (\theta_i-1) e_i e_i^t \text{ respects \ref{Ass=theta} (A4)},
\end{eqnarray*} 
where $\theta_1>\theta_2>...>\theta_k$.
We define
\begin{eqnarray*}
&&\hat{\Sigma}_{P_{s,r}}=P_{s,r}^{1/2} W P_{s,r}^{1/2},\\
&&\hat{\Sigma}_{P_k}= P_k^{1/2} W P_k^{1/2}.
\end{eqnarray*}
Moreover, for $s,r=1,2,...,k$ with $s\neq r$, we define
\begin{eqnarray*}
\hat{u}_{P_{s,r},1}, \hat{\theta}_{P_{s,r},1} &\text{ s.t. }& \hat{\Sigma}_{P_{s,r}} \hat{u}_{P_{s,r},1} =\hat{\theta}_{P_{s,r},1} \hat{u}_{P_{s,r},1},  \\
\hat{u}_{P_k,s},\hat{\theta}_{P_k,s} &\text{ s.t. }& \hat{\Sigma}_{P_k} \hat{u}_{P_k,s}= \hat{\theta}_{P_k,s}\hat{u}_{P_k,s},  
\end{eqnarray*}
where $\hat{\theta}_{P_{s,r},1}= \hat{\lambda}_{\hat{\Sigma}_{P_{s,r},1}}$ and
 $\hat{\theta}_{P_k,s}=\hat{\lambda}_{\hat{\Sigma}_{P_k},s}$. \\
Assuming the convention
\begin{eqnarray*}
\text{For $s=1,2,...,k$ and $i=1,2,...,s$, } \hat{u}_{P_s,i,i}>0\,,
\end{eqnarray*}
leads to  
\begin{eqnarray*}
\boxed{\sum_{\underset{i\neq s,r}{i=1}}^m \hat{u}_{P_{s,r},1,i} \hat{u}_{P_{s,r},2,i}=\sum_{i=k+1}^m \hat{u}_{P_k,s,i} \hat{u}_{P_k,r,i} +O_p\left( \frac{1}{ \sqrt{\theta_s \theta_r} m} \right).}
\end{eqnarray*}
\end{Th}
(Proof in appendix \ref{appendixproof}.)

\subsection{Component distribution Theorem}

\begin{Th}\label{Thcomponentdistribution}
Suppose Assumption \ref{Ass=matrice} holds with canonical $P$ and \ref{Ass=theta} (A4). We define:
\begin{eqnarray*}
U&=&
\begin{pmatrix}
\hat{u}_{P_k,1}^t\\
\hat{u}_{P_k,2}^t\\
\vdots \\
\hat{u}_{P_k,m}^t
\end{pmatrix}=\begin{pmatrix}
\hat{u}_{P_k,1:k,1:k} & \hat{u}_{P_k,1:k,k+1:m} \\ 
\hat{u}_{P_k,k+1:m,1:k} & \hat{u}_{P_k,k+1:m,k+1:m}.
\end{pmatrix}
\end{eqnarray*}
To simplify the result we use the sign convention, 
\begin{eqnarray*}
\text{For $s=1,2,...,k$ and $i=1,2,...,s$, } \hat{u}_{P_s,i,i}>0.
\end{eqnarray*}
\begin{enumerate}
\item Without loss of generality on the $k$ first components, the $k^{\text{th}}$ element of the first eigenvector is \\
\scalebox{0.85}{
\begin{minipage}{1\textwidth}
\begin{eqnarray*}
\hat{u}_{P_k,1,k}&=& \frac{\sqrt{\theta_k}\theta_1}{|\theta_k-\theta_1|} \hat{u}_{P_{k-1},1,k} +O_p\left(\frac{\min(\theta_1,\theta_k)}{\theta_1^{1/2}\theta_k^{1/2}m} \right)+O_p\left(\frac{1}{\sqrt{\theta_1 \theta_k m}} \right)\\
&=&\frac{\theta_1 \sqrt{\theta_k}}{|\theta_k-\theta_1|} \frac{1}{m} \sqrt{1-\hat{\alpha}_1^2} \ Z +O_p\left(\frac{\min(\theta_1,\theta_k)}{\theta_1^{1/2}\theta_k^{1/2}m} \right)+O_p\left(\frac{1}{\sqrt{\theta_1 \theta_k m}} \right)\\
&=&\frac{\sqrt{\theta_1\theta_k}}{|\theta_k-\theta_1|} \frac{1}{\sqrt{m}} \sqrt{M_2-1} \ Z +O_p\left(\frac{\min(\theta_1,\theta_k)}{\theta_1^{1/2}\theta_k^{1/2}m} \right)+O_p\left(\frac{1}{ \sqrt{\theta_1 \theta_k m}} \right), 
\end{eqnarray*}
\end{minipage}}\\

where $Z$ is a standard normal and $M_2=\frac{1}{m}\sum_{i=1}^m \hat{\lambda}_{W,i}^2$ is obtained by conditioning on the spectrum. 
\begin{itemize}
\item Thus, knowing the spectrum and  assuming $\theta_1,\theta_k \rightarrow \infty$,
\begin{eqnarray*}
\hat{u}_{P_k,1,k} \overset{Asy}{\sim} \Normal\left(0,\frac{\theta_1 \theta_k}{|\theta_1-\theta_k|} \frac{M_2-1}{m} \right).
\end{eqnarray*}
\item If $\theta_k$ is finite, 
\begin{eqnarray*}
\hat{u}_{P_k,1,k} =O_p\left( \frac{1}{\sqrt{\theta_1 m}} \right).
\end{eqnarray*}
\end{itemize}
\noindent This result holds for any components $\hat{u}_{P_k,s,t}$ where $s \neq t \in \lbrace1,2,...,k \rbrace $.

\begin{Rem}\ \\
The sign of $\hat{u}_{P_{k},1,k}$ obtained by the construction using Theorem \ref{Theoremcaraceigenstructure} is always positive. By convention $(\hat{u}_{P_{k},i,i}>0,$ for $i=1,2,...,k)$, we multiply by $\text{sign}\left(\hat{u}_{P_k,1,1} \right)$ obtained in the construction. Thus,  the remark of Theorem \ref{Theoremcaraceigenstructure} describes the sign of the component assuming the convention. 
\begin{equation*}
\resizebox{.9\hsize}{!}{$ P\left\lbrace\text{sign}\left( \hat{u}_{P_{k},1,k} \right) = 
\text{sign}\left( \left( \hat{\theta}_{P_k,1} -\hat{\theta}_{P_{k-1},1}\right)\hat{u}_{P_{k-1},1,k} \hat{u}_{P_{k-1},1,1} \right)
\right\rbrace  = 1 +O \left(\frac{1}{m} \right).$}
\end{equation*}
\end{Rem}

\item For $s=1,...,k$, the vector $\frac{\hat{u}_{s,k+1:m}}{\sqrt{1-\hat{\alpha}_s^2}}$, where $\hat{\alpha}_s^2=\sum_{i=1}^k \hat{u}_{i,s}^2$, is unit invariant by rotation. Moreover, for $j>k$,
\begin{eqnarray*}
\hat{u}_{j,s} \sim \Normal\left(0,\frac{1-\alpha_s^2}{m }\right),
\end{eqnarray*}
where $\alpha_s^2$ is the limit of $\hat{\alpha}_s^2$.\\
Finally, the columns of $U^t[k+1:m,k+1:m]$ are invariant by rotation.

\item 
Assuming $P_k=\I_m + \sum_{i=1}^k (\theta_i-1) \epsilon_i \epsilon_i^t$ is such that 
\begin{eqnarray*}
&&\theta_1,\theta_2,...,\theta_{k_1} \text{ are proportional, and}\\
&&\theta_{k_1+1},\theta_{k_1+2},...,\theta_{k} \text{ are proportional},
\end{eqnarray*}
then\\
\scalebox{0.9}{
\begin{minipage}{1\textwidth}
\begin{eqnarray*}
\sum \hat{u}_{k+1:m,1}^2 &<& \sum \hat{u}_{k+1:m,1:k_{1}}^2\\
&\sim & {\rm{RV}}\left( O\left( \frac{1}{\theta_1 } \right), O \left( \frac{1}{\theta_1^2 m}\right)\right)+ O_p\left( \frac{\min(\theta_1,\theta_k)}{\max(\theta_1,\theta_k) m}\right).
\end{eqnarray*}
\end{minipage}}\\

If $P$ satisfies Assumption \ref{Ass=theta}(A4) with $ \min\left(\frac{\theta_1}{\theta_k},\frac{\theta_k}{\theta_1} \right) \rightarrow 0$, then
\begin{eqnarray*}
\sum \hat{u}_{k+1:m,1}^2 
&\sim & {\rm{RV}}\left( O\left( \frac{1}{\theta_1 } \right), O \left( \frac{1}{\theta_1^2 m}\right)\right)+ O_p\left( \frac{1}{\theta_1 m}\right).
\end{eqnarray*}

\end{enumerate}
\end{Th}
(Proof in appendix \ref{appendixproof}.)

\subsection{Invariant Double Angle Theorem}
Finally, using the previous Theorem, we can prove the Invariant Theorem of the double angle.

\begin{Co}
\label{ThInvariantdouble}  
Suppose $W_X$ and $W_Y$ satisfies Assumption \ref{Ass=matrice} and 
\begin{eqnarray*}
&&\tilde{P}_s=\I_m+(\theta_s-1) e_s e_s^t, \text{ for } s=1,2,...,k,\\
&&P_k=\I_m+\sum_{i=1}^k (\theta_i-1) e_i e_i^t \text{ respects \ref{Ass=theta} (A4)},
\end{eqnarray*} 
where $\theta_1>\theta_2>...>\theta_k$.
We define
\begin{eqnarray*}
&&\hat{\Sigma}_{X,\tilde{P}_s}=\tilde{P}_s^{1/2} W_X \tilde{P}_s^{1/2} \text{ and } \hat{\Sigma}_{X,\tilde{P}_s}= \tilde{P}_s^{1/2} W_Y \tilde{P}_s^{1/2},\\
&&\hat{\Sigma}_{X,P_k}=P_k^{1/2} W_X P_k^{1/2} \text{ and } \hat{\Sigma}_{Y,P_k}= P_k^{1/2} W_Y P_k^{1/2}.
\end{eqnarray*}
For $s=1,...,k$, we define
\begin{eqnarray*}
\hat{u}_{\hat{\Sigma}_{X,\tilde{P}_s},1}, \hat{\theta}_{\hat{\Sigma}_{X,\tilde{P}_s},1} &\text{ s.t. }& \hat{\Sigma}_{X,\tilde{P}_s} \hat{u}_{\hat{\Sigma}_{X,\tilde{P}_s},1} =\hat{\theta}_{\hat{\Sigma}_{X,\tilde{P}_s},1} \hat{u}_{\hat{\Sigma}_{X,\tilde{P}_s},1},\\
\hat{u}_{\hat{\Sigma}_{X,P_k},s}, \hat{\theta}_{\hat{\Sigma}_{X,P_k},s} &\text{ s.t. }& \hat{\Sigma}_{X,P_k} \hat{u}_{\hat{\Sigma}_{X,P_k},s} =\hat{\theta}_{\hat{\Sigma}_{X,P_k},s} \hat{u}_{\hat{\Sigma}_{X,P_k},s}, 
\end{eqnarray*}
where $\hat{\theta}_{\hat{\Sigma}_{X,\tilde{P}_s},1}= \hat{\lambda}_{\hat{\Sigma}_{X,\tilde{P}_s,1}}$ and
 $\hat{\theta}_{\hat{\Sigma}_{X,P_k},s}=\hat{\lambda}_{\hat{\Sigma}_{X,P_k},s}$. The statistics of the group $Y$ are defined in analogous manner.\\
\noindent Then,

\begin{align*}
\Aboxed{\left\langle \hat{u}_{\hat{\Sigma}_{X,\tilde{P}_s},1},\hat{u}_{\hat{\Sigma}_{Y,\tilde{P}_s},1} \right\rangle^2 &= \ \sum_{i=1}^{k} \left\langle \hat{u}_{\hat{\Sigma}_{X,P_k},s},\hat{u}_{\hat{\Sigma}_{Y,P_k},s} \right\rangle^2 +O_p\left( \frac{1}{\theta_s m} \right)}\\
&= \ \sum_{i=1}^{k+\epsilon} \left\langle \hat{u}_{\hat{\Sigma}_{X,P_k},s},\hat{u}_{\hat{\Sigma}_{Y,P_k},i} \right\rangle^2 +O_p\left( \frac{1}{\theta_s m} \right),
\end{align*}
where $\epsilon$ is a small integer.
\begin{Rem}\
\begin{enumerate}
\item The procedure of the proof shows an interesting invariant:\\
Assuming the sign convention $\hat{u}_{P_s,i,i}>0$ for $s=1,2,...,k$ and $i=1,2,...,s$,
\begin{eqnarray*}
 \sum_{i=k+1}^m \hat{u}_{P_k,1,i}\hat{\hat{u}}_{P_k,1,i} =  \sum_{i=k}^m \hat{u}_{P_{k-1},1,i}\hat{\hat{u}}_{P_{k-1},1,i}+ O_p\left( \frac{1}{\theta_1 m}\right).
\end{eqnarray*}
\item The distribution of $\left\langle \hat{u}_{\hat{\Sigma}_{X,P_1},1},\hat{u}_{\hat{\Sigma}_{Y,P_1},1} \right\rangle^2$ is computed in Theorem \ref{jointdistribution}.
\item An error of $\epsilon$ principal components does not affect the asymptotic distribution of the general double angle. This property allows us to construct a robust test. 
\end{enumerate}
\end{Rem}
\end{Co}
(Proof in appendix \ref{appendixproof}.)

\section{Tools for the proofs}
In this section we present intermediary results necessary to prove the main theorems of this paper. 
\subsection{Characterization of the eigenstructure}
The next theorem concerns eigenvalues and eigenvectors. In order to show the result for $u_1$, without loss of generality we use the following condition for the other eigenvalues.
\begin{Not}\  \label{Not=thetadifferent}
Usually we assume $\theta_1>\theta_2>...>\theta_k$ such that $\hat{\theta}_{P_k,s}$, the $s^{\rm th}$ largest eigenvalue of $\hat{\Sigma}_{P_k}$ corresponds to $\theta_s$.\\
We can relax the strict ordering $\theta_1>\theta2>...>\theta_k$ in the following manner. The order of $\theta_s$ in the eigenvalues $\theta_1,\theta_2,...,\theta_t$, $t\geqslant s$ is $\text{rank}_t(\theta_s)=r_{t,s}$. Assuming a perturbation $P_t$, $\theta_s$ corresponds to the $r_{t,s}^{\text{th}}$ largest eigenvalue of $\hat{\Sigma}_{P_t}$. In order to use simple notation, we again call this corresponding estimated eigenvalue, $\hat{\theta}_{P_t,s}$. \\
We also change the notation for the eigenvector. For $i=1,2,...,t$, $\hat{u}_{P_t,s}$ is the eigenvector corresponding to $\hat{\theta}_{P_t,s}$.
\end{Not}

\begin{Th} \label{Theoremcaraceigenstructure}
Using the same notation as in the Invariant Theorem (\ref{Invariantth}, \ref{Thinvarianteigenvalue}) and under Assumption \ref{Ass=matrice} and \ref{Ass=theta}(A4), we can compute the eigenvalues and the components  of interest of the eigenvector of $\hat{\Sigma}_{\P_k}$. Using assumption \ref{Ass=matrice}, we can without loss of generality suppose the canonical form for the perturbation $P_k$.
\begin{itemize}
\item Eigenvalues :\\
\scalebox{0.75}{
\begin{minipage}{1\textwidth}
\begin{eqnarray*}
\underbrace{\sum_{i=k}^m \frac{\hat{\lambda}_{P_{k-1},i}}{\hat{\theta}_{P_k,s}-\hat{\lambda}_{P_{k-1},i}}  \hat{u}_{P_{k-1},i,k}^2}_{(a) O_p\left( \frac{1}{\theta_s} \right)} + \underbrace{ \frac{\hat{\theta}_{P_{k-1},s}}{\hat{\theta}_{P_k,s}-\hat{\theta}_{P_{k-1},s}} \hat{u}_{P_{k-1},s,k}^2}_{(b) \overset{\scalebox{0.5}{order}}{\sim} \left( \frac{\theta_k-\theta_s}{ \theta_s \theta_k} \right)}+ \underbrace{\sum_{\underset{i\neq s}{i=1}}^{k-1} \frac{\hat{\theta}_{P_{k-1},i}}{\hat{\theta}_{P_k,s}-\hat{\theta}_{P_{k-1},i}} \hat{u}_{P_{k-1},i,k}^2}_{(c) O_p\left(\frac{1}{\theta_s m} \right)} = \frac{1}{\theta_k-1},
\end{eqnarray*}
\end{minipage}}\\
for $s=1,2,...,k$.
\begin{Rem}\
If we do not assume canonical perturbations, then the formula is longer but the structure remains essentially the same. Assuming Condition \ref{Ass=matrice} to hold, leads to matrices that are invariant under rotations. Elementary linear algebra methods extend the result to any perturbation.
\end{Rem}

\item Eigenvectors:\\
We define $\tilde{u}_{P_k,i}$ such that $ W P_k \tilde{u}_{P_k,i} = \hat{\theta}_{P_k,i} \tilde{u}_{P_k,i}$ and  $\hat{u}_{P_k,i}$ such that \linebreak $ P_k^{1/2} W P_k^{1/2} \hat{u}_{P_k,i} = \hat{\theta}_{P_k,i} \hat{u}_{P_k,i}$. To simplify notation we assume that $\theta_i$ corresponds to $\hat{\theta}_{P_k,i}$. This notation is explained in \ref{Not=thetadifferent} and allows us without loss of generality to describe only the eigenvector $\hat{u}_{P_k,1}$.\\ 
\scalebox{0.58}{
\begin{minipage}{1\textwidth}
\begin{eqnarray*} 
&& \hspace{-0.5cm} \left\langle \tilde{u}_{P_k,1},e_1 \right\rangle^2 \\
  &&  =\frac{\left( 
\overbrace{ \sum_{i=k}^m \frac{\hat{\lambda}_{P_{k-1},i}}{\hat{\theta}_{P_k,1}-\hat{\lambda}_{P_{k-1},i}} \hat{u}_{P_{k-1},i,1} \hat{u}_{P_{k-1},i,k}}^{(a) O_p\left( \frac{1}{\theta_1^{3/2} \sqrt{m}}\right)} + 
\overbrace{ \frac{\hat{\theta}_{P_{k-1},1}}{\hat{\theta}_{P_{k},1}-\hat{\theta}_{P_{k-1},1}}  \hat{u}_{P_{k-1},1,1} \hat{u}_{P_{k-1},1,k}}^{(b) \overset{\scalebox{0.5}{order}}{\sim} \ \frac{\sqrt{\theta_1 m}}{\min\left( \theta_1, \theta_k \right)} } +  
\overbrace{ \sum_{i=2}^{k-1} \frac{\hat{\theta}_{P_{k-1},i}}{\hat{\theta}_{P_k,1}-\hat{\theta}_{P_{k-1},i}} \hat{u}_{P_{k-1},i,1} \hat{u}_{P_{k-1},i,k} }^{(c) O_p\left( \frac{ 1 }{\theta_1^{1/2} m}\right)} \right)^2}
{\underbrace{\sum_{i=k}^m \frac{\hat{\lambda}_{P_{k-1},i}^2}{(\hat{\theta}_{P_k,1}-\hat{\lambda}_{P_{k-1},i})^2} \hat{u}_{P_{k-1},i,k}^2}_{(d) O_p\left( \frac{1}{\theta_1^2}\right) } + 
\underbrace{\frac{\hat{\theta}_{P_{k-1},1}^2}{(\hat{\theta}_{P_{k},1}-\hat{\theta}_{P_{k-1},1})^2}   \hat{u}_{P_{k-1},1,k}^2 }_{(e)\overset{\scalebox{0.5}{order}}{\sim} \ \frac{\theta_1 m}{ \min \left( \theta_1,\theta_k \right)^2}  }+  
\underbrace{\sum_{i=2}^{k-1} \frac{\hat{\theta}_{P_{k-1},i}^2}{(\hat{\theta}_{P_k,1}-\hat{\theta}_{P_{k-1},i})^2}  \hat{u}_{P_{k-1},i,k}^2}_{(f) O_p\left( \frac{1}{\theta_1 m}\right) }}, \hspace{20cm}\\
&& \hspace{-0.5cm} \left\langle \tilde{u}_{P_k,1},e_k \right\rangle^2 = \frac{1}{D_1 (\theta_k-1)^2} (g),\hspace{20cm} \\
 && \hspace{-0.5cm} \left\langle \tilde{u}_{P_k,1},e_s \right\rangle^2\hspace{-1.5cm}  \hspace{20cm} \\
&&\hspace{0.5cm}  = \frac{1}{D_1}  \left( 
\overbrace{\sum_{i=k}^m \frac{\hat{\lambda}_{P_{k-1},i}}{\hat{\theta}_{P_k,1}-\hat{\lambda}_{P_{k-1},i}} \hat{u}_{P_{k-1},i,s} \hat{u}_{P_{k-1},i,k}}^{(h) O_p\left( \frac{1}{\theta_s^{1/2} \theta_1 \sqrt{m}}\right)} 
+ 
\overbrace{\frac{\hat{\theta}_{P_{k-1},1}}{\hat{\theta}_{P_{k},1}-\hat{\theta}_{P_{k-1},1}}  \hat{u}_{P_{k-1},1,s} \hat{u}_{P_{k-1},1,k}}^{(i) \ \overset{\scalebox{0.5}{order}}{\sim} \frac{\min\left( \theta_1 , \theta_s \right)}{\sqrt{\theta_s}\min\left( \theta_1 , \theta_k \right)}}\right. \\
&&\hspace{3cm}\left.+  
\overbrace{\sum_{i=2,\neq s}^{k-1} \frac{\hat{\theta}_{P_{k-1},i}}{\hat{\theta}_{P_k,1}-\hat{\theta}_{P_{k-1},i}} \hat{u}_{P_{k-1},i,s} \hat{u}_{P_{k-1},i,k} }^{(j) O_p\left( \underset{i \neq 1,s}{\max}\left( \frac{\min\left( \theta_1,\theta_i \right)\min\left( \theta_s, \theta_i \right)}{\sqrt{\theta_s} \theta_1 \theta_i \sqrt{m}} \right)\right)}
+  
\overbrace{ \frac{\hat{\theta}_{P_{k-1},s}}{\hat{\theta}_{P_k,1}-\hat{\theta}_{P_{k-1},s}} \hat{u}_{P_{k-1},s,s} \hat{u}_{P_{k-1},s,k} }^{(k) O_p\left( \frac{\min\left( \theta_1 , \theta_s \right)}{\sqrt{\theta_s} \theta_1 \sqrt{m}} \right)} 
\right)^2 .
\end{eqnarray*} 
\end{minipage}}\\
\noindent Finally, 
$$\hat{u}_{P_k,1}= \frac{\left( \tilde{u}_{P_k,1,1},\tilde{u}_{P_k,1,2},...,\sqrt{\theta_k} \tilde{u}_{P_k,1,k},...,\tilde{u}_{P_k,} \right)}{\underbrace{ \sqrt{1+ \left(\theta_k-1 \right) \tilde{u}_{P_k1,k}^2}}_{1+O_p\left(\frac{\min(\theta_1,\theta_k)}{\max(\theta_1,\theta_k)m} \right)}},$$
where $\sqrt{1+ \left(\theta-1 \right) \tilde{u}_{P_k1,k}^2}$ is the norm of $P_k^{1/2}\tilde{u}_{P_k,1}$ that we will call $N_1$.

\begin{Rem}\ 
\begin{enumerate}
\item By construction, the sign of $\hat{u}_{P_k,1,k}$ is always positive. This is, however, not the case of $\hat{u}_{P_{k-1},i,i}$.  We can show that:\\
\scalebox{0.85}{
\begin{minipage}{1\textwidth}
\begin{eqnarray*}
P\left\lbrace {\rm sign} \left( \hat{u}_{P_k,1,1} \right) ={\rm sign} \left( \left(\hat{\theta}_{P_k,1}-\hat{\theta}_{P_{k-1},1}\right) \hat{u}_{P_{k-1},1,1} \hat{u}_{P_{k-1},1,k}\right) \right\rbrace \underset{m \rightarrow \infty}{\rightarrow} 1.
\end{eqnarray*}
\end{minipage}}\\

Moreover, the convergence to $1$ is of order $1/m$. If $\theta_1$ tends to infinity, then 
\begin{eqnarray*}
P\left\lbrace {\rm sign} \left( \hat{u}_{P_k,1,1} \right) ={\rm sign} \left( \left(\theta_1-\theta_k\right)  \hat{u}_{P_{k-1},1,k}\right) \right\rbrace \underset{m,\theta_1 \rightarrow \infty}{\rightarrow} 1.
\end{eqnarray*}
Thus, if we use a convention such as  $\text{sign}\left(\hat{u}_{P_{k},i,i}\right) >0 $ for $i=1,...,k-1$, 
then the sign of $\hat{u}_{P_k,1,k}$ is distributed as a Bernoulli with parameter 1/2.

\item Without loss of generality, the other eigenvectors $\hat{u}_{P_k,r}$ for $r=1,2,...,k-1$ can be computed by the same formula thanks to the notation linking the estimated eigenvector to the eigenvalue $theta_i$.\\
This formula does, however, not work for the vector $\hat{u}_{P_k,k}$. Applying a different order of perturbation shows that similar formulas exist for $\hat{u}_{P_k,k}$. (If the perturbation in $e_1$ is applied at the end for example.) \\
This observation leads to a problem in the proofs of the Dot Product Theorems \ref{ThDotproduct} and \ref{ThInvariantdot}. Deeper investigations are necessary to understand the two eigenvectors when $k=2$.
\begin{eqnarray*}
D_2&=&\underbrace{\sum_{i=2}^m \frac{\hat{\lambda}_{P_{1},i}^2}{(\hat{\theta}_{P_2,2}-\hat{\lambda}_{P_{1},i})^2} \hat{u}_{P_{1},i,2}^2}_{ O_p\left( \frac{1}{\theta_2^2}\right) } + 
\underbrace{\frac{\hat{\theta}_{P_{1},1}^2}{(\hat{\theta}_{P_{2},2}-\hat{\theta}_{P_{1},1})^2}   \hat{u}_{P_{1},1,2}^2 }_{O_p\left( \frac{\theta_1}{(\theta_2-\theta_1)^2 m}\right) },\\
N_2^2&=&1+\frac{1}{(\theta_2-1)D_2},\\
N_2 D_2 &=& D_2 + \frac{1}{\theta_2-1} \\
&=& \frac{1}{\theta_2-1}+O_p\left( \frac{1}{\theta_2^2} \right)+O_p\left( \frac{\theta_1}{(\theta_2-\theta_1) m} \right).
\end{eqnarray*}
Furthermore, the theorem requires investigation of the $m-k$ noisy components of the eigenvectors. For $r=1,2$ and $s=3,4,...,m$,
\begin{eqnarray*}
\hat{u}_{P_2,r,s}= \frac{\sum_{i=1}^m \frac{\hat{\lambda}_{P_{1},i}}{\hat{\theta}_{P_2,r}-\hat{\lambda}_{P_{1},i}} \hat{u}_{P_{1},i,s} \hat{u}_{P_{1},i,2}}{\sqrt{D_r} N_r}.
\end{eqnarray*}
The estimations using this last formula are difficult. It is beneficial to look at 
$$\hat{u}_{P_2,1,t}/\sqrt{\sum_{s=3}^m \hat{u}_{P_2,1,s}^2} \text{ and } \hat{u}_{P_2,2,t}/\sqrt{\sum_{s=3}^m \hat{u}_{P_2,2,s}^2}$$
 for $t=3,4,...,m$.
\item If the perturbation is not canonical, then we can apply a rotation $U$, such that $U u_s = \epsilon_s$, and replace $\hat{u}_{P_{k-1},i}$ by $U^t \hat{u}_{P_{k-1},i}$. Then, $\left\langle \tilde{u}_{P_k,1},e_s \right\rangle^2$ is replaced by $\left\langle \tilde{u}_{P_k,1},u_s \right\rangle^2$.
\end{enumerate}
\end{Rem}
\end{itemize}

\end{Th}
(Proof in appendix \ref{appendixproof}.)

\subsection{Double dot product}
\begin{Th}
 \label{Thdoubledot} 
Suppose $W_X$ and $W_Y$ satisfies Assumption \ref{Ass=matrice} and 
$P_k=\I_m+\sum_{i=1}^k (\theta_i-1) e_i e_i^t$ satisfies \ref{Ass=theta} (A4),
 where $\theta_1>\theta_2>...>\theta_k$.
We set
\begin{eqnarray*}
&&\hat{\Sigma}_{X}=\hat{\Sigma}_{X,P_k}=P_k^{1/2} W_X P_k^{1/2} \text{ and } \hat{\Sigma}_{Y,P_k}= P_k^{1/2} W_Y P_k^{1/2}.
\end{eqnarray*}
and for $s=1,...,k$,
\begin{eqnarray*}
\hat{u}_{\hat{\Sigma}_{X},s}, \hat{\theta}_{\hat{\Sigma}_{X},s} &\text{ s.t. }& \hat{\Sigma}_{X} \hat{u}_{\hat{\Sigma}_{X},s} =\hat{\theta}_{\hat{\Sigma}_{X},s} \hat{u}_{\hat{\Sigma}_{X},s},\\
\hat{u}_{\hat{\Sigma}_{Y},s}, \hat{\theta}_{\hat{\Sigma}_{Y},s} &\text{ s.t. }& \hat{\Sigma}_{Y} \hat{u}_{\hat{\Sigma}_{Y},s} =\hat{\theta}_{\hat{\Sigma}_{Y},s} \hat{u}_{\hat{\Sigma}_{Y},s}, 
\end{eqnarray*}
where $\hat{\theta}_{\hat{\Sigma}_{Y},s}=\hat{\lambda}_{\hat{\Sigma}_{Y},s}$ and
 $\hat{\theta}_{\hat{\Sigma}_{X},s}=\hat{\lambda}_{\hat{\Sigma}_{X},s}$.  To simplify the result we assume the sign convention: 
\begin{eqnarray*}
\text{For $s=1,2,...,k$ and $i=1,2,...,s$, } \hat{u}_{\hat{\Sigma}_X,i,i}>0, \ \hat{u}_{\hat{\Sigma}_Y,i,i}>0.
\end{eqnarray*}
Finally, we define
\begin{eqnarray*}
\tilde{u}_s = \hat{U}_X^t \hat{\hat{u}}_{\hat{\Sigma}_{Y},s},
\end{eqnarray*}
where,
\begin{eqnarray*}
\hat{U}_X=\left(
v_1,
v_2, 
\cdots,
v_m
\right)=\left(
\hat{u}_{\hat{\Sigma}_{X},1},
\hat{u}_{\hat{\Sigma}_{X},2}, 
\cdots 
\hat{u}_{\hat{\Sigma}_{X},k}, 
v_{k+1},
v_{k+2},
\cdots,
v_m
\right),
\end{eqnarray*}
where the vectors $v_{k+1},...,v_m$ are chosen such that the matrix $\hat{U}_X$ is orthonormal. Then, 
\begin{itemize}
\item If $\theta_j,\theta_t \rightarrow \infty$: \\
\scalebox{0.77}{
\begin{minipage}{1\textwidth}
\begin{eqnarray*}
\sum_{i=k+1}^m \tilde{u}_{j,i} \tilde{u}_{t,i}&=& \sum_{i=k+1}^m \hat{u}_{\hat{\Sigma}_{Y},j,i} \hat{u}_{\hat{\Sigma}_{Y},t,i} +\sum_{i=k+1}^m \hat{u}_{\hat{\Sigma}_{X},j,i} \hat{u}_{\hat{\Sigma}_{X},t,i} -\sum_{i=k+1}^m \hat{u}_{\hat{\Sigma}_{X},j,i} \hat{u}_{\hat{\Sigma}_{Y},t,i} \\
&& \hspace{2cm} - \sum_{i=k+1}^m \hat{u}_{\hat{\Sigma}_{Y},j,i} \hat{u}_{\hat{\Sigma}_{X},t,i}- \left( \hat{u}_{\hat{\Sigma}_{X},t,j}+\hat{u}_{\hat{\Sigma}_{Y},j,t} \right) \left( \hat{\alpha}^2_{\hat{\Sigma}_{X},j} - \hat{\alpha}^2_{\hat{\Sigma}_{X},t} \right)\\
&& \hspace{2cm}+ O_p\left( \frac{1}{\theta_1 m} \right)+O_p\left( \frac{1}{\theta_1^2  \sqrt{m}} \right),
\end{eqnarray*}
\end{minipage}}\\

where $\hat{\alpha}^2_{\hat{\Sigma}_{X},t}= \sum_{i=1}^k \hat{u}_{\hat{\Sigma}_{X},t,i}^2$. 
\item If $\theta_t$ is finite: 
\begin{eqnarray*}
\sum_{i=k+1}^m \tilde{u}_{j,i} \tilde{u}_{t,i}=O_p\left( \frac{1}{\sqrt{m}\sqrt{\theta_1}} \right).
\end{eqnarray*}
\end{itemize}

\noindent Moreover, for $s=1,...,k$, $t=2,...,k$ and $j=k+1,...m$,
\begin{eqnarray*}
&&\sum_{i=1}^k \tilde{u}_{s,i}^2 = \sum_{i=1}^k \left\langle \hat{u}_{\hat{\Sigma}_{X},i},\hat{u}_{\hat{\Sigma}_{Y},s} \right\rangle^2,\\
&& \tilde{u}_{s,s}=\hat{u}_{\hat{\Sigma}_{X},s,s}\hat{u}_{\hat{\Sigma}_{Y},s,s}+O_p\left( \frac{1}{m} \right) +O_p\left( \frac{1}{\theta_s^{1/2}m^{1/2}} \right),\\
 &&\tilde{u}_{s,t}=\hat{u}_{\hat{\Sigma}_{X},t,s}+\hat{u}_{\hat{\Sigma}_{X},s,t}+O_p\left( \frac{\sqrt{\min(\theta_s,\theta_t)}}{m\sqrt{\max(\theta_s,\theta_t)}} \right)+O_p\left( \frac{1}{ \theta_t m^{1/2}} \right),\\
  &&\tilde{u}_{t,s}=O_p\left( \frac{\sqrt{\min(\theta_s,\theta_t)}}{m\sqrt{\max(\theta_s,\theta_t)}} \right)+O_p\left( \frac{1}{ \theta_s m^{1/2}} \right),\\
&& \tilde{u}_{s,j}=\hat{u}_{\hat{\Sigma}_{Y},s,j}-\hat{u}_{\hat{\Sigma}_{X},s,j} \left\langle \hat{u}_{\hat{\Sigma}_{Y},j},\hat{u}_{\hat{\Sigma}_{X},j} \right\rangle+ O_p\left( \frac{1}{\theta_s^{1/2} m } \right).
\end{eqnarray*}
\end{Th}
(Proof in appendix \ref{appendixproof}.)

\subsection{Lemmas for Invariant Dot product Theorem}
This section introduces a lemma used in the proof of the Dot Product Theorem \ref{ThDotproduct}. 

\begin{Lem}\label{LemstatW}
Assuming $W$ and $\hat{\Sigma}_{P_1}$ as in Theorem \ref{ThDotproduct}, then by construction of the eigenvectors using Theorem \ref{Theoremcaraceigenstructure},\\
\scalebox{0.68}{
\begin{minipage}{1\textwidth}
\begin{eqnarray*}
&\hat{u}_{P_1,1,2} &=\frac{W_{1,2}}{\sqrt{\theta_1} W_{1,1}} -\frac{W_{1,2}}{\theta_1^{3/2}}\left( -1/2+3/2 M_2 \right)+ \frac{\left(W^2\right)_{1,2}}{\theta_1^{3/2}} + O_p\left( \frac{1}{\theta_1^{3/2}m} \right)+ O_p\left( \frac{1}{\theta_1^{5/2}m^{1/2}} \right)\\
&&= \frac{W_{1,2}}{\sqrt{\theta_1}}+O_p\left( \frac{1}{\theta_1^{1/2} m} \right)+O_p\left( \frac{1}{\theta_1^{3/2} m^{1/2}} \right),\\
&\sum_{i=2}^m \hat{\lambda}_{P_1,i}^2 \hat{u}_{P_1,i,2}^2 &=W_{2,2}+O_p\left( \frac{1}{m}\right),\\
&\sum_{i=2}^m \hat{\lambda}_{P_1,i} \hat{u}_{P_1,i,1}\hat{u}_{P_1,i,2}&=W_{1,2} \frac{  M_2}{\sqrt{\theta_1}} -\left(W^2\right)_{1,2}\frac{1}{\sqrt{\theta_1}}
+ O_p\left( \frac{1}{\theta_1^{1/2}m} \right)+ O_p\left( \frac{1}{\theta_1^{3/2}m^{1/2}} \right).
\end{eqnarray*}
\end{minipage}}

\begin{Rem}\ \\
Because the perturbation is of order $1$, the two sign conventions defined in \ref{Not=Theorem} are the same.
\end{Rem}

\end{Lem}
(Proof in appendix \ref{appendixproof}.)

\section{Conclusion}
In this paper we extend results of \cite{mainarticle} to perturbation of order $k>1$. Theorem \ref{jointdistribution} provides all the background results needed to build powerful test. The approach contains two deficiencies:
\begin{itemize}
\item We cannot treat the case with equal perturbing eigenvalues, $\theta_1=\theta_2$. Indeed, all our theorems always assume different eigenvalues. In the case of equality, the procedures do not stay conservative.
\item The distribution of the data before the perturbation is applied are assumed to be invariant under rotation. If we relax this assumption, then our procedure are no longer necessarily conservative.
\end{itemize} 

In future work we will present a procedure based on the residual spikes introduced in \cite{mainarticle} for perturbations of order $1$. These statistics seems to capture the differences between two populations very effectively and the problem of equal eigenvalues of the perturbation does not affect these tests. Relaxing the hypotheses of invariance under rotation still influences the properties of these alternative tests, but have a lesser impact.

\pagebreak

\appendix
\section{Statistical applications of Random matrix theory:\\ comparison of two populations II,\\
Supplement} \label{appendixproof}

\subsection{Introduction} 

This appendix contains the supplemental material presenting the proofs of the theorems and lemmas of the paper. These results are first introduced with the same notation as in the main paper and directly proved. Because some assumptions are used in the proofs, we also introduce the notation, some definitions and some assumptions. 

\subsection{Notations, Definitions, Assumptions and Previous Theorems}
As presented in \cite{mainarticle} we use the following notation.

\begin{Not}\ \label{AANot=Theorem}\\
Although we use a precise notation to enunciate the theorems, the proofs rely on a simpler notation when no confusion is possible. This difference is always specified at the beginning of a proof. 
\begin{itemize}
\item If $W$ is a symmetric random matrix, we denote by $\left(\hat{\lambda}_{W,i}, \hat{u}_{W,i} \right)$ its $i^{\rm th}$ eigenvalue and eigenvector.
\item A finite perturbation of order $k$ is denoted by $P_k= \I_m  +\sum_{i=1}^k (\theta_i-1) u_i u_i^t \in \mathbb{R}^{m \times m}$ with $u_1,u_2,...,u_k \in \mathbb{R}^{m \times m}$ orthonormal vectors.
\item We denote by $W \in  \mathbb{R}^{m\times m}$ random matrix that is invariant under rotation as defined in Assumption \ref{AAAss=matrice}. Moreover, the estimated covariance matrix is $\hat{\Sigma}=P_k^{1/2} W P_k^{1/2}$. \\
When comparing two groups, we use $W_X$, $W_Y$ and $\hat{\Sigma}_X$, $\hat{\Sigma}_Y$.
\item When we consider only one group, $\hat{\Sigma}_{P_r}=P_r^{1/2} W P_r^{1/2}$ is a perturbation of order $r\leq k$ of the matrix $W$ and
\begin{itemize} 
\item $\hat{u}_{P_r,i}$ is its $i^{\text{th}}$ eigenvector. When $r=k$ we use the simpler notation $\hat{u}_{i}=\hat{u}_{P_k,i}$.
\item  $\hat{u}_{P_r,i,j}$ is the $j^{\rm th}$ component of the $i^{\rm th}$ eigenvector.
\item $\hat{\lambda}_{P_r,i}$ is its $i^{\text{th}}$ eigenvalue. If $\theta_1>\theta_2>...>\theta_r$, then for $i=1,2,...,r$ we use also the notation $\hat{\theta}_{P_r,i}=\hat{\lambda}_{P_r,i}$. We call these eigenvalues the spikes. When $r=k$, we use the simpler notation $\hat{\theta}_i=\hat{\theta}_{P_k,i}$.
\item $\hat{\alpha}_{P_r,i}^2=\sum_{j=1}^r \left\langle \hat{u}_{P_r,i},u_j \right\rangle^2$ is called the \textbf{general angle}.
\end{itemize}  
With this notation, we have $\hat{\Sigma}=\hat{\Sigma}_{P_k}=P_k^{1/2} W P_k^{1/2}$.
\item When we look at two groups $X$ and $Y$, we use a notation similar to the above. The perturbation of order $r$ of the matrices $W_X$ and $W_Y$ are $\hat{\Sigma}_{X,P_r}=P_r^{1/2} W_X P_r^{1/2}$ and $\hat{\Sigma}_{Y,P_r}=P_r^{1/2} W_Y P_r^{1/2}$, respectively. Then, we define for the group $\hat{\Sigma}_{X,P_r}$ (and similarly for $\hat{\Sigma}_{Y,P_r}$):
\begin{itemize}
\item $\hat{u}_{\hat{\Sigma}_{X,P_r},i}$ is its $i^{\text{th}}$ eigenvector. When $r=k$ we use the simpler notation $\hat{u}_{X,i}=\hat{u}_{\hat{\Sigma}_{X,P_k},i}$.
\item  $\hat{u}_{\hat{\Sigma}_{X,P_r},i,j}$ is the $j^{\rm th}$ component of the $i^{\rm th}$ eigenvector.
\item $\hat{\lambda}_{\hat{\Sigma}_{X,P_r},i}$ is its $i^{\text{th}}$ eigenvalue. If $\theta_1>\theta_2>...>\theta_r$, then for $i=1,2,...,r$ we use the notation $\hat{\theta}_{\hat{\Sigma}_{X,P_r},i}=\hat{\lambda}_{\hat{\Sigma}_{X,P_r},i}$. When $r=k$, we use the simpler notation $\hat{\theta}_{X,i}=\hat{\theta}_{\hat{\Sigma}_{X,P_k},i}$.
\item $\hat{\alpha}_{\hat{\Sigma}_{X,P_r},i}^2=\sum_{j=1}^r \left\langle \hat{u}_{\hat{\Sigma}_{X,P_r},i},u_j \right\rangle^2$. 
\item $\hat{\alpha}_{X,Y,P_r,i}^2=\sum_{j=1}^r \left\langle \hat{u}_{\hat{\Sigma}_{X,P_r},i},\hat{u}_{\hat{\Sigma}_{Y,P_r},j} \right\rangle^2$ is the \textbf{double angle} and, when no confusion is possible, we use the simpler notation $\hat{\alpha}_{P_r,i}^2$. When this simpler notation is used, it is stated explicitly.
\end{itemize} 
\item The theorems can assume a sign convention 
\begin{eqnarray*}
\hat{u}_{P_s,i,i}>0, \text{ for $s=1,2,...,k$ and $i=1,2,...,s$,}
\end{eqnarray*}
 as in Theorem \ref{AAThInvariantdot} or \ref{AAThcomponentdistribution}. On the other hand, some theorems assume the convention 
\begin{eqnarray*}
\hat{u}_{P_s,i,s}>0, \text{ for $s=1,2,...,k$ and $i=1,2,...,s$,}
\end{eqnarray*}
as in Theorem \ref{AATheoremcaraceigenstructure}.\\
Other theorems are not affected by this convention and do not specify it. Nevertheless, the convention will be given in the proofs when confusion is possible.
\item We define the function $M_{s_1,s_2,X}(\rho_X)$, $M_{s_1,s_2,Y}(\rho_Y)$ and $M_{s_1,s_2}(\rho_X,\rho_Y)$ as 
\begin{eqnarray*}
M_{s_1,s_2,X}(\rho_X)&=&\frac{1}{m} \sum_{i=1}^m \frac{\hat{\lambda}_{W_X,i}^{s_1}}{\left( \rho_X -\hat{\lambda}_{W_X,i} \right)^2},\\
M_{s_1,s_2,Y}(\rho_Y)&=&\frac{1}{m} \sum_{i=1}^m \frac{\hat{\lambda}_{W_Y,i}^{s_1}}{\left( \rho_Y -\hat{\lambda}_{W_Y,i} \right)^2},\\
M_{s_1,s_2}(\rho_X,\rho_Y)&=& \frac{M_{s_1,s_2,X}(\rho_X)+M_{s_1,s_2,Y}(\rho_Y)}{2}.
\end{eqnarray*}
In particular, when $s_2=0$, we use $M_{s_1,X}=M_{s_1,0,X}$. When we only study one group, we use the simpler notation $M_{s_1,s_2}(\rho)$ when no confusion is possible.
\item We use two transforms inspired by the T-transform:
\begin{itemize}
\item $T_{W,u}(z)= \sum_{i=1}^m \frac{\hat{\lambda}_{W,i}}{z-\hat{\lambda}_{W,i}} \left\langle \hat{u}_{W,i},u \right\rangle^2$ is the T-transform in direction $u$ using the random matrix $W$.
\item $\hat{T}_{\hat{\Sigma}_X}(z)= \frac{1}{m}\sum_{i=k+1}^m \frac{\hat{\lambda}_{\hat{\Sigma}_X,i}}{z-\hat{\lambda}_{\hat{\Sigma}_X,i}}$,
 and  $\hat{T}_{W_X}(z)= \frac{1}{m} \sum_{i=1}^m \frac{\hat{\lambda}_{W_X,i}}{z-\hat{\lambda}_{W_X,i}}$ are the estimated T-transforms using $\hat{\Sigma}_X$ and $W$ respectively.
\end{itemize}
\item In some theorems we use the notation $\overset{\scalebox{0.5}{order}}{\sim}$ to describe the order size in probability of a positive random variable. For example, $X_m \overset{\scalebox{0.5}{order}}{\sim} 1/m$ if
$\frac{X_m}{1/m}$ tends to a random variable $X$ independent of $m$, with $P\left\lbrace |X|>\epsilon_j \right\rbrace \overset{j\rightarrow \infty}{\longrightarrow} 1$ for any sequences $\epsilon_j$ tending to $0$.
\end{itemize}
\end{Not}

We recall the assumptions of the main paper.
\begin{Ass}\label{AAAss=matrice} 
Let $W_X$ and $W_Y$ be such that
\begin{eqnarray*}
W_X=O_X \Lambda_X O_X \text{ and } W_Y=O_Y \Lambda_Y O_Y,
\end{eqnarray*}
where 
\begin{eqnarray*}
&&O_X, O_Y  \text{ are unit orthonormal invariant and independent random matrices,}\\
&&\Lambda_X,\Lambda_Y \text{ are diagonal bounded matrices and independent of } O_X, O_Y,\\
&&\Tr\left( W_X \right)=1 \text{ and } \Tr\left( W_Y \right)=1.
\end{eqnarray*}

\noindent Assume $P_X= \I_m + \sum_{i=1}^k(\theta_{X,i}-1) e_i e_i^t$ and $P_Y= \I_m+ \sum_{i=1}^k(\theta_{Y,i}-1) e_i e_i^t$. Then 
\begin{eqnarray*}
\hat{\Sigma}_X=P_X^{1/2} W_X P_X^{1/2} \text{ and } \hat{\Sigma}_Y=P_Y^{1/2} W_Y P_Y^{1/2}.
\end{eqnarray*}
\end{Ass}

\begin{Ass}\label{AAAss=theta}
\begin{itemize} \
\item[(A1)] $\frac{\theta}{\sqrt{m} }\rightarrow \infty.$
\item[(A2)] $\theta \rightarrow \infty.$
\item[(A3)] $\theta_i=p_i \theta$, where $p_i$ is fixed and different from $1$.
\item[(A4)] For  $i=1,...,k_\infty, \ \theta_i=p_i \theta$, $\theta \rightarrow \infty$  according to  (A1) or (A2),\\
For $i=k_\infty+1,...,k, \ \theta_i=p_i \theta_0$.\\
For all $i \neq j$, $p_i \not= p_j$.
\end{itemize}
\end{Ass}

We recall the definitions.
\begin{Def} \label{AADef:detectable} 
\begin{enumerate}
\item We assume that a perturbation $P=\I_m+ (\theta-1) u u^t$ is \textbf{detectable} in $\hat{\Sigma}=P^{1/2} W P^{1/2}$ if the perturbation creates a largest isolated eigenvalue, $\hat{\theta}$.
\item We say that a finite perturbation of order $k$ is \textbf{detectable} if it creates $k$ largest eigenvalues well separated from the spectrum of $W$.
\end{enumerate}
\end{Def}
\begin{Def} \label{AADef:bloc} 
The perturbation $P_k=\I_m+ \sum_{i=1}^k (\theta_i-1) u_i u_i^t$ is in two blocs if,
\begin{itemize}
\item For  $i=1,...,k_\infty, \ \theta_i=p_i \theta$, $\theta \rightarrow \infty$ for fixed $p_1>p_2>...>p_{k_\infty}$.
\item For $i=k_\infty+1,...,k, \ \theta_i=p_i \theta_0$ for fixed $p_{k_\infty+1}>p_{k_\infty+2}>...>p_{k}$.
\end{itemize}
\end{Def}
\begin{Def} \label{AADef=unbiased} 
Suppose $\hat{\Sigma}$ satisfies Assumption \ref{AAAss=matrice}.\\
The \textbf{unbiased estimator of} $\theta$ is defined as 
$$ \hat{\hat{\theta}}=1+\frac{1}{\frac{1}{m-k} \sum_{i=k+1}^{m} \frac{\hat{\lambda}_{\hat{\Sigma},i}}{\hat{\theta}-\hat{\lambda}_{\hat{\Sigma},i}}},$$
where $\hat{\lambda}_{\hat{\Sigma},i}$ is the $i^{\text{th}}$ of $\hat{\Sigma}$.\\
Suppose that $\hat{\theta}$ and $\hat{u}_i$ are the $i^{\text{th}}$ eigenvalue and eigenvector of $\hat{\Sigma}$, the \textbf{filtered estimated covariance matrix} is defined as 
$$\hat{\hat{\Sigma}}= \I_m+\sum_{i=1}^k (\hat{\hat{\theta}}_i-1) \hat{u}_i \hat{u}_i^t.$$
\end{Def}
\begin{Def} \label{AADef:Invariant}
Let $W$ be a random matrix. Moreover, let $P_1=\I_m+(\theta_1-1) u_1 u_1^t$ and $P_k=\I_m+\sum_{i=1}^k (\theta_i-1)u_i u_i^t$ be perturbations of order $1$ and $k$, respectively. We say that a statistic $T\left( W_m,P_1\right)$ is \textbf{invariant} with respect to $k$, if $T\left( W_m,P_k\right)$ is such that 
\begin{eqnarray*}
\resizebox{.9\hsize}{!}{$T\left( W_m,P_k\right) = T\left( W_m,P_1\right) + \epsilon_m, \text{ where }
\max\left(\frac{\epsilon_m}{\E\left[ T\left( W,P_1\right)\right]},\frac{\epsilon_m^2}{\var\left( T\left( W,P_1\right) \right)}\right) \rightarrow 0.$}
\end{eqnarray*}
\end{Def}
We recall the main Theorems of \cite{mainarticle} in a lighter form.
\begin{Th} \label{AAjointdistribution}
Suppose $W_X$ and $W_Y$ satisfy \ref{AAAss=matrice} with $P=P_X=P_Y$, a detectable perturbation of order $k=1$. Moreover, we assume as known, \linebreak $S_{W_X}=\left\lbrace \hat{\lambda}_{W_X,1},\hat{\lambda}_{W_X,2},...,\hat{\lambda}_{W_X,m} \right\rbrace$ and $S_{W_Y}=\left\lbrace \hat{\lambda}_{W_Y,1},\hat{\lambda}_{W_Y,2},...,\hat{\lambda}_{W_Y,m} \right\rbrace$, the eigenvalues of $W_X$ and $W_Y$. We defined
\begin{eqnarray*}
\hat{\Sigma}_X&=&P^{1/2} W_X P^{1/2},\\
\hat{\Sigma}_Y&=&P^{1/2} W_Y P^{1/2},\\
P&=&\I_m+(\theta-1) u u^t,
\end{eqnarray*}
where $u$ is fixed. We construct the unbiased estimators of $\theta$,
\begin{eqnarray*}
\hat{\hat{\theta}}_X \ \bigg| \ \frac{1}{\hat{\hat{\theta}}_X-1}= \frac{1}{m}\sum_{i=1}^m \frac{\hat{\lambda}_{W_X,i}}{\hat{\theta}_X-\hat{\lambda}_{W_X,i}} &\text{ and }& \hat{\hat{\theta}}_Y \ \bigg| \ \frac{1}{\hat{\hat{\theta}}_Y-1}= \frac{1}{m}\sum_{i=1}^m \frac{\hat{\lambda}_{W_Y,i}}{\hat{\theta}_Y-\hat{\lambda}_{W_Y,i}}
\end{eqnarray*}
where $\hat{\theta}_X=\hat{\lambda}_{\hat{\Sigma}_X,1}$ and $\hat{\theta}_Y=\hat{\lambda}_{\hat{\Sigma}_Y,1}$ are the largest eigenvalues of $\hat{\Sigma}_X$ and $\hat{\Sigma}_Y$ with corresponding eigenvectors $\hat{u}_X=\hat{u}_{\hat{\Sigma}_X,1}$ and $\hat{u}_Y=\hat{u}_{\hat{\Sigma}_Y,1}$. \\
Using this notation and assuming a convergence rate of $\left(\hat{\theta}_X,\hat{\theta}_Y\right)$ to $\left(\rho_X,\rho_Y\right)$ in $O_p\left(\theta/\sqrt{m}\right)$ with $\E\left[ \hat{\theta}_X \right]=\rho_X+o\left( \frac{\theta}{\sqrt{m}} \right)$ and $\E\left[ \hat{\theta}_Y \right]=\rho_X+o\left( \frac{\theta}{\sqrt{m}} \right)$, we have
\begin{equation*}
\resizebox{.95\hsize}{!}{$ \left. \begin{pmatrix}
\hat{\hat{\theta}}_X  \\ 
\hat{\hat{\theta}}_Y  \\ 
\left\langle \hat{u}_X,\hat{u}_Y \right\rangle^2
\end{pmatrix} \right| S_{W_X}, S_{W_Y} \sim \Normal \left(
\begin{pmatrix}
\theta \\
\theta  \\ 
\alpha_{X,Y}^2
\end{pmatrix}
,\frac{1}{m}\begin{pmatrix}
\sigma_{\theta,X}^2 & 0 & \sigma_{\theta,\alpha^2,X} \\ 
0 & \sigma_{\theta,Y}^2  &  \sigma_{\theta,\alpha^2,Y} \\
\sigma_{\theta,\alpha^2,X} & \sigma_{\theta,\alpha^2,Y} & \sigma_{\alpha^2,X,Y}^2
\end{pmatrix}
 \right)+\begin{pmatrix}
o_p\left(\frac{\theta}{\sqrt{m}}\right) \\ 
o_p\left(\frac{\theta}{\sqrt{m}}\right) \\ 
o_p\left( \frac{1}{\theta \sqrt{m}}\right)
\end{pmatrix}.$}
\end{equation*}

\noindent Here, all the parameters depend on 
\begin{eqnarray*}
 \resizebox{.95\hsize}{!}{$ M_{s,r,X}(\rho_X)= \frac{1}{m} \sum_{i=1}^m \frac{\hat{\lambda}_{W_X,i}^s}{(\rho_X-\hat{\lambda}_{W_X,i})^r} \text{ and } M_{s,r,Y}(\rho_Y)=\frac{1}{m} \sum_{i=1}^m \frac{\hat{\lambda}_{W_X,i}^s}{(\rho_Y-\hat{\lambda}_{W_X,i})^r} .$}
\end{eqnarray*}
\end{Th}

\begin{Th} 
\label{AAconvergence}
In this theorem, $P=\I_m+ (\theta-1) u u^t$ is a finite perturbation of order $1$. Suppose $W$ is a symmetric matrix with eigenvalues $\hat{\lambda}_{W,i}\geqslant 0$ and eigenvectors $\hat{u}_{W,i}$ for $i=1,2,...,m$. The perturbation of $W$ by $P$ leads to $\hat{\Sigma}=P^{1/2}WP^{1/2}$. \\
For $i=1,2,...,m$, we define $\tilde{u}_{\hat{\Sigma},i}$ and $\hat{\lambda}_{\hat{\Sigma},i}$ such that 
$$ W P \tilde{u}_{\hat{\Sigma},i} = \hat{\lambda}_{\hat{\Sigma},i} \tilde{u}_{\hat{\Sigma},i},$$ 
and the usual $\hat{u}_{\hat{\Sigma},i}$ such that if $\hat{\Sigma}=P^{1/2}WP^{1/2}$, then
$$\hat{\Sigma} \hat{u}_{\hat{\Sigma},i} = P^{1/2} W P^{1/2} \hat{u}_{\hat{\Sigma},i} = \hat{\lambda}_{\hat{\Sigma},i} \hat{u}_{\hat{\Sigma},i}.$$
\begin{itemize}
\item The eigenvalues $\hat{\lambda}_{\hat{\Sigma},s}$ are such that for $s=1,2,...,m$,
\begin{eqnarray*}
\sum_{i=1}^m \frac{\hat{\lambda}_{W,i}}{\hat{\lambda}_{\hat{\Sigma},s}-\hat{\lambda}_{W,i}}  \left\langle \hat{u}_{W,i},u \right\rangle^2 = \frac{1}{\theta_k-1}.
\end{eqnarray*}
\item The eigenvectors $\tilde{u}_{\hat{\Sigma},s}$ are such that 
\begin{eqnarray*} 
&&  \left\langle \tilde{u}_{\hat{\Sigma},s},v \right\rangle^2  =\frac{\left( 
 \sum_{i=1}^m \frac{\hat{\lambda}_{W,i}}{\hat{\lambda}_{\hat{\Sigma},s}-\hat{\lambda}_{W,i}} \left\langle \hat{u}_{W,i},v \right\rangle \left\langle \hat{u}_{W,i},u \right\rangle  \right)^2}{
\sum_{i=1}^m \frac{\hat{\lambda}_{W,i}^2}{(\hat{\lambda}_{\hat{\Sigma},s}-\hat{\lambda}_{W,i})^2} \left\langle \hat{u}_{W,i},u \right\rangle^2 }. 
\end{eqnarray*}
In particular if $v=u$,
\begin{eqnarray*}
 \left\langle \tilde{u}_{\hat{\Sigma},s},u \right\rangle^2 &=& \frac{1}{\left(\theta_k-1 \right)^2 \left(\sum_{i=1}^m \frac{\hat{\lambda}_{W,i}^2}{(\hat{\lambda}_{\hat{\Sigma},s}-\hat{\lambda}_{W,i})^2} \left\langle \hat{u}_{W,i},u \right\rangle ^2 \right)}.
\end{eqnarray*}
Moreover, 
\begin{eqnarray*}
\hat{u}_{\hat{\Sigma},s}= \frac{P^{1/2} \tilde{u}_{\hat{\Sigma},s}} {\sqrt{1+ \left(\theta-1 \right) \left\langle \tilde{u}_{\hat{\Sigma},s},u \right\rangle^2}}.
\end{eqnarray*}
\end{itemize}
\end{Th}

\begin{Th}
 \label{AAThunitstatistic} \ \\
Let $W$ be a random matrix with spectrum $S_W=\left\lbrace\hat{\lambda}_{W,1},\hat{\lambda}_{W,2},...,\hat{\lambda}_{W,m} \right\rbrace$ normalized to have a trace of $1$. We denote by $u_{p_1}$ and $u_{p_2}$, two orthonormal invariant random vectors of size $m$ and independent of the eigenvalues of $W$. We set 
\begin{eqnarray*}
\vec{B}_m\left(\rho,\vec{s},\vec{r},\vec{p} \right)=\sqrt{m} \left( \begin{pmatrix}
\sum_{i=1}^m \frac{\hat{\lambda}_{W,i}^{s_1}}{\left(\rho-\hat{\lambda}_{W,i}\right)^{s_2}} {u}_{p_1,i} {u}_{p_2,i} \\ 
\sum_{i=1}^m \frac{\hat{\lambda}_{W,i}^{r_1}}{\left(\rho-\hat{\lambda}_{W,i}\right)^{r_2}} {u}_{p_1,i} {u}_{p_2,i} 
\end{pmatrix} - \begin{pmatrix}
M_{s_1,s_2} \\ 
M_{r_1,r_2}
\end{pmatrix} \mathbf{1}_{p_1=p_2}\right),
\end{eqnarray*}
where $\vec{s}=\left(s_1,s_2 \right)$, $\vec{r}=\left(r_1,r_2 \right)$ and $\vec{p}=\left(p_1,p_2 \right)$ with indices $1 \leqslant p_1 \leqslant p_2 \leqslant m$ and $s_1,s_2,r_1,r_2 \in \mathbf{N}$.\\
If $p=p_1=p_2$, we have
\begin{equation*}
\resizebox{.95\hsize}{!}{$\left. \vec{B}_m\left(\rho,\vec{s},\vec{r},\vec{p} \right) \right| S_W \sim {\Normal}\left(\vec{0},  \begin{pmatrix}
2 \left(M_{2 s_1, 2 s_2}-M_{s_1,s_1}^2 \right) & 2 \left( M_{s_1+r_1,s_2+r_2}-M_{s_1,s_2} M_{r_1,r_2} \right) \\ 
2 \left( M_{s_1+r_1,s_2+r_2}-M_{s_1,s_2} M_{r_1,r_2} \right) & 2 \left(M_{2 r_1, 2 r_2}-M_{r_1,r_1}^2 \right)
\end{pmatrix} \right)+ o_{p;m}(1),$}
\end{equation*}
where $M_{s,r}=M_{s,r}(\rho)=\frac{1}{m}\sum_{i=1}^m \frac{\hat{\lambda}_{W,i}^s}{\left(\rho-\hat{\lambda}_{W,i}\right)^r} $.\\
Moreover, for $p_1 \neq p_2$,
\begin{equation*}
\resizebox{.95\hsize}{!}{$\left.\vec{B}_m\left(\rho,\vec{s},\vec{r},\vec{p} \right) \right| S_W \sim {\Normal}\left(\vec{0},  \begin{pmatrix}
 M_{2 s_1, 2 s_2}-M_{s_1,s_1}^2  &   M_{s_1+r_1,s_2+r_2}-M_{s_1,s_2} M_{r_1,r_2}  \\ 
  M_{s_1+r_1,s_2+r_2}-M_{s_1,s_2} M_{r_1,r_2}  & M_{2 r_1, 2 r_2}-M_{r_1,r_1}^2 
\end{pmatrix} \right)+ o_{p;m}(1).$}
\end{equation*}
In particular, with the notation $M_{s,0}=M_s=\frac{1}{m}\sum_{i=1}^m \hat{\lambda}_{W,i}^s $,
\begin{equation*}
\resizebox{.95\hsize}{!}{$\left. \sqrt{m} \left( \begin{pmatrix}
\sum_{i=1}^m \hat{\lambda}_{W,i} {u}_{p,i}^2 \\ 
\sum_{i=1}^m \hat{\lambda}_{W,i}^2 {u}_{p,i}^2
\end{pmatrix} - \begin{pmatrix}
1 \\ 
M_2
\end{pmatrix}\right) \right| S_W \sim \Normal\left(\vec{0},  \begin{pmatrix}
2 \left(M_2-1 \right) & 2 \left( M_3-M_2\right) \\ 
2 \left( M_3-M_2 \right) & 2 \left(M_4-M_2^2 \right)
\end{pmatrix} \right)+ o_{p;m}(1),$}
\end{equation*}
and 
\begin{equation*}
\resizebox{.95\hsize}{!}{$\left.\sqrt{m} \left( \begin{pmatrix}
\sum_{i=1}^m \hat{\lambda}_{W,i} {u}_{p_1,i} {u}_{p_2,i} \\ 
\sum_{i=1}^m \hat{\lambda}_{W,i}^2 {u}_{p_1,i} {u}_{p_2,i}
\end{pmatrix} - \begin{pmatrix}
0 \\ 
0
\end{pmatrix}\right) \right| S_W \sim \Normal \left(\vec{0},  \begin{pmatrix}
M_2-1  &   M_3-M_2 \\ 
  M_3-M_2  &  M_4-M_2^2 
\end{pmatrix} \right) + o_{p;m}(1).$}
\end{equation*} 
Finally if we look at $K$ bivariate normal random variables :
\begin{equation*}
\resizebox{.95\hsize}{!}{$\mathbf{B}_{m}\left(\vec{\rho},\mathbf{{s}},\mathbf{{r}},\mathbf{{p}} \right)=\left(\vec{B}_m\left(\rho_1,\vec{s}_1,\vec{r}_1,\vec{p}_1 \right), 
\vec{B}_m\left(\rho_2,\vec{s}_2,\vec{r}_2,\vec{p}_2 \right),
...,
\vec{B}_m\left(\rho_K,\vec{s}_K,\vec{r}_K,\vec{p}_K \right) \right),$}
\end{equation*}
where $,\vec{p}_i \neq\vec{p}_j$ if $i\neq j$. Then, conditioning on the spectrum $S_W$, \linebreak $\mathbf{\vec{B}}_{m}\left(\vec{\rho},\mathbf{{s}},\mathbf{{r}},\mathbf{{p}} \right)$ tends to a multivariate Normal. Moreover, all the bivariate elements are asymptotically independent.
\end{Th}

\subsection{Main Theorems} \label{AAsec:Theorems}
In this section we present and prove the theorems and lemmas of this paper.

\subsubsection{Invariant Eigenvalue Theorem}
\begin{Th} \label{AAThinvarianteigenvalue} Suppose that $W$ satisfies Assumption \ref{AAAss=matrice} and
\begin{eqnarray*}
&&\tilde{P}_s=\I_m+(\theta_s-1) e_s e_s^t, \text{ for } s=1,2,...,k,\\
&&P_k=\I_m+\sum_{i=1}^k (\theta_i-1) e_i e_i^t \text{ satisfies Assumption \ref{AAAss=theta} (A4),}
\end{eqnarray*} 
where $\theta_1>\theta_2>...>\theta_k$.
We define
$$\hat{\Sigma}_{\tilde{P}_s}=\tilde{P}_s^{1/2} W \tilde{P}_s^{1/2}, \text{ and }
\hat{\Sigma}_{P_k}= P_k^{1/2} W P_k^{1/2}\,.$$
Moreover, for $s=1,2,...,k$, we define
\begin{eqnarray*}
\hat{u}_{\tilde{P}_s,1}, \hat{\theta}_{\tilde{P}_s,1} &\text{ s.t. }& \hat{\Sigma}_{\tilde{P}_s} \hat{u}_{\tilde{P}_s,1} =\hat{\theta}_{\tilde{P}_s,1} \hat{u}_{\tilde{P}_s,1},  \\
\hat{u}_{P_k,s},\hat{\theta}_{P_k,s} &\text{ s.t. }& \hat{\Sigma}_{P_k} \hat{u}_{P_k,s}= \hat{\theta}_{P_k,s}\hat{u}_{P_k,s},  
\end{eqnarray*}
where $\hat{\theta}_{\tilde{P}_s,1}= \hat{\lambda}_{\hat{\Sigma}_{\tilde{P}_s,1}}$ and $\hat{\theta}_{P_k,s}=\hat{\lambda}_{\hat{\Sigma}_{P_k},s}$. The following results hold:
\begin{enumerate}
\item  For $s>1$,
$$ \boxed{\hat{\theta}_{P_{k},s} - \hat{\theta}_{\tilde{P}_{s},1} \overset{\scalebox{0.5}{order}}{\sim} \frac{\theta_s}{m}} $$
and
$$ \boxed{\hat{\theta}_{P_{k},1} - \hat{\theta}_{\tilde{P}_{1},1} \overset{\scalebox{0.5}{order}}{\sim} \frac{\theta_2}{m},},$$
where $\overset{\scalebox{0.5}{order}}{\sim}$ is the order size in probability. The distribution of $ \hat{\theta}_{P_{k},s} $ is therefore asymptotically the same as the distribution of $\hat{\theta}_{\tilde{P}_{s},1}$ studied in Theorem \ref{AAjointdistribution}.
\item More precisely, we define for $r,s \in \left\lbrace 1,2,...,k \right\rbrace$ with $r \neq s $ ,
$$ P_{-r}=\I_m + \sum_{\underset{i\neq r}{i=1}}^k \left( \theta_i-1 \right) e_i e_i^t. $$
\begin{itemize}
\item If $\theta_s> \theta_r$, then 
\begin{equation*}
 \resizebox{.8\hsize}{!}{$\hat{\theta}_{P_{k},s}-\hat{\theta}_{P_{-r},s}
= 
-\frac{\hat{\theta}_{P_{-r},s} \hat{\theta}_{P_{k},s} (\theta_r-1)}{\theta_r-1 -\hat{\theta}_{P_{k},s}} \hat{u}_{P_{-r},s,r}^2 + O_p \left(\frac{1}{m} \right)+ O_p \left(\frac{\theta_r}{m^{3/2}} \right) .$}
\end{equation*}
\item If $\theta_s < \theta_r$, then 
\begin{equation*}
\resizebox{.8\hsize}{!}{$ \hat{\theta}_{P_{k},s}-\hat{\theta}_{P_{-r},s-1}
= 
-\frac{\hat{\theta}_{P_{-r},s-1} \hat{\theta}_{P_{k},s} (\theta_r-1)}{\theta_r-1 -\hat{\theta}_{P_{k},s}} \hat{u}_{P_{-r},s-1,r}^2 + O_p \left(\frac{1}{m} \right)+ O_p \left(\frac{\theta_s}{m^{3/2}} \right).$}
\end{equation*}
\end{itemize}

\end{enumerate} 
 \end{Th}

\subsubsection{Invariant Angle Theorem}\label{AAInvariantanglesection}

\begin{Th}\label{AAInvariantth}
Using the same notation as Theorem \ref{AAThinvarianteigenvalue}, we have
\begin{enumerate}
\item The general angle is invariant in the sense of Definition \ref{AADef:Invariant},
\begin{eqnarray*}
&& \boxed{\sum_{i=1}^k \hat{u}_{P_k,s,i}^2=\hat{u}_{\tilde{P}_s,1,s}^2 + O_p\left( \frac{1}{\theta_s m} \right).}
\end{eqnarray*}
Therefore,  the distribution of $\sum_{i=1}^k \hat{u}_{P_k,s,i}^2 $ is asymptotically the same as the distribution of $\hat{u}_{\tilde{P}_s,1,s}^2$ studied in Theorem \ref{AAjointdistribution}.
\item Moreover, 
$$ \hat{u}_{P_k,s,s}^2 = \hat{u}_{\tilde{P}_s,1,s}^2 + O_p\left( \frac{1}{m} \right). $$
\end{enumerate}
\end{Th}

\subsubsection{Asymptotic distribution of the dot product}

\begin{Th}
\label{AAThDotproduct} 
Suppose that $W$ satisfies Assumption \ref{AAAss=matrice} and $P_2=\I_m+\sum_{i=1}^2 (\theta_i-1) e_i e_i^t$ with $\theta_1>\theta_2$.
Let
\begin{eqnarray*}
&&\hat{\Sigma}_{P_2}= P_2^{1/2} W P_2^{1/2} \text{ and } \hat{\Sigma}_{P_1}= P_1^{1/2} W P_1^{1/2}.
\end{eqnarray*}
Moreover, for $s,k=1,2$ and $s \leqslant k$, let
\begin{eqnarray*}
\hat{u}_{P_k,s},\hat{\theta}_{P_k,s} &\text{ s.t. }& \hat{\Sigma}_{P_k} \hat{u}_{P_k,s}= \hat{\theta}_{P_k,s}\hat{u}_{P_k,s},  
\end{eqnarray*}
where $\hat{\theta}_{P_k,s}=\hat{\lambda}_{\hat{\Sigma}_{P_k},s}$. Finally, assume that for $s=1,2,...,k$ and $i=1,2,...,s$,  $\hat{u}_{P_s,i,i}>0$. Then we have the following
\begin{enumerate}
\item If the Assumptions \ref{AAAss=theta} (A2) and (A3) $(\theta_i=p_i\theta \rightarrow \infty)$ hold, then

\scalebox{0.73}{
\begin{minipage}{1\textwidth}
\begin{eqnarray*}
\sum_{s=3}^m  \hat{u}_{P_2,1,s}  \hat{u}_{P_2,2,s} &=&\hat{u}_{P_2,1,2} \left(  \frac{ 1}{\theta_1} - \frac{1}{\theta_2}\right) -
 \frac{1}{\theta_2^{1/2}} \sum_{j>1 }^m \hat{\lambda}_{P_{1},j} \hat{u}_{P_{1},j,1} \hat{u}_{P_{1},j,2} \\
 &&\hspace{1cm}+ O_p\left(\frac{1}{\theta_1^{1/2}\theta_2^{1/2} m} \right)+ O_p\left(\frac{1}{\theta_1^{1/2} \theta_2^{3/2} m^{1/2}} \right)\\
  &=&\frac{-\left(1+M_2 \right) W_{1,2} +\left(W^2\right)_{1,2}}{\sqrt{\theta_1 \theta_2}}+ O_p\left(\frac{1}{\theta_1^{1/2}\theta_2^{1/2} m} \right)+ O_p\left(\frac{1}{\theta_1^{1/2} \theta_2^{3/2} m^{1/2}} \right).
\end{eqnarray*}
\end{minipage}}\\

Thus,  we can approximate the distribution conditional on the spectrum of $W$, \\
\scalebox{0.75}{
\begin{minipage}{1\textwidth}
\begin{eqnarray*}
 \sum_{s=3}^m  \hat{u}_{P_2,1,s}  \hat{u}_{P_k,2,s} &\sim & \Normal \left(0, \frac{\left(1+M_2 \right)^2(M_2-1)+\left(M_4-(M_2)^2\right)-2 \left(1+M_2\right)\left(M_3-M_2\right)}{\theta_1 \theta_2 m} \right)\\
&&\hspace{1cm} + O_p\left(\frac{1}{\theta_1^{1/2}\theta_2^{1/2} m} \right)+ O_p\left(\frac{1}{\theta_1^{1/2} \theta_2^{3/2} m^{1/2}} \right). 
\end{eqnarray*}
\end{minipage}} \vspace{0.2cm}

\item If $\theta_2$ is finite, then 
\begin{eqnarray*}
\sum_{s=3}^m  \hat{u}_{P_2,1,s}  \hat{u}_{P_2,2,s} &=& O_p\left( \frac{1}{\sqrt{\theta_1 m}}\right).
\end{eqnarray*}
\end{enumerate}
\begin{Rem}\ 
From the above, we can easily show that\\
\scalebox{0.9}{
\begin{minipage}{1\textwidth}
\begin{eqnarray*}
&&\hspace{-1cm} \hat{u}_{P_2,1,2}  \left(  \frac{ 1}{\theta_1} - \frac{1}{\theta_2}\right)\delta +\sum_{s=3}^m  \hat{u}_{P_2,1,s}  \hat{u}_{P_2,2,s}\\
 &&= \frac{-\left(\delta +M_2 \right) W_{1,2} +\left(W^2\right)_{1,2}}{\sqrt{\theta_1 \theta_2}}+O_p\left(\frac{1}{\theta m} \right)+O_p\left(\frac{1}{\theta^2 m^{1/2}} \right)\\
&&\sim  \Normal \left(0, \frac{\left(\delta+M_2 \right)^2(M_2-1)+\left(M_4-(M_2)^2\right)-2 \left(\delta +M_2 \right)\left(M_3-M_2\right)}{\theta_1 \theta_2 m} \right)\\
&&\hspace{1cm} +O_p\left(\frac{1}{\theta m} \right)+O_p\left(\frac{1}{\theta^2 m^{1/2}} \right).
\end{eqnarray*}
\end{minipage}} \vspace{0.2cm}
\end{Rem}
\end{Th}

\subsubsection{Invariant Dot Product Theorem}
\begin{Th}\label{AAThInvariantdot} 
Suppose that $W$ satisfies Assumption \ref{AAAss=matrice} and
\begin{eqnarray*}
&&P_{s,r}=\I_m+\sum_{i=s,r}^2 (\theta_i-1) e_i e_i^t\\
&&P_k=\I_m+\sum_{i=1}^k (\theta_i-1) e_i e_i^t \text{ verifies \ref{AAAss=theta} (A4)},
\end{eqnarray*} 
where $\theta_1>\theta_2>...>\theta_k$.
We define
\begin{eqnarray*}
&&\hat{\Sigma}_{P_{s,r}}=P_{s,r}^{1/2} W P_{s,r}^{1/2},\\
&&\hat{\Sigma}_{P_k}= P_k^{1/2} W P_k^{1/2}.
\end{eqnarray*}
Moreover, for $s,r=1,2,...,k$ with $s\neq r$, we define
\begin{eqnarray*}
\hat{u}_{P_{s,r},1}, \hat{\theta}_{P_{s,r},1} &\text{ s.t. }& \hat{\Sigma}_{P_{s,r}} \hat{u}_{P_{s,r},1} =\hat{\theta}_{P_{s,r},1} \hat{u}_{P_{s,r},1},  \\
\hat{u}_{P_k,s},\hat{\theta}_{P_k,s} &\text{ s.t. }& \hat{\Sigma}_{P_k} \hat{u}_{P_k,s}= \hat{\theta}_{P_k,s}\hat{u}_{P_k,s},  
\end{eqnarray*}
where $\hat{\theta}_{P_{s,r},1}= \hat{\lambda}_{\hat{\Sigma}_{P_{s,r},1}}$ and
 $\hat{\theta}_{P_k,s}=\hat{\lambda}_{\hat{\Sigma}_{P_k},s}$. \\
If
\begin{eqnarray*}
\text{For $s=1,2,...,k$ and $i=1,2,...,s$, } \hat{u}_{P_s,i,i}>0\,,
\end{eqnarray*}
then 
\begin{eqnarray*}
\boxed{\sum_{\underset{i\neq s,r}{i=1}}^m \hat{u}_{P_{s,r},1,i} \hat{u}_{P_{s,r},2,i}=\sum_{i=k+1}^m \hat{u}_{P_k,s,i} \hat{u}_{P_k,r,i} +O_p\left( \frac{1}{ \sqrt{\theta_s \theta_r} m} \right).}
\end{eqnarray*}
\end{Th}

\subsubsection{Component distribution Theorem}
\begin{Th}\label{AAThcomponentdistribution}
Suppose Assumption \ref{AAAss=matrice} holds with canonical $P$ and \ref{AAAss=theta} (A4). We define:
\begin{eqnarray*}
U&=&
\begin{pmatrix}
\hat{u}_{P_k,1}^t\\
\hat{u}_{P_k,2}^t\\
\vdots \\
\hat{u}_{P_k,m}^t
\end{pmatrix}=\begin{pmatrix}
\hat{u}_{P_k,1:k,1:k} & \hat{u}_{P_k,1:k,k+1:m} \\ 
\hat{u}_{P_k,k+1:m,1:k} & \hat{u}_{P_k,k+1:m,k+1:m}.
\end{pmatrix}
\end{eqnarray*}
To simplify the result we assume the sign convention, 
\begin{eqnarray*}
\text{For $s=1,2,...,k$ and $i=1,2,...,s$, } \hat{u}_{P_s,i,i}>0.
\end{eqnarray*}
\begin{enumerate}
\item Without loss of generality on the $k$ first components, the $k^{\text{th}}$ element of the first eigenvector is \\
\scalebox{0.85}{
\begin{minipage}{1\textwidth}
\begin{eqnarray*}
\hat{u}_{P_k,1,k}&=& \frac{\sqrt{\theta_k}\theta_1}{|\theta_k-\theta_1|} \hat{u}_{P_{k-1},1,k} +O_p\left(\frac{\min(\theta_1,\theta_k)}{\theta_1^{1/2}\theta_k^{1/2}m} \right)+O_p\left(\frac{1}{\sqrt{\theta_1 \theta_k m}} \right)\\
&=&\frac{\theta_1 \sqrt{\theta_k}}{|\theta_k-\theta_1|} \frac{1}{m} \sqrt{1-\hat{\alpha}_1^2} \ Z +O_p\left(\frac{\min(\theta_1,\theta_k)}{\theta_1^{1/2}\theta_k^{1/2}m} \right)+O_p\left(\frac{1}{\sqrt{\theta_1 \theta_k m}} \right)\\
&=&\frac{\sqrt{\theta_1\theta_k}}{|\theta_k-\theta_1|} \frac{1}{\sqrt{m}} \sqrt{M_2-1} \ Z +O_p\left(\frac{\min(\theta_1,\theta_k)}{\theta_1^{1/2}\theta_k^{1/2}m} \right)+O_p\left(\frac{1}{ \sqrt{\theta_1 \theta_k m}} \right), 
\end{eqnarray*}
\end{minipage}}\\

where $Z$ is a standard normal and $M_2=\frac{1}{m}\sum_{i=1}^m \hat{\lambda}_{W,i}^2$ is obtained by conditioning on the spectrum. 
\begin{itemize}
\item Thus, knowing the spectrum and  assuming $\theta_1,\theta_k \rightarrow \infty$,
\begin{eqnarray*}
\hat{u}_{P_k,1,k} \overset{Asy}{\sim} \Normal\left(0,\frac{\theta_1 \theta_k}{|\theta_1-\theta_k|} \frac{M_2-1}{m} \right).
\end{eqnarray*}
\item If $\theta_k$ is finite, 
\begin{eqnarray*}
\hat{u}_{P_k,1,k} =O_p\left( \frac{1}{\sqrt{\theta_1 m}} \right).
\end{eqnarray*}
\end{itemize}
\noindent This result holds for any components $\hat{u}_{P_k,s,t}$ where $s \neq t \in \lbrace1,2,...,k \rbrace $.

\item For $s=1,...,k$, the vector $\frac{\hat{u}_{s,k+1:m}}{\sqrt{1-\hat{\alpha}_s^2}}$, where $\hat{\alpha}_s^2=\sum_{i=1}^k \hat{u}_{i,s}^2$, is unit invariant under rotation. Moreover, for $j>k$,
\begin{eqnarray*}
\hat{u}_{j,s} \sim \Normal\left(0,\frac{1-\alpha_s^2}{m }\right),
\end{eqnarray*}
where $\alpha_s^2$ is the limit of $\hat{\alpha}_s^2$.\\
Moreover, the columns of $U^t[k+1:m,k+1:m]$ are rotation invariant.

\item  Assuming $P_k=\I_m + \sum_{i=1}^k (\theta_i-1) \epsilon_i \epsilon_i^t$ is such that 
\begin{eqnarray*}
&&\theta_1,\theta_2,...,\theta_{k_1} \text{ are proportional, and}\\
&&\theta_{k_1+1},\theta_{k_1+2},...,\theta_{k} \text{ are proportional},
\end{eqnarray*}
then\\
\scalebox{0.9}{
\begin{minipage}{1\textwidth}
\begin{eqnarray*}
\sum \hat{u}_{k+1:m,1}^2 &<& \sum \hat{u}_{k+1:m,1:k_{1}}^2\\
&\sim & {\rm{RV}}\left( O\left( \frac{1}{\theta_1 } \right), O \left( \frac{1}{\theta_1^2 m}\right)\right)+ O_p\left( \frac{\min(\theta_1,\theta_k)}{\max(\theta_1,\theta_k) m}\right).
\end{eqnarray*}
\end{minipage}}\\

Moreover, if $P$ satisfies Assumption \ref{AAAss=theta}(A4) with $ \min\left(\frac{\theta_1}{\theta_k},\frac{\theta_k}{\theta_1} \right) \rightarrow 0$, then
\begin{eqnarray*}
\sum \hat{u}_{k+1:m,1}^2 
&\sim & {\rm{RV}}\left( O\left( \frac{1}{\theta_1 } \right), O \left( \frac{1}{\theta_1^2 m}\right)\right)+ O_p\left( \frac{1}{\theta_1 m}\right).
\end{eqnarray*}

\end{enumerate}
\end{Th}

\subsubsection{Invariant Double Angle Theorem}
\begin{Co}
\label{AAThInvariantdouble}  
Suppose $W_X$ and $W_Y$ satisfies Assumption \ref{AAAss=matrice} and 
\begin{eqnarray*}
&&\tilde{P}_s=\I_m+(\theta_s-1) e_s e_s^t, \text{ for } s=1,2,...,k,\\
&&P_k=\I_m+\sum_{i=1}^k (\theta_i-1) e_i e_i^t \text{ respects \ref{AAAss=theta} (A4)},
\end{eqnarray*} 
where $\theta_1>\theta_2>...>\theta_k$.
We define
\begin{eqnarray*}
&&\hat{\Sigma}_{X,\tilde{P}_s}=\tilde{P}_s^{1/2} W_X \tilde{P}_s^{1/2} \text{ and } \hat{\Sigma}_{X,\tilde{P}_s}= \tilde{P}_s^{1/2} W_Y \tilde{P}_s^{1/2},\\
&&\hat{\Sigma}_{X,P_k}=P_k^{1/2} W_X P_k^{1/2} \text{ and } \hat{\Sigma}_{Y,P_k}= P_k^{1/2} W_Y P_k^{1/2}.
\end{eqnarray*}
Moreover, for $s=1,...,k$, we define
\begin{eqnarray*}
\hat{u}_{\hat{\Sigma}_{X,\tilde{P}_s},1}, \hat{\theta}_{\hat{\Sigma}_{X,\tilde{P}_s},1} &\text{ s.t. }& \hat{\Sigma}_{X,\tilde{P}_s} \hat{u}_{\hat{\Sigma}_{X,\tilde{P}_s},1} =\hat{\theta}_{\hat{\Sigma}_{X,\tilde{P}_s},1} \hat{u}_{\hat{\Sigma}_{X,\tilde{P}_s},1},\\
\hat{u}_{\hat{\Sigma}_{X,P_k},s}, \hat{\theta}_{\hat{\Sigma}_{X,P_k},s} &\text{ s.t. }& \hat{\Sigma}_{X,P_k} \hat{u}_{\hat{\Sigma}_{X,P_k},s} =\hat{\theta}_{\hat{\Sigma}_{X,P_k},s} \hat{u}_{\hat{\Sigma}_{X,P_k},s}, 
\end{eqnarray*}
where $\hat{\theta}_{\hat{\Sigma}_{X,\tilde{P}_s},1}= \hat{\lambda}_{\hat{\Sigma}_{X,\tilde{P}_s,1}}$ and
 $\hat{\theta}_{\hat{\Sigma}_{X,P_k},s}=\hat{\lambda}_{\hat{\Sigma}_{X,P_k},s}$. The statistics of the group $Y$ are defined in analogous manner.\\
\noindent Then,

\begin{align*}
\Aboxed{\left\langle \hat{u}_{\hat{\Sigma}_{X,\tilde{P}_s},1},\hat{u}_{\hat{\Sigma}_{Y,\tilde{P}_s},1} \right\rangle^2 &= \ \sum_{i=1}^{k} \left\langle \hat{u}_{\hat{\Sigma}_{X,P_k},s},\hat{u}_{\hat{\Sigma}_{Y,P_k},s} \right\rangle^2 +O_p\left( \frac{1}{\theta_s m} \right)}\\
&= \ \sum_{i=1}^{k+\epsilon} \left\langle \hat{u}_{\hat{\Sigma}_{X,P_k},s},\hat{u}_{\hat{\Sigma}_{Y,P_k},i} \right\rangle^2 +O_p\left( \frac{1}{\theta_s m} \right),
\end{align*}
where $\epsilon$ is a small integer.
\end{Co}

\subsection{Tools for the proofs}

\subsubsection{Characterization of the eigenstructure}
\begin{Th} \label{AATheoremcaraceigenstructure}
Using the same notation as in the Invariant Theorem (\ref{AAInvariantth}, \ref{AAThinvarianteigenvalue}) and under Assumption \ref{AAAss=matrice} and \ref{AAAss=theta}(A4), we can compute the eigenvalues and the components  of interest of the eigenvector of $\hat{\Sigma}_{\P_k}$. Using these conditions, we can without loss of generality suppose the canonical form for the perturbation $P_k$.
\begin{itemize}
\item Eigenvalues :\\
\scalebox{0.75}{
\begin{minipage}{1\textwidth}
\begin{eqnarray*}
\underbrace{\sum_{i=k}^m \frac{\hat{\lambda}_{P_{k-1},i}}{\hat{\theta}_{P_k,s}-\hat{\lambda}_{P_{k-1},i}}  \hat{u}_{P_{k-1},i,k}^2}_{(a) O_p\left( \frac{1}{\theta_s} \right)} + \underbrace{ \frac{\hat{\theta}_{P_{k-1},s}}{\hat{\theta}_{P_k,s}-\hat{\theta}_{P_{k-1},s}} \hat{u}_{P_{k-1},s,k}^2}_{(b) \overset{\scalebox{0.5}{order}}{\sim} \left( \frac{\theta_k-\theta_s}{ \theta_s \theta_k} \right)}+ \underbrace{\sum_{\underset{i\neq s}{i=1}}^{k-1} \frac{\hat{\theta}_{P_{k-1},i}}{\hat{\theta}_{P_k,s}-\hat{\theta}_{P_{k-1},i}} \hat{u}_{P_{k-1},i,k}^2}_{(c) O_p\left(\frac{1}{\theta_s m} \right)} = \frac{1}{\theta_k-1},
\end{eqnarray*}
\end{minipage}}\\
for $s=1,2,...,k$.
\begin{Rem}\
Without the canonical form for the perturbations, the formula is longer but the structure remains essentially the same. Elementary linear algebra methods extend the result from rotationally invariant matrices to arbitrary perturbations.
\end{Rem}

\item Eigenvectors :\\
We define $\tilde{u}_{P_k,i}$ such that $ W P_k \tilde{u}_{P_k,i} = \hat{\theta}_{P_k,i} \tilde{u}_{P_k,i}$ and  $\hat{u}_{P_k,i}$ such that \linebreak $ P_k^{1/2} W P_k^{1/2} \hat{u}_{P_k,i} = \hat{\theta}_{P_k,i} \hat{u}_{P_k,i}$. To simplify notation we assume that $\theta_i$ corresponds to $\hat{\theta}_{P_k,i}$. This notation is explained in \ref{AANot=thetadifferent} and allows for a more efficient description of the first $k$ eigenvectors.\\ 
\scalebox{0.58}{
\begin{minipage}{1\textwidth}
\begin{eqnarray*} 
&& \hspace{-0.5cm} \left\langle \tilde{u}_{P_k,1},e_1 \right\rangle^2 \\
  &&  =\frac{\left( 
\overbrace{ \sum_{i=k}^m \frac{\hat{\lambda}_{P_{k-1},i}}{\hat{\theta}_{P_k,1}-\hat{\lambda}_{P_{k-1},i}} \hat{u}_{P_{k-1},i,1} \hat{u}_{P_{k-1},i,k}}^{(a) O_p\left( \frac{1}{\theta_1^{3/2} \sqrt{m}}\right)} + 
\overbrace{ \frac{\hat{\theta}_{P_{k-1},1}}{\hat{\theta}_{P_{k},1}-\hat{\theta}_{P_{k-1},1}}  \hat{u}_{P_{k-1},1,1} \hat{u}_{P_{k-1},1,k}}^{(b) \overset{\scalebox{0.5}{order}}{\sim} \ \frac{\sqrt{\theta_1 m}}{\min\left( \theta_1, \theta_k \right)} } +  
\overbrace{ \sum_{i=2}^{k-1} \frac{\hat{\theta}_{P_{k-1},i}}{\hat{\theta}_{P_k,1}-\hat{\theta}_{P_{k-1},i}} \hat{u}_{P_{k-1},i,1} \hat{u}_{P_{k-1},i,k} }^{(c) O_p\left( \frac{ 1 }{\theta_1^{1/2} m}\right)} \right)^2}
{\underbrace{\sum_{i=k}^m \frac{\hat{\lambda}_{P_{k-1},i}^2}{(\hat{\theta}_{P_k,1}-\hat{\lambda}_{P_{k-1},i})^2} \hat{u}_{P_{k-1},i,k}^2}_{(d) O_p\left( \frac{1}{\theta_1^2}\right) } + 
\underbrace{\frac{\hat{\theta}_{P_{k-1},1}^2}{(\hat{\theta}_{P_{k},1}-\hat{\theta}_{P_{k-1},1})^2}   \hat{u}_{P_{k-1},1,k}^2 }_{(e)\overset{\scalebox{0.5}{order}}{\sim} \ \frac{\theta_1 m}{ \min \left( \theta_1,\theta_k \right)^2}  }+  
\underbrace{\sum_{i=2}^{k-1} \frac{\hat{\theta}_{P_{k-1},i}^2}{(\hat{\theta}_{P_k,1}-\hat{\theta}_{P_{k-1},i})^2}  \hat{u}_{P_{k-1},i,k}^2}_{(f) O_p\left( \frac{1}{\theta_1 m}\right) }}, \hspace{20cm}\\
&& \hspace{-0.5cm} \left\langle \tilde{u}_{P_k,1},e_k \right\rangle^2 = \frac{1}{D_1 (\theta_k-1)^2} (g),\hspace{20cm} \\
 && \hspace{-0.5cm} \left\langle \tilde{u}_{P_k,1},e_s \right\rangle^2\hspace{-1.5cm}  \hspace{20cm} \\
&&\hspace{0.5cm}  = \frac{1}{D_1}  \left( 
\overbrace{\sum_{i=k}^m \frac{\hat{\lambda}_{P_{k-1},i}}{\hat{\theta}_{P_k,1}-\hat{\lambda}_{P_{k-1},i}} \hat{u}_{P_{k-1},i,s} \hat{u}_{P_{k-1},i,k}}^{(h) O_p\left( \frac{1}{\theta_s^{1/2} \theta_1 \sqrt{m}}\right)} 
+ 
\overbrace{\frac{\hat{\theta}_{P_{k-1},1}}{\hat{\theta}_{P_{k},1}-\hat{\theta}_{P_{k-1},1}}  \hat{u}_{P_{k-1},1,s} \hat{u}_{P_{k-1},1,k}}^{(i) \ \overset{\scalebox{0.5}{order}}{\sim} \frac{\min\left( \theta_1 , \theta_s \right)}{\sqrt{\theta_s}\min\left( \theta_1 , \theta_k \right)}}\right. \\
&&\hspace{3cm}\left.+  
\overbrace{\sum_{i=2,\neq s}^{k-1} \frac{\hat{\theta}_{P_{k-1},i}}{\hat{\theta}_{P_k,1}-\hat{\theta}_{P_{k-1},i}} \hat{u}_{P_{k-1},i,s} \hat{u}_{P_{k-1},i,k} }^{(j) O_p\left( \underset{i \neq 1,s}{\max}\left( \frac{\min\left( \theta_1,\theta_i \right)\min\left( \theta_s, \theta_i \right)}{\sqrt{\theta_s} \theta_1 \theta_i \sqrt{m}} \right)\right)}
+  
\overbrace{ \frac{\hat{\theta}_{P_{k-1},s}}{\hat{\theta}_{P_k,1}-\hat{\theta}_{P_{k-1},s}} \hat{u}_{P_{k-1},s,s} \hat{u}_{P_{k-1},s,k} }^{(k) O_p\left( \frac{\min\left( \theta_1 , \theta_s \right)}{\sqrt{\theta_s} \theta_1 \sqrt{m}} \right)} 
\right)^2 .
\end{eqnarray*} 
\end{minipage}}\\
\noindent Finally, 
$$\hat{u}_{P_k,1}= \frac{\left( \tilde{u}_{P_k,1,1},\tilde{u}_{P_k,1,2},...,\sqrt{\theta_k} \tilde{u}_{P_k,1,k},...,\tilde{u}_{P_k,} \right)}{\underbrace{ \sqrt{1+ \left(\theta_k-1 \right) \tilde{u}_{P_k1,k}^2}}_{1+O_p\left(\frac{\min(\theta_1,\theta_k)}{\max(\theta_1,\theta_k)m} \right)}},$$
where $\sqrt{1+ \left(\theta-1 \right) \tilde{u}_{P_k1,k}^2}$ is the norm of $P_k^{1/2}\tilde{u}_{P_k,1}$ that we will call $N_1$.
\begin{Rem}\ 
\begin{enumerate}
\item By construction, the sign of $\hat{u}_{P_k,1,k}$ is always positive. This is, however, not the case of $\hat{u}_{P_{k-1},i,i}$.  We can show that\\
\scalebox{0.85}{
\begin{minipage}{1\textwidth}
\begin{eqnarray*}
P\left\lbrace {\rm sign} \left( \hat{u}_{P_k,1,1} \right) ={\rm sign} \left( \left(\hat{\theta}_{P_k,1}-\hat{\theta}_{P_{k-1},1}\right) \hat{u}_{P_{k-1},1,1} \hat{u}_{P_{k-1},1,k}\right) \right\rbrace \underset{m \rightarrow \infty}{\rightarrow} 1.
\end{eqnarray*}
\end{minipage}}\\

Moreover, the convergence to $1$ is of order $1/m$. If $\theta_1$ tends to infinity, then 
\begin{eqnarray*}
P\left\lbrace {\rm sign} \left( \hat{u}_{P_k,1,1} \right) ={\rm sign} \left( \left(\theta_1-\theta_k\right)  \hat{u}_{P_{k-1},1,k}\right) \right\rbrace \underset{m,\theta_1 \rightarrow \infty}{\rightarrow} 1.
\end{eqnarray*}
Thus, if we use a convention such as  $\text{sign}\left(\hat{u}_{P_{k},i,i}\right) >0 $ for $i=1,...,k-1$, 
then the sign of $\hat{u}_{P_k,1,k}$ is distributed as a Bernoulli with parameter 1/2.

\item Without loss of generality, the other eigenvector $\hat{u}_{P_k,r}$ for $r=1,2,...,k-1$ can be computed by the same formula thanks to the notation linking the estimated eigenvector $\hat{u}_{P_k,r}$ to the eigenvalue $\theta_r$.\\ 
However, the formula does not work for the vector $\hat{u}_{P_k,k}$. Indeed it allows to express the $k-1$ eigenvectors, $\hat{u}_{P_{k},s}$ for $s=1,2,...,k-1$, as a function of the $\hat{u}_{P_{k-1},i}$, $\hat{\lambda}_{P_{k-1},\tilde{i}}$ , $\hat{\theta}_{P_{k-1},s}$ and $\hat{\theta}_{P_{k},s}$ for $i=1,2...,m$, $\tilde{i}=k,k+1...,m$, $s=1,2...,k-1$. Applying the perturbation in a different fashion shows that similar formulas do exist for $\hat{u}_{P_k,k}$. (If by permuting the indices $k$ and $1$, the perturbation in $e_1$ is applied at the end, for example.) However the eigenstructure of this last vector will not be expressed in function of the same random variables.\\
This observation exhibits a problem in the proofs of the Dot Product Theorems \ref{AAThDotproduct} and \ref{AAThInvariantdot}. Deeper investigations are necessary to understand the two eigenvectors when $k=2$ and express both $\hat{u}_{P_2,1}$ and $\hat{u}_{P_2,2}$ as a function of $\hat{u}_{P_{1},i}$, $\hat{\lambda}_{P_{1},\tilde{i}}$ , $\hat{\theta}_{P_{1},1}$ and $\hat{\theta}_{P_{2},2}$ for $1=1,2.,,,,m$.
\begin{eqnarray*}
D_2&=&\underbrace{\sum_{i=2}^m \frac{\hat{\lambda}_{P_{1},i}^2}{(\hat{\theta}_{P_2,2}-\hat{\lambda}_{P_{1},i})^2} \hat{u}_{P_{1},i,2}^2}_{ O_p\left( \frac{1}{\theta_2^2}\right) } + 
\underbrace{\frac{\hat{\theta}_{P_{1},1}^2}{(\hat{\theta}_{P_{2},2}-\hat{\theta}_{P_{1},1})^2}   \hat{u}_{P_{1},1,2}^2 }_{O_p\left( \frac{\theta_1}{(\theta_2-\theta_1)^2 m}\right) },\\
N_2^2&=&1+\frac{1}{(\theta_2-1)D_2},\\
N_2 D_2 &=& D_2 + \frac{1}{\theta_2-1} \\
&=& \frac{1}{\theta_2-1}+O_p\left( \frac{1}{\theta_2^2} \right)+O_p\left( \frac{\theta_1}{(\theta_2-\theta_1) m} \right).
\end{eqnarray*}
Furthermore, the theorem must investigate the $m-k$ noisy components of the eigenvectors. For $r=1,2$ and $s=3,4,...,m$,
\begin{eqnarray*}
\hat{u}_{P_2,r,s}= \frac{\sum_{i=1}^m \frac{\hat{\lambda}_{P_{1},i}}{\hat{\theta}_{P_2,r}-\hat{\lambda}_{P_{1},i}} \hat{u}_{P_{1},i,s} \hat{u}_{P_{1},i,2}}{\sqrt{D_r} N_r}.
\end{eqnarray*}
The estimations using this last formula are difficult. When we investigate these components, it is profitable to look at 
$$\hat{u}_{P_2,1,t}/\sqrt{\sum_{s=3}^m \hat{u}_{P_2,1,s}^2} \text{ and } \hat{u}_{P_2,2,t}/\sqrt{\sum_{s=3}^m \hat{u}_{P_2,2,s}^2}$$
 for $t=3,4,...,m$.
\item If the perturbation is not canonical, then we can apply a rotation $U$, such that $U u_s = \epsilon_s$, and replace $\hat{u}_{P_{k-1},i}$ by $U^t \hat{u}_{P_{k-1},i}$. Then, $\left\langle \tilde{u}_{P_k,1},e_s \right\rangle^2$ is replaced by $\left\langle \tilde{u}_{P_k,1},u_s \right\rangle^2$.
\end{enumerate}
\end{Rem}
\end{itemize}

\end{Th}

\subsubsection{Double dot product}
\begin{Th}
 \label{AAThdoubledot} 
Suppose $W_X$ and $W_Y$ satisfy Assumption \ref{AAAss=matrice} and 
$P_k=\I_m+\sum_{i=1}^k (\theta_i-1) e_i e_i^t$ satisfies \ref{AAAss=theta} (A4),
where $\theta_1>\theta_2>...>\theta_k$. We set
\begin{eqnarray*}
&&\hat{\Sigma}_{X}=\hat{\Sigma}_{X,P_k}=P_k^{1/2} W_X P_k^{1/2} \text{ and } \hat{\Sigma}_{Y,P_k}= P_k^{1/2} W_Y P_k^{1/2}.
\end{eqnarray*}
and for $s=1,...,k$,
\begin{eqnarray*}
\hat{u}_{\hat{\Sigma}_{X},s}, \hat{\theta}_{\hat{\Sigma}_{X},s} &\text{ s.t. }& \hat{\Sigma}_{X} \hat{u}_{\hat{\Sigma}_{X},s} =\hat{\theta}_{\hat{\Sigma}_{X},s} \hat{u}_{\hat{\Sigma}_{X},s},\\
\hat{u}_{\hat{\Sigma}_{Y},s}, \hat{\theta}_{\hat{\Sigma}_{Y},s} &\text{ s.t. }& \hat{\Sigma}_{Y} \hat{u}_{\hat{\Sigma}_{Y},s} =\hat{\theta}_{\hat{\Sigma}_{Y},s} \hat{u}_{\hat{\Sigma}_{Y},s}, 
\end{eqnarray*}
where $\hat{\theta}_{\hat{\Sigma}_{Y},s}=\hat{\lambda}_{\hat{\Sigma}_{Y},s}$ and
 $\hat{\theta}_{\hat{\Sigma}_{X},s}=\hat{\lambda}_{\hat{\Sigma}_{X},s}$.  To simplify the result we assume the sign convention: 
\begin{eqnarray*}
\text{For $s=1,2,...,k$ and $i=1,2,...,s$, } \hat{u}_{\hat{\Sigma}_X,i,i}>0, \ \hat{u}_{\hat{\Sigma}_Y,i,i}>0.
\end{eqnarray*}
Finally, we define
\begin{eqnarray*}
\tilde{u}_s = \hat{U}_X^t \hat{\hat{u}}_{\hat{\Sigma}_{Y},s},
\end{eqnarray*}
where,
\begin{eqnarray*}
\hat{U}_X=\left(
v_1,
v_2, 
\cdots,
v_m
\right)=\left(
\hat{u}_{\hat{\Sigma}_{X},1},
\hat{u}_{\hat{\Sigma}_{X},2}, 
\cdots 
\hat{u}_{\hat{\Sigma}_{X},k}, 
v_{k+1},
v_{k+2},
\cdots,
v_m
\right),
\end{eqnarray*}
where the vectors $v_{k+1},...,v_m$ are chosen such that the matrix $\hat{U}_X$ is orthonormal. Then, 
\begin{itemize}
\item If $\theta_j,\theta_t \rightarrow \infty$: \\
\scalebox{0.77}{
\begin{minipage}{1\textwidth}
\begin{eqnarray*}
\sum_{i=k+1}^m \tilde{u}_{j,i} \tilde{u}_{t,i}&=& \sum_{i=k+1}^m \hat{u}_{\hat{\Sigma}_{Y},j,i} \hat{u}_{\hat{\Sigma}_{Y},t,i} +\sum_{i=k+1}^m \hat{u}_{\hat{\Sigma}_{X},j,i} \hat{u}_{\hat{\Sigma}_{X},t,i} -\sum_{i=k+1}^m \hat{u}_{\hat{\Sigma}_{X},j,i} \hat{u}_{\hat{\Sigma}_{Y},t,i} \\
&& \hspace{2cm} - \sum_{i=k+1}^m \hat{u}_{\hat{\Sigma}_{Y},j,i} \hat{u}_{\hat{\Sigma}_{X},t,i}- \left( \hat{u}_{\hat{\Sigma}_{X},t,j}+\hat{u}_{\hat{\Sigma}_{Y},j,t} \right) \left( \hat{\alpha}^2_{\hat{\Sigma}_{X},j} - \hat{\alpha}^2_{\hat{\Sigma}_{X},t} \right)\\
&& \hspace{2cm}+ O_p\left( \frac{1}{\theta_1 m} \right)+O_p\left( \frac{1}{\theta_1^2  \sqrt{m}} \right),
\end{eqnarray*}
\end{minipage}}\\

where $\hat{\alpha}^2_{\hat{\Sigma}_{X},t}= \sum_{i=1}^k \hat{u}_{\hat{\Sigma}_{X},t,i}^2$. 
\item If $\theta_t$ is finite, then 
\begin{eqnarray*}
\sum_{i=k+1}^m \tilde{u}_{j,i} \tilde{u}_{t,i}=O_p\left( \frac{1}{\sqrt{m}\sqrt{\theta_1}} \right).
\end{eqnarray*}
\end{itemize}

\noindent Moreover, for $s=1,...,k$, $t=2,...,k$ and $j=k+1,...m$,
\begin{eqnarray*}
&&\sum_{i=1}^k \tilde{u}_{s,i}^2 = \sum_{i=1}^k \left\langle \hat{u}_{\hat{\Sigma}_{X},i},\hat{u}_{\hat{\Sigma}_{Y},s} \right\rangle^2,\\
&& \tilde{u}_{s,s}=\hat{u}_{\hat{\Sigma}_{X},s,s}\hat{u}_{\hat{\Sigma}_{Y},s,s}+O_p\left( \frac{1}{m} \right) +O_p\left( \frac{1}{\theta_s^{1/2}m^{1/2}} \right),\\
 &&\tilde{u}_{s,t}=\hat{u}_{\hat{\Sigma}_{X},t,s}+\hat{u}_{\hat{\Sigma}_{X},s,t}+O_p\left( \frac{\sqrt{\min(\theta_s,\theta_t)}}{m\sqrt{\max(\theta_s,\theta_t)}} \right)+O_p\left( \frac{1}{ \theta_t m^{1/2}} \right),\\
  &&\tilde{u}_{t,s}=O_p\left( \frac{\sqrt{\min(\theta_s,\theta_t)}}{m\sqrt{\max(\theta_s,\theta_t)}} \right)+O_p\left( \frac{1}{ \theta_s m^{1/2}} \right),\\
&& \tilde{u}_{s,j}=\hat{u}_{\hat{\Sigma}_{Y},s,j}-\hat{u}_{\hat{\Sigma}_{X},s,j} \left\langle \hat{u}_{\hat{\Sigma}_{Y},j},\hat{u}_{\hat{\Sigma}_{X},j} \right\rangle+ O_p\left( \frac{1}{\theta_s^{1/2} m } \right).
\end{eqnarray*}
\end{Th}

\subsubsection{Lemmas}
\begin{Lem}\label{AALemstatW}
Suppose $W$ and $\hat{\Sigma}_{P_1}$ are as in Theorem \ref{AAThDotproduct}, then by construction of the eigenvectors using Theorem \ref{AATheoremcaraceigenstructure},\\
\scalebox{0.68}{
\begin{minipage}{1\textwidth}
\begin{eqnarray*}
&\hat{u}_{P_1,1,2} &=\frac{W_{1,2}}{\sqrt{\theta_1} W_{1,1}} -\frac{W_{1,2}}{\theta_1^{3/2}}\left( -1/2+3/2 M_2 \right)+ \frac{\left(W^2\right)_{1,2}}{\theta_1^{3/2}} + O_p\left( \frac{1}{\theta_1^{3/2}m} \right)+ O_p\left( \frac{1}{\theta_1^{5/2}m^{1/2}} \right)\\
&&= \frac{W_{1,2}}{\sqrt{\theta_1}}+O_p\left( \frac{1}{\theta_1^{1/2} m} \right)+O_p\left( \frac{1}{\theta_1^{3/2} m^{1/2}} \right),\\
&\sum_{i=2}^m \hat{\lambda}_{P_1,i}^2 \hat{u}_{P_1,i,2}^2 &=W_{2,2}+O_p\left( \frac{1}{m}\right),\\
&\sum_{i=2}^m \hat{\lambda}_{P_1,i} \hat{u}_{P_1,i,1}\hat{u}_{P_1,i,2}&=W_{1,2} \frac{  M_2}{\sqrt{\theta_1}} -\left(W^2\right)_{1,2}\frac{1}{\sqrt{\theta_1}}
+ O_p\left( \frac{1}{\theta_1^{1/2}m} \right)+ O_p\left( \frac{1}{\theta_1^{3/2}m^{1/2}} \right).
\end{eqnarray*}
\end{minipage}}

\end{Lem}

\begin{Lem} \label{AALemreducedim} 
Suppose $w_1,...,w_k \in \mathbb{R}^m$ and $\lambda_1,...,\lambda_k \in \mathbb{R}^*$, then if the function $\lambda()$ provides non-trivial eigenvalues,
\begin{eqnarray*}
\lambda\Bigg( \sum_{i=1}^k \lambda_i w_i w_i^t \Bigg) = \lambda\Bigg( H \Bigg) ,
\end{eqnarray*}
where\\
\scalebox{0.9}{
\begin{minipage}{1\textwidth}
\begin{eqnarray*}
H=\begin{pmatrix}
{\lambda}_1 & \sqrt{{\lambda}_1 {\lambda}_2} \left\langle w_1,w_2 \right\rangle & \sqrt{{\lambda}_1 {\lambda}_3} \left\langle w_1,w_3 \right\rangle & \cdots & \sqrt{{\lambda}_k {\lambda}_2} \left\langle w_1,w_k \right\rangle \\ 
\sqrt{{\lambda}_2 {\lambda}_1} \left\langle w_2,w_1 \right\rangle  & {\lambda}_2 & \sqrt{{\lambda}_2 {\lambda}_3} \left\langle w_2,w_3 \right\rangle & \cdots  & \sqrt{{\lambda}_2 {\lambda}_k} \left\langle w_2,w_k \right\rangle \\ 
\sqrt{{\lambda}_3 {\lambda}_1} \left\langle w_3,w_1 \right\rangle &\sqrt{{\lambda}_3 {\lambda}_2} \left\langle w_3,w_2 \right\rangle  & {\lambda}_3 & \cdots  & \sqrt{{\lambda}_3 {\lambda}_k} \left\langle w_3,w_k \right\rangle \\ 
\vdots & \vdots & \ddots & \ddots &  \vdots \\ 
 \sqrt{{\lambda}_k {\lambda}_1} \left\langle w_k,w_1 \right\rangle & \sqrt{{\lambda}_k {\lambda}_2} \left\langle w_k,w_2 \right\rangle & \sqrt{{\lambda}_k {\lambda}_3} \left\langle w_k,w_3 \right\rangle & \cdots  & {\lambda}_k  \\ 
\end{pmatrix} .
\end{eqnarray*}
\end{minipage}}
\end{Lem}

\subsection{Proofs}

\subsubsection{Invariant proofs}
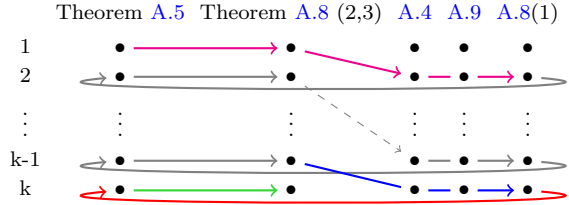
\begin{figure}[htbp]

\begin{center} 
\begin{tikzpicture}
  \matrix (magic) [matrix of nodes]
  { \  & Theorem \ref{AAInvariantth} & Theorem \ref{AAThcomponentdistribution} (2,3)    &   \ref{AAThinvarianteigenvalue} & \ref{AATheoremcaraceigenstructure} & \ref{AAThcomponentdistribution}(1)  \\
   1  & $\bullet$   & $\bullet$    & $\bullet$ &   $\bullet$&  $\bullet$   \\
   2  & $\bullet$   & $\bullet$ & $\bullet$ &    $\bullet$&  $\bullet$\\
   \vdots & \vdots   & \vdots  & \vdots & \vdots & \vdots \\
    k-1  & $\bullet$  & $\bullet$  & $\bullet$  &  $\bullet$ &  $\bullet$ \\
   k  & $\bullet$ & $\bullet$   &  $\bullet$ & $\bullet$  &  $\bullet$ \\
};
  \draw[thick,magenta,->] (magic-2-2) -- (magic-2-3) ;
   \draw[thick,magenta,->] (magic-2-3) -- (magic-3-4);
   \draw[thick,magenta,->] (magic-3-4) -- (magic-3-5) -- (magic-3-6)   ;
    
\draw[thick,gray,->] (magic-3-6) to [out=-5,in=-175,looseness=1](magic-3-2);
\draw[thick,gray,->] (magic-3-2) -- (magic-3-3);
\draw[dashed,gray,->] (magic-3-3) -- (magic-5-4);
\draw[thick,gray,->]  (magic-5-4) -- (magic-5-5) -- (magic-5-6);
\draw[thick,gray,->] (magic-5-6) to [out=-5,in=-175,looseness=1](magic-5-2);
\draw[thick,gray,->] (magic-5-2) -- (magic-5-3) ;

\draw[thick,green!80!black!80,->] (magic-6-2) -- (magic-6-3)  ;
\draw[thick,blue,->]  (magic-5-3) -- (magic-6-4)  -- (magic-6-5) -- (magic-6-6);
\draw[thick,red,->] (magic-6-6) to [out=-5,in=-175,looseness=1](magic-6-2);
\end{tikzpicture}
\end{center}
\caption{Procedure used in the proof.} \label{AAfig=proof}
\end{figure}

In this section, we prove some invariance results by induction. The procedure is summarized in Figure \ref{AAfig=proof}. First we initialize the induction (in pink). Then, the induction assumes the proven results in the grey part and proves the blue, red and green parts. 

\paragraph{Pink}
First, we discuss the initialization part.
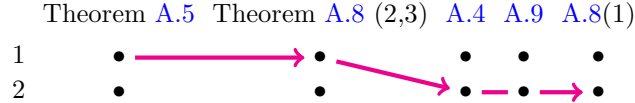
\begin{center}
\begin{tikzpicture}
  \matrix (magic) [matrix of nodes]
  { \  & Theorem \ref{AAInvariantth} & Theorem \ref{AAThcomponentdistribution} (2,3)    &   \ref{AAThinvarianteigenvalue} & \ref{AATheoremcaraceigenstructure} & \ref{AAThcomponentdistribution}(1)  \\
   1  & $\bullet$   & $\bullet$   & $\bullet$ &  $\bullet$   & $\bullet$\\
   2  & $\bullet$   & $\bullet$ & $\bullet$ &  $\bullet$&  $\bullet$\\
};
   \draw[ultra thick,magenta,->] (magic-2-2) -- (magic-2-3) ;
   \draw[ultra thick,magenta,->] (magic-2-3) -- (magic-3-4);
   \draw[ultra thick,magenta,->] (magic-3-4) -- (magic-3-5) -- (magic-3-6)   ;
\end{tikzpicture}
\end{center}

\noindent The Invariant Theorem \ref{AAInvariantth} is trivially true for perturbations of order $k=1$. 
\begin{proof}{\textbf{Theorem \ref{AAThcomponentdistribution} (2,3), $\mathbf{k=1}$}}  \label{AAproofThcomponentdistribution23k1}

\noindent In the following picture we can assume the first result for $k=1$ is proven.
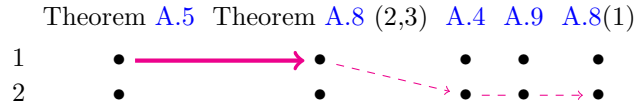
\begin{center}
\begin{tikzpicture}
  \matrix (magic) [matrix of nodes]
  { \  & Theorem \ref{AAInvariantth} & Theorem \ref{AAThcomponentdistribution} (2,3)    &   \ref{AAThinvarianteigenvalue} & \ref{AATheoremcaraceigenstructure} & \ref{AAThcomponentdistribution}(1)  \\
   1  & $\bullet$   & $\bullet$   & $\bullet$ &   $\bullet$   & $\bullet$\\
   2  & $\bullet$   & $\bullet$ & $\bullet$ &  $\bullet$&  $\bullet$\\
};
   \draw[ultra thick,magenta,->] (magic-2-2) -- (magic-2-3) ;
   \draw[dashed,magenta,->] (magic-2-3) -- (magic-3-4);
   \draw[dashed,magenta,->] (magic-3-4) -- (magic-3-5) -- (magic-3-6)   ;
\end{tikzpicture}
\end{center}
We define
\begin{eqnarray*}
U&=&
\begin{pmatrix}
\hat{u}_{P_k,1}^t\\
\hat{u}_{P_k,2}^t\\
\vdots \\
\hat{u}_{P_k,m}^t
\end{pmatrix}=\begin{pmatrix}
\hat{u}_{P_k,1:k,1:k} & \hat{u}_{P_k,1:k,k+1:m} \\ 
\hat{u}_{P_k,k+1:m,1:k} & \hat{u}_{P_k,k+1:m,k+1:m}.
\end{pmatrix}
\end{eqnarray*}
and
$$ O_1= \begin{pmatrix}
\rm{I}_1 & 0 \\ 
0 & O_{m-1}
\end{pmatrix},$$
where $O_{m-1}$ is a rotation matrix.
\begin{enumerate}
\item[2.] Assuming a canonical $P_1=\I_m+(\theta_1-1) e_1 e_1^t$, we know that $\hat{\Sigma} \sim P_1^{1/2} W P_1^{1/2}$ and $O_1 \hat{\Sigma} O_1^t$ follow the same distribution under Assumption \ref{AAAss=matrice}. Although the eigenvectors change, they still follow the same distribution, $O_1 U^t \sim U^t$. Therefore,  $\hat{u}_{i,(k+1):m}$ is rotationally invariant and $\corr \left( \hat{u}_{i,j_1}, \hat{u}_{i,j_2} \right)=\delta_{j_1}(j_2)$. \\
We can show that knowing the first line of the matrix, then  \linebreak
$\hat{u}_{P_1,i,2:m}/ ||\hat{u}_{P_1,i,2:m}||$ is unit uniform for $i=1,2,...,m$. Therefore,  these statistics are independent (not jointly) of the first line.\\
Uniformity of $\hat{u}_{P_1,i,2:m}$ implies for $s=2,3,...,m$, 
\begin{eqnarray*}
\sqrt{m} \frac{\hat{u}_{P_1,1,s}}{ ||\hat{u}_{P_1,1,2:m}||}= \sqrt{m} \frac{\hat{u}_{P_1,1,s}}{ \sqrt{1-\hat{\alpha}^2_{P_1,1}}}\sim \Normal\left(0,1\right)+o_p \left( 1 \right).
\end{eqnarray*}
By Slutsky's Theorem and the distribution of the angle for $k=1$, Theorem \ref{AAjointdistribution},
 \begin{eqnarray*}
\hat{u}_{P_1,1,s} \sim \Normal \left( 0,\frac{1-\alpha_1^2}{m}\right)+o_p \left( \frac{1}{\sqrt{m}} \right),
\end{eqnarray*}
where $\alpha_1^2$ is the limit of the angle and can be approximated by $1-\frac{M_2-1}{\theta_1}+O_p\left( \frac{1}{\theta_1^2}\right)<1$.

\item[3.] 
Using the distribution of $\hat{\alpha}^2_{P_1,1,1}$ given in Theorem \ref{AAjointdistribution},
 $$\sum \hat{u}_{2:m,1}^2=1-\hat{\alpha}^2_{P_1,1,1} \sim {\rm{RV}}\left( O\left( \frac{1}{\theta_1 } \right), O \left( \frac{1}{\theta_1^2 m}\right)\right).$$
\end{enumerate}
\end{proof}

\noindent Then, we prove the Invariant Angle Theorem for the eigenvalues, Theorem \ref{AAThinvarianteigenvalue} for k=2.

\begin{proof}{\textbf{Theorem \ref{AAThinvarianteigenvalue}, $\mathbf{k=2}$}}  \label{AAproofThinvarianteigenvaluek2}
\noindent We prove the theorem for $k=2$. In the following picture we can assume the grey results as proven.
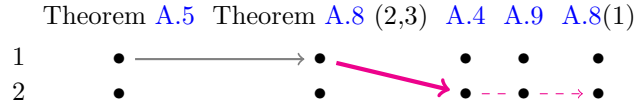
\begin{center}
\begin{tikzpicture}
  \matrix (magic) [matrix of nodes]
  { \  & Theorem \ref{AAInvariantth} & Theorem \ref{AAThcomponentdistribution} (2,3)    &   \ref{AAThinvarianteigenvalue} & \ref{AATheoremcaraceigenstructure} & \ref{AAThcomponentdistribution}(1)  \\
   1  & $\bullet$   & $\bullet$    &   $\bullet$&  $\bullet$   & $\bullet$\\
   2  & $\bullet$   & $\bullet$  &  $\bullet$ &  $\bullet$&  $\bullet$\\
};
   \draw[thick,gray,->] (magic-2-2) -- (magic-2-3) ;
   \draw[ultra thick,magenta,->] (magic-2-3) -- (magic-3-4);
   \draw[dashed,magenta,->] (magic-3-4) -- (magic-3-5) -- (magic-3-6)   ;
\end{tikzpicture}
\end{center}

\noindent Without loss of generality, we only prove the invariance of $\hat{\theta}_{P_{1},1}$. For simplicity, we assume $\theta_1> \theta_2$, but this assumption is only used to simplify notation. Each step can be done assuming $\theta_1 < \theta_2$. Using Theorem \ref{AAconvergence} and the canonical perturbation $\tilde{P_2}$ lead to
\begin{eqnarray*}
&& \sum_{i=2}^m \frac{\hat{\lambda}_{P_1,i}}{\hat{\theta}_{P_{2},1}-\hat{\lambda}_{P_1,i}} \hat{u}_{P_{1},i,2}^2 
+ \frac{\hat{\theta}_{P_{1},1}}{\hat{\theta}_{P_{2},1}-\hat{\theta}_{P_{1},1}} \hat{u}_{P_{1},1,2}^2 =\frac{1}{\theta_2-1}.
\end{eqnarray*}
Therefore,  
\begin{eqnarray*}
\frac{\hat{\theta}_{P_{1},1}}{\hat{\theta}_{P_{2},1}-\hat{\theta}_{P_{1},1}} \hat{u}_{P_{1},1,2}^2 &=& -\sum_{i=2}^m \frac{\hat{\lambda}_{P_1,i}}{\hat{\theta}_{P_{2},1}-\hat{\lambda}_{P_1,i}} \hat{u}_{P_{1},i,2}^2 + \frac{1}{\theta_2-1}  \\
&=& - \frac{1}{\hat{\theta}_{P_{2},1}}\sum_{i=2}^m \hat{\lambda}_{P_1,i} \hat{u}_{P_{1},i,2}^2 + \frac{1}{\theta_2-1} + O_p\left( \frac{1}{\theta_1^2} \right) \\
&\overset{1*}{=}& - \frac{1}{\hat{\theta}_{P_{2},1}} \left( 1+O_p\left(\frac{1}{\sqrt{m}}\right) \right) + \frac{1}{\theta_2-1} + O_p\left( \frac{1}{\theta_1^2} \right) \\
&=& -\frac{\theta_2-1 -\hat{\theta}_{P_{2},1}}{\hat{\theta}_{P_{2},1} (\theta_2-1)} + O_p\left(\frac{1}{\theta_1^2} \right)+ O_p\left(\frac{1}{\theta_1 \sqrt{m}} \right),
\end{eqnarray*}
where 1* is true because
\begin{eqnarray*}
\sum_{i=2}^m \hat{\lambda}_{P_1,i} \hat{u}_{P_{1},i,2}^2 &=& \sum_{i=1}^m \hat{\lambda}_{P_1,i} \hat{u}_{P_{1},1,2}^2 -  \hat{\theta}_{P_{1},1} \hat{u}_{P_{1},1,2}^2\\
&=& \hat{\Sigma}_{P_1,2,2}-\hat{\theta}_{P_{1},1} \hat{u}_{P_{1},1,2}^2\\
&=& W_{2,2}-\hat{\theta}_{P_{1},1} \hat{u}_{P_{1},1,2}^2\\
&=& 1+O_p\left(\frac{1}{\sqrt{m}}\right).
\end{eqnarray*}
The last line is obtained using the fact that the canonical perturbation $P_1$ does not affect $ W_{2:m,2:m}$. Moreover, $W$ satisfies Assumption \ref{AAAss=matrice} and thus $W_{2,2}=1+O_p\left(1/\sqrt{m}\right)$. On the other hand the second term $\hat{\theta}_{P_{1},1} \hat{u}_{P_{1},1,2}^2=O_p\left(1/m\right) $ by Theorem \ref{AAThcomponentdistribution}(2) for $k=1$.\\

\noindent By Theorem \ref{AAThcomponentdistribution}(2),\\ 
\scalebox{0.77}{
\begin{minipage}{1\textwidth}
\begin{eqnarray*}
 \left(1+ O_p\left(\frac{\theta_2}{\theta_1 (\theta_2-\theta_1)} \right)+ O_p\left(\frac{\theta_2}{ \sqrt{m} (\theta_2-\theta_1 )} \right)\right)  \left(\hat{\theta}_{P_{2},1}- \hat{\theta}_{P_{1},1}\right) &=& 
-\frac{\hat{\theta}_{P_{1},1} \hat{\theta}_{P_{2},1} \left(\theta_2-1 \right)}{\theta_2-1 -\hat{\theta}_{P_{2},1}} \hat{u}_{P_{1},1,2}^2\\
&=& O_p\left( \frac{\theta_1 \theta_2}{m(\theta_2-\theta_1)} \right).
\end{eqnarray*}
\end{minipage}}\\

We note that even without Assumption \ref{AAAss=theta}(A4), we have
\begin{eqnarray*}
&&\hat{\theta}_{P_{2},1}- \hat{\theta}_{P_{1},1} \overset{\scalebox{0.5}{order}}{\sim} \frac{\min\left(\theta_1,\theta_2 \right)}{m}.
\end{eqnarray*}
More precisely we can write\\ 
\scalebox{0.95}{
\begin{minipage}{1\textwidth}
\begin{eqnarray*}
 \hat{\theta}_{P_{2},1}- \hat{\theta}_{P_{1},1} &=& 
-\frac{\hat{\theta}_{P_{1},1} \hat{\theta}_{P_{2},1} \left(\theta_2-1 \right)}{\theta_2-1 -\hat{\theta}_{P_{2},1}} \hat{u}_{P_{1},1,2}^2 + O_p \left(\frac{1}{m} \right)+ O_p \left(\frac{\min(\theta_1,\theta_2)}{m^{3/2}} \right).
\end{eqnarray*}
\end{minipage}}\\

\noindent Each step of the above computation can be done for $\hat{\theta}_{\tilde{P}_{2},1}-\hat{\theta}_{P_2,2}$. Therefore,  for $s\neq t \in \lbrace 1,2 \rbrace$ we obtain the general result. \\ 
\scalebox{0.8}{
\begin{minipage}{1\textwidth}
\begin{eqnarray*}
 \left(1+ O_p\left(\frac{\theta_t}{\theta_s (\theta_t-\theta_s)} \right)+ O_p\left(\frac{\theta_t}{ \sqrt{m} (\theta_t-\theta_s )} \right)\right)  \left( \hat{\theta}_{P_{2},s}-\hat{\theta}_{\tilde{P}_{s},1}\right) &=& 
-\frac{\hat{\theta}_{\tilde{P}_{s},1} \hat{\theta}_{P_{2},s} \left(\theta_t-1 \right)}{\theta_t-1 -\hat{\theta}_{P_{2},s}} \hat{u}_{\tilde{P}_{s},1,t}^2\\
&\overset{\scalebox{0.5}{order}}{\sim} &  \frac{\theta_1 \theta_2}{m(\theta_2-\theta_1)} .
\end{eqnarray*}
\end{minipage}}\\
This leads to \\ 
\scalebox{0.95}{
\begin{minipage}{1\textwidth}
\begin{eqnarray*}
\hat{\theta}_{P_{2},s} -\hat{\theta}_{\tilde{P}_{s},1}
&=& 
-\frac{\hat{\theta}_{\tilde{P}_{s},1} \hat{\theta}_{P_{2},s} \left(\theta_t-1 \right)}{\theta_t-1 -\hat{\theta}_{P_{2},s}} \hat{u}_{\tilde{P}_{s},1,t}^2 + O_p \left(\frac{1}{m} \right)+ O_p \left(\frac{\min(\theta_1,\theta_2)}{m^{3/2}} \right)\\
&\overset{\scalebox{0.5}{order}}{\sim} & \frac{\min\left(\theta_1,\theta_2 \right)}{m} .
\end{eqnarray*}
\end{minipage}}

\end{proof}

\begin{proof}{\textbf{Theorem \ref{AATheoremcaraceigenstructure}, \label{AAproofTheoremcaraceigenstructurek2} and \ref{AAThcomponentdistribution}(1), $\mathbf{k=2}$}}\label{AAproofThcomponentdistribution1k2} 
\noindent We prove the theorems for $k=2$. In the following picture we can assume the grey results as proven.
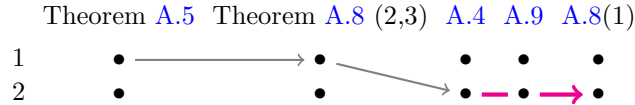
\begin{center}
\begin{tikzpicture}
  \matrix (magic) [matrix of nodes]
  { \  & Theorem \ref{AAInvariantth} & Theorem \ref{AAThcomponentdistribution} (2,3)    &   \ref{AAThinvarianteigenvalue} & \ref{AATheoremcaraceigenstructure} & \ref{AAThcomponentdistribution}(1)  \\
   1  & $\bullet$   & $\bullet$    &   $\bullet$&  $\bullet$   & $\bullet$\\
   2  & $\bullet$   & $\bullet$  &  $\bullet$ &  $\bullet$&  $\bullet$\\
};
   \draw[thick,gray,->] (magic-2-2) -- (magic-2-3) ;
   \draw[thick,gray,->] (magic-2-3) -- (magic-3-4);
   \draw[ultra thick,magenta,->] (magic-3-4) -- (magic-3-5) -- (magic-3-6)   ;
\end{tikzpicture}
\end{center}
These proofs are exactly the same when the perturbation is of order $k$. Thus,  we will do it only once in pages \pageref{AAproofThcomponentdistribution1k} and \pageref{AAproofTheoremcaraceigenstructurek}. As we will see, the proofs of these theorems uses only the grey results and the proof of Theorem \ref{AAThcomponentdistribution}(1) for $k$ uses Theorem \ref{AATheoremcaraceigenstructure} for $k$. Moreover, although the proof of Theorem \ref{AATheoremcaraceigenstructure} for $k>2$ uses Theorem \ref{AAThcomponentdistribution}(1) for $k-1$, the initializing part $k=2$ does not need Theorem \ref{AAThcomponentdistribution}(1).
\end{proof}

\paragraph{Blue}
In this section, we assume all the results for $k-1$. These results appear in grey in the following picture. We want to prove  \ref{AAThinvarianteigenvalue}, \ref{AATheoremcaraceigenstructure} and \ref{AAThcomponentdistribution}(1) for $k$.
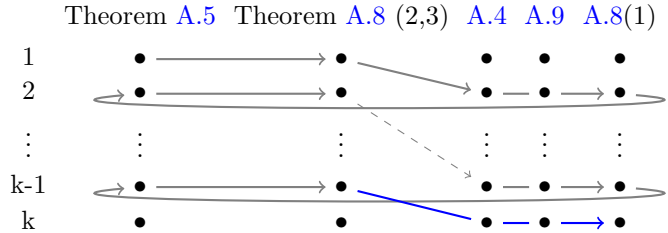
\begin{center}
\begin{tikzpicture}
  \matrix (magic) [matrix of nodes]
  { \ & Theorem \ref{AAInvariantth} & Theorem \ref{AAThcomponentdistribution} (2,3)    &   \ref{AAThinvarianteigenvalue} & \ref{AATheoremcaraceigenstructure} & \ref{AAThcomponentdistribution}(1) \\
   1  & $\bullet$   & $\bullet$    & $\bullet$ &  $\bullet$   & $\bullet$\\
   2  & $\bullet$   & $\bullet$ & $\bullet$  &  $\bullet$&  $\bullet$\\
   \vdots & \vdots   & \vdots  & \vdots & \vdots  & \vdots\\
    k-1  & $\bullet$  & $\bullet$  & $\bullet$  &  $\bullet$ &  $\bullet$ \\
   k  & $\bullet$ & $\bullet$   &  $\bullet$ &  $\bullet$ &  $\bullet$ \\
};
  \draw[thick,gray,->] (magic-2-2) -- (magic-2-3) ;
   \draw[thick,gray,->] (magic-2-3) -- (magic-3-4);
   \draw[thick,gray,->] (magic-3-4) -- (magic-3-5) -- (magic-3-6)   ;
    
\draw[thick,gray,->] (magic-3-6) to [out=-5,in=-175,looseness=1](magic-3-2);
\draw[thick,gray,->] (magic-3-2) -- (magic-3-3);
\draw[dashed,gray,->] (magic-3-3) -- (magic-5-4);
\draw[thick,gray,->]  (magic-5-4) -- (magic-5-5) -- (magic-5-6);
\draw[thick,gray,->] (magic-5-6) to [out=-5,in=-175,looseness=1](magic-5-2);
\draw[thick,gray,->] (magic-5-2) -- (magic-5-3) ;
\draw[thick,blue,->]  (magic-5-3) -- (magic-6-4)  -- (magic-6-5) -- (magic-6-6);
\end{tikzpicture}
\end{center}

\noindent First, we prove the Invariant Eigenvalue Theorem.
\begin{proof}{\textbf{Theorem \ref{AAThinvarianteigenvalue}}}  \label{AAproofThinvarianteigenvaluek} 
\noindent We can assume the grey results in the following picture as proven.
\begin{center}
\begin{tikzpicture}
  \matrix (magic) [matrix of nodes]
  { \ & Theorem \ref{AAInvariantth} & Theorem \ref{AAThcomponentdistribution} (2,3)    &   \ref{AAThinvarianteigenvalue} & \ref{AATheoremcaraceigenstructure} & \ref{AAThcomponentdistribution}(1) \\
   1  & $\bullet$   & $\bullet$    & $\bullet$ &  $\bullet$   & $\bullet$\\
   2  & $\bullet$   & $\bullet$ & $\bullet$  &  $\bullet$&  $\bullet$\\
   \vdots & \vdots   & \vdots  & \vdots & \vdots  & \vdots\\
    k-1  & $\bullet$  & $\bullet$  & $\bullet$  &  $\bullet$ &  $\bullet$ \\
   k  & $\bullet$ & $\bullet$   &  $\bullet$ &  $\bullet$ &  $\bullet$ \\
};
  \draw[thick,gray,->] (magic-2-2) -- (magic-2-3) ;
   \draw[thick,gray,->] (magic-2-3) -- (magic-3-4);
   \draw[thick,gray,->] (magic-3-4) -- (magic-3-5) -- (magic-3-6) ;
    
\draw[thick,gray,->] (magic-3-6) to [out=-5,in=-175,looseness=1](magic-3-2);
\draw[thick,gray,->] (magic-3-2) -- (magic-3-3);
\draw[dashed,gray,->] (magic-3-3) -- (magic-5-4);
\draw[thick,gray,->]  (magic-5-4) -- (magic-5-5) -- (magic-5-6);
\draw[thick,gray,->] (magic-5-6) to [out=-5,in=-175,looseness=1](magic-5-2);
\draw[thick,gray,->] (magic-5-2) -- (magic-5-3) ;
\draw[thick,blue,->]  (magic-5-3) -- (magic-6-4)  ;
\draw[dashed,blue,->]   (magic-6-4)  -- (magic-6-5) -- (magic-6-6);
\end{tikzpicture}
\end{center}
The proof for $k$ is the same as the proof for $k=2$ with a small negligible error. We present the proof for $\hat{\theta}_{P_{k},s}-\hat{\theta}_{P_{k-1},s}$, where the last added perturbation is of order $\theta_k$ and $\theta_1>\theta_2>...>\theta_{k}$. Similar computations can be done to demonstrate the result when the last added perturbation is of order $\theta_r$, $r \neq s$. \\
By using Theorem \ref{AAconvergence}, \ref{AAThcomponentdistribution}(1) for $k-1$ and using the fact that the $p_i$ are different in Assumption \ref{AAAss=theta}(A4),\\ 
\scalebox{0.77}{
\begin{minipage}{1\textwidth}
\begin{eqnarray*}
&&\hspace{-0.7cm} \sum_{i=k}^m \frac{\hat{\lambda}_{P_{k-1},i}}{\hat{\theta}_{P_{k},s}-\hat{\lambda}_{P_{k-1},i}} \hat{u}_{P_{k-1},i,k}^2 
+ \frac{\hat{\theta}_{P_{k-1},s}}{\hat{\theta}_{P_{k},s}-\hat{\theta}_{P_{k-1},s}} \hat{u}_{P_{k-1},s,k}^2 + \sum_{\underset{i \neq s}{i=1}}^k  \frac{\hat{\theta}_{P_{k-1},i}}{\hat{\theta}_{P_{k},s}-\hat{\theta}_{P_{k-1},i}} \hat{u}_{P_{k-1},i,k}^2 =\frac{1}{\theta_k-1}\\ 
&&\hspace{-0.7cm}  \Rightarrow \sum_{i=k-1}^m \frac{\hat{\lambda}_{P_{k-1},i}}{\hat{\theta}_{P_{k},s}-\hat{\lambda}_{P_{k-1},i}} \hat{u}_{P_{k-1},i,k}^2 
+ \frac{\hat{\theta}_{P_{k-1},s}}{\hat{\theta}_{P_{k},s}-\hat{\theta}_{P_{k-1},s}} \hat{u}_{P_{k-1},s,k}^2 +O_p\left( \frac{1}{m \underset{\underset{i \neq s}{i=1,2,3,...,k-1}}{\min}\left(\theta_s-\theta_i \right)} \right) =\frac{1}{\theta_k-1}
\end{eqnarray*}
\end{minipage}}\\

Therefore,  \\ 
\scalebox{0.8}{
\begin{minipage}{1\textwidth}
\begin{eqnarray*}
\frac{\hat{\theta}_{P_{k-1},s}}{\hat{\theta}_{P_{k},s}-\hat{\theta}_{P_{k-1},s}} \hat{u}_{P_{k-1},s,k}^2 &=& -\sum_{i=k}^m \frac{\hat{\lambda}_{P_{k-1},i}}{\hat{\theta}_{P_{k},s}-\hat{\lambda}_{P_{k-1},i}} \hat{u}_{P_{k-1},i,k}^2 + \frac{1}{\theta_k-1} +O_p\left( \frac{1}{m \theta_s} \right) \\
&\overset{1*}{=}& - \frac{1}{\hat{\theta}_{P_{k},s}} \left( 1+ O_p\left( \frac{1}{\sqrt{m}}\right) \right) + \frac{1}{\theta_k-1}  +O_p\left( \frac{1}{m \theta_s} \right) + O_p\left( \frac{1}{\theta_s^2} \right) \\
&=& -\frac{\theta_k -1 -\hat{\theta}_{P_{k},s}}{\hat{\theta}_{P_{k},s} \left(\theta_k-1 \right)} +O_p\left( \frac{1}{m \theta_s} \right) +O_p\left( \frac{1}{\sqrt{m} \theta_s} \right) + O_p\left( \frac{1}{\theta_s^2} \right) ,
\end{eqnarray*}
\end{minipage}}\\

where 1* is true because
\begin{eqnarray*}
\sum_{i=k}^m \hat{\lambda}_{P_{k-1},i} \hat{u}_{P_{k-1},i,k}^2 &=& \sum_{i=1}^m \hat{\lambda}_{P_{k-1},i} \hat{u}_{P_{k-1},i,k}^2 -  \sum_{i=1}^{k-1}\hat{\theta}_{P_{k-1},i} \hat{u}_{P_{k-1},i,k}^2\\
&=& \hat{\Sigma}_{P_{k-1},k,k}-  \sum_{i=1}^{k-1}\hat{\theta}_{P_{k-1},i} \hat{u}_{P_{k-1},i,k}^2\\
&=& W_{k,k}-  \sum_{i=1}^{k-1}\hat{\theta}_{P_{k-1},i} \hat{u}_{P_{k-1},i,k}^2\\
&=&1+O_p\left(\frac{1}{\sqrt{m}}\right).
\end{eqnarray*}
The last line is obtained because the canonical perturbation $P_{k-1}$ does not affect $ W_{k:m,k:m}$. Moreover, $W$ satisfies Assumption \ref{AAAss=matrice}; therefore, $W_{k,k}=1+O_p\left(1/\sqrt{m}\right)$. On the other hand, the second term $ \sum_{i=1}^{k-1}\hat{\theta}_{P_{k-1},i} \hat{u}_{P_{k-1},i,k}^2=O_p\left(1/m\right) $ by Theorem \ref{AAThcomponentdistribution}(2) for $k-1$.\\

\noindent Thus, by Theorem \ref{AAThcomponentdistribution}(2) for $k-1$,\\ 
\scalebox{0.77}{
\begin{minipage}{1\textwidth}
\begin{eqnarray*}
&&\hspace{-0.5cm} \left(1 +O_p\left( \frac{\theta_k}{\theta_s(\theta_k-\theta_s) } \right) + O_p\left( \frac{\theta_k}{\sqrt{m}(\theta_k-\theta_s)} \right) + O_p\left( \frac{\theta_s \theta_k}{m(\theta_k-\theta_s)(\theta_{k-1}-\theta_s)} \right)\right) \left(\hat{\theta}_{P_{k},s}- \hat{\theta}_{P_{k-1},s}\right) \\
&&\hspace{1cm}= 
-\frac{\hat{\theta}_{P_{k-1},s} \hat{\theta}_{P_{k},s} (\theta_k-1)}{\theta_k-1 -\hat{\theta}_{P_{k},s}} \hat{u}_{P_{k-1},s,k}^2\\
&&\hspace{1cm}= O_p\left( \frac{\theta_s \theta_k}{m(\theta_k-\theta_s)} \right).
\end{eqnarray*}
\end{minipage}}\\

and
\begin{eqnarray*}
&&\hat{\theta}_{P_{k},s}- \hat{\theta}_{P_{k-1},s} \overset{\scalebox{0.5}{order}}{\sim} \frac{\min\left(\theta_s, \theta_k \right)}{m}.
\end{eqnarray*}
More precisely we can write\\ 
\scalebox{0.85}{
\begin{minipage}{1\textwidth}
\begin{eqnarray*}
 \hat{\theta}_{P_{k},s}- \hat{\theta}_{P_{k-1},s} &=& 
-\frac{\hat{\theta}_{P_{k-1},s} \hat{\theta}_{P_{k},s} (\theta_k-1)}{\theta_k-1 -\hat{\theta}_{P_{k},s}} \hat{u}_{P_{k-1},s,k}^2 + O_p \left(\frac{1}{m} \right)+ O_p \left(\frac{\min(\theta_s,\theta_k)}{m^{3/2}} \right).
\end{eqnarray*}
\end{minipage}}\\

The $\min$ function can be simplified in our case $\theta_k< \theta_s$; however the above notation is more easily generalized.\\
\noindent Each step of the computation can be done assuming that the last applied perturbation is $\theta_r$ instead of $\theta_k$ for $r=1,2,...,k$. Moreover, in this case, similar computations lead to $\hat{\theta}_{P_k,s}-\hat{\theta}_{P_{-r},s}$ where $s=1,2,...,k$, $s\neq r$. We use the notation 
\begin{eqnarray*}
&&P_{-r}= \I_m+\sum_{\underset{i\neq r}{i=1}}^k (\theta_i-1) e_i e_i^t.
\end{eqnarray*}
Therefore,  for $s\neq r \in \lbrace 1,2,...,k \rbrace$ we obtain the general result.
\begin{itemize}
\item If $\theta_s> \theta_r$, then \\ 
\scalebox{0.9}{
\begin{minipage}{1\textwidth}
\begin{eqnarray*}
 \hat{\theta}_{P_{k},s}-\hat{\theta}_{P_{-r},s}
&=& 
-\frac{\hat{\theta}_{P_{-r},s} \hat{\theta}_{P_{k},s} (\theta_r-1)}{\theta_r-1 -\hat{\theta}_{P_{k},s}} \hat{u}_{P_{-r},s,r}^2 + O_p \left(\frac{1}{m} \right)+ O_p \left(\frac{\theta_r}{m^{3/2}} \right) \\
&\overset{\scalebox{0.5}{order}}{\sim} & \frac{\theta_r}{m}.
\end{eqnarray*}
\end{minipage}}\\

\item If $\theta_s < \theta_r$, then \\ 
\scalebox{0.8}{
\begin{minipage}{1\textwidth}
\begin{eqnarray*}
 \hat{\theta}_{P_{k},s}-\hat{\theta}_{P_{-r},s-1}
&=& 
-\frac{\hat{\theta}_{P_{-r},s-1} \hat{\theta}_{P_{k},s} (\theta_r-1)}{\theta_r-1 -\hat{\theta}_{P_{k},s}} \hat{u}_{P_{-r},s-1,r}^2 + O_p \left(\frac{1}{m} \right)+ O_p \left(\frac{\theta_s}{m^{3/2}} \right) \\
&\overset{\scalebox{0.5}{order}}{\sim} & \frac{\theta_s}{m}.
\end{eqnarray*}
\end{minipage}}\\

\noindent Finally, we obtain for $s>1$,
$$ \hat{\theta}_{P_{k},s} - \hat{\theta}_{\tilde{P}_{s},1} \overset{\scalebox{0.5}{order}}{\sim} \frac{\theta_s}{m} $$
and for $s=1$,
$$ \hat{\theta}_{P_{k},1} - \hat{\theta}_{\tilde{P}_{1},1} \overset{\scalebox{0.5}{order}}{\sim} \frac{\theta_2}{m} .$$
\end{itemize}
\end{proof}

\noindent Next, we prove the characterization of eigenvalues and eigenvectors.
\begin{proof}{\textbf{Theorem \ref{AATheoremcaraceigenstructure}}}  \label{AAproofTheoremcaraceigenstructurek}
\noindent To obtain the result we can assume the grey results in the following picture as proven.
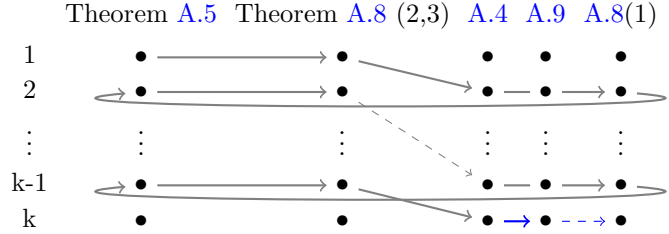
\begin{center}
\begin{tikzpicture}
  \matrix (magic) [matrix of nodes]
 { \ & Theorem \ref{AAInvariantth} & Theorem \ref{AAThcomponentdistribution} (2,3)    &   \ref{AAThinvarianteigenvalue} & \ref{AATheoremcaraceigenstructure} & \ref{AAThcomponentdistribution}(1) \\
   1  & $\bullet$   & $\bullet$    & $\bullet$ &  $\bullet$   & $\bullet$\\
   2  & $\bullet$   & $\bullet$ & $\bullet$  &  $\bullet$&  $\bullet$\\
   \vdots & \vdots   & \vdots  & \vdots & \vdots  & \vdots\\
    k-1  & $\bullet$  & $\bullet$  & $\bullet$  &  $\bullet$ &  $\bullet$ \\
   k  & $\bullet$ & $\bullet$   &  $\bullet$ &  $\bullet$ &  $\bullet$ \\
};
  \draw[thick,gray,->] (magic-2-2) -- (magic-2-3) ;
   \draw[thick,gray,->] (magic-2-3) -- (magic-3-4);
   \draw[thick,gray,->] (magic-3-4) -- (magic-3-5) -- (magic-3-6)   ;
    
\draw[thick,gray,->] (magic-3-6) to [out=-5,in=-175,looseness=1](magic-3-2);
\draw[thick,gray,->] (magic-3-2) -- (magic-3-3);
\draw[dashed,gray,->] (magic-3-3) -- (magic-5-4);
\draw[thick,gray,->]  (magic-5-4) -- (magic-5-5) -- (magic-5-6) ;
\draw[thick,gray,->] (magic-5-6) to [out=-5,in=-175,looseness=1](magic-5-2);
\draw[thick,gray,->] (magic-5-2) -- (magic-5-3) ;
\draw[thick,gray,->] (magic-5-3) -- (magic-6-4) ;
\draw[thick,blue,->]  (magic-6-4) -- (magic-6-5)  ;
\draw[dashed,blue,->]    (magic-6-5) -- (magic-6-6);
\end{tikzpicture}
\end{center}
\noindent The initialisation of the induction, $k=2$, follows from Theorem \ref{AAconvergence}, \ref{AAThinvarianteigenvalue} for $k=2$ and \ref{AAThcomponentdistribution} for $k=1$. Therefore,  we directly prove the result for $k$. \\
Assumption \ref{AAAss=theta}(A4) implies that we have two groups of eigenvalues composing the perturbation. The first group is finite with bounded eigenvalues and the second group has proportional eigenvalues tending to infinity.\\
In order to do a general proof we need to discuss the notation.
\begin{Not}\  \label{AANot=thetadifferent}
\begin{itemize}
\item Usually we assume $\theta_1>\theta_2>...>\theta_k$ such that $\hat{\theta}_{P_k,s}$, the $s^{\rm th}$ largest eigenvalue of $\hat{\Sigma}_{P_k}$ corresponds to $\theta_s$.\\
In this proof we relax the order $\theta_1>\theta2>...>\theta_k$ to do a general proof. The order of $\theta_s$ among the eigenvalues $\theta_1,\theta_2,...,\theta_t$, $t\geqslant s$ is $\text{rank}_t(\theta_s)=r_{t,s}$. Therefore,  assuming a perturbation $P_t$, $\theta_s$ corresponds to the $r_{t,s}^{\text{th}}$ largest eigenvalue of $\hat{\Sigma}_{P_t}$. In order to use simple notation, we again call this corresponding estimated eigenvalue, $\hat{\theta}_{P_r,s}$. \\
Moreover, we change the notation for the eigenvector. In this theorem, for $i=1,2,...,r$, $\hat{u}_{P_r,s}$ is the eigenvector corresponding to $\hat{\theta}_{P_r,s}$.
\item We assume two groups of eigenvalues of size $k_1$ and $k-k_1$ such that these groups satisfy Assumption \ref{AAAss=theta}(A4). Moreover, $\theta_1$ is supposed to be in the first group. We say that the groups are of order $\theta_1$ and $\theta$, respectively, such that only one of them tends to infinity.
\end{itemize}
\end{Not}
\noindent Using this new notation we can without loss of generality construct the proof for $\hat{u}_{P_k,1}$. Note that $\theta_1$ is no longer the largest eigenvalue.

\begin{enumerate}
\item[(a),(h)]  By Cauchy-Schwarz and using $\rho_1= \E\left[\hat{\theta}_{P_k,1}\right]$, \\
\scalebox{0.75}{
\begin{minipage}{1\textwidth}
$$\left| \sum_{i=k}^m \frac{\hat{\lambda}_{P_{k-1},i}}{\hat{\theta}_{P_k,1}-\hat{\lambda}_{P_{k-1},i}} \hat{u}_{P_{k-1},i,1} \hat{u}_{P_{k-1},i,k} \right| \leqslant \sum_{i=k}^m \frac{\hat{\lambda}_{P_{k-1},i}}{\rho_1-\hat{\lambda}_{P_{k-1},i}} \hat{u}_{P_{k-1},i,1} \hat{u}_{P_{k-1},i,k}+O_p\left( \frac{1}{\theta_1^{3/2}m^{1/2}} \right).$$
\vspace{0.1cm}
\end{minipage}}\\
Some prerequisite results are easily proven using theorems for $k-1$: \\ 
\scalebox{0.85}{
\begin{minipage}{1\textwidth}
\begin{eqnarray*}
&& \hat{u}_{P_{k-1},i,1} \text{ s.t. } \sum_{i=k}^m \hat{u}_{P_{k-1},i,1}^2 = O_p\left( \frac{1}{\theta_1} \right) \text{ (By Theorem \ref{AAThcomponentdistribution} part 3)},\\
&& \hat{u}_{P_{k-1},i,k} \sim \rm{RV}\left(0,1/m \right),\\
&&\E\left[\hat{\lambda}_{P_{k-1},i}\hat{u}_{P_{k-1},i,1}  \hat{u}_{P_{k-1},i,k}  \right]=0, \text{ (By invariance under rotation),}\\
&& \var \left( \sum_{i=k}^m  \hat{u}_{P_{k-1},i,1} \hat{u}_{P_{k-1},i,k} \right)
=\var \left( \sum_{i=1}^{k-1} \hat{u}_{P_{k-1},i,1} \hat{u}_{P_{k-1},i,k} \right)
= O_p\left( \frac{1}{\theta_1 m} \right).
\end{eqnarray*}
\end{minipage}}\\

This leads to \\ 
\scalebox{0.8}{
\begin{minipage}{1\textwidth}
\begin{eqnarray*}
 \sum_{i=k}^m  \var \left( \hat{u}_{P_{k-1},i,1} \hat{u}_{P_{k-1},i,k} \right) &=& \sum_{i=k}^m \E \left[ \hat{u}_{P_{k-1},i,1}^2 \hat{u}_{P_{k-1},i,k}^2 \right]\\
 &=& \frac{1}{m-k+1} \sum_{i=k}^m \E \left[ \hat{u}_{P_{k-1},i,1}^2 \left(1- \sum_{s=1}^{k-1} \hat{u}_{P_{k-1},i,s}^2 \right) \right]\\
 &=&O_p\left(\frac{1}{\theta_1 m} \right).
\end{eqnarray*}
\end{minipage}}\\

In order to obtain the order of the size, we use the last part of Theorem \ref{AAThcomponentdistribution}. Either the perturbation in direction $e_1$ is finite and the result follows directly, or the perturbation tends to infinity and we can separate the perturbations into two groups, one finite and the other one tending to infinity. The last result of  Theorem \ref{AAThcomponentdistribution} gives the required estimate.\\
\scalebox{0.75}{
\begin{minipage}{1\textwidth}
\begin{eqnarray*}
&&\hspace{-0.5cm}\var \left( \sum_{i=k}^m \frac{\hat{\lambda}_{P_{k-1},i}}{\hat{\theta}_{P_k,1}-\hat{\lambda}_{P_{k-1},i}} \hat{u}_{P_{k-1},i,1} \hat{u}_{P_{k-1},i,k} \right)\\
&&\hspace{0.5cm}=
 \sum_{i=k}^m  \var \left( \frac{\hat{\lambda}_{P_{k-1},i}}{\hat{\theta}_{P_k,1}-\hat{\lambda}_{P_{k-1},i}} \hat{u}_{P_{k-1},i,1} \hat{u}_{P_{k-1},i,k} \right)\\
 &&\hspace{1.5cm}+
 \sum_{i\neq j =k}^m \cov\left(\frac{\hat{\lambda}_{P_{k-1},i}}{\hat{\theta}_{P_k,1}-\hat{\lambda}_{P_{k-1},i}}\hat{u}_{P_{k-1},i,1} \hat{u}_{P_{k-1},i,k},\frac{\hat{\lambda}_{P_{k-1},j}}{\hat{\theta}_{P_k,1}-\hat{\lambda}_{P_{k-1},j}} \hat{u}_{P_{k-1},j,1} \hat{u}_{P_{k-1},j,k} \right)\\
&&\hspace{0.5cm}=A+B.
\end{eqnarray*}
\end{minipage}}\\

The parts $A$ and $B$ are studied separately. By Assumption \ref{AAAss=matrice}, $\hat{\lambda}_{P_{k-1},k}$ is bounded by a constant $\lambda$. 
\begin{eqnarray*}
A&=&\sum_{i=k}^m  \var \left(  \frac{\hat{\lambda}_{P_{k-1},i}}{\rho_1-\hat{\lambda}_{P_{k-1},i}} \hat{u}_{P_{k-1},i,1} \hat{u}_{P_{k-1},i,k} \right)\\
&=&\sum_{i=k}^m  \E \left[  \left(\frac{\hat{\lambda}_{P_{k-1},i}}{\rho_1-\hat{\lambda}_{P_{k-1},i}}\right)^2 \hat{u}_{P_{k-1},i,1}^2 \hat{u}_{P_{k-1},i,k}^2 \right]\\
&\leqslant& \left(\frac{\lambda}{\rho_1 -\lambda}\right)^2 \sum_{i=k}^m  \E \left[   \hat{u}_{P_{k-1},i,1}^2 \hat{u}_{P_{k-1},i,k}^2 \right]\\
&=& \left(\frac{\lambda}{\rho_1 -\lambda}\right)^2 \sum_{i=k}^m  \var \left( \hat{u}_{P_{k-1},i,1} \hat{u}_{P_{k-1},i,k} \right)\\
&=&O\left(\frac{1}{\theta_1^3 m} \right).
\end{eqnarray*}\\ 
\scalebox{0.65}{
\begin{minipage}{1\textwidth}
\begin{eqnarray*}
|B|&=& \left| \sum_{i\neq j =k}^m \cov\left(\frac{\hat{\lambda}_{P_{k-1},i}}{\rho_1-\hat{\lambda}_{P_{k-1},i}} \hat{u}_{P_{k-1},i,1} \hat{u}_{P_{k-1},i,k},\frac{\hat{\lambda}_{P_{k-1},j}}{\rho_1-\hat{\lambda}_{P_{k-1},j}}  \hat{u}_{P_{k-1},j,1} \hat{u}_{P_{k-1},j,k} \right)\right|\\
&=& \left| \sum_{i\neq j =k}^m \left( \E\left[\frac{\hat{\lambda}_{P_{k-1},i}}{\rho_1-\hat{\lambda}_{P_{k-1},i}} \frac{\hat{\lambda}_{P_{k-1},j}}{\rho_1-\hat{\lambda}_{P_{k-1},j}}\hat{u}_{P_{k-1},i,1}  \hat{u}_{P_{k-1},j,1} \hat{u}_{P_{k-1},i,k} \hat{u}_{P_{k-1},j,k} \right]-0 \right)\right|\\
&=& \left|  \sum_{i\neq j =k}^m \frac{1}{m-k+1}  \E_p\left[\frac{\hat{\lambda}_{P_{k-1},i}}{\rho_1-\hat{\lambda}_{P_{k-1},i}} \frac{\hat{\lambda}_{P_{k-1},j}}{\rho_1-\hat{\lambda}_{P_{k-1},j}}\hat{u}_{P_{k-1},i,1}  \hat{u}_{P_{k-1},j,1} \sum_{r =k}^m \hat{u}_{P_{k-1},i,r} \hat{u}_{P_{k-1},j,r} \right]\right| \\
&=& \left|  \sum_{i\neq j =k}^m \frac{1}{m-k+1}  \E\left[\frac{\hat{\lambda}_{P_{k-1},i}}{\rho_1-\hat{\lambda}_{P_{k-1},i}} \frac{\hat{\lambda}_{P_{k-1},j}}{\rho_1-\hat{\lambda}_{P_{k-1},j}}\hat{u}_{P_{k-1},i,1}  \hat{u}_{P_{k-1},j,1} \sum_{r =1}^{k-1} \hat{u}_{P_{k-1},i,r} \hat{u}_{P_{k-1},j,r} \right]\right| \\
&\leqslant&    \frac{\left(\frac{\lambda}{\rho_1 -\lambda}\right)^2}{m-k} \sum_{r =1}^{k-1}   \E\left[  \sum_{i\neq j =k}^m \left|\hat{u}_{P_{k-1},i,1}  \hat{u}_{P_{k-1},j,1}  \hat{u}_{P_{k-1},i,r} \hat{u}_{P_{k-1},j,r}  \right| \right] \\
&\leqslant &  O\left( \frac{1}{\theta_1^2 m} \right) \sum_{r =1}^{k-1}   \E_p\left[  \left(\sum_{i =k}^m | \hat{u}_{P_{k-1},i,1}\hat{u}_{P_{k-1},i,r}| \right)^2  \right] \\
&\leqslant & O\left( \frac{1}{\theta_1^2 m} \right)  \sum_{r =1}^{k-1}    \E_p\left[  \left( \sum_{i =k}^m  \hat{u}_{P_{k-1},i,1}^2 \right) \left(\sum_{i =k}^m\hat{u}_{P_{k-1},i,r}^2 \right) \right] \\
&=&O\left( \frac{1}{\theta_1^3 m} \right).
\end{eqnarray*}
\end{minipage}}\\

Thus 
\begin{eqnarray*}
\var \left( \sum_{i=k}^m \frac{\hat{\lambda}_{P_{k-1},i}}{\rho_1-\hat{\lambda}_{P_{k-1},i}} \hat{u}_{P_{k-1},i,1} \hat{u}_{P_{k-1},i,k} \right) = O_p\left( \frac{1}{\theta_1^3 m} \right).
\end{eqnarray*}
Therefore,  because the expectation is $0$ by invariance under rotation,
\begin{eqnarray*}
\sum_{i=k}^m \frac{\hat{\lambda}_{P_{k-1},i}}{\hat{\theta}_{P_k,1}-\hat{\lambda}_{P_{k-1},i}} \hat{u}_{P_{k-1},i,1} \hat{u}_{P_{k-1},i,k}
&=& O_p\left(\frac{1}{\theta_1^{3/2} \sqrt{m}} \right) .
\end{eqnarray*}

\item[(b)] We study
$$ \frac{\hat{\theta}_{P_{k-1},1}}{\hat{\theta}_{P_{k},1}-\hat{\theta}_{P_{k-1},1}}  \hat{u}_{P_{k-1},1,1} \hat{u}_{P_{k-1},1,k}. $$
By Theorem \ref{AAInvariantth} and Theorem \ref{AAThcomponentdistribution} for $k-1$, we obtain \\ 
\scalebox{0.8}{
\begin{minipage}{1\textwidth}
\begin{eqnarray*}
&&\text{If }\theta_1 >D, \text{ for $D>0$ fixed }, \exists \  d(D) \text{ such that, }  1> |\hat{u}_{P_{k-1},1,1}| > d(D)  >0 \\
&&\hspace{2cm} \text{ with probability tending to } 1, \\
&&\hat{u}_{P_{k-1},1,k} \overset{\scalebox{0.5}{order}}{\sim} \frac{1}{\sqrt{\theta_1 m}} 
 .
\end{eqnarray*}
\end{minipage}}\\

We see thanks to Theorem \ref{AAThinvarianteigenvalue} for $k$ that 
\begin{eqnarray*}
\frac{\hat{\theta}_{P_{k-1},1}}{\hat{\theta}_{P_{k},1}-\hat{\theta}_{P_{k-1},1}} \overset{\scalebox{0.5}{order}}{\sim} \frac{\theta_1 m}{\min\left(\theta_1,\theta_k \right)}.
\end{eqnarray*}
The result is straightforward.
\item[(c)] We study 
$$\sum_{i=2}^{k-1} \frac{\hat{\theta}_{P_{k-1},i}}{\hat{\theta}_{P_k,1}-\hat{\theta}_{P_{k-1},i}} \hat{u}_{P_{k-1},i,1} \hat{u}_{P_{k-1},i,k}.$$
By Theorem \ref{AAThcomponentdistribution}, 
\begin{eqnarray*}
\hat{u}_{P_{k-1},i,1} &=& O_p\left(\frac{\sqrt{\theta_1 \theta_i}}{\left( \theta_1-\theta_i \right) \sqrt{m}} \right), \\
\hat{u}_{P_{k-1},i,k} &=& O_p \left( \frac{1}{\sqrt{\theta_i} \sqrt{m}}\right).
\end{eqnarray*}
Therefore,  
\begin{eqnarray*}
\frac{\hat{\theta}_{P_{k-1},i}}{\hat{\theta}_{P_k,1}-\hat{\theta}_{P_{k-1},i}} \hat{u}_{P_{k-1},i,1} \hat{u}_{P_{k-1},i,k}&=&  O_p\left(\frac{\sqrt{\theta_1} \theta_i}{\left( \theta_1-\theta_i \right)^2 \sqrt{m}} \right).
\end{eqnarray*}
Studying the different possibilities for $\theta_i$ and  $\theta_1$ leads to the desired result.
\item[(d)] We consider 
$$\sum_{i=k}^m \frac{\hat{\lambda}_{P_{k-1},i}^2}{(\hat{\theta}_{P_k,1}-\hat{\lambda}_{P_{k-1},i})^2} \hat{u}_{P_{k-1},i,k}^2.$$
A straightforward computation leads to \\ 
\scalebox{0.9}{
\begin{minipage}{1\textwidth}
\begin{eqnarray*}
\sum_{i=k}^m \frac{\hat{\lambda}_{P_{k-1},i}^2}{(\hat{\theta}_{P_k,1}-\hat{\lambda}_{P_{k-1},i})^2} \hat{u}_{P_{k-1},i,k}^2 \leqslant    O_p\left(\frac{1}{\theta_1^2} \right) \left(1- \sum_{i=1}^{k-1} \hat{u}_{P_{k-1},i,k}^2 \right)=O_p\left(\frac{1}{\theta_1^2} \right).
\end{eqnarray*}
\end{minipage}}\\

\item[(e)] We study 
$$\frac{\hat{\theta}_{P_{k-1},1}^2}{(\hat{\theta}_{P_{k},1}-\hat{\theta}_{P_{k-1},1})^2}   \hat{u}_{P_{k-1},1,k}^2. $$
By Theorems \ref{AAThinvarianteigenvalue} and \ref{AAThcomponentdistribution},
\begin{eqnarray*}
&&\frac{\hat{\theta}_{P_{k-1},1}^2}{(\hat{\theta}_{P_{k},1}-\hat{\theta}_{P_{k-1},1})^2} \overset{\scalebox{0.5}{order}}{\sim} \frac{\theta_1^2 m^2}{\min\left(\theta_1,\theta_k \right)^2},\\
 &&\hat{u}_{P_{k-1},1,k}^2 \overset{\scalebox{0.5}{order}}{\sim} \frac{1}{m \theta_1} .
\end{eqnarray*}
The result is straightforward.

\item[(f)] We study 
$$\sum_{i=2}^{k-1} \frac{\hat{\theta}_{P_{k-1},i}^2}{(\hat{\theta}_{P_k,1}-\hat{\theta}_{P_{k-1},i})^2}  \hat{u}_{P_{k-1},i,k}^2.$$
By Theorem \ref{AAThcomponentdistribution} and \ref{AAThinvarianteigenvalue}, 
\begin{eqnarray*}
&&\hat{u}_{P_{k-1},i,k}^2=O_p\left( \frac{1}{m \theta_i} \right).\\
\end{eqnarray*}
Then, 
\begin{eqnarray*}
\frac{\hat{\theta}_{P_{k-1},i}^2}{(\hat{\theta}_{P_k,1}-\hat{\theta}_{P_{k-1},i})^2}  \hat{u}_{P_{k-1},i,k}^2 = O_p\left( \frac{\theta_i}{\left( \theta_1 - \theta_i \right)^2 m} \right).
\end{eqnarray*}
Studying the different possibilities for $\theta_i$ and  $\theta_1$ leads to the result.
\item[(g)] The result is obtained directly from Theorem \ref{AAconvergence}.
\item[(h)] The same proof as in (a) leads to the result.
\item[(i)]  We study 
$$\frac{\hat{\theta}_{P_{k-1},1}}{\hat{\theta}_{P_{k},1}-\hat{\theta}_{P_{k-1},1}}  \hat{u}_{P_{k-1},1,s} \hat{u}_{P_{k-1},1,k},$$
for $s=2,...,k-1$.

Using Theorems \ref{AAThinvarianteigenvalue} and \ref{AAThcomponentdistribution} we have
\begin{eqnarray*}
&&  \hat{u}_{P_{k-1},1,s} \hat{u}_{P_{k-1},1,k} 
\overset{\scalebox{0.5}{order}}{\sim} \frac{\min\left( \theta_1,\theta_s \right)}{\theta_1 \sqrt{\theta_s} m},\\
&& \hat{\theta}_{P_{k},1}-\hat{\theta}_{P_{k-1},1} \overset{\scalebox{0.5}{order}}{\sim}  \frac{\min \left( \theta_1, \theta_k \right)}{m} .
\end{eqnarray*}
The result follows directly.

\item[(j)] We consider 
$$\sum_{i=2,\neq s}^{k-1} \frac{\hat{\theta}_{P_{k-1},i}}{\hat{\theta}_{P_k,1}-\hat{\theta}_{P_{k-1},i}} \hat{u}_{P_{k-1},i,s} \hat{u}_{P_{k-1},i,k}.$$
Using the Theorems \ref{AAThinvarianteigenvalue} and \ref{AAThcomponentdistribution}  the result is straightforward as for part (c).

\item[(k)] We study 
$$\frac{\hat{\theta}_{P_{k-1},s}}{\hat{\theta}_{P_k,1}-\hat{\theta}_{P_{k-1},s}} \hat{u}_{P_{k-1},s,s} \hat{u}_{P_{k-1},s,k}.$$
By Theorem \ref{AAInvariantth} and  Lemma \ref{AAThcomponentdistribution}, \\ 
\scalebox{0.85}{
\begin{minipage}{1\textwidth}
\begin{eqnarray*}
&&\text{If }\theta_s >D, \text{ for $D>0$ fixed }, \exists \  d(D) \text{ such that, }  1> |\hat{u}_{P_{k-1},s,s}| > d(D)  >0\\
&&\hspace{2 cm}\text{ with probability tending to } 1, \\
&&\hat{u}_{P_{k-1},s,k} = O_p\left(\frac{1}{\sqrt{\theta_s m}}\right) 
 .
\end{eqnarray*}
\end{minipage}}\\
The result follows.
\end{enumerate}
The link between $\hat{u}_{P_k,1,1} $ and $\tilde{u}_{P_k,1,1} $ is obtained by basic notions of linear algebra and similar estimations of the norm. \\
\noindent We now prove the first point of the remark.
\begin{enumerate}
\item First, we study ${\rm sign} \left( \hat{u}_{P_k,1,1} \right)$ by investigating $ \tilde{u}_{P_k,1,1}$ which was defined in the statement of the theorem. Then, by construction, the results hold for  $\hat{u}_{P_k,1,1}$ because we just rescale $\tilde{u}_{P_k,1}$ to obtain $\hat{u}_{P_k,1}$. The theorem says  \\ 
\scalebox{0.6}{
\begin{minipage}{1\textwidth}
\begin{eqnarray*}  
&&\hspace{-1cm} {\rm sign} \left( \tilde{u}_{P_k,1,1}\right) \\
 && \hspace{-1cm} = {\rm sign} \left( \frac{
\overbrace{ \sum_{i=k}^m \frac{\hat{\lambda}_{P_{k-1},i}}{\hat{\theta}_{P_k,1}-\hat{\lambda}_{P_{k-1},i}} \hat{u}_{P_{k-1},i,1} \hat{u}_{P_{k-1},i,k}}^{(a) O_p\left( \frac{1}{\theta_1^{3/2} \sqrt{m}}\right)} + 
\overbrace{ \frac{\hat{\theta}_{P_{k-1},1}}{\hat{\theta}_{P_{k},1}-\hat{\theta}_{P_{k-1},1}}  \hat{u}_{P_{k-1},1,1} \hat{u}_{P_{k-1},1,k}}^{(b) \overset{\scalebox{0.5}{order}}{\sim} \ \frac{\sqrt{\theta_1 m}}{\min\left( \theta_1, \theta_k \right)} } +  
\overbrace{ \sum_{i=2}^{k-1} \frac{\hat{\theta}_{P_{k-1},i}}{\hat{\theta}_{P_k,1}-\hat{\theta}_{P_{k-1},i}} \hat{u}_{P_{k-1},i,1} \hat{u}_{P_{k-1},i,k} }^{(c) O_p\left( \frac{ 1 }{\theta_1^{1/2} m}\right)} }
{\sqrt{\underbrace{\sum_{i=k}^m \frac{\hat{\lambda}_{P_{k-1},i}^2}{(\hat{\theta}_{P_k,1}-\hat{\lambda}_{P_{k-1},i})^2} \hat{u}_{P_{k-1},i,k}^2}_{(d) O_p\left( \frac{1}{\theta_1^2}\right)} + 
\underbrace{\frac{\hat{\theta}_{P_{k-1},1}^2}{(\hat{\theta}_{P_{k},1}-\hat{\theta}_{P_{k-1},1})^2}   \hat{u}_{P_{k-1},1,k}^2 }_{(e)\overset{\scalebox{0.5}{order}}{\sim} \ \frac{\theta_1 m}{ \min \left( \theta_1,\theta_k \right)^2}}+  
\underbrace{\sum_{i=2}^{k-1} \frac{\hat{\theta}_{P_{k-1},i}^2}{(\hat{\theta}_{P_k,1}-\hat{\theta}_{P_{k-1},i})^2}  \hat{u}_{P_{k-1},i,k}^2}_{(f) O_p\left( \frac{1}{\theta_1 m}\right) }}}\right). \hspace{20cm}
\end{eqnarray*}
\end{minipage}}\\

\noindent The first convergence is directly obtained from \\ 
\scalebox{0.7}{
\begin{minipage}{1\textwidth}
\begin{eqnarray*}
{\rm sign} \left( \tilde{u}_{P_k,1,1}\right) &=& {\rm sign} \left(  \frac{\hat{\theta}_{P_{k-1},1}}{\hat{\theta}_{P_{k},1}-\hat{\theta}_{P_{k-1},1}}  \hat{u}_{P_{k-1},1,1} \hat{u}_{P_{k-1},1,k} + O_p\left( \frac{1}{\theta_1^{1/2}m} \right)+O_p\left( \frac{1}{\theta_1^{3/2} \sqrt{m}}\right) \right).
\end{eqnarray*}
\end{minipage}}\\

Using Theorem \ref{AAThinvarianteigenvalue} and assuming $m$ and $\theta_1$ sufficiently large lead to
\begin{eqnarray*}
{\rm sign} \left( \tilde{u}_{P_k,1,1}\right) &=& {\rm sign} \left( \left(\hat{\theta}_{P_{k},1}-\hat{\theta}_{P_{k-1},1}\right)  \hat{u}_{P_{k-1},1,k}  \right)\\
&=&{\rm sign} \left( \left(\theta_1-\theta_k\right)  \hat{u}_{P_{k-1},1,k}\right).
\end{eqnarray*}

\item The second remark supposes a perturbation of order $k=2$. We already know the behaviour of the first eigenvector. In order to obtain results for the second vector, we need to replace $\hat{\theta}_{P_{k},1}$ by $\hat{\theta}_{P_{k},2}$ in the formula and the order size changes. Similar arguments as above lead to the result.
\end{enumerate}

\end{proof}
\noindent For the last part of the proof of the blue part in the Figure \ref{AAfig=proof}, we study the first point of the component Theorem \ref{AAThcomponentdistribution}.
\begin{proof}{\textbf{Theorem \ref{AAThcomponentdistribution}}}  \label{AAproofThcomponentdistribution1k}
\noindent To prove this result we can assume the grey results in the following picture as proven.
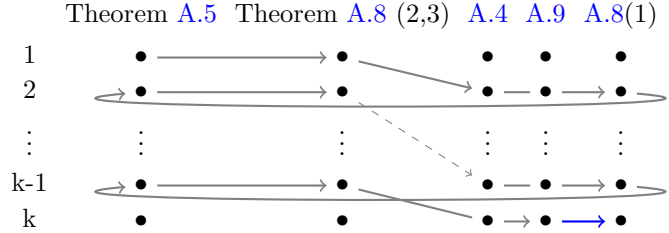
\begin{center}
\begin{tikzpicture}
  \matrix (magic) [matrix of nodes]
 { \ & Theorem \ref{AAInvariantth} & Theorem \ref{AAThcomponentdistribution} (2,3)    &   \ref{AAThinvarianteigenvalue} & \ref{AATheoremcaraceigenstructure} & \ref{AAThcomponentdistribution}(1) \\
   1  & $\bullet$   & $\bullet$  &  $\bullet$   &  $\bullet$   & $\bullet$\\
   2  & $\bullet$   & $\bullet$ &  $\bullet$  &  $\bullet$&  $\bullet$\\
   \vdots & \vdots   & \vdots & \vdots  & \vdots  & \vdots\\
    k-1  & $\bullet$  & $\bullet$ &  $\bullet$   &  $\bullet$ &  $\bullet$ \\
   k  & $\bullet$ & $\bullet$  &  $\bullet$  &  $\bullet$ &  $\bullet$ \\
};
  \draw[thick,gray,->] (magic-2-2) -- (magic-2-3) ;
   \draw[thick,gray,->] (magic-2-3) -- (magic-3-4);
   \draw[thick,gray,->] (magic-3-4) -- (magic-3-5) -- (magic-3-6)  ;
    \draw[thick,gray,->] (magic-3-6) to [out=-5,in=-175,looseness=1](magic-3-2);
\draw[thick,gray,->] (magic-3-2) -- (magic-3-3);
\draw[dashed,gray,->] (magic-3-3) -- (magic-5-4);
\draw[thick,gray,->]  (magic-5-4) -- (magic-5-5) -- (magic-5-6);
\draw[thick,gray,->] (magic-5-6) to [out=-5,in=-175,looseness=1](magic-5-2);
\draw[thick,gray,->] (magic-5-2) -- (magic-5-3) ;
\draw[thick,gray,->] (magic-5-3) -- (magic-6-4)-- (magic-6-5);
\draw[thick,blue,->]  (magic-6-5)-- (magic-6-6);
\end{tikzpicture}
\end{center}
This proof computes $\hat{u}_{P_k,1,k}$, but the method can be used to study any components $\hat{u}_{P_k,s,t}$ where $s \neq t \in \lbrace1,2,...,k \rbrace $. In order to extend it we must use Notation \ref{AANot=thetadifferent}. First we assume the convention of Theorem \ref{AATheoremcaraceigenstructure}, $\hat{u}_{P_k,1,k}>0$.\\ 
\scalebox{0.7}{
\begin{minipage}{1\textwidth} \begin{eqnarray*} 
\left\langle \tilde{u}_{P_k,1},e_k \right\rangle \hspace{-1cm}&&\\
 &=& \frac{\sum_{i=k}^m \frac{\hat{\lambda}_{P_{k-1},i}}{\hat{\theta}_{P_k,1}-\hat{\lambda}_{P_{k-1},i}} \hat{u}_{P_{k-1},i,k}^2 +
\frac{\hat{\theta}_{P_{k-1},1}}{\hat{\theta}_{P_{k},1}-\hat{\theta}_{P_{k-1},1}}  \hat{u}_{P_{k-1},1,k}^2  + 
\sum_{i=2}^{k-1} \frac{\hat{\theta}_{P_{k-1},i}}{\hat{\theta}_{P_k,1}-\hat{\theta}_{P_{k-1},i}} \hat{u}_{P_{k-1},i,k}^2  }
{\sqrt{\underbrace{ \sum_{i=k}^m \frac{\hat{\lambda}_{P_{k-1},i}^2}{(\hat{\theta}_{P_k,1}-\hat{\lambda}_{P_{k-1},i})^2} \hat{u}_{P_{k-1},i,k}^2}_{O_p\left(\frac{1}{\theta_1^2}\right)} +\underbrace{ 
\frac{\hat{\theta}_{P_{k-1},1}^2}{(\hat{\theta}_{P_{k},1}-\hat{\theta}_{P_{k-1},1})^2}   \hat{u}_{P_{k-1},1,k}^2 }_{\ \overset{\scalebox{0.5}{order}}{\sim} \frac{\theta_1 m}{\min\left(\theta_1,\theta_k \right)^2}}+  
\underbrace{\sum_{i=2}^{k-1} \frac{\hat{\theta}_{P_{k-1},i}^2}{(\hat{\theta}_{P_k,1}-\hat{\theta}_{P_{k-1},i})^2}  \hat{u}_{P_{k-1},i,k}^2}_{O_p\left(\frac{1}{\theta_1 m}\right)}}}\\
&=& \frac{1}{\theta_{k}-1} \frac{1}{\sqrt{\frac{\hat{\theta}_{P_{k-1},1}^2}{(\hat{\theta}_{P_{k},1}-\hat{\theta}_{P_{k-1},1})^2}   \hat{u}_{P_{k-1},1,k}^2}}+ O_p\left(\frac{\min(\theta_1,\theta_k)^{3}}{\theta_k \theta_1^{7/2} m^{3/2}} \right)+ O_p\left(\frac{\min(\theta_1,\theta_k)^{3}}{\theta_k \theta_1^{5/2} m^{5/2}} \right)\\
&=&\frac{1}{\theta_{k}-1}\frac{|\hat{\theta}_{P_{k},1}-\hat{\theta}_{P_{k-1},1}|}{|\hat{\theta}_{P_{k-1},1}|  |\hat{u}_{P_{k-1},1,k}|}   +O_p\left(\frac{\min(\theta_1,\theta_k)^{3}}{\theta_k \theta_1^{7/2} m^{3/2}} \right)+ O_p\left(\frac{\min(\theta_1,\theta_k)^{3}}{\theta_k \theta_1^{5/2} m^{5/2}} \right).\\
\end{eqnarray*}
\end{minipage}}\\

Then, $ \hat{u}_{P_k,1} = P_k^{1/2} \tilde{u}_{P_k,1}/N_1$, and
\begin{eqnarray*}
N_1^2&=& \sum_{i=1}^{k-1}   \tilde{u}_{P_k,i}^2+ \theta_k \tilde{u}_{P_k,k}^2 +\sum_{i=k+1}^m \tilde{u}_{P_k,i}^2\\
&=& 1+(\theta-1)\tilde{u}_{P_k,k}^2.
\end{eqnarray*}
We also know by Theorem \ref{AAThinvarianteigenvalue} that \\ 
\scalebox{0.95}{
\begin{minipage}{1\textwidth} 
\begin{eqnarray*}
 \hat{\theta}_{P_{k},1}- \hat{\theta}_{P_{k-1},1} &=& 
-\frac{\hat{\theta}_{P_{k-1},1} \hat{\theta}_{P_{k},1} (\theta_k-1)}{\theta_k-1 -\hat{\theta}_{P_{k},1}} \hat{u}_{P_{k-1},1,k}^2 + O_p \left(\frac{1}{m} \right)+ O_p \left(\frac{\theta_1}{m^{3/2}} \right).
\end{eqnarray*}
\end{minipage}}\\

and
$$ \hat{\theta}_{P_{k-1},1}- \hat{\theta}_{P_{k},1} = O_p \left(\frac{\min(\theta_1,\theta_k)}{m} \right) 
.$$
Therefore,  Theorem \ref{AATheoremcaraceigenstructure} and  \ref{AAThinvarianteigenvalue} for $k$ leads to \\ 
\scalebox{0.75}{
\begin{minipage}{1\textwidth}
\begin{eqnarray*} 
\hat{u}_{P_k,1,k} &=&\frac{\left\langle \tilde{u}_{P_k,1},e_k \right\rangle}{{\rm Norm}} \sqrt{\theta_k} \\
&=&\left(\underbrace{\frac{1}{\theta_{k}-1}\frac{|\hat{\theta}_{P_{k},1}-\hat{\theta}_{P_{k-1},1}|}{|\hat{\theta}_{P_{k-1},1}|  |\hat{u}_{P_{k-1},1,k}|}}_{\ \overset{\scalebox{0.5}{order}}{\sim} \frac{\min(\theta_1,\theta_k)}{\theta_k \theta_1^{1/2} m^{1/2}} } +  O_p\left(\frac{\min(\theta_1,\theta_k)^{3}}{\theta_k \theta_1^{7/2} m^{3/2}} \right)+ O_p\left(\frac{\min(\theta_1,\theta_k)^{3}}{\theta_k \theta_1^{5/2} m^{5/2}} \right)\right) \frac{\sqrt{\theta_k}}{1+O_p\left(\frac{1}{m} \right)}\\
&=&\frac{1}{\sqrt{\theta_{k}}}\frac{|\hat{\theta}_{P_{k},1}-\hat{\theta}_{P_{k-1},1}|}{|\hat{\theta}_{P_{k-1},1}|  |\hat{u}_{P_{k-1},1,k}|}
+ O_p\left( \frac{\min(\theta_1,\theta_k)}{\theta_k^{1/2} \theta_1^{1/2} m^{3/2}} \right)
\\
&=& 
 \frac{\left| -\frac{\hat{\theta}_{P_{k-1},1} \hat{\theta}_{P_{k},1} (\theta_k-1)}{\theta_k-1 -\hat{\theta}_{P_{k},1}} \hat{u}_{P_{k-1},1,k}^2 + O_p \left(\frac{1}{m} \right)+ O_p \left(\frac{\min(\theta_1,\theta_k)}{m^{3/2}} \right) \right|}{\sqrt{\theta_{k}}|\hat{\theta}_{P_{k-1},1}|  |\hat{u}_{P_{k-1},1,k}|}  
+ O_p\left( \frac{\min(\theta_1,\theta_k)}{\theta_k^{1/2} \theta_1^{1/2} m^{3/2}} \right)\\
&=&\frac{\sqrt{\theta_k} \hat{\theta}_{P_{k},1} }{|\theta_k-\hat{\theta}_{P_{k},1}|} |\hat{u}_{P_{k-1},1,k}|+O_p\left(\frac{\min(\theta_1,\theta_k)}{\theta_1^{1/2}\theta_k^{1/2}m} \right)+O_p\left(\frac{1}{ \theta_1^{1/2} \theta_k^{1/2} m^{1/2}} \right).
\end{eqnarray*}
\end{minipage}}\\

\noindent Note that the sign is always positive! We can use the Remark of Theorem \ref{AATheoremcaraceigenstructure} and set $\hat{u}_{P_s,i,i} >0$ for $ s=1,2,...,k$ and $i=1,2,...,s$. Then, the previous result becomes more convenient:\\
\noindent Under the sign condition for the eigenvector,
\begin{eqnarray*} 
\hat{u}_{P_k,1,k} &=&\frac{\sqrt{\theta_k} \theta_1 }{\theta_1-\theta_k} \hat{u}_{P_{k-1},1,k}+O_p\left(\frac{\min(\theta_1,\theta_k)}{\theta_1^{1/2}\theta_k^{1/2}m} \right)+O_p\left(\frac{1}{ \theta_1^{1/2} \theta_k^{1/2} m^{1/2}} \right).
\end{eqnarray*}
\noindent Therefore,  we directly obtain the distribution when $\theta_1,\theta_k \rightarrow \infty$. Using
\begin{eqnarray*}
\hat{\alpha}^2_{P_{k-1},1} &=& \sum_{i=1}^{k-1} \left\langle \hat{u}_{P_{k-1},1},\epsilon_i \right\rangle^2,\\
\hat{\alpha}^2_{P_{k-1},1}&=&\alpha_1^2+O_p\left( \frac{1}{\theta_1 m} \right)=1- \frac{M_2-1}{\theta_1}+O_p\left( \frac{1}{\theta_1^2} \right)+O_p\left( \frac{1}{\theta_1 \sqrt{m}} \right)
\end{eqnarray*}
and the second part of this Theorem \ref{AAThcomponentdistribution} for $k-1$,
\begin{eqnarray*}
\hat{u}_{P_{k-1},1,k}|  \hat{\alpha}^2_{P_{k-1},1} &\overset{Asy}{\sim}& \Normal \left(0, \frac{1-\hat{\alpha}^2_{P_{k-1},1}}{m} \right),
\end{eqnarray*}
gives \\ 
\scalebox{0.77}{
\begin{minipage}{1\textwidth}
\begin{eqnarray*}
&&\hat{u}_{P_{k},1,k}| \hat{\alpha}^2_{P_{k},1}\sim
\Normal \left(0, \frac{\theta_k \theta_1^2 }{(\theta_k-\theta_1)^2}\frac{\hat{\alpha}^2_{P_{k-1},1}-1}{m} \right)+O_p\left(\frac{\min(\theta_1,\theta_k)}{\theta_1^{1/2}\theta_k^{1/2}m} \right)+O_p\left(\frac{1}{ \theta_1^{1/2} \theta_k^{1/2} m^{1/2}} \right)
\end{eqnarray*}
\end{minipage}}\\

and  \\ 
\scalebox{0.85}{
\begin{minipage}{1\textwidth}
\begin{eqnarray*}
&&\hat{u}_{P_{k},1,k}\sim
\Normal \left(0, \frac{\theta_k \theta_1 }{(\theta_k-\theta_1)^2}\frac{M_{2}-1}{m} \right)+O_p\left(\frac{\min(\theta_1,\theta_k)}{\theta_1^{1/2}\theta_k^{1/2}m} \right)+O_p\left(\frac{1}{ \theta_1^{1/2} \theta_k^{1/2} m^{1/2}} \right).
\end{eqnarray*}
\end{minipage}}\\

\noindent Finally, we extend this result to small eigenvalues,
\begin{eqnarray*}
&&\text{If } \theta_1 \rightarrow \infty \text{ and } \theta_k \text{ is finite, then } \hat{u}_{P_{k},1,k}=O_p\left( \frac{1}{\sqrt{\theta_1 m}} \right),\\
&&\text{If } \theta_1 \text{ and } \theta_k \text{ are finite, then } \hat{u}_{P_{k},1,k}=O_p\left( \frac{1}{\sqrt{m}} \right).
\end{eqnarray*}

\end{proof}

\paragraph{Red}
By induction we show the part of the Invariant Theorem \ref{AAInvariantth} shown in red in the picture. We assume the truth of the grey theorems.
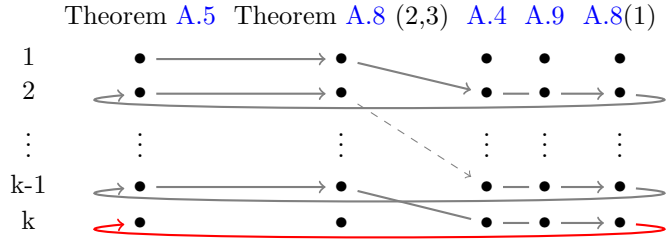
\begin{center}
\begin{tikzpicture}
  \matrix (magic) [matrix of nodes]
  { \  & Theorem \ref{AAInvariantth} & Theorem \ref{AAThcomponentdistribution} (2,3)    &   \ref{AAThinvarianteigenvalue} & \ref{AATheoremcaraceigenstructure} & \ref{AAThcomponentdistribution}(1)\\
   1  & $\bullet$   & $\bullet$    & $\bullet$ &  $\bullet$   & $\bullet$\\
   2  & $\bullet$   & $\bullet$ & $\bullet$  &  $\bullet$&  $\bullet$\\
   \vdots & \vdots   & \vdots  & \vdots & \vdots  & \vdots\\
    k-1  & $\bullet$  & $\bullet$  & $\bullet$  &  $\bullet$ &  $\bullet$ \\
   k  & $\bullet$ & $\bullet$   &  $\bullet$  &  $\bullet$ &  $\bullet$ \\
};
  \draw[thick,gray,->] (magic-2-2) -- (magic-2-3) ;
   \draw[thick,gray,->] (magic-2-3) -- (magic-3-4);
   \draw[thick,gray,->] (magic-3-4) -- (magic-3-5) -- (magic-3-6)  ;
    
\draw[thick,gray,->] (magic-3-6) to [out=-5,in=-175,looseness=1](magic-3-2);
\draw[thick,gray,->] (magic-3-2) -- (magic-3-3);
\draw[dashed,gray,->] (magic-3-3) -- (magic-5-4);
\draw[thick,gray,->]  (magic-5-4) -- (magic-5-5) -- (magic-5-6) ;
\draw[thick,gray,->] (magic-5-6) to [out=-5,in=-175,looseness=1](magic-5-2);
\draw[thick,gray,->] (magic-5-2) -- (magic-5-3) ;

\draw[thick,gray,->]  (magic-5-3) -- (magic-6-4)  -- (magic-6-5) -- (magic-6-6);
\draw[thick,red,->] (magic-6-6) to [out=-5,in=-175,looseness=1](magic-6-2);
\end{tikzpicture}
\end{center}

\begin{proof}{\textbf{Theorem \ref{AAInvariantth}}}  \label{AAproofInvariantth}
\noindent We assume the induction hypotheses and prove the result for $k$. The idea is to use Theorem \ref{AATheoremcaraceigenstructure} to simplify the $k$ first entries of the eigenvector $\tilde{u}_{P_k,1}$. Then, we show that 
$$\tilde{F}_{P_{k}}^2= \sum_{i=k+1}^m \tilde{u}_{P_{k},1,i}^2= \sum_{i=k}^m \hat{u}_{P_{k-1},1,i}^2 + O_p\left( \frac{1}{m \theta_1} \right) = \hat{F}_{P_{k-1}}^2+ O_p\left( \frac{1}{m \theta_1} \right).$$
Finally, we easily prove 
$$ \hat{F}_{P_{k}}^2= \hat{F}_{P_{k-1}}^2+O_p\left( \frac{1}{m \theta_1} \right).$$

\begin{Rem}\ \\ 
The following proof studies $\sum_{i=k+1}^m \tilde{u}_{P_{k},1,i}^2$ with $\theta_1>\theta_2>...>\theta_k$. However, the proof is easily extended to $\sum_{i=k+1}^m \tilde{u}_{P_{k},s,i}^2$ for $s=1,2,...,k$ and $\theta_s>\theta_k$. Finally, the proof is also valid for $\theta_s>\theta_k$ with more elaborate notation as in \ref{AANot=thetadifferent}. In order to simplify the two expansions for the reader, we will not further reduce values such as $\min(\theta_1,\theta_i)$.
\end{Rem}

\begin{enumerate}
\item[\textbf{A:}] First, we investigate 
\begin{eqnarray*}
\tilde{\Sigma}_{P_k}=\hat{\Sigma}_{P_{k-1}} P_k
\end{eqnarray*}
using Theorem \ref{AATheoremcaraceigenstructure} and \ref{AAconvergence}.
The eigenvectors of $\tilde{\Sigma}_{P_k}$ are
\begin{eqnarray*}
\tilde{u}_{P_k,i}=\frac{( \hat{\theta}_{P_k,i} \I_m - \hat{\Sigma}_{P_{k-1}})^{-1} \hat{\Sigma}_{P_{k-1}} \epsilon_k  }{ \sqrt{e_k^t \hat{\Sigma}_{P_{k-1}} ( \hat{\theta}_{P_k,i}\rm{I} - \hat{\Sigma}_{P_{k-1}})^{-2} \hat{\Sigma}_{P_{k-1}} \epsilon_k}}.
\end{eqnarray*}
We then have \\ 
\scalebox{0.62}{
\begin{minipage}{1\textwidth} \begin{eqnarray*} 
\left\langle \tilde{u}_{P_k,1},e_s \right\rangle^2 \hspace{-1.3cm} &&\\
 &&\hspace{-1cm}= \frac{\left( \sum_{i=k}^m \frac{\hat{\lambda}_{P_{k-1},i}}{\hat{\theta}_{P_k,1}-\hat{\lambda}_{P_{k-1},i}} \hat{u}_{P_{k-1},i,s} \hat{u}_{P_{k-1},i,k} + 
\frac{\hat{\theta}_{P_{k-1},1}}{\hat{\theta}_{P_{k},1}-\hat{\theta}_{P_{k-1},1}}  \hat{u}_{P_{k-1},1,s} \hat{u}_{P_{k-1},1,k} +  
\sum_{i=2}^{k-1} \frac{\hat{\theta}_{P_{k-1},i}}{\hat{\theta}_{P_k,1}-\hat{\theta}_{P_{k-1},i}} \hat{u}_{P_{k-1},i,s} \hat{u}_{P_{k-1},i,k} \right)^2 }
{\sum_{i=k}^m \frac{\hat{\lambda}_{P_{k-1},i}^2}{(\hat{\theta}_{P_k,1}-\hat{\lambda}_{P_{k-1},i})^2} \hat{u}_{P_{k-1},i,k}^2 + 
\frac{\hat{\theta}_{P_{k-1},1}^2}{(\hat{\theta}_{P_{k},1}-\hat{\theta}_{P_{k-1},1})^2}   \hat{u}_{P_{k-1},1,k}^2 +  
\sum_{i=2}^{k-1} \frac{\hat{\theta}_{P_{k-1},i}^2}{(\hat{\theta}_{P_k,1}-\hat{\theta}_{P_{k-1},i})^2}  \hat{u}_{P_{k-1},i,k}^2}\\
&&\hspace{-1cm}= \frac{\left(A_{1,s,k:m}+A_{1,s,1}+A_{1,s,2:k-1} \right)^2}{D_{1,k:m}+D_{1,1}+D_{1,2:k-1}}\\
&&\hspace{-1cm}= \frac{A_{1,s}^2}{D_{1}}. \hspace{25cm}
\end{eqnarray*}
\end{minipage}}\\

The size of each element of the equation can be estimated by Theorem \ref{AATheoremcaraceigenstructure}.

\item[\textbf{B:}] We investigate the norm of the noisy part of the eigenvector. Let
\begin{eqnarray*}
\tilde{F}_{P_k}^2=\sum_{i=k+1}^m \tilde{u}_{P_{k},1,i}^2=1-\sum_{i=1}^k \tilde{u}_{P_{k},1,i}^2
=1-\frac{\sum_{i=s}^{k}A_{1,s}^2}{D_1}.
\end{eqnarray*}
We want to show that $\tilde{F}_{P_{k}}^2 \approx \hat{F}_{P_{k-1}}^2$ using Theorem \ref{AATheoremcaraceigenstructure}.\\
\noindent First, we approximate $A_{1,s}$, $A_{1,s}^2$ and $D_1$: \\ 
\scalebox{0.84}{
\begin{minipage}{1\textwidth}
\begin{eqnarray*}
&&A_{1,1}=\overbrace{ \frac{\hat{\theta}_{P_{k-1},1}}{\hat{\theta}_{P_{k},1}-\hat{\theta}_{P_{k-1},1}}  \hat{u}_{P_{k-1},1,1} \hat{u}_{P_{k-1},1,k}}^{ O_p\left( \frac{\theta_1^{1/2} m^{1/2}}{\min(\theta_1,\theta_k)} \right)}+\overbrace{ \sum_{i=2}^{k-1} \frac{\hat{\theta}_{P_{k-1},i}}{\hat{\theta}_{P_k,1}-\hat{\theta}_{P_{k-1},i}} \hat{u}_{P_{k-1},i,1} \hat{u}_{P_{k-1},i,k} }^{ O_p\left( \frac{1}{\theta_1^{1/2} m}\right)}\\
&&\hspace{3cm}+O_p\left(\frac{1}{\theta_1^{3/2}m^{1/2}} \right), \hspace{30cm}
\end{eqnarray*}
\end{minipage}} \\ 
\scalebox{0.8}{
\begin{minipage}{1\textwidth}
\begin{eqnarray*}
&&A_{1,s}=\overbrace{\frac{\hat{\theta}_{P_{k-1},1}}{\hat{\theta}_{P_{k},1}-\hat{\theta}_{P_{k-1},1}}  \hat{u}_{P_{k-1},1,s} \hat{u}_{P_{k-1},1,k}}^{ O_p\left( \frac{\min(\theta_1,\theta_s)}{\theta_s^{1/2} \min(\theta_1,\theta_k)} \right)} 
+  
\overbrace{\sum_{i=2,\neq s}^{k-1} \frac{\hat{\theta}_{P_{k-1},i}}{\hat{\theta}_{P_k,1}-\hat{\theta}_{P_{k-1},i}} \hat{u}_{P_{k-1},i,s} \hat{u}_{P_{k-1},i,k} }^{O_p\left( \underset{i\neq 1,s,k}{\max}\left( \frac{\min(\theta_1, \theta_i) \min(\theta_s,\theta_i)}{\theta_s^{1/2}\theta_1 \theta_i m^{1/2}}\right)\right)}\\
&&\hspace{3cm}+ 
\overbrace{ \frac{\hat{\theta}_{P_{k-1},s}}{\hat{\theta}_{P_k,1}-\hat{\theta}_{P_{k-1},s}} \hat{u}_{P_{k-1},s,s} \hat{u}_{P_{k-1},s,k} }^{O_p\left( \frac{\min(\theta_1,\theta_s)}{\theta_s^{1/2}  \theta_1 m^{1/2}}\right)} +O_p\left(\frac{1}{\theta_s^{1/2} \theta_1 m^{1/2}} \right),\hspace{30cm}
\end{eqnarray*}
\end{minipage}}\\
\begin{eqnarray*}
&&A_{1,k}=\frac{1}{\theta_k-1}, \hspace{30cm}
\end{eqnarray*} \\ 
\scalebox{0.8}{
\begin{minipage}{1\textwidth}
\begin{eqnarray*}
&&D_1= 
\overbrace{\frac{\hat{\theta}_{P_{k-1},1}^2}{(\hat{\theta}_{P_{k},1}-\hat{\theta}_{P_{k-1},1})^2}   \hat{u}_{P_{k-1},1,k}^2 }^{O_p\left( \frac{\theta_1 m}{\min(\theta_1,\theta_k)^2}\right) }+  
\overbrace{\sum_{i=2}^{k-1} \frac{\hat{\theta}_{P_{k-1},i}^2}{(\hat{\theta}_{P_k,1}-\hat{\theta}_{P_{k-1},i})^2}  \hat{u}_{P_{k-1},i,k}^2}^{O_p\left( \frac{1}{\theta_1 m}\right) }+O_p\left( \frac{1}{\theta_1^2}\right) \\
&&\hspace{0.6cm}=O_p\left(\frac{\theta_1 m}{\min(\theta_1,\theta_k)^2}\right), \hspace{30cm}
\end{eqnarray*}
\end{minipage}}\\
\begin{eqnarray*}
&&A_{1,k}^2=\frac{1}{(\theta_k-1)^2}, \hspace{30cm}
\end{eqnarray*}
\scalebox{0.78}{
\begin{minipage}{1\textwidth}
\begin{eqnarray*}
&& A_{1,s}^2=
\overbrace{\sum_{i=1}^{k-1} \frac{\hat{\theta}_{P_{k-1},i}^2}{(\hat{\theta}_{P_k,1}-\hat{\theta}_{P_{k-1},i})^2} \hat{u}_{P_{k-1},i,s}^2 \hat{u}_{P_{k-1},i,k}^2 }^{A_{1,s,1}} 
  \\
&&\hspace{1.5cm} +\underbrace{2\sum_{i=1}^{k-1} \sum_{j>i}^{k-1} \frac{\hat{\theta}_{P_{k-1},i}\hat{\theta}_{P_{k-1},j}}{(\hat{\theta}_{P_k,1}-\hat{\theta}_{P_{k-1},i})(\hat{\theta}_{P_k,1}-\hat{\theta}_{P_{k-1},j})} \hat{u}_{P_{k-1},i,s} \hat{u}_{P_{k-1},i,k} \hat{u}_{P_{k-1},j,s} \hat{u}_{P_{k-1},j,k}}_{A_{1,s,2}}\\
&&\hspace{1.5cm} + O_p\left( \frac{1}{\theta_1 \min(\theta_1,\theta_k)} \right).\hspace{30cm}
\end{eqnarray*}
\end{minipage}}\\

Further  investigations allow us to estimate the $ \sum_{s=2}^{k-1} A_{1,s}^2$:
\begin{eqnarray*}
\sum_{s=1}^{k-1} A_{1,s}^2 &=& \sum_{s=1}^{k-1} A_{1,s,1}+\sum_{s=2}^{k-1} A_{1,s,2}+O_p\left( \frac{1}{\theta_1 \min(\theta_1,\theta_k)} \right), \hspace{30cm}
\end{eqnarray*} 
\scalebox{0.7}{
\begin{minipage}{1\textwidth}
\begin{eqnarray*}
\sum_{s=1}^{k-1} A_{1,s,1}&=&\sum_{s=1}^{k-1} \sum_{i=1}^{k-1} \frac{\hat{\theta}_{P_{k-1},i}^2}{(\hat{\theta}_{P_k,1}-\hat{\theta}_{P_{k-1},i})^2} \hat{u}_{P_{k-1},i,s}^2 \hat{u}_{P_{k-1},i,k}^2\\
&=& \sum_{i=1}^{k-1} \left( \sum_{s=1}^{k-1} \hat{u}_{P_{k-1},i,s}^2\right) \frac{\hat{\theta}_{P_{k-1},i}^2}{(\hat{\theta}_{P_k,1}-\hat{\theta}_{P_{k-1},i})^2}  \hat{u}_{P_{k-1},i,k}^2 \\
&=& \left( \sum_{s=1}^{k-1} \hat{u}_{P_{k-1},1,s}^2\right) D_1
+
\underbrace{\sum_{i=1}^{k-1} \underbrace{\left( \sum_{s=1}^{k-1} \hat{u}_{P_{k-1},i,s}^2-\sum_{s=1}^{k-1} \hat{u}_{P_{k-1},1,s}^2\right)}_{O_p\left(\frac{1}{\min(\theta_1,\theta_i)}\right) \text{ by induction}} \underbrace{\frac{\hat{\theta}_{P_{k-1},i}^2}{(\hat{\theta}_{P_k,1}-\hat{\theta}_{P_{k-1},i})^2}  \hat{u}_{P_{k-1},i,k}^2}_{O_p\left( \frac{\min(\theta_1,\theta_i)}{\theta_1(\theta_1-\theta_i)m} \right)}}_{O_p\left(\underset{i=2,...,k-1}{\max}\frac{1}{\theta_1 (\theta_1-\theta_i)m} \right)}, \hspace{30cm}
\end{eqnarray*}
\end{minipage}} \\ 
\scalebox{0.7}{
\begin{minipage}{1\textwidth}
\begin{eqnarray*} 
\sum_{s=1}^{k-1} A_{1,s,2} \hspace{-1.5cm} &&\\
&=& 2\sum_{s=1}^{k-1} \sum_{i=1}^{k-1} \sum_{j>i}^{k-1} \frac{\hat{\theta}_{P_{k-1},i}\hat{\theta}_{P_{k-1},j}}{(\hat{\theta}_{P_k,1}-\hat{\theta}_{P_{k-1},i})(\hat{\theta}_{P_k,1}-\hat{\theta}_{P_{k-1},j})} \hat{u}_{P_{k-1},i,s} \hat{u}_{P_{k-1},i,k} \hat{u}_{P_{k-1},j,s} \hat{u}_{P_{k-1},j,k}\\
&=&2 \sum_{i=1}^{k-1} \sum_{j>i}^{k-1} \frac{\hat{\theta}_{P_{k-1},i}\hat{\theta}_{P_{k-1},j}}{(\hat{\theta}_{P_k,1}-\hat{\theta}_{P_{k-1},i})(\hat{\theta}_{P_k,1}-\hat{\theta}_{P_{k-1},j})} \hat{u}_{P_{k-1},i,k}  \hat{u}_{P_{k-1},j,k} \left(\sum_{s=1}^{k-1}  \hat{u}_{P_{k-1},i,s}\hat{u}_{P_{k-1},j,s} \right)\\
&=&  2 \sum_{i=1}^{k-1} \sum_{j>i}^{k-1} \underbrace{\frac{\hat{\theta}_{P_{k-1},i}\hat{\theta}_{P_{k-1},j}}{(\hat{\theta}_{P_k,1}-\hat{\theta}_{P_{k-1},i})(\hat{\theta}_{P_k,1}-\hat{\theta}_{P_{k-1},j})} \hat{u}_{P_{k-1},i,k}  \hat{u}_{P_{k-1},j,k}}_{\text{If }i=1, \ O_p\left( \frac{\min(\theta_1,\theta_j)}{\min(\theta_1,\theta_k)\theta_1^{1/2} \theta_j^{1/2} } \right) \text{ and if } i>1, \ O_p\left( \frac{\min(\theta_i,\theta_1) \min(\theta_j,\theta_1)}{\theta_1^2 \theta_i^{1/2} \theta_j^{1/2} m} \right)}
 \underbrace{\left(- \sum_{s=k}^{m}  \hat{u}_{P_{k-1},i,s}\hat{u}_{P_{k-1},j,s} \right)}_{O_p\left( \frac{1}{\theta_i^{1/2} \theta_j^{1/2}}\right)}\\
 &=&O_p\left(\underset{j=2,...,k-1}{\max}\frac{1}{\min(\theta_1,\theta_k) \max(\theta_1,\theta_j)}\right).\hspace{30cm}
\end{eqnarray*}
\end{minipage}}\\

Thus,
\begin{eqnarray*}
\tilde{F}_{P_k}^2&=&1-\frac{1}{D_1 (\theta_k-1)^2}-\frac{\sum_{i=s}^{k-1}A_{1,s}^2}{D_1}\\
&=&1-O_p\left(\frac{1}{\theta_1 m }\right)-  \sum_{s=1}^{k-1} \hat{u}_{P_{k-1},1,s}^2  +O_p\left( \frac{\theta_1}{\max(\theta_1,\theta_k)^2 m}\right)\\
&=&1-  \sum_{s=1}^{k-1} \hat{u}_{P_{k-1},1,s}^2  +O_p\left( \frac{1}{\theta_1 m}\right)\\
&=&\hat{F}_{P_{k-1}}+O_p\left( \frac{1}{\theta_1 m}\right).
\end{eqnarray*}

\item[\textbf{C:}] The result is already demonstrated for the eigenvector of 
\begin{eqnarray*}
\tilde{\Sigma}_{P_k}=\hat{\Sigma}_{P_{k-1}} P_k.
\end{eqnarray*}
Now, we need to extend this to 
\begin{eqnarray*}
\hat{\Sigma}_{P_k}=P_k^{1/2} \hat{\Sigma}_{P_{k-1}} P_k^{1/2}.
\end{eqnarray*}
The link between the eigenvectors is 
\begin{eqnarray*}
\hat{u}_{P_{k-1},1}&=&\frac{P_k^{1/2} \tilde{u}_{P_{k},1}}{ \sqrt{{\rm Norm}^2}}\\
{\rm Norm}^2 &=& \sum_{i=1}^{k-1} \tilde{u}_{P_{k},1,i}^2 + \theta_k  \tilde{u}_{P_{k},1,k}^2 +\sum_{i=k+1}^{m} \tilde{u}_{P_{k},1,i}^2\\
&=&1+\underbrace{(\theta_k-1)  \overbrace{\tilde{u}_{P_{k},1,k}^2}^{O_p\left(\frac{\theta_1}{\max(\theta_1,\theta_k)^2 m}\right)}}_{O_p\left(\frac{1}{m}\right)}
\end{eqnarray*}
Using the induction hypothesis, the result is true for $k-1$; therefore, by Theorem \ref{AAjointdistribution}, 
$$\hat{F}_{P_{k-1}}={\rm RV} \left( O_p\left( \frac{1}{\theta_1} \right),O_p\left( \frac{1}{\theta_1^2 m} \right) \right).$$ 
Then,
\begin{eqnarray*}
\hat{F}_{P_k}^2&=&\sum_{i=k+1}^m  \hat{u}_{P_{k-1},1,i}^2 \\
&=&\frac{1}{\rm Norm^2}\sum_{i=k+1}^m  \tilde{u}_{P_{k-1},1,i}^2\\
&=&\frac{1}{1+O_p\left( \frac{1}{ m} \right)}\tilde{F}_{P_k}^2 \\
&=&\frac{1}{1+O_p\left( \frac{1}{ m} \right)} \left( \hat{F}_{P_{k-1}}+O_p\left( \frac{1}{\theta_1 m}\right) \right) \\
&=&  \hat{F}_{P_{k-1}} +O_p\left( \frac{1}{\theta_1 m}\right).
\end{eqnarray*}
This last equation concludes the proof by induction
\begin{eqnarray*}
\sum_{i=1}^k \hat{u}_{P_k,1,i}^2 =\hat{u}_{P_1,1,1}^2+O_p\left( \frac{1}{\theta_1 m}\right).
\end{eqnarray*}
\end{enumerate}
\end{proof}

\paragraph{Green}
In this section we want to prove the green part in the following picture.
\noindent In order to prove Theorem \ref{AAThcomponentdistribution} (2, 3) for $k$, we only assume the truth of the grey results in the picture.

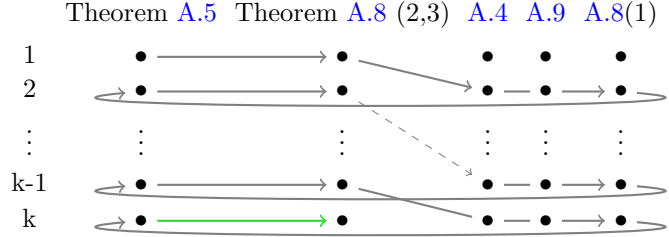
\begin{center}
\begin{tikzpicture}
  \matrix (magic) [matrix of nodes]
  { \  & Theorem \ref{AAInvariantth} & Theorem \ref{AAThcomponentdistribution} (2,3)    &   \ref{AAThinvarianteigenvalue} & \ref{AATheoremcaraceigenstructure} & \ref{AAThcomponentdistribution}(1)\\
   1  & $\bullet$   & $\bullet$    & $\bullet$ &  $\bullet$   & $\bullet$\\
   2  & $\bullet$   & $\bullet$ & $\bullet$  &  $\bullet$&  $\bullet$\\
   \vdots & \vdots   & \vdots  & \vdots & \vdots  & \vdots\\
    k-1  & $\bullet$  & $\bullet$  & $\bullet$  &  $\bullet$ &  $\bullet$ \\
   k  & $\bullet$ & $\bullet$   &  $\bullet$  &  $\bullet$ &  $\bullet$ \\
};
  \draw[thick,gray,->] (magic-2-2) -- (magic-2-3) ;
   \draw[thick,gray,->] (magic-2-3) -- (magic-3-4);
   \draw[thick,gray,->] (magic-3-4) -- (magic-3-5) -- (magic-3-6)   ;
    
\draw[thick,gray,->] (magic-3-6) to [out=-5,in=-175,looseness=1](magic-3-2);
\draw[thick,gray,->] (magic-3-2) -- (magic-3-3);
\draw[dashed,gray,->] (magic-3-3) -- (magic-5-4);
\draw[thick,gray,->]  (magic-5-4) -- (magic-5-5) -- (magic-5-6) ;
\draw[thick,gray,->] (magic-5-6) to [out=-5,in=-175,looseness=1](magic-5-2);
\draw[thick,gray,->] (magic-5-2) -- (magic-5-3) ;

\draw[thick,green!80!black!80,->] (magic-6-2) -- (magic-6-3)  ;
\draw[thick,gray,->]  (magic-5-3) -- (magic-6-4)  -- (magic-6-5) -- (magic-6-6);
\draw[thick,gray,->] (magic-6-6) to [out=-5,in=-175,looseness=1](magic-6-2);
\end{tikzpicture}
\end{center}

\begin{proof}{\textbf{Theorem \ref{AAThcomponentdistribution} (2,3)}}  \label{AAproofThcomponentdistribution23k}
\noindent To prove this theorem for $k$ we use the same procedure as for $k=1$.\\
Let 
\begin{eqnarray*}
U&=&
\begin{pmatrix}
\hat{u}_{P_k,1}^t\\
\hat{u}_{P_k,2}^t\\
\vdots \\
\hat{u}_{P_k,m}^t
\end{pmatrix}=\begin{pmatrix}
\hat{u}_{P_k,1:k,1:k} & \hat{u}_{P_k,1:k,k+1:m} \\ 
\hat{u}_{P_k,k+1:m,1:k} & \hat{u}_{P_k,k+1:m,k+1:m}.
\end{pmatrix}
\end{eqnarray*}
and
$$ O_k= \begin{pmatrix}
\rm{I}_k & 0 \\ 
0 & O_{m-k}
\end{pmatrix},$$
where $O_{m-k}$ is Haar invariant.
\begin{itemize}
\item[2.] When $P_k$ is canonical, we know that $\hat{\Sigma} \sim P_k^{1/2} W P_k^{1/2}$ and $O_k \hat{\Sigma} O_k^t$ follow the same distribution under Assumption \ref{AAAss=matrice}. Therefore,  $\hat{u}_{i,k+1:m}$ is rotationally invariant and $\corr \left( \hat{u}_{i,j_1}, \hat{u}_{i,j_2} \right)=\delta_{j_1}(j_2)$.
Knowing \linebreak $\hat{u}_{P_k,1:m,1:k}$, we can show that  
$\hat{u}_{P_k,i,k+1:m}/ ||\hat{u}_{P_k,i,k+1:m}||$ is uniform for $i=1,2,...,m$. Therefore, these statistics are independent (not jointly) of $\hat{u}_{P_k,1:m,1:k}$. 
Uniformity of $\hat{u}_{P_k,r,k+1:m}$ implies that, for $s=k+1,...,m$ and $r=1,2,...,k$,
\begin{eqnarray*}
\sqrt{m} \frac{\hat{u}_{P_k,r,s}}{ ||\hat{u}_{P_k,r,(k+1):m}||}= \sqrt{m} \frac{\hat{u}_{P_k,r,s}}{ \sqrt{1-\hat{\alpha}^2_{P_k,r}}}\sim \Normal\left(0,1\right)+o_p \left( 1 \right),
\end{eqnarray*}
where
\begin{eqnarray*}
\hat{\alpha}^2_{P_k,r}=\sum_{i=1}^k \left\langle \hat{u}_{P_k,r},\epsilon_i \right\rangle^2.
\end{eqnarray*}
By Slutsky's Theorem and the Invariant Angle Theorem \ref{AAInvariantth} for $k$,
 \begin{eqnarray*}
\hat{u}_{P_k,r,s} \sim \Normal \left( 0,\frac{1-\alpha_r^2}{m}\right)+o_p \left( \frac{1}{\sqrt{m}} \right),
\end{eqnarray*}
where $\alpha_r^2= \underset{m \rightarrow \infty}{\lim} \hat{\alpha}^2_{P_k,r} = 1-\frac{M_2-1}{\theta_r}+O_p\left( \frac{1}{\theta^2}\right)<1$.

\item[3.] Next, we estimate the order of $\sum \hat{u}_{k+1:m,1}^2$.\\
Without loss of generality we assume that the perturbation 
\begin{eqnarray*}
P_k=\I_m + \sum_{i=1}^k (\theta_i-1) \epsilon_i \epsilon_i^t
\end{eqnarray*}
verifies Assumption \ref{AAAss=theta}(A4) and is such that 
\begin{eqnarray*}
&&\theta_1,\theta_2,...,\theta_{k_1} \text{ are proportional,}\\
&&\theta_{k_1+1},\theta_{k_1+2},...,\theta_{k} \text{ are proportional}.
\end{eqnarray*}
Then by Theorem \ref{AAInvariantth} and \ref{AAThcomponentdistribution} Part 1 for perturbations of order $k$, \\ 
\scalebox{0.9}{
\begin{minipage}{1\textwidth}
\begin{eqnarray*}
  \sum \hat{u}_{k+1:m,1:k_{1}}^2
 &=& \sum \hat{u}_{k_1+1:m,1:k_{1}}^2-  \sum \hat{u}_{k_1+1:k,1:k_{1}}^2 \\
&=& \sum \hat{u}_{1:k_{1},k_1+1:m}^2  + O_p\left( \frac{\min(\theta_1,\theta_k)}{\max(\theta_1,\theta_k) m}\right)\\
&=& \sum \hat{u}_{1:k_{1},k+1:m}^2  + O_p\left( \frac{\min(\theta_1,\theta_k)}{\max(\theta_1,\theta_k) m}\right)\\
&\sim & {\rm{RV}}\left( O\left( \frac{1}{\theta_1 } \right), O \left( \frac{1}{\theta_1^2 m}\right)\right)+ O_p\left( \frac{\min(\theta_1,\theta_k)}{\max(\theta_1,\theta_k) m}\right).
\end{eqnarray*}
\end{minipage}}\\

The result is straightforward. 

\end{itemize}
\end{proof}

\subsubsection{Dot product distribution and perquisite Lemma}

In this section we prove the results concerning the partial dot product between two estimated eigenvectors. First, we show a useful small Lemma. Then, we investigate its distribution when $k=2$. Finally, we prove the invariance to increasing $k$.

\paragraph{Prerequisite}
\begin{proof}{\textbf{Theorem \ref{AALemstatW}}} \label{AAproofLemstatW}
\noindent The proofs of the three results use Theorem \ref{AATheoremcaraceigenstructure}. \\
First, we recall that 
\begin{eqnarray*}
\sum_{i=2}^m \hat{\lambda}_{P_1,i}^2 \hat{u}_{P_1,i,2}^2 &=& \hat{\Sigma}^2_{P_1,2,2}- \hat{\theta}_{P_1,1}^2 \hat{u}_{P_1,1,2}^2,\\
\sum_{i=2}^m \hat{\lambda}_{P_1,i} \hat{u}_{P_1,i,1}\hat{u}_{P_1,i,2} &=& \hat{\Sigma}_{P_1,1,2}- \hat{\theta}_{P_1,1} \hat{u}_{P_1,1,1}\hat{u}_{P_1,1,2}.
\end{eqnarray*}
Moreover, if $\tilde{P}_1= \left( \sqrt{\theta_1}-1 \right) e_1 e_1^t$, then
\begin{eqnarray*}
\hat{\Sigma}_{P_1}&=& W+W \tilde{P}_1+ \tilde{P}_1 W + \tilde{P}_1 W \tilde{P}_1,\\
\hat{\Sigma}_{P_1,1,2}&=& W_{1,2} \sqrt{\theta_1},\\
\left(\hat{\Sigma}_{P_1}^2\right)_{2,2}&=& \left( W+W \tilde{P}_1+ \tilde{P}_1 W + \tilde{P}_1 W \tilde{P}_1 \right)^2 [2,2]\\
&=& \left(W^2\right)_{2,2}+(\theta-1) \left(W_{1,2}\right)^2,
\end{eqnarray*}
where $A[2,2]$ is the entry $A_{2,2}$ of the matrix $A$.\\
In order to prove the formulas, we need some estimations of 
\begin{eqnarray*}
\hat{u}_{P_1,1,1}^2, \ \hat{\theta}_{P_1,1}^2 \text{ and } \frac{\hat{u}_{P_1,1,2}}{\sqrt{1-\hat{u}_{P_1,1,1}^2}}.
\end{eqnarray*}
A more precise estimation of $\hat{u}_{P_1,1,1}^2$ leads to \\ 
\scalebox{0.74}{
\begin{minipage}{1\textwidth}
\begin{eqnarray*}
\hat{u}_{P_1,1,1}^2&=& 1-\frac{\left(W^2\right)_{1,1}-\left(W_{1,1}\right)^2}{\theta_1\left(W_{1,1}\right)^2}+\frac{1+\frac{3 \left( \left(W^2\right)_{1,1}\right)^2}{\left(W_{1,1}\right)^4}-\frac{2 \left(W^2\right)_{1,1}}{\left(W_{1,1}\right)^2}-\frac{2 \left(W^3\right)_{1,1}}{\left(W_{1,1}\right)^3}}{\theta_1^2}+O_p\left(\frac{1}{\theta_1^3} \right),\\
\hat{u}_{P_1,1,1}&=&1-\frac{\left(W^2\right)_{1,1}-\left(W_{1,1}\right)^2}{2 \theta_1 \left(W_{1,1}\right)^2}+\frac{1+\frac{3 \left(\left(W^2\right)_{1,1}\right)^2}{\left(W_{1,1}\right)^4}-\frac{2 \left(W^2\right)_{1,1}}{\left(W_{1,1}\right)^2}-\frac{2 \left(W^3\right)_{1,1}}{\left(W_{1,1}\right)^3}}{2 \theta_1^2}+O_p\left(\frac{1}{\theta_1^3} \right),\\
\sqrt{1-\hat{u}_{P_1,1,1}^2}&=&\frac{1}{\sqrt{\theta_1}}\left( \frac{\sqrt{\left(W^2\right)_{1,1}-\left(W_{1,1}\right)^2}}{W_{1,1}}-\frac{W_{1,1} \left( 1+\frac{3 \left(\left(W^2\right)_{1,1}\right)^2}{\left(W_{1,1}\right)^4}-\frac{2 \left(W^2\right)_{1,1}}{\left(W_{1,1}\right)^2}-\frac{2 \left(W^3\right)_{1,1}}{\left(W_{1,1}\right)^3} \right)}{2 \theta \sqrt{\left(W^2\right)_{1,1}-\left(W_{1,1}\right)^2}} \right) \\
&&\hspace{1cm}+O_p\left(\frac{1}{\theta_1^{5/2}} \right).
\end{eqnarray*} 
\end{minipage}}\\

Then, we estimate $\hat{\theta}_{P_1,1}^2$, \small
\begin{eqnarray*}
\frac{1}{\theta_1-1} &=& \sum_{i=1}^m \frac{\hat{\lambda}_{W,i}}{\hat{\theta}_{P_1,1}-\hat{\lambda}_{W,i}} \hat{u}_{W,i,1}^2\\
&=& \frac{W_{1,1}}{\hat{\theta}_{P_1,1}} + \frac{\left(W^2\right)_{1,1}}{\hat{\theta}_{P_1,1}^2} + O_p\left( \frac{1}{\theta_1^3} \right)\\
&&\hspace{-2cm} \Rightarrow \hat{\theta}_{P_1,1} = \theta_1 W_{1,1}+\frac{\left(W^2\right)_{1,1}-\left(W_{1,1}\right)^2}{W_{1,1}}+O_p\left(\frac{1}{\theta_1}\right)\\
&&\hspace{-2cm} \Rightarrow \hat{\theta}_{P_1,1}^2 = \theta_1^2 W_{1,1}^2+2 \theta_1 \left( \left(W^2\right)_{1,1}-\left(W_{1,1}\right)^2\right) +O_p\left(1\right).
\end{eqnarray*} 
Finally, we estimate the rescaled component, \\ 
\scalebox{0.85}{
\begin{minipage}{1\textwidth}
\begin{eqnarray*}
\frac{\hat{u}_{P_1,1,2}}{\sqrt{1-\hat{u}_{P_1,1,1}^2}}&=&\frac{\sum_{i=1}^m \frac{\hat{\lambda}_{W,i}}{\hat{\theta}_{P_1,1}-\hat{\lambda}_{W,i}} \hat{u}_{W,i,1}\hat{u}_{W,i,2}}{\sqrt{\sum_{s=2}^m \left( \sum_{i=1}^m \frac{\hat{\lambda}_{W,i}}{\hat{\theta}_{P_1,1}-\hat{\lambda}_{W,i}} \hat{u}_{W,i,1}\hat{u}_{W,i,s} \right)^2}}\\
&=& \frac{\frac{1}{\hat{\theta}_{P_1,1}} W_{1,2}+ \frac{1}{\hat{\theta}_{P_1,1}^2} \left(W^2\right)_{1,2}+ O_p\left(\frac{1}{\theta^3 \sqrt{m}} \right)}{\sqrt{\sum_{s=2}^m \left( \frac{1}{\hat{\theta}_{P_1,1}} W_{1,s}+ \frac{1}{\hat{\theta}_{P_1,1}^2} \left(W^2\right)_{1,s}+ O_p\left(\frac{1}{\theta_1^3 \sqrt{m}} \right) \right)^2}}\\
&=& \frac{ W_{1,2}+ \frac{1}{\hat{\theta}_{P_1,1}} \left(W^2\right)_{1,2}+ O_p\left(\frac{1}{\theta_1^2 \sqrt{m}} \right)}{\sqrt{\sum_{s=2}^m \left(  \left(W_{1,s}\right)^2+ 2 \frac{1}{\hat{\theta}_{P_1,1}}W_{1,s}  \left(W^2\right)_{1,s}+ O_p\left(\frac{1}{\theta_1^2 m} \right) \right)}}\\
&=& \frac{ W_{1,2}+ \frac{1}{\hat{\theta}_{P_1,1}} \left(W^2\right)_{1,2}+ O_p\left(\frac{1}{\theta_1^2 \sqrt{m}} \right)}{\sqrt{\left(W^2\right)_{1,1}-\left(W_{1,1}\right)^2+ \frac{1}{\hat{\theta}_{P_1,1}} \left[ \left(W^3\right)_{1,1}-W_{1,1} \left(W^2\right)_{1,1} \right]+ O_p\left(\frac{1}{\theta_1^2 m} \right) }}\\
&=&W_{1,2} \left(\frac{1}{\sqrt{\left(W^2\right)_{1,1}-\left(W_{1,1}\right)^2}} - \frac{\left(W^3\right)_{1,1}-W_{1,1} \left(W^2\right)_{1,1}}{\left( \left(W^2\right)_{1,1}-\left(W_{1,1}\right)^2 \right)^{3/2} \hat{\theta}_{P_1,1}} \right)\\
&&\hspace{1cm} +  \frac{\left(W^2\right)_{1,2}}{\sqrt{\left(W^2\right)_{1,1}-\left(W_{1,1}\right)^2} \hat{\theta}_{P_1,1}} +O_p\left( \frac{1}{\theta_1^2 \sqrt{m}} \right).
\end{eqnarray*}
\end{minipage}}\\

Using this estimation, the three formulas are easily proven.\\
We start with the first formula:\\ 
\scalebox{0.8}{
\begin{minipage}{1\textwidth}
\begin{eqnarray*}
\hat{u}_{P_1,1,2}&=& \frac{\hat{u}_{P_1,1,2}}{\sqrt{1-\hat{u}_{P_1,1,1}^2}} \sqrt{1-\hat{u}_{P_1,1,1}^2}\\
&=& \frac{W_{1,2}}{\sqrt{\theta_1} W_{1,1}} -\frac{W_{1,2}}{\theta_1^{3/2}}\left( -1/2+3/2 M_2 \right)+ \frac{\left(W^2\right)_{1,2}}{\theta_1^{3/2}} + O_p\left( \frac{1}{\theta_1^{3/2}m} \right)+ O_p\left( \frac{1}{\theta_1^{5/2}m^{1/2}} \right).
\end{eqnarray*}
\end{minipage}}\\

Then, the second formula:\\ 
\scalebox{0.8}{
\begin{minipage}{1\textwidth}
\begin{eqnarray*}
\sum_{i=2}^m \hat{\lambda}_{P_1,i}^2 \hat{u}_{P_1,i,2}^2-\left(W^2\right)_{2,2} &=& (\theta_1-1) \left(W_{1,2}\right)^2- \hat{\theta}_{P_1,1}^2 \hat{u}_{P_1,1,2}^2\\
&=&(\theta_1-1) \left(W_{1,2}\right)^2  - \left(  \theta_1^2 \left(W_{1,1}\right)^2+2 \theta_1 \left( \left(W^2\right)_{1,1}-\left(W_{1,1}\right)^2\right) +O_p\left(1\right) \right)\\
&&\hspace{1cm} \left( \frac{\left(W_{1,2}\right)^2}{\theta_1 \left(W_{1,1}\right)^2}+O_p\left( \frac{1}{\theta^2 m} \right) \right)\\
&=&  - \left(W_{1,2}\right)^2-2  \frac{\left(W_{1,2}\right)^2 \left( \left(W^2\right)_{1,1}-\left(W_{1,1}\right)^2\right)}{\left(W_{1,1}\right)^2}+O_p\left(\frac{1}{m} \right)\\
&=&O_p\left(\frac{1}{m} \right).
\end{eqnarray*}
\end{minipage}}\\

Finally, some computations lead to the last formula,\\ 
\scalebox{0.8}{
\begin{minipage}{1\textwidth}
\begin{eqnarray*}
&&\hat{\theta}_{P_1,1} \hat{u}_{P_1,1,1}\hat{u}_{P_1,1,2}\\
&& \hspace{1cm}=
 \left( \theta_1 W_{1,1}+\frac{\left(W^2\right)_{1,1}-\left(W_{1,1}\right)^2}{W_{1,1}}\right) \left( 1-\frac{\left(W^2\right)_{1,1}-\left(W_{1,1}\right)^2}{2 \theta_1 \left(W_{1,1}\right)^2}\right) \\
&&\hspace{1.5cm} \left( \frac{W_{1,2}}{\sqrt{\theta_1} W_{1,1}} -\frac{W_{1,2}}{\theta_1^{3/2}}\left(-1/2+3/2 M_2 \right)+\frac{\left(W^2\right)_{1,2}}{\theta_1^{3/2}}\right) + O_p\left( \frac{1}{\theta_1^{1/2}m} \right)+ O_p\left( \frac{1}{\theta_1^{3/2}m^{1/2}} \right) \\
&&\hspace{1cm}=  W_{1,2} \left( \sqrt{\theta_1} -  \frac{  M_2}{\sqrt{\theta_1}} \right) + \left(W^2\right)_{1,2}\frac{1}{\sqrt{\theta_1}}
+ O_p\left( \frac{1}{\theta_1^{1/2}m} \right)+ O_p\left( \frac{1}{\theta_1^{3/2}m^{1/2}} \right).
\end{eqnarray*}
\end{minipage}}\\

Therefore,\\ 
\scalebox{0.9}{
\begin{minipage}{1\textwidth}
\begin{eqnarray*}
\sum_{i=2}^m \hat{\lambda}_{P_1,i} \hat{u}_{P_1,i,1}\hat{u}_{P_1,i,2} &=& \sqrt{\theta_1} W_{1,2} - \hat{\theta}_{P_1,1} \hat{u}_{P_1,1,1}\hat{u}_{P_1,1,2}\\
&=& W_{1,2} \frac{  M_2}{\sqrt{\theta_1}} -\left(W^2\right)_{1,2}\frac{1}{\sqrt{\theta_1}}
+ O_p\left( \frac{1}{\theta_1^{1/2}m} \right)+ O_p\left( \frac{1}{\theta_1^{3/2}m^{1/2}} \right).
\end{eqnarray*}
\end{minipage}}

\end{proof}

\paragraph{Distribution}
\noindent (Page \pageref{AAThDotproduct})
\begin{proof}{\textbf{Theorem \ref{AAThDotproduct}}}\label{AAproofThDotproduct} 
\noindent We begin this proof with a remark about the sign convention. This Theorem assumes $\hat{u}_{P_s,i,i}>0$ for $s=1,2,...,k$ and $i=1,2,...,s$. The Theorem \ref{AATheoremcaraceigenstructure}, however, constructs the eigenvectors of random matrices with another sign convention,
\begin{eqnarray*}
\hat{u}_{P_s,i,s}>0, \text{ for $s=1,2,...,k$ and $i=1,2,...,s$.}
\end{eqnarray*}
We will use the same notation for both and invite the reader to be aware of the following. The parts $\textbf{A}$ and $\textbf{B}$ use the convention of Theorem \ref{AATheoremcaraceigenstructure}. This changes in the end of part $\textbf{B}$.   Finally part $\textbf{C}$ uses the convention of this theorem.\\
The first part, $\textbf{A}$, expresses the components of an eigenvector using Theorem \ref{AATheoremcaraceigenstructure}. The second part, $\textbf{B}$, expresses the dot product of $\hat{\Sigma}_{P_2}$ with the eigenstructure of $\hat{\Sigma}_{P_1}$. Finally, with the previous part leading to a nice formula, we investigate in $\textbf{C}$ the distribution of this statistic.\\
We will often replace $\hat{\theta}_{P_1,1}$ by $ \hat{\lambda}_{P_1,1}$ to simplify computations. 

\paragraph*{A: }
For $t=1,2$, we study the expression:
\begin{eqnarray*}
\tilde{{u}}_{P_2,t,s} &=&  \frac{\sum_{i=1}^m \frac{\hat{\lambda}_{P_{1},i}}{\hat{\theta}_{P_2,t}-\hat{\lambda}_{P_{1},i}} \hat{u}_{P_{1},i,s} \hat{u}_{P_{1},i,2}}{\sqrt{\sum_{i=1}^m \frac{\hat{\lambda}_{P_{1},i}^2}{\left(\hat{\theta}_{P_2,t}-\hat{\lambda}_{P_{1},i}\right)^2} \hat{u}_{P_{1},i,2}^2}}
= \frac{\sum_{i=1}^m \frac{\hat{\lambda}_{P_{1},i}}{\hat{\theta}_{P_2,t}-\hat{\lambda}_{P_{1},i}} \hat{u}_{P_{1},i,s} \hat{u}_{P_{1},i,2}}{\sqrt{D_t}},
\end{eqnarray*}
where by Theorem \ref{AATheoremcaraceigenstructure} and assuming $\theta_1>\theta_2$,
\begin{eqnarray*}
D_1&=&\sum_{i=2}^m \frac{\hat{\lambda}_{P_{1},i}^2}{(\hat{\theta}_{P_2,1}-\hat{\lambda}_{P_{1},i})^2} \hat{u}_{P_{1},i,2}^2 + 
\frac{\hat{\theta}_{P_{1},1}^2}{(\hat{\theta}_{P_{2},1}-\hat{\theta}_{P_{1},1})^2}   \hat{u}_{P_{1},1,2}^2 \\
&=& \underbrace{\frac{\hat{\theta}_{P_{1},1}^2}{(\hat{\theta}_{P_{2},1}-\hat{\theta}_{P_{1},1})^2}   \hat{u}_{P_{1},1,2}^2}_{\ \overset{\scalebox{0.5}{order}}{\sim}  \frac{\theta_1 m}{\theta_2^2} } +O_p\left( \frac{1}{\theta_1^2} \right),\\
D_2&=&\sum_{i=2}^m \frac{\hat{\lambda}_{P_{1},i}^2}{(\hat{\theta}_{P_2,2}-\hat{\lambda}_{P_{1},i})^2} \hat{u}_{P_{1},i,2}^2 + 
\frac{\hat{\theta}_{P_{1},1}^2}{(\hat{\theta}_{P_{2},2}-\hat{\theta}_{P_{1},1})^2}   \hat{u}_{P_{1},1,2}^2 \\
&=&O_p\left( \frac{1}{\theta_2^2} \right)+O_p\left( \frac{1}{\theta_1  m} \right).
\end{eqnarray*}
By Theorem \ref{AATheoremcaraceigenstructure}, $\hat{{u}}_{P_2,t,s}=\frac{P_2^{1/2} \tilde{{u}}_{P_2,t}}{{N_t}} $, where
\begin{eqnarray*}
N_t^2&=& \tilde{{u}}_{P_2,t,1}^2 + 
\sum_{i=3}^{m} \tilde{{u}}_{P_2,t,i}^2+\tilde{{u}}_{P_2,t,2}^2 \theta_2\\
&=& 1+ (\theta_2-1)  \tilde{{u}}_{P_2,t,2}^2\\
&=& 1+ \frac{1}{(\theta_2-1) D_t}.
\end{eqnarray*}
Then,
\begin{eqnarray*}
N_t^2 D_t &=& D_t+ \frac{1}{(\theta_2-1) },\\
N_1^2 D_1 &=& \frac{\hat{\theta}_{P_{1},1}^2}{(\hat{\theta}_{P_{2},1}-\hat{\theta}_{P_{1},1})^2}   \hat{u}_{P_{1},1,2}^2 +O_p\left( \frac{1}{\theta_2}\right),\\
N_2^2 D_2 &=&\frac{1}{(\theta_2-1) }+O_p\left( \frac{1}{\theta_2^2} \right)+O_p\left( \frac{1}{\theta_1 m} \right).
\end{eqnarray*}
Therefore,
\begin{eqnarray*}
\frac{1}{N_1 \sqrt{D_1}}&=& \frac{|\hat{\theta}_{P_{2},1}-\hat{\theta}_{P_{1},1}|}{\hat{\theta}_{P_{1},1}  |\hat{u}_{P_{1},1,2}|}   +O_p\left( \frac{\theta_2^{2}}{\theta_1^{3/2} m^{3/2}} \right) \\
&=& O_p\left( \frac{\theta_2}{\theta_1^{1/2}m^{1/2}}\right),\\
\frac{1}{N_2 \sqrt{D_2}}&=& \sqrt{\theta_2-1} +O_p\left( \frac{1}{\theta_2^{1/2}} \right)+O_p\left( \frac{\theta_2^{3/2}}{\theta_1 m} \right).\\
\end{eqnarray*}

\paragraph*{B: }
\noindent We are now in a position to investigate:
\begin{eqnarray*}
\sum_{s=3}^m  \hat{u}_{P_k,1,s}  \hat{u}_{P_k,2,s}
\end{eqnarray*}
First,
\begin{eqnarray*}
 \hat{u}_{P_2,1,s}  \hat{u}_{P_2,2,s}&=& \frac{
 \sum_{i,j=1}^m \frac{\hat{\lambda}_{P_{1},i}\hat{\lambda}_{P_{1},j}}{(\hat{\theta}_{P_2,1}-\hat{\lambda}_{P_{1},i})(\hat{\theta}_{P_2,2}-\hat{\lambda}_{P_{1},i})} \hat{u}_{P_{1},i,s}\hat{u}_{P_{1},j,s} \hat{u}_{P_{1},i,2} \hat{u}_{P_{1},j,2}
 }{\sqrt{D_1 D_2}N_1 N_2 }.\\
\end{eqnarray*}
Then,\\ 
\scalebox{0.75}{
\begin{minipage}{1\textwidth}
\begin{eqnarray*}
 &&\sum_{s=3}^m \hat{u}_{P_2,1,s}  \hat{u}_{P_2,2,s}\\
 &&\hspace{1cm}= \frac{\sum_{s=3}^m
 \sum_{i,j=1}^m \frac{\hat{\lambda}_{P_{1},i}\hat{\lambda}_{P_{1},j}}{(\hat{\theta}_{P_2,1}-\hat{\lambda}_{P_{1},i})(\hat{\theta}_{P_2,2}-\hat{\lambda}_{P_{1},j})} \hat{u}_{P_{1},i,s}\hat{u}_{P_{1},j,s} \hat{u}_{P_{1},i,2} \hat{u}_{P_{1},j,2}
 }{\sqrt{D_1 D_2}N_1 N_2 }\\
 &&\hspace{1cm}= \frac{1}{\sqrt{D_1 D_2}N_1 N_2 } \left( 
 \sum_{i,j=1, i\not=j }^m \frac{\hat{\lambda}_{P_{1},i}\hat{\lambda}_{P_{1},j}}{(\hat{\theta}_{P_2,1}-\hat{\lambda}_{P_{1},i})(\hat{\theta}_{P_2,2}-\hat{\lambda}_{P_{1},j})}  \hat{u}_{P_{1},i,2} \hat{u}_{P_{1},j,2} \left( \sum_{s=3}^m\hat{u}_{P_{1},i,s}\hat{u}_{P_{1},j,s} \right)\right. \\
 &&\hspace{1.5cm}+\left. 
 \sum_{i=1}^m \frac{\hat{\lambda}_{P_{1},i}^2}{(\hat{\theta}_{P_2,1}-\hat{\lambda}_{P_{1},i})(\hat{\theta}_{P_2,2}-\hat{\lambda}_{P_{1},i})}  \hat{u}_{P_{1},i,2}^2 \left( \sum_{s=3}^m\hat{u}_{P_{1},i,s}^2 \right)
 \right)\\
 &&\hspace{1cm}= \frac{1}{\sqrt{D_1 D_2}N_1 N_2 } \\
 &&\hspace{1.5cm} \left( -\underbrace{
 \sum_{i,j=1, i\not=j }^m \frac{\hat{\lambda}_{P_{1},i}\hat{\lambda}_{P_{1},j}}{(\hat{\theta}_{P_2,1}-\hat{\lambda}_{P_{1},i})(\hat{\theta}_{P_2,2}-\hat{\lambda}_{P_{1},j})}  \left( \hat{u}_{P_{1},i,2}^2 \hat{u}_{P_{1},j,2}^2+\hat{u}_{P_{1},i,1}\hat{u}_{P_{1},j,1}\hat{u}_{P_{1},i,2} \hat{u}_{P_{1},j,2}\right)}_{\text{Part 2}}\right. \\
 &&\hspace{2cm}+\left. \underbrace{
 \sum_{i=1}^m \frac{\hat{\lambda}_{P_{1},i}^2}{(\hat{\theta}_{P_2,1}-\hat{\lambda}_{P_{1},i})(\hat{\theta}_{P_2,2}-\hat{\lambda}_{P_{1},i})}  \hat{u}_{P_{1},i,2}^2 \left( 1-\hat{u}_{P_{1},i,1}^2-\hat{u}_{P_{1},i,2}^2\right)}_{\text{Part 1}} 
 \right).
\end{eqnarray*}
\end{minipage}}\\

Using part $ \mathbf{A}$ gives,
\begin{eqnarray*}
\frac{1}{\sqrt{D_1 D_2}N_1 N_2 }= O_p\left( \frac{\theta_2^{3/2}}{\theta_1^{1/2} m^{1/2}} \right).
\end{eqnarray*}
Next, we consider the sum of Part 1 and Part 2 in the above equation and neglect terms smaller than $O_p\left( \frac{1}{\theta_2^2} \right)$. (If at least one term is of order $\frac{1}{\theta_2^2}$.) 
 \begin{itemize}
 \item[] \textbf{Part 1: } We decompose the sum of Part 1 into $i=1$ and  $i>1$. Then, using Theorems \ref{AATheoremcaraceigenstructure}, \ref{AAThcomponentdistribution} and \ref{AAThinvarianteigenvalue}, each term can be estimated.
 \begin{itemize}
 \item[1.1) ]i=1 : \\ 
\scalebox{0.7}{
\begin{minipage}{1\textwidth}
 \begin{eqnarray*}
&& \hspace{-1.5cm}\frac{\hat{\lambda}_{P_{1},1}^2}{(\hat{\theta}_{P_2,1}-\hat{\lambda}_{P_{1},1})(\hat{\theta}_{P_2,2}-\hat{\lambda}_{P_{1},1})}  \hat{u}_{P_{1},1,2}^2 \left( 1-\hat{u}_{P_{1},1,1}^2-\hat{u}_{P_{1},1,2}^2\right)\\
&& \hspace{-1cm}= \underbrace{ \mathbf{\frac{\hat{\theta}_{P_{1},1}^2}{(\hat{\theta}_{P_2,1}-\hat{\theta}_{P_{1},1})(\hat{\theta}_{P_2,2}-\hat{\theta}_{P_{1},1})}  \hat{u}_{P_{1},1,2}^2 \left( 1-\hat{\alpha}_{P_1,1}^2\right)}}_{O_p\left(\frac{1}{\theta_1 \theta_2 }\right)} - \underbrace{\frac{\hat{\theta}_{P_{1},1}^2}{(\hat{\theta}_{P_2,1}-\hat{\theta}_{P_{1},1})(\hat{\theta}_{P_2,2}-\hat{\theta}_{P_{1},1})}  \hat{u}_{P_{1},1,2}^4 }_{O_p\left(\frac{1}{\theta_1\theta_2 m}\right)}.
 \end{eqnarray*}
 \end{minipage}}\\
 
 \item[1.2) ]i$>$1 : 
\begin{itemize}
\item First, we show a small non-optimal result
 \begin{eqnarray*}
 \sum_{i=2}^m \hat{u}_{P_{1},i,2}^2\hat{u}_{P_{1},i,1}^2= O_p\left( \frac{1}{\theta_1 m^{1/2}} \right).
 \end{eqnarray*} 
 \noindent We easily obtain this result by using inequalities on the sums,
  \begin{eqnarray*}
\sum_{i=2}^m \hat{u}_{P_{1},i,2}^2\hat{u}_{P_{1},i,1}^2& \leqslant & \left( \sum_{i=2}^m \hat{u}_{P_{1},i,2}^4\right)^{1/2} \left( \sum_{i=2}^m \hat{u}_{P_{1},i,1}^4\right)^{1/2}\\
&=& O_p\left(\frac{1}{\theta_1 m^{1/2}} \right).
 \end{eqnarray*}
By Theorem \ref{AAThcomponentdistribution} Part 3, $\sum_{i=2}^m \hat{u}_{P_{1},i,1}^4=O_p\left(\frac{1}{\theta_1^2} \right)$, and the estimation $\sum_{i=2}^m \hat{u}_{P_{1},i,2}^4 =O_p\left( 1/\sqrt{m} \right)$ holds by the spherical property. Indeed, because
 $\hat{u}_{P_{1},i,2:m}$ is invariant by rotation, then \linebreak  $\hat{u}_{P_{1},i,2:m} / ||\hat{u}_{P_{1},i,2:m}|| $ is uniform. Therefore,\\ 
\scalebox{0.8}{
\begin{minipage}{1\textwidth}
 \begin{eqnarray*}
 \E\left[\frac{\hat{u}_{P_{1},i,2}^4}{||\hat{u}_{P_{1},i,2:m}||^4} \right] &=& O_p\left(\frac{1}{m^2} \right) \text{ and }
 \E\left[\frac{ \hat{u}_{P_{1},i,2}^8}{||\hat{u}_{P_{1},i,2:m}||^8} \right] = O_p\left(\frac{1}{m^4} \right).
\end{eqnarray*} 
\end{minipage}}\\
 
We see that $\hat{u}_{P_{1},i,2}^4 \sim {\rm RV}\left( O\left(\frac{1}{m^2} \right) ,O\left(\frac{1}{m^4} \right)  \right)$. Finally, summing the random variables leads to
\begin{eqnarray*}
\E\left[\sum_{i=2}^m \hat{u}_{P_{1},i,2}^4 \right]&=& O_p\left(\frac{1}{m}\right), \\
\var\left(\sum_{i=2}^m \hat{u}_{P_{1},i,2}^4 \right)&=& O_p\left(\frac{1}{m^2} \right).
\end{eqnarray*}

\item We can finally estimate the sum of interest:\\ 
\scalebox{0.71}{
\begin{minipage}{1\textwidth}
 \begin{eqnarray*}
 &&\hspace{-1.6cm}\sum_{i=2}^m \frac{\hat{\lambda}_{P_{1},i}^2}{(\hat{\theta}_{P_2,1}-\hat{\lambda}_{P_{1},i})(\hat{\theta}_{P_2,2}-\hat{\lambda}_{P_{1},i})}  \hat{u}_{P_{1},i,2}^2 \left( 1-\hat{u}_{P_{1},i,1}^2-\hat{u}_{P_{1},i,2}^2\right)\\
 &&\hspace{-1.4cm}=\frac{1}{\hat{\theta}_{P_2,1}\hat{\theta}_{P_2,2}}
 \sum_{i=2}^m \hat{\lambda}_{P_{1},i}^2 \hat{u}_{P_{1},i,2}^2 \left( 1-\hat{u}_{P_{1},i,1}^2-\hat{u}_{P_{1},i,2}^2\right)+O_p\left( \frac{1}{\theta_1 \theta_2^2 } \right)\\
 && \hspace{-1.4cm}= \frac{1}{\hat{\theta}_{P_2,1}\hat{\theta}_{P_2,2}} \left(
 \sum_{i=2}^m \hat{\lambda}_{P_{1},i}^2 \hat{u}_{P_{1},i,2}^2 + O_p\left(1 \right) \sum_{i=2}^m \hat{u}_{P_{1},i,2}^2\hat{u}_{P_{1},i,1}^2+ O_p\left( 1 \right) \sum_{i=2}^m \hat{u}_{P_{1},i,2}^4
  \right)+ O_p\left( \frac{1}{\theta_1 \theta_2^2} \right)\\
 && \hspace{-1.4cm}= \frac{1}{\hat{\theta}_{P_2,1}\hat{\theta}_{P_2,2}} 
 \sum_{i=2}^m \hat{\lambda}_{P_{1},i}^2 \hat{u}_{P_{1},i,2}^2 + O_p\left(\frac{1}{\theta_1 \theta_2} \right) \sum_{i=2}^m \hat{u}_{P_{1},i,2}^2\hat{u}_{P_{1},i,1}^2\\
 &&\hspace{1cm}+ O_p\left(\frac{1}{\theta_1 \theta_2} \right) \sum_{i=2}^m \hat{u}_{P_{1},i,2}^4
  + O_p\left( \frac{1}{\theta_1 \theta_2^2}\right) \\
   && \hspace{-1.4cm}=  \underbrace{\mathbf{\frac{1}{\hat{\theta}_{P_2,1}\hat{\theta}_{P_2,2}} 
 \sum_{i=2}^m \hat{\lambda}_{P_{1},i}^2 \hat{u}_{P_{1},i,2}^2}}_{O_p\left(\frac{1}{ \theta_1 \theta_2} \right)} + O_p\left(\frac{1}{\theta_1^2 \theta_2 m^{1/2} } \right) + O_p\left(\frac{1}{\theta_1 \theta_2 m } \right) 
  + O_p\left( \frac{1}{\theta_1 \theta_2^2}\right).
 \end{eqnarray*}
 \end{minipage}}
\end{itemize}

 \end{itemize}
 \item[] \textbf{Part 2: } As for the previous part, we divide this term.
  \begin{itemize}
  \item[2.1)] $\sum_{i,j=1, i\not=j }^m \frac{\hat{\lambda}_{P_{1},i}\hat{\lambda}_{P_{1},j}}{(\hat{\theta}_{P_2,1}-\hat{\lambda}_{P_{1},i})(\hat{\theta}_{P_2,2}-\hat{\lambda}_{P_{1},j})}  \hat{u}_{P_{1},i,1}\hat{u}_{P_{1},j,1}\hat{u}_{P_{1},i,2} \hat{u}_{P_{1},j,2}$.
  \begin{itemize}
  \item[2.1.1) ]i$=$1,j$>$1 : We want to prove\\ 
\scalebox{0.77}{
\begin{minipage}{1\textwidth}
  \begin{eqnarray*}
  \mathbf{\frac{\hat{\theta}_{P_{1},1}}{(\hat{\theta}_{P_2,1}-\hat{\theta}_{P_{1},1})}  \hat{u}_{P_{1},1,1} \hat{u}_{P_{1},1,2}\sum_{j>1 }^m \frac{\hat{\lambda}_{P_{1},j}}{(\hat{\theta}_{P_2,2}-\hat{\lambda}_{P_{1},j})}  \hat{u}_{P_{1},j,1} \hat{u}_{P_{1},j,2}}=O_p\left( \frac{1}{\theta_2^2} \right).
\end{eqnarray*}
\end{minipage}}\\
   
The order size follows from Theorems \ref{AATheoremcaraceigenstructure}, \ref{AAThcomponentdistribution} and \ref{AAThinvarianteigenvalue},\\ 
\scalebox{0.8}{
\begin{minipage}{1\textwidth}
\begin{eqnarray*}
  \underbrace{\frac{\hat{\theta}_{P_{1},1}}{(\hat{\theta}_{P_2,1}-\hat{\theta}_{P_{1},1})} \hat{u}_{P_{1},1,1}\hat{u}_{P_{1},1,2}}_{O_p\left( \frac{\theta_1^{1/2}m^{1/2}}{\theta_2 } \right)}\underbrace{\sum_{j>1 }^m \frac{\hat{\lambda}_{P_{1},j}}{(\hat{\theta}_{P_2,2}-\hat{\lambda}_{P_{1},j})}  \hat{u}_{P_{1},j,1} \hat{u}_{P_{1},j,2}}_{O_p\left( \frac{1}{\theta_1^{1/2} \theta_2 m^{1/2}} \right)}=O_p\left( \frac{1}{\theta_2^2} \right).
\end{eqnarray*}
\end{minipage}} 
\begin{Rem}\ \\
The Theorem \ref{AATheoremcaraceigenstructure} estimates the order size of the second term for $\hat{\theta}_{P_2,1}$. This same proof is still valid in this new case.
\end{Rem}

\item[2.1.2) ]i$>$1,j$=$1 : Using the fact that $\hat{\lambda}_{P_{1},i}$ is bounded for $i>1$, we find that\\ 
\scalebox{0.74}{
\begin{minipage}{1\textwidth}
\begin{eqnarray*}
\underbrace{\frac{\hat{\theta}_{P_{1},1}}{(\hat{\theta}_{P_2,2}-\hat{\theta}_{P_{1},1})}}_{O_p(1)} \underbrace{ \hat{u}_{P_{1},1,1}}_{O_p(1)} \underbrace{\hat{u}_{P_{1},1,2}}_{O_p\left( \frac{1}{m^{1/2}\theta_1^{1/2}} \right)} \underbrace{\sum_{i>1 }^m \frac{\hat{\lambda}_{P_{1},i}}{(\hat{\theta}_{P_2,1}-\hat{\lambda}_{P_{1},i})}  \hat{u}_{P_{1},i,1} \hat{u}_{P_{1},i,2}}_{O_p\left(  \frac{1}{\theta_1^{3/2}  m^{1/2}} \right)}=O_p\left( \frac{1}{\theta_1^2  m }\right).
\end{eqnarray*}
\end{minipage}}\\

\item[2.1.3) ]i$>$1,j$>$1,i$\not=$j :\\ 
\scalebox{0.7}{
\begin{minipage}{1\textwidth}
\begin{eqnarray*}
&&\hspace{-1.5cm}\left| \sum_{i,j>1,i \not= j }^m \frac{\hat{\lambda}_{P_{1},i}\hat{\lambda}_{P_{1},j}}{(\hat{\theta}_{P_2,1}-\hat{\lambda}_{P_{1},i})(\hat{\theta}_{P_2,2}-\hat{\lambda}_{P_{1},j})}  \hat{u}_{P_{1},i,1}\hat{u}_{P_{1},j,1}\hat{u}_{P_{1},i,2} \hat{u}_{P_{1},j,2}\right| \hspace{30cm}\\
&&\hspace{-0.8cm} \leqslant
\left( \frac{1}{\hat{\theta}_{P_2,1} \hat{\theta}_{P_2,2}}+ O_p\left( \frac{1}{\theta_1 \theta_2^2} \right) \right)
 \left(\sum_{i>1}^m \hat{\lambda}_{P_{1},i} |\hat{u}_{P_{1},i,1}||\hat{u}_{P_{1},i,2}|\right) 
\left(\sum_{j>1}^m \hat{\lambda}_{P_{1},j} |\hat{u}_{P_{1},j,1}| |\hat{u}_{P_{1},j,2}|\right)\\
&&\hspace{-0.8cm} \leqslant \left( \frac{1}{\hat{\theta}_{P_2,1} \hat{\theta}_{P_2,2}}+ O_p\left( \frac{1}{\theta_1 \theta_2^2} \right) \right) 
\hat{\lambda}_{\max}^2 \left(\sum_{i>1}^m \hat{u}_{P_{1},i,1}^2\right)  \left(\sum_{i>1}^m \hat{u}_{P_{1},i,2}^2 \right) \\
&&\hspace{-0.8cm} \leqslant \left( \frac{1}{\hat{\theta}_{P_2,1} \hat{\theta}_{P_2,2}}+ O_p\left( \frac{1}{\theta_1 \theta_2^2} \right) \right) 
\hat{\lambda}_{\max}^2 \left(1-\hat{\alpha}_{P_1,1}^2\right)  \left(1- \hat{u}_{P_{1},1,2}^2 \right) \\
&&\hspace{-0.8cm}= O_p\left( \frac{1}{\theta_1^2 \theta_2} \right).
\end{eqnarray*}
\end{minipage}}\\

  \end{itemize}
  \item[2.2)] Here,\\ 
\scalebox{0.7}{
\begin{minipage}{1\textwidth}
  \begin{eqnarray*}
  &&\hspace{-1.5cm} \sum_{i,j=1, i\not=j }^m \frac{\hat{\lambda}_{P_{1},i}\hat{\lambda}_{P_{1},j}}{(\hat{\theta}_{P_2,1}-\hat{\lambda}_{P_{1},i})(\hat{\theta}_{P_2,2}-\hat{\lambda}_{P_{1},j})}   \hat{u}_{P_{1},i,2}^2 \hat{u}_{P_{1},j,2}^2 \hspace{20cm}\\
  &&\hspace{-1.2cm} =  \sum_{i,j=1}^m \frac{\hat{\lambda}_{P_{1},i}\hat{\lambda}_{P_{1},j}}{(\hat{\theta}_{P_2,1}-\hat{\lambda}_{P_{1},i})(\hat{\theta}_{P_2,2}-\hat{\lambda}_{P_{1},j})}   \hat{u}_{P_{1},i,2}^2 \hat{u}_{P_{1},j,2}^2- 
  \sum_{i=1}^m \frac{\hat{\lambda}_{P_{1},i}^2}{(\hat{\theta}_{P_2,1}-\hat{\lambda}_{P_{1},i})(\hat{\theta}_{P_2,2}-\hat{\lambda}_{P_{1},i})}   \hat{u}_{P_{1},i,2}^4\\
  &&\hspace{-1.2cm} =  \mathbf{\frac{1}{(\theta_2-1)^2}}
  -\underbrace{\frac{\hat{\theta}_{P_{1},1}^2}{(\hat{\theta}_{P_2,1}-\hat{\theta}_{P_{1},1})(\hat{\theta}_{P_2,2}-\hat{\theta}_{P_{1},1})}   \hat{u}_{P_{1},1,2}^4}_{O_p\left( \frac{1}{\theta_1 \theta_2 m} \right)}
  -\underbrace{\sum_{i=2}^m \frac{\hat{\lambda}_{P_{1},i}^2}{(\hat{\theta}_{P_2,1}-\hat{\lambda}_{P_{1},i})(\hat{\theta}_{P_2,2}-\hat{\lambda}_{P_{1},i})}   \hat{u}_{P_{1},i,2}^4}_{O_p\left(\frac{1}{\theta_1 \theta_2 m}\right)}.
  \end{eqnarray*}
  \end{minipage}}\\
  
  \end{itemize}
 \end{itemize}

Combining the two parts leads to\\ 
\scalebox{0.8}{
\begin{minipage}{1\textwidth}
\begin{eqnarray*}
&&\hspace{-1cm} \sum_{s=3}^m \hat{u}_{P_2,1,s}  \hat{u}_{P_2,2,s}\\  
&&\hspace{-0.5cm}= \frac{1}{\sqrt{D_1 D_2}N_1 N_2 } \left(\frac{\hat{\theta}_{P_{1},1}^2}{(\hat{\theta}_{P_2,1}-\hat{\theta}_{P_{1},1})(\hat{\theta}_{P_2,2}-\hat{\theta}_{P_{1},1})}  \hat{u}_{P_{1},1,2}^2 \left( 1-\hat{\alpha}_{P_1,1}^2\right)
+ \frac{1}{\hat{\theta}_{P_2,1}\hat{\theta}_{P_2,2}} 
 \sum_{i=2}^m \hat{\lambda}_{P_{1},i}^2 \hat{u}_{P_{1},i,2}^2 \right.\\
&&\hspace{1.5cm} \left.-
 \frac{\hat{\theta}_{P_{1},1}}{(\hat{\theta}_{P_2,1}-\hat{\theta}_{P_{1},1})}  \hat{u}_{P_{1},1,1} \hat{u}_{P_{1},1,2}\frac{1}{\hat{\theta}_{P_2,2}}\sum_{j>1 }^m \hat{\lambda}_{P_{1},j} \hat{u}_{P_{1},j,1} \hat{u}_{P_{1},j,2}
 -\frac{1}{(\theta_2-1)^2}\right)\\
 &&\hspace{1.5cm}+O_p\left(\frac{\theta_2^{1/2}}{\theta_1^{3/2} m^{3/2}}\right)+O_p\left(\frac{1}{\theta_1^{3/2} \theta_2^{1/2} m^{1/2}} \right)\\
 &&\hspace{-0.5cm}= \frac{|\hat{\theta}_{P_{2},1}-\hat{\theta}_{P_{1},1}|}{\hat{\theta}_{P_{1},1}  |\hat{u}_{P_{1},1,2}|} \sqrt{\theta_2-1}
  \left(\frac{\hat{\theta}_{P_{1},1}^2}{(\hat{\theta}_{P_2,1}-\hat{\theta}_{P_{1},1})(\hat{\theta}_{P_2,2}-\hat{\theta}_{P_{1},1})}  \hat{u}_{P_{1},1,2}^2 \left( 1-\hat{\alpha}_{P_1,1}^2\right)\right.\\
&&\hspace{1.5cm}+ \frac{1}{\hat{\theta}_{P_2,1}\hat{\theta}_{P_2,2}} 
 \sum_{i=2}^m \hat{\lambda}_{P_{1},i}^2 \hat{u}_{P_{1},i,2}^2 \\
&&\hspace{1.5cm}\left. -
 \frac{\hat{\theta}_{P_{1},1}}{(\hat{\theta}_{P_2,1}-\hat{\theta}_{P_{1},1})}  \hat{u}_{P_{1},1,1} \hat{u}_{P_{1},1,2}\frac{1}{\hat{\theta}_{P_2,2}}\sum_{j>1 }^m \hat{\lambda}_{P_{1},j} \hat{u}_{P_{1},j,1} \hat{u}_{P_{1},j,2}-\frac{1}{(\theta_2-1)^2}\right)\\
 &&\hspace{1.5cm}+O_p\left(\frac{\theta_2^{1/2}}{\theta_1^{3/2} m^{3/2}}\right)+O_p\left(\frac{1}{\theta_1^{3/2} \theta_2^{1/2} m^{1/2}} \right)+O_p\left(\frac{1}{\theta_1^{1/2} \theta_2^{5/2} m^{1/2}} \right).
  \end{eqnarray*}
  \end{minipage}}\\
  
\noindent In this second part we simplify the terms using Theorems \ref{AAThcomponentdistribution} and \ref{AAThinvarianteigenvalue}:
  $$\hat{u}_{P_k,1,k}= \frac{\sqrt{\theta_k}\theta_1}{|\theta_k-\theta_1|} |\hat{u}_{P_{k-1},1,k}|+O_p\left(\frac{\min(\theta_1,\theta_k)}{\theta_1^{1/2}\theta_k^{1/2}m} \right)+O_p\left(\frac{1}{ \theta_1^{1/2} \theta_2^{1/2}m^{1/2}} \right)$$
  and 
  $$  \hat{\theta}_{P_{2},1}-\hat{\theta}_{P_{1},1} =-\frac{\hat{\theta}_{P_{1},1} \hat{\theta}_{P_{2},1} (\theta_2-1)}{\theta_2-1 -\hat{\theta}_{P_{2},1}} \hat{u}_{P_{1},1,2}^2+O_p\left(\frac{\theta_2}{m^{3/2}} \right) +O_p\left( \frac{1}{m}\right).$$
   Recall that without the convention $\hat{u}_{P_k,1,1}>0$, by construction we have $\hat{u}_{P_k,1,k}>0$. Because $\theta_1>\theta_2$,
  \begin{itemize}
  \item{$P1:$ } \\ 
\scalebox{0.85}{
\begin{minipage}{1\textwidth} 
  \begin{eqnarray*}
   P_1&=&\frac{|\hat{\theta}_{P_{2},1}-\hat{\theta}_{P_{1},1}|}{\hat{\theta}_{P_{1},1}  |\hat{u}_{P_{1},1,2}|} \sqrt{\theta_2-1}  \frac{\hat{\theta}_{P_{1},1}^2}{(\hat{\theta}_{P_2,1}-\hat{\theta}_{P_{1},1})(\hat{\theta}_{P_2,2}-\hat{\theta}_{P_{1},1})}  \hat{u}_{P_{1},1,2}^2 \left( 1-\hat{\alpha}_{P_1,1}^2\right) \\
 &\overset{Asy}{=}& \sqrt{\theta_2-1}  \frac{\hat{\theta}_{P_{1},1}}{\hat{\theta}_{P_2,2}-\hat{\theta}_{P_{1},1}}  |\hat{u}_{P_{1},1,2}| \left( 1-\hat{\alpha}_{P_1,1}^2\right) \\
 &\overset{Asy}{=} & -\hat{u}_{P_2,1,2}\left( 1-\hat{\alpha}_{P_1,1}^2\right) + O_p\left(\frac{\theta_2^{1/2}}{\theta_1^{3/2} m} \right)+ O_p\left(\frac{1}{\theta_1^{3/2} \theta_2^{1/2} m^{1/2}} \right).
  \end{eqnarray*}
  \end{minipage}}\\
  
We use the notation $\overset{Asy}{=}$ because the probability that the sign is wrong tends to $0$ in $1/m$ when $\theta_1$ tends to infinity. Moreover, when $\theta_1$ is finite, the order size is $1/\sqrt{m}$.
   \item{$P2:$ } \\ 
\scalebox{0.85}{
\begin{minipage}{1\textwidth}
  \begin{eqnarray*}
  P_2&=&\frac{|\hat{\theta}_{P_{2},1}-\hat{\theta}_{P_{1},1}|}{\hat{\theta}_{P_{1},1}  |\hat{u}_{P_{1},1,2}|} \sqrt{\theta_2-1}  
  \frac{1}{\hat{\theta}_{P_2,1}\hat{\theta}_{P_2,2}} 
 \sum_{i=2}^m \hat{\lambda}_{P_{1},i}^2 \hat{u}_{P_{1},i,2}^2\\
 &=&\frac{\hat{\theta}_{P_{1},1} \hat{\theta}_{P_{2},1} \theta_2}{|\theta_2 -\hat{\theta}_{P_{2},1}|} \hat{u}_{P_{1},1,2}^2 \frac{1}{\hat{\theta}_{P_{1},1}  |\hat{u}_{P_{1},1,2}|} \sqrt{\theta_2-1}  
  \frac{1}{\hat{\theta}_{P_2,1}\hat{\theta}_{P_2,2}} 
 \sum_{i=2}^m \hat{\lambda}_{P_{1},i}^2 \hat{u}_{P_{1},i,2}^2 \\
 &&\hspace{1cm}+ O_p\left(\frac{\theta_2^{1/2}}{\theta_1^{3/2} m} \right)+ O_p\left(\frac{1}{\theta_1^{3/2} \theta_2^{1/2} m^{1/2}} \right)\\
 &=&\frac{|\hat{u}_{P_{1},1,2}|  \sqrt{\theta_2-1}   }{|\theta_2 -\hat{\theta}_{P_{2},1}|} 
 \sum_{i=2}^m \hat{\lambda}_{P_{1},i}^2 \hat{u}_{P_{1},i,2}^2+ O_p\left(\frac{\theta_2^{1/2}}{\theta_1^{3/2} m} \right)+ O_p\left(\frac{1}{\theta_1^{3/2} \theta_2^{1/2} m^{1/2}} \right)\\
 &=&\frac{|\hat{u}_{P_{2},1,2}| }{\theta_1} 
 \sum_{i=2}^m \hat{\lambda}_{P_{1},i}^2 \hat{u}_{P_{1},i,2}^2 + O_p\left(\frac{\theta_2^{1/2}}{\theta_1^{3/2} m} \right)+ O_p\left(\frac{1}{\theta_1^{3/2} \theta_2^{1/2} m^{1/2}} \right).\\
  \end{eqnarray*}
  \end{minipage}}\\
 \item{$P3:$ }  Using Lemma \ref{AALemstatW},\\ 
\scalebox{0.73}{
\begin{minipage}{1\textwidth}
  \begin{eqnarray*}
 P_3&=&\frac{|\hat{\theta}_{P_{2},1}-\hat{\theta}_{P_{1},1}|}{\hat{\theta}_{P_{1},1}  |\hat{u}_{P_{1},1,2}|} \sqrt{\theta_2-1}  \frac{\hat{\theta}_{P_{1},1}}{(\hat{\theta}_{P_2,1}-\hat{\theta}_{P_{1},1})}  \hat{u}_{P_{1},1,1} \hat{u}_{P_{1},1,2}\frac{1}{\hat{\theta}_{P_2,2}}\sum_{j>1 }^m \hat{\lambda}_{P_{1},j} \hat{u}_{P_{1},j,1} \hat{u}_{P_{1},j,2}\\
 &\overset{Asy}{=}& \text{sign}\left( \hat{u}_{P_{2},1,1} \right) \frac{1}{\theta_2^{1/2}} \sum_{j>1 }^m \hat{\lambda}_{P_{1},j} \hat{u}_{P_{1},j,1} \hat{u}_{P_{1},j,2} +O_p\left( \frac{1}{\theta_1^{1/2} \theta_2^{3/2} m^{1/2}} \right)+O_p\left( \frac{1}{\theta_1^{1/2} \theta_2^{1/2} m} \right),
  \end{eqnarray*}
  \end{minipage}}\\
  
  where the sign equality is obtained by the remark of Theorem \ref{AATheoremcaraceigenstructure} and tends to be correct in $1/m$.
  \item{$P4:$ }  
  \begin{eqnarray*}
 P_4&=& \frac{1}{\sqrt{D_1}}\sqrt{\theta_2-1} \frac{1}{(\theta_2-1)^2}\\
&=&\frac{\sqrt{\theta_2-1} }{\theta_2-1}\tilde{u}_{P_2,1,2}\\
&=&\frac{1}{\theta_2-1}\hat{u}_{P_2,1,2}\\
&=&\frac{1}{\theta_2}\hat{u}_{P_2,1,2}+O_p\left(\frac{1}{\theta_1^{1/2} \theta_2^{3/2} m^{1/2}}\right)
  \end{eqnarray*}
  \end{itemize}
  
\noindent By construction we know that $\hat{u}_{P_2,1,2}>0$, but this is not the case for $\hat{u}_{P_2,1,1}$. We will correct this later, but first combine $P_1+P_2-P_4$ to obtain\\ 
\scalebox{0.92}{
\begin{minipage}{1\textwidth}
 \begin{eqnarray*}
 P_1+P_2-P_4 &\overset{Asy}{=}& \hat{u}_{P_2,1,2} \left( -\left( 1-\hat{\alpha}_{P_1,1}^2\right) + \frac{ \sum_{i=2}^m \hat{\lambda}_{P_{1},i}^2 \hat{u}_{P_{1},i,2}^2}{\theta_1} - \frac{1}{\theta_2}\right) \\
 &&\hspace{1cm}+ O_p\left(\frac{\theta_2^{1/2}}{\theta_1^{3/2} m} \right)+ O_p\left(\frac{1}{\theta_1^{1/2} \theta_2^{3/2} m^{1/2}} \right)\\
 &\overset{Asy}{=}&\hat{u}_{P_2,1,2} \left( -\frac{ \sum_{i=1}^m \hat{\lambda}_{W,i}^2\hat{u}_{W,i,1}^2-1}{\theta_1} + \frac{ \sum_{i=2}^m \hat{\lambda}_{P_{1},i}^2 \hat{u}_{P_{1},i,2}^2}{\theta_1} - \frac{1}{\theta_2}\right) \\
 &&\hspace{1cm}+ O_p\left(\frac{\theta_2^{1/2}}{\theta_1^{3/2} m} \right)+ O_p\left(\frac{1}{\theta_1^{1/2} \theta_2^{3/2} m^{1/2}} \right)\\
  &\overset{Asy}{=}&\hat{u}_{P_2,1,2} \left(  \frac{ 1}{\theta_1} - \frac{1}{\theta_2}\right) + O_p\left(\frac{\theta_2^{1/2}}{\theta_1^{3/2} m} \right)+ O_p\left(\frac{1}{\theta_1^{1/2} \theta_2^{3/2} m^{1/2}} \right).
\end{eqnarray*}
\end{minipage}}\\
  
Indeed Lemma \ref{AALemstatW} shows that 
\begin{eqnarray*}
&&\sum_{i=1}^m \hat{\lambda}_{W,i}^2\hat{u}_{W,i,1}^2= \left(W^2\right)_{1,1},\\
&&\sum_{i=2}^m \hat{\lambda}_{P_{1},i}^2\hat{u}_{P_{1},i,2}^2= \left(W^2\right)_{2,2}+O_p\left( \frac{1}{m} \right).
\end{eqnarray*}
The result follows by invariance of $W^2$ under rotation.\\

\noindent Finally, we combine the different parts \\ 
\scalebox{0.85}{
\begin{minipage}{1\textwidth}
\begin{eqnarray*}
P_1+P_2-P_3-P_4&\overset{Asy}{=}& \hat{u}_{P_2,1,2} \left(  \frac{ 1}{\theta_1} - \frac{1}{\theta_2}\right) -
\text{sign}\left( \hat{u}_{P_{2},1,1} \right) \frac{1}{\theta_2^{1/2}} \sum_{j>1 }^m \hat{\lambda}_{P_{1},j} \hat{u}_{P_{1},j,1} \hat{u}_{P_{1},j,2} \\
&&\hspace{1cm} + O_p\left(\frac{1}{\theta_1^{1/2}\theta_2^{1/2} m} \right)+ O_p\left(\frac{1}{\theta_1^{1/2} \theta_2^{3/2} m^{1/2}} \right),
\end{eqnarray*}
\end{minipage}}\\

where the asymptotic equality is discussed in Remark \ref{AAREMasy}.\\
We change the convention of the sign such that $\hat{u}_{P_{2},i,i}>0$, $i=1,2$. Therefore,  we multiply by $\text{sign}\left( \hat{u}_{P_{2},1,1} \right)$. With this convention $\hat{u}_{P_2,1,2}$ is no longer strictly positive. Nevertheless, we keep using the same notation. 
\begin{eqnarray*}
P_1+P_2-P_3-P_4&=& \hat{u}_{P_2,1,2} \left(  \frac{ 1}{\theta_1} - \frac{1}{\theta_2}\right) -
 \frac{1}{\theta_2^{1/2}} \sum_{j>1 }^m \hat{\lambda}_{P_{1},j} \hat{u}_{P_{1},j,1} \hat{u}_{P_{1},j,2} \\
&& \hspace{1cm}+ O_p\left(\frac{1}{\theta_1^{1/2}\theta_2^{1/2} m} \right)+ O_p\left(\frac{1}{\theta_1^{1/2} \theta_2^{3/2} m^{1/2}} \right).
\end{eqnarray*}

\begin{Rem}\ \label{AAREMasy} \\
First, we recall that the $O$ errors are in probability and take care of this possible fluctuation with probability tending to $0$.\\
The simplification of $P_1+P_2-P_3-P_4$ is possible thanks to the remark of Theorem \ref{AATheoremcaraceigenstructure} showing that the signs are correct with probability tending to $1$ in $1/m$ when $\theta_2$ is large. In particular, there is a probability of order $1/m$ to have an error of size $O_p\left(\frac{1}{\theta_1^{1/2} \theta_2^{1/2} m^{1/2}} \right)$. Luckily this rare error will not affect the moment estimation of the statistic. \\
Then, when $\theta_2$ is finite, the formula just provides order size.
\end{Rem} 

\noindent This estimation concludes part $\mathbf{B}$. 

\paragraph*{C: }
In this section we express \\ 
\scalebox{0.8}{
\begin{minipage}{1\textwidth}
\begin{equation*}
\hat{u}_{P_2,1,2} \left(  \frac{ 1}{\theta_1} - \frac{1}{\theta_2}\right)\delta  -
 \frac{1}{\theta_2^{1/2}} \sum_{j>1 }^m \hat{\lambda}_{P_{1},j} \hat{u}_{P_{1},j,1} \hat{u}_{P_{1},j,2} + O_p\left(\frac{1}{\theta_1^{1/2}\theta_2^{1/2} m} \right)+ O_p\left(\frac{1}{\theta_1^{1/2} \theta_2^{3/2} m^{1/2}} \right)
 \end{equation*}
 \end{minipage}}\\
 
as a function of the unit statistic defined in Theorem \ref{AAThunitstatistic}. 
Using Theorem \ref{AAThcomponentdistribution} and Lemma \ref{AALemstatW} leads to the following estimations,\\ 
\scalebox{0.8}{
\begin{minipage}{1\textwidth}
\begin{eqnarray*}
&\hat{u}_{P_2,1,2}&= \frac{\sqrt{\theta_2}\theta_1}{|\theta_2-\theta_1|} \hat{u}_{P_{1},1,2} +O_p\left(\frac{\theta_2^{1/2}}{\theta_1^{1/2} m} \right)+O_p\left(\frac{1}{\theta_1^{1/2} \theta_2^{1/2} m^{1/2}} \right),\\
&\hat{u}_{P_1,1,2} &=\frac{W_{1,2}}{\sqrt{\theta_1} }+ O_p\left( \frac{1}{\theta_1^{3/2}m^{1/2}} \right)+ O_p\left( \frac{1}{\theta_1^{1/2}m} \right),\\
&\sum_{i=2}^m \hat{\lambda}_{P_1,i} \hat{u}_{P_1,i,1}\hat{u}_{P_1,i,2}&=W_{1,2} \frac{  M_2}{\sqrt{\theta_1}} -\left(W^2\right)_{1,2}\frac{1}{\sqrt{\theta_1}}
+ O_p\left( \frac{1}{\theta_1^{1/2}m} \right)+ O_p\left( \frac{1}{\theta_1^{3/2}m^{1/2}} \right).
\end{eqnarray*}
\end{minipage}}\\

Based on this, we can show that \\ 
\scalebox{0.85}{
\begin{minipage}{1\textwidth}
\begin{eqnarray*}
&&\hat{u}_{P_2,1,2} \left(  \frac{ 1}{\theta_1} - \frac{1}{\theta_2}\right) \delta -
 \frac{1}{\theta_2^{1/2}} \sum_{j>1 }^m \hat{\lambda}_{P_{1},j} \hat{u}_{P_{1},j,1} \hat{u}_{P_{1},j,2} \\
 &&\hspace{1cm}= \frac{-\left(\delta +M_2 \right) W_{1,2} +\left(W^2\right)_{1,2}}{\sqrt{\theta_1 \theta_2}}+ O_p\left(\frac{1}{\theta_1^{1/2}\theta_2^{1/2} m} \right)+ O_p\left(\frac{1}{\theta_1^{1/2} \theta_2^{3/2} m^{1/2}} \right).
\end{eqnarray*} 
\end{minipage}}\\

 The result is straightforward using a delta method and Theorem \ref{AAThunitstatistic}.
\end{proof}

\subsubsection{Invariant Dot Product}

\noindent (Page \pageref{AAThInvariantdot})

\begin{proof}{\textbf{Theorem \ref{AAThInvariantdot}}}\label{AAproofThInvariantdot}  
\noindent We start this proof with two important remarks.
\begin{itemize}
\item  This proof will assume the sign convention of Theorem \ref{AATheoremcaraceigenstructure}. We will correct for this at the end of the proof.
\item We use the notation \ref{AANot=thetadifferent} to prove the result based only on  $\theta_1 > \theta_2$ and relaxing the order of the other eigenvalues. This notation permutes the estimated eigenvalues and their eigenvectors, but the reader can also read this proof as if  $\theta_1>\theta_2>...>\theta_k$ and realize that the notation allows for this generalisation. Moreover, we add the notation $\hat{\lambda}_{P_{r},i}=\hat{\theta}_{P_r,i}$ for $i=1,2,...,r$ in order to simplify formulas.
\end{itemize}
Theorem \ref{AATheoremcaraceigenstructure} leads to \\ 
\scalebox{0.72}{
\begin{minipage}{1\textwidth}
\begin{eqnarray*}
\hat{u}_{P_k,1,s} \hat{u}_{P_k,2,s}= \frac{1}{\sqrt{D_1 D_2}N_1N_2}\sum_{i,j} \frac{\hat{\lambda}_{P_{k-1},i}\hat{\lambda}_{P_{k-1},j}}{\left(\hat{\theta}_{P_k,1}-\hat{\lambda}_{P_{k-1},i}\right)\left(\hat{\theta}_{P_k,2}-\hat{\lambda}_{P_{k-1},j}\right)} \hat{u}_{P_{k-1},i,k} \hat{u}_{P_{k-1},j,k} \hat{u}_{P_{k-1},i,s} \hat{u}_{P_{k-1},j,s},
\end{eqnarray*}
\end{minipage}}\\

where $N_1$ and $N_2$ are scalars such that the vectors are of unit length. It then follows that\\ 
\scalebox{0.74}{
\begin{minipage}{1\textwidth}
\begin{eqnarray*}
\sum_{s=k+1}^m \hat{u}_{P_k,1,s} \hat{u}_{P_k,2,s}&=& \\
&&\hspace{-4cm} \frac{1}{\sqrt{D_1 D_2}N_1N_2} 
 \underbrace{\left( \sum_{i \not = j} \frac{\hat{\lambda}_{P_{k-1},i}\hat{\lambda}_{P_{k-1},j}}{\left(\hat{\theta}_{P_k,1}-\hat{\lambda}_{P_{k-1},i}\right)\left(\hat{\theta}_{P_k,2}-\hat{\lambda}_{P_{k-1},j}\right)}  \hat{u}_{P_{k-1},i,k} \hat{u}_{P_{k-1},j,k} \left(- \sum_{r=1}^k \hat{u}_{P_{k-1},i,r} \hat{u}_{P_{k-1},j,r} \right) \right.}_{\text{Part 2}}\\
&&\hspace{0.8cm}\underbrace{\left. \sum_{i=1}^m \frac{\hat{\lambda}_{P_{k-1},i}^2 }{\left(\hat{\theta}_{P_k,1}-\hat{\lambda}_{P_{k-1},i}\right)\left(\hat{\theta}_{P_k,2}-\hat{\lambda}_{P_{k-1},i}\right)}  \hat{u}_{P_{k-1},i,k}^2 \left(1- \sum_{r=1}^k \hat{u}_{P_{k-1},i,r}^2 \right)\right)}_{\text{Part 1}}.
\end{eqnarray*}
\end{minipage}}\\

\noindent First we will study Part 1 and Part 2 in \textbf{A}. Then in \textbf{B}, we will show 
\begin{eqnarray*}
 \frac{1}{\sqrt{D_1 D_2}N_1N_2}=O_p\left( \frac{\min(\theta_1,\theta_k) \min(\theta_2,\theta_k)}{\theta_1^{1/2} \theta_2^{1/2} m} \right).
\end{eqnarray*}
Finally, in part \textbf{C}, we combine \textbf{A} and \textbf{B} to conclude the proof.
\begin{itemize}
\item[\textbf{A:}] 
Assuming the previous estimation, we can neglect all the terms of order $o_p\left(\frac{\sqrt{m}}{\min(\theta_1,\theta_k) \min(\theta_2,\theta_k)} \right)$ in Part 1 and 2. The order sizes of the elements are obtained using Theorems \ref{AAThcomponentdistribution}, \ref{AAjointdistribution}, \ref{AAThinvarianteigenvalue}, \ref{AATheoremcaraceigenstructure}, the Invariant Angle Theorem \ref{AAInvariantth}, the Dot Product Theorem \ref{AAThDotproduct} and its Invariant Theorem \ref{AAThInvariantdot}.

\paragraph*{Part 1 :} We will show that we can neglect this entire part.
\begin{itemize}
\item[1.1)]  $i=1$ : Assuming without loss of generality that $\theta_1<\theta_2$ leads to \\ 
\scalebox{0.69}{
\begin{minipage}{1\textwidth}
\begin{eqnarray*}
&&\hspace{-1cm}\frac{\hat{\theta}_{P_{k-1},1}^2 }{\left(\hat{\theta}_{P_k,1}-\hat{\theta}_{P_{k-1},1}\right)\left(\hat{\theta}_{P_k,2}-\hat{\theta}_{P_{k-1},1}\right)}  \hat{u}_{P_{k-1},1,k}^2 \left(1- \sum_{r=1}^k \hat{u}_{P_{k-1},1,r}^2 \right)\\
&& \hspace{1cm}= 
\underbrace{\frac{\hat{\theta}_{P_{k-1},1}^2 }{\left(\hat{\theta}_{P_k,1}-\hat{\theta}_{P_{k-1},1}\right)\left(\hat{\theta}_{P_k,2}-\hat{\theta}_{P_{k-1},1}\right)} \hat{u}_{P_{k-1},1,k}^2}_{O_p\left( \frac{1}{\min(\theta_1,\theta_k)} \right)} \left(\underbrace{1- \sum_{r=1}^{k-1} \hat{u}_{P_{k-1},1,r}^2 }_{=(1-\hat{\alpha}_{P{k-1},1}^2)=O_p\left( \frac{1}{\theta_1}\right)}- \underbrace{\hat{u}_{P_{k-1},1,k}^2}_{O_p\left( \frac{1}{\theta_1 m} \right)}\right)\\
&& \hspace{1cm}= O_p\left( \frac{1}{\theta_1 \min(\theta_1,\theta_k)} \right).
\end{eqnarray*}
\end{minipage}}\\

\item[1.2)] $i=2$ :\\ 
\scalebox{0.67}{
\begin{minipage}{1\textwidth}
\begin{eqnarray*}
&&\frac{\hat{\theta}_{P_{k-1},2}^2 }{\left(\hat{\theta}_{P_k,1}-\hat{\theta}_{P_{k-1},2}\right)\left(\hat{\theta}_{P_k,2}-\hat{\theta}_{P_{k-1},2}\right)}  \hat{u}_{P_{k-1},2,k}^2 \left(1- \sum_{r=1}^k \hat{u}_{P_{k-1},2,r}^2 \right)= O_p\left( \frac{1}{\theta_1 \min(\theta_2,\theta_k) } \right).
\end{eqnarray*}
\end{minipage}}\\

\item[1.3)] $i=3,...,k-1$ : \\ 
\scalebox{0.65}{
\begin{minipage}{1\textwidth}
\begin{eqnarray*}
&&\hspace{-0.5cm}\frac{\hat{\theta}_{P_{k-1},i}^2 }{\left(\hat{\theta}_{P_k,1}-\hat{\theta}_{P_{k-1},i}\right)\left(\hat{\theta}_{P_k,2}-\hat{\theta}_{P_{k-1},i}\right)}  \hat{u}_{P_{k-1},i,k}^2 \left(1- \sum_{r=1}^k \hat{u}_{P_{k-1},i,r}^2 \right)= O_p\left( \frac{1}{\max(\theta_1,\theta_i) \max(\theta_2,\theta_i) m} \right).
\end{eqnarray*}
\end{minipage}}\\

\item[1.4)] $i \geqslant k$ : \\ 
\scalebox{0.71}{
\begin{minipage}{1\textwidth}
\begin{eqnarray*}
&&\frac{\hat{\lambda}_{P_{k-1},i}^2 }{\left(\hat{\theta}_{P_k,1}-\hat{\lambda}_{P_{k-1},i}\right)\left(\hat{\theta}_{P_k,2}-\hat{\lambda}_{P_{k-1},i}\right)}  \hat{u}_{P_{k-1},i,k}^2 \left(1- \sum_{r=1}^k \hat{u}_{P_{k-1},i,r}^2 \right)= O_p\left( \frac{1}{\theta_1 \theta_2 m} \right)\\
&&\Rightarrow \sum_{i=k}^m \frac{\hat{\lambda}_{P_{k-1},i}^2 }{\left(\hat{\theta}_{P_k,1}-\hat{\lambda}_{P_{k-1},i}\right)\left(\hat{\theta}_{P_k,2}-\hat{\lambda}_{P_{k-1},i}\right)}  \hat{u}_{P_{k-1},i,k}^2 \left(1- \sum_{r=1}^k \hat{u}_{P_{k-1},i,r}^2 \right)= O_p\left( \frac{1}{\theta_1 \theta_2 } \right).
\end{eqnarray*}
\end{minipage}}
\end{itemize}
\paragraph*{Part 2 :} The second part is trickier but, again, many elements can be neglected.
\item[2.1)] $i \not = j \geqslant k$ :  By the previous part, if $i=j\geqslant k$, then the sum is $O_p\left( \frac{1}{\theta_1 \theta_2} \right)$.\\ 
\scalebox{0.66}{
\begin{minipage}{1\textwidth}
\begin{eqnarray*} 
&&\hspace{-0.4cm} \left| \sum_{i \not = j \geqslant k} \frac{\hat{\lambda}_{P_{k-1},i}\hat{\lambda}_{P_{k-1},j}}{\left(\hat{\theta}_{P_k,1}-\hat{\lambda}_{P_{k-1},i}\right)\left(\hat{\theta}_{P_k,2}-\hat{\lambda}_{P_{k-1},j}\right)}  \hat{u}_{P_{k-1},i,k} \hat{u}_{P_{k-1},j,k} \left(- \sum_{r=1}^k \hat{u}_{P_{k-1},i,r} \hat{u}_{P_{k-1},j,r} \right)+O_p\left(\frac{1}{\theta_1 \theta_2} \right)  \right|\\
&&\hspace{1cm} = \left| \sum_{i ,j \geqslant k} \frac{\hat{\lambda}_{P_{k-1},i}\hat{\lambda}_{P_{k-1},j}}{\left(\hat{\theta}_{P_k,1}-\hat{\lambda}_{P_{k-1},i}\right)\left(\hat{\theta}_{P_k,2}-\hat{\lambda}_{P_{k-1},j}\right)}  \hat{u}_{P_{k-1},i,k} \hat{u}_{P_{k-1},j,k} \left(- \sum_{r=1}^k \hat{u}_{P_{k-1},i,r} \hat{u}_{P_{k-1},j,r} \right) \right|\\
&&\hspace{1cm} = \left| \sum_{r=1}^k \sum_{i ,j \geqslant k} \frac{\hat{\lambda}_{P_{k-1},i}\hat{\lambda}_{P_{k-1},j}}{\left(\hat{\theta}_{P_k,1}-\hat{\lambda}_{P_{k-1},i}\right)\left(\hat{\theta}_{P_k,2}-\hat{\lambda}_{P_{k-1},j}\right)}  \hat{u}_{P_{k-1},i,k} \hat{u}_{P_{k-1},j,k} \left(- \hat{u}_{P_{k-1},i,r} \hat{u}_{P_{k-1},j,r} \right) \right|\\
&&\hspace{1cm}  \leqslant O_p\left(1 \right) \times\sum_{r=1}^k \frac{1}{\hat{\theta}_{P_k,2} \hat{\theta}_{P_k,1}}  \left( \sum_{i  \geqslant k} \left|\hat{\lambda}_{P_{k-1},i}
  \hat{u}_{P_{k-1},i,k}   \hat{u}_{P_{k-1},i,r} \right|\right)^2 
  \\
  &&\hspace{1cm} \leqslant O_p\left(1 \right) \times \sum_{r=1}^k \frac{1}{\hat{\theta}_{P_k,2} \hat{\theta}_{P_k,1}}  \left( \sum_{i  \geqslant k} \hat{\lambda}_{P_{k-1},i}^2
  \hat{u}_{P_{k-1},i,k}^2\right)   \left( \sum_{i  \geqslant k}\hat{u}_{P_{k-1},i,r}^2  \right)  \\
  &&\hspace{1cm} \leqslant O_p\left(1 \right) \times \sum_{r=1}^k \frac{\lambda_{\max}^2}{\hat{\theta}_{P_k,2} \hat{\theta}_{P_k,1}}  \left( \sum_{i  \geqslant k} 
  \hat{u}_{P_{k-1},i,k}^2\right)   \left( \sum_{i  \geqslant k}\hat{u}_{P_{k-1},i,r}^2  \right)  \\
&&\hspace{1cm} = O_p\left( \frac{1}{\theta_1 \theta_2}  \right).
\end{eqnarray*}
\end{minipage}}\\

\item[2.2)] $i=2,...,k-1$, $j\geqslant k$ : 
\begin{itemize}
\item[2.2.1)] $r=1,...,k-1$: \\ 
\scalebox{0.68}{
\begin{minipage}{1\textwidth}
\begin{eqnarray*} 
&&\hspace{-0.5cm}\left| \sum_{j \geqslant k} \frac{\hat{\lambda}_{P_{k-1},i}\hat{\lambda}_{P_{k-1},j}}{\left(\hat{\theta}_{P_k,1}-\hat{\lambda}_{P_{k-1},i}\right)\left(\hat{\theta}_{P_k,2}-\hat{\lambda}_{P_{k-1},j}\right)}  \hat{u}_{P_{k-1},i,k} \hat{u}_{P_{k-1},j,k} \left(- \sum_{r=2}^{k-1} \hat{u}_{P_{k-1},i,r} \hat{u}_{P_{k-1},j,r} \right) \right| \\
&&\hspace{-0.2cm} \leqslant \sum_{j \geqslant k} \frac{\hat{\lambda}_{P_{k-1},i}\hat{\lambda}_{P_{k-1},j}}{\left(\hat{\theta}_{P_k,1}-\hat{\lambda}_{P_{k-1},i}\right)\left(\hat{\theta}_{P_k,2}-\hat{\lambda}_{P_{k-1},j}\right)}  \left|\hat{u}_{P_{k-1},i,k}\right| \left|\hat{u}_{P_{k-1},j,k}\right| \left( \sum_{r=2}^{k-1} \left|\hat{u}_{P_{k-1},i,r}\right| \left|\hat{u}_{P_{k-1},j,r}\right| \right) \\
 &&\hspace{-0.2cm} \leqslant O_p\left( \frac{1}{\theta_1 \theta_2} \right) \sum_{r=2}^{k-1} \hat{\lambda}_{P_{k-1},i}  \underbrace{\left|\hat{u}_{P_{k-1},i,k} \hat{u}_{P_{k-1},i,r}\right|}_{O_p\left( \frac{1}{ m^{1/2}} \right) }  \underbrace{\sum_{j=k}^m \hat{\lambda}_{P_{k-1},j}
 \left|\hat{u}_{P_{k-1},j,k} \hat{u}_{P_{k-1},j,r}\right|}_{O_p\left( m^{1/2} \right)}\\
 &&\hspace{-0.2cm} = O_p\left( \frac{1}{\theta_1 \theta_2} \right).
 \end{eqnarray*} 
 \end{minipage}}\\
 
 The size could be improved; however, this estimation is enough to justify neglecting the term.
\item [2.2.2)] $r=k$: \\ 
\scalebox{0.76}{
\begin{minipage}{1\textwidth} \begin{eqnarray*}
&&\left| \sum_{j \geqslant k} \frac{\hat{\lambda}_{P_{k-1},i}\hat{\lambda}_{P_{k-1},j}}{\left(\hat{\theta}_{P_k,1}-\hat{\lambda}_{P_{k-1},i}\right)\left(\hat{\theta}_{P_k,2}-\hat{\lambda}_{P_{k-1},j}\right)}  \hat{u}_{P_{k-1},i,k} \hat{u}_{P_{k-1},j,k} \left(-\hat{u}_{P_{k-1},i,k} \hat{u}_{P_{k-1},j,k} \right) \right|\\
&&\hspace{1cm} \leqslant O_p\left( \frac{1}{\theta_1 \theta_2} \right)  \hat{\lambda}_{P_{k-1},i}  \underbrace{\hat{u}_{P_{k-1},i,k}^2}_{O_p\left( \frac{1}{\theta_i m} \right)}  \underbrace{\sum_{j=k}^m \hat{\lambda}_{P_{k-1},j}
 \hat{u}_{P_{k-1},j,k}^2 }_{O_p\left( 1 \right)}\\
 &&\hspace{1cm} = O_p\left( \frac{1}{\theta_1 \theta_2 m} \right).
 \end{eqnarray*}
 \end{minipage}}
\end{itemize} 
\item[2.3)]  $i=1$, $j\geqslant k$ : 
\begin{itemize}
\item[2.3.1)] $r=2,3,4,...,k-1$: \\ 
\scalebox{0.65}{
\begin{minipage}{1\textwidth}
\begin{eqnarray*}
&&\hspace{-1cm}\left| \sum_{j \geqslant k} \frac{\hat{\theta}_{P_{k-1},1}\hat{\lambda}_{P_{k-1},j}}{\left(\hat{\theta}_{P_k,1}-\hat{\theta}_{P_{k-1},1}\right)\left(\hat{\theta}_{P_k,2}-\hat{\lambda}_{P_{k-1},j}\right)}  \hat{u}_{P_{k-1},1,k} \hat{u}_{P_{k-1},j,k} \left(- \sum_{r=2}^{k-1} \hat{u}_{P_{k-1},1,r} \hat{u}_{P_{k-1},j,r} \right) \right|\\
 &&\hspace{-0.3cm} = O_p\left( \frac{\theta_1 m}{\theta_2 \min(\theta_1,\theta_k)} \right) \sum_{r=2}^{k-1}  \underbrace{\left|\hat{u}_{P_{k-1},1,k} \hat{u}_{P_{k-1},1,r}\right|}_{O_p\left( \frac{\min(\theta_1,\theta_r)^{1/2}}{\theta_1^{1/2} \max(\theta_1,\theta_r)^{1/2} m} \right)}  \sum_{j=k}^m \hat{\lambda}_{P_{k-1},j}
 \left| \hat{u}_{P_{k-1},j,k} \hat{u}_{P_{k-1},j,r}\right|\\
 &&\hspace{-0.3cm} \leqslant \underset{r=2,...,k-1}{\max} \left( O_p\left( \frac{\theta_1^{1/2}\min(\theta_1,\theta_r)^{1/2}}{\theta_2 \max(\theta_1,\theta_r)^{1/2}\min(\theta_1,\theta_k)} \right)
\left(\sum_{j=k}^m \hat{\lambda}_{P_{k-1},j}^2
 \hat{u}_{P_{k-1},j,k}^2 \right)^{1/2}
 \left(\sum_{j=k}^m \hat{u}_{P_{k-1},j,r}^2\right)^{1/2}\right)\\
  &&\hspace{-0.3cm} \leqslant O_p\left( \frac{1}{\theta_1 \theta_2} \right).
 \end{eqnarray*} 
 \end{minipage}}
\item [2.3.2)] $r=k$:  \\ 
\scalebox{0.71}{
\begin{minipage}{1\textwidth}
\begin{eqnarray*}
&&\left| \sum_{j \geqslant k} \frac{\hat{\theta}_{P_{k-1},1}\hat{\lambda}_{P_{k-1},j}}{\left(\hat{\theta}_{P_k,1}-\hat{\theta}_{P_{k-1},1}\right)\left(\hat{\theta}_{P_k,2}-\hat{\lambda}_{P_{k-1},j}\right)}  \hat{u}_{P_{k-1},1,k} \hat{u}_{P_{k-1},j,k} \left(- \hat{u}_{P_{k-1},1,k} \hat{u}_{P_{k-1},j,k} \right) \right|\\
&&\hspace{1cm} =O_p\left( \frac{1}{\theta_2 \min(\theta_1,\theta_k)} \right).
 \end{eqnarray*} 
 \end{minipage}}
 \item[2.3.3)] $r=1$: We use Theorem \ref{AATheoremcaraceigenstructure} part (b) and (h). \\ 
\scalebox{0.71}{
\begin{minipage}{1\textwidth} 
 \begin{eqnarray*}
&&\left| \sum_{j \geqslant k} \frac{\hat{\theta}_{P_{k-1},1}\hat{\lambda}_{P_{k-1},j}}{\left(\hat{\theta}_{P_k,1}-\hat{\theta}_{P_{k-1},1}\right)\left(\hat{\theta}_{P_k,2}-\hat{\lambda}_{P_{k-1},j}\right)}  \hat{u}_{P_{k-1},1,k} \hat{u}_{P_{k-1},j,k} \left(- \hat{u}_{P_{k-1},1,1} \hat{u}_{P_{k-1},j,1} \right) \right|\\
&&\hspace{2cm} =\underbrace{\left| \frac{\hat{\theta}_{P_{k-1},1}}{\hat{\theta}_{P_k,1}-\hat{\theta}_{P_{k-1},1}}  \hat{u}_{P_{k-1},1,k} \hat{u}_{P_{k-1},1,1} \right|}_{O_p\left(\frac{\theta_1^{1/2}m^{1/2}}{\min(\theta_1,\theta_k)} \right)} \underbrace{\left| \sum_{j \geqslant k} \frac{\hat{\lambda}_{P_{k-1},j}}{\hat{\theta}_{P_k,2}-\hat{\lambda}_{P_{k-1},j}}  \hat{u}_{P_{k-1},j,k}  \hat{u}_{P_{k-1},j,1}\right| }_{O_p\left( \frac{1}{\theta_2 \theta_1^{1/2} m^{1/2}} \right)} \\
&&\hspace{2cm} = O_p\left(\frac{1}{\theta_2 \min(\theta_1,\theta_k)} \right).
  \end{eqnarray*}
  \end{minipage}}
\end{itemize} 
\item[2.4)]  $j<k$, $i\geqslant k$ : As in 2.2 and 2.3, we can show that this part is of order $O_p\left(\frac{1}{\theta_1 \min(\theta_2,\theta_k)} \right)$.
\item[2.5)] $i,j<k$ 
\begin{itemize}
\item[2.5.1)] $i,j<k$, $i \not = 1$, $j \not = 2$: \\ 
\scalebox{0.69}{
\begin{minipage}{1\textwidth}
\begin{eqnarray*}
&&\underbrace{\frac{\hat{\theta}_{P_{k-1},i}\hat{\theta}_{P_{k-1},j}}{\left(\hat{\theta}_{P_k,1}-\hat{\theta}_{P_{k-1},i}\right)\left(\hat{\theta}_{P_k,2}-\hat{\theta}_{P_{k-1},j}\right)} \hat{u}_{P_{k-1},i,k} \hat{u}_{P_{k-1},j,k}}_{O_p\left( \frac{\theta_i^{1/2} \theta_j^{1/2}}{\max(\theta_1,\theta_i) \max(\theta_2,\theta_j) m} \right)} \underbrace{\left(- \sum_{r=1}^k \hat{u}_{P_{k-1},i,r} \hat{u}_{P_{k-1},j,r} \right)}_{O_p\left(\frac{1}{\theta_i^{1/2} \theta_j^{1/2} m^{1/2}}\right) \text{ (by induction on $k-1$)} }\\
&&\hspace{1cm}= O_p\left( \frac{1}{\theta_1 \theta_2 m^{3/2}} \right).
\end{eqnarray*}
\end{minipage}}

\item[2.5.2)] $i=1$, $j=3,4,...,k-1$:  \\ 
\scalebox{0.67}{
\begin{minipage}{1\textwidth}
\begin{eqnarray*}
&&\underbrace{\frac{\hat{\theta}_{P_{k-1},1}\hat{\theta}_{P_{k-1},j}}{\left(\hat{\theta}_{P_k,1}-\hat{\theta}_{P_{k-1},1}\right)\left(\hat{\theta}_{P_k,2}-\hat{\theta}_{P_{k-1},j}\right)} \hat{u}_{P_{k-1},1,k} \hat{u}_{P_{k-1},j,k}}_{O_p\left( \frac{\theta_1^{1/2} \theta_j^{1/2}}{\min(\theta_1,\theta_k) \max(\theta_2,\theta_j)} \right)} \underbrace{\left(- \sum_{r=1}^k \hat{u}_{P_{k-1},1,r} \hat{u}_{P_{k-1},j,r} \right)}_{O_p\left(\frac{1}{\theta_1^{1/2} \theta_j^{1/2} m^{1/2}}\right) \text{ (by induction on $k-1$)} }\\
&&\hspace{1cm}= O_p\left( \frac{1}{\theta_2 \min(\theta_1,\theta_k) m^{1/2}} \right).
\end{eqnarray*}
\end{minipage}}

\item[2.5.3)] $j=2$, $i=3,4,...,k-1$: By similar simplifications as 2.5.2,\\ 
\scalebox{0.67}{
\begin{minipage}{1\textwidth}
\begin{eqnarray*}
&&\frac{\hat{\theta}_{P_{k-1},i}\hat{\theta}_{P_{k-1},2}}{\left(\hat{\theta}_{P_k,1}-\hat{\theta}_{P_{k-1},i}\right)\left(\hat{\theta}_{P_k,2}-\hat{\theta}_{P_{k-1},2}\right)} \hat{u}_{P_{k-1},i,k} \hat{u}_{P_{k-1},2,k}\left(- \sum_{r=1}^k \hat{u}_{P_{k-1},i,r} \hat{u}_{P_{k-1},2,r} \right)\\
&&\hspace{1cm}= O_p\left( \frac{1}{\theta_1 \min(\theta_2,\theta_k) m^{1/2}} \right).
\end{eqnarray*}
\end{minipage}}

\item[2.5.4)] $i=1$, $j=2$ : \\ 
\scalebox{0.67}{
\begin{minipage}{1\textwidth}
\begin{eqnarray*}
&&\underbrace{\frac{\hat{\theta}_{P_{k-1},1}\hat{\theta}_{P_{k-1},2}}{\left(\hat{\theta}_{P_k,1}-\hat{\theta}_{P_{k-1},1}\right)\left(\hat{\theta}_{P_k,2}-\hat{\theta}_{P_{k-1},2}\right)} \hat{u}_{P_{k-1},1,k} \hat{u}_{P_{k-1},2,k}}_{\overset{\scalebox{0.5}{order}}{\sim} \left( \frac{\theta_1^{1/2} \theta_2^{1/2} m}{\min(\theta_1,\theta_k) \min(\theta_2,\theta_k)} \right)} \underbrace{\left(- \sum_{r=1}^k \hat{u}_{P_{k-1},1,r} \hat{u}_{P_{k-1},2,r} \right)}_{O_p\left(\frac{1}{\theta_1^{1/2} \theta_2^{1/2} m^{1/2}}\right) \text{ (by induction on $k-1$)} }\\
&&\hspace{1cm}= \mathbf{O_p\left( \frac{ m^{1/2}}{\min(\theta_1,\theta_k) \min(\theta_2,\theta_k)} \right)}.
\end{eqnarray*}
\end{minipage}}\\

\noindent This term cannot be neglected and its estimation is presented in $\textbf{C}$.
\end{itemize}

\noindent Finally, \\ 
\scalebox{0.67}{
\begin{minipage}{1\textwidth}
\begin{eqnarray*}
\sum_{s=k+1}^m \hat{u}_{P_k,1,s} \hat{u}_{P_k,2,s}&=&\\
&&\hspace{-3.5cm} \frac{1}{\sqrt{D_1 D_2}N_1N_2} \left( 
 \frac{\hat{\theta}_{P_{k-1},1}\hat{\theta}_{P_{k-1},2}}{\left(\hat{\theta}_{P_k,1}-\hat{\theta}_{P_{k-1},1}\right)\left(\hat{\theta}_{P_k,2}-\hat{\theta}_{P_{k-1},2}\right)}  \hat{u}_{P_{k-1},1,k} \hat{u}_{P_{k-1},2,k} \left(- \sum_{r=1}^k \hat{u}_{P_{k-1},1,r} \hat{u}_{P_{k-1},2,r} \right) \right)\\
 && + O_p\left(\frac{1}{\theta_1^{1/2} \theta_2^{1/2} m}\right) .
\end{eqnarray*}
\end{minipage}}

\item[\textbf{B:}] In this paragraph we investigate $ \frac{1}{\sqrt{D_1 D_2}N_1N_2}$. \\ 
\scalebox{0.69}{
\begin{minipage}{1\textwidth}
\begin{eqnarray*}
D_1&=&\underbrace{\sum_{i=k}^m \frac{\hat{\lambda}_{P_{k-1},i}^2}{(\hat{\theta}_{P_k,1}-\hat{\lambda}_{P_{k-1},i})^2} \hat{u}_{P_{k-1},i,k}^2}_{ O_p\left( \frac{1}{\theta_1^2}\right) } + 
\underbrace{\frac{\hat{\theta}_{P_{k-1},1}^2}{(\hat{\theta}_{P_{k},1}-\hat{\theta}_{P_{k-1},1})^2}   \hat{u}_{P_{k-1},1,k}^2 }_{ \overset{\scalebox{0.5}{order}}{\sim}  \frac{\theta_1 m}{\min(\theta_1,\theta_k)^2} }+  
\underbrace{\sum_{i=2}^{k-1} \frac{\hat{\theta}_{P_{k-1},i}^2}{(\hat{\theta}_{P_k,1}-\hat{\theta}_{P_{k-1},i})^2}  \hat{u}_{P_{k-1},i,k}^2}_{ O_p\left( \frac{1}{\theta_1 m}\right) }\\
&=& \frac{\hat{\theta}_{P_{k-1},1}^2}{(\hat{\theta}_{P_{k},1}-\hat{\theta}_{P_{k-1},1})^2}   \hat{u}_{P_{k-1},1,k}^2 + O_p\left( \frac{1}{\theta_1 m} \right)+ O_p\left( \frac{1}{\theta_1^2} \right)\\
&=&O_p\left( \frac{\theta_1 m}{\min(\theta_1,\theta_k)^2} \right).
\end{eqnarray*}
\end{minipage}} \\

\noindent Because $\hat{{u}}_{P_k,t,s}=\frac{P_k^{1/2} \tilde{{u}}_{P_k,t}}{ N_t} $, it follows that \small
\begin{eqnarray*}
N_1^2&=& 
\sum_{i\not = k}^{m} \tilde{{u}}_{P_k,1,i}^2+\tilde{{u}}_{P_k,1,k}^2 \theta_k \\
&=& 1+ (\theta_k-1)  \tilde{{u}}_{P_2,1,2}^2\\
&=& 1+ \frac{1}{(\theta_k-1) D_1}\\
&=& 1+ O_p\left( \frac{\min(\theta_1,\theta_k)}{\max(\theta_1,\theta_k)m} \right).
\end{eqnarray*}\normalsize
We easily obtain  \small
\begin{eqnarray*}
\frac{1}{N_1 \sqrt{D_1}}&=& \frac{|\hat{\theta}_{P_{k},1}-\hat{\theta}_{P_{k-1},1}|}{\hat{\theta}_{P_{k-1},1}  |\hat{u}_{P_{k-1},1,k}|}   +O_p\left( \frac{\min(\theta_1,\theta_k)}{ \theta_1^{1/2} m^{3/2}} \right) \\
&= & O_p\left( \frac{\min(\theta_1,\theta_k)}{\theta_1^{1/2}m^{1/2}} \right)
\end{eqnarray*}\normalsize
and \\ 
\scalebox{0.9}{
\begin{minipage}{1\textwidth}
\begin{eqnarray*}
\frac{1}{N_1 N_2 \sqrt{D_1 D_2 }}&=& \frac{|\hat{\theta}_{P_{k},1}-\hat{\theta}_{P_{k-1},1}|}{\hat{\theta}_{P_{k-1},1}  |\hat{u}_{P_{k-1},1,k}|}  \frac{|\hat{\theta}_{P_{k},2}-\hat{\theta}_{P_{k-1},2}|}{\hat{\theta}_{P_{k-1},2}  |\hat{u}_{P_{k-1},2,k}|}   +O_p\left( \frac{\min(\theta_1,\theta_k)}{\theta_1^{1/2} \theta_2^{1/2}m^2} \right) \\
&=&O_p\left( \frac{\min(\theta_1,\theta_k) \min(\theta_2,\theta_k)}{\theta_1^{1/2} \theta_2^{1/2} m} \right).
\end{eqnarray*}
\end{minipage}}
\item[\textbf{C:}]
From \textbf{A} and \textbf{B}, we conclude using Theorem \ref{AAThinvarianteigenvalue},\\ 
\scalebox{0.65}{
\begin{minipage}{1\textwidth}
\begin{eqnarray*}
&&\sum_{s=k+1}^m \hat{u}_{P_k,1,s} \hat{u}_{P_k,2,s} \hspace{40cm}\\
&&\hspace{0.5cm} =\frac{\frac{\hat{\theta}_{P_{k-1},1}\hat{\theta}_{P_{k-1},2}}{\left(\hat{\theta}_{P_k,1}-\hat{\theta}_{P_{k-1},1}\right)\left(\hat{\theta}_{P_k,2}-\hat{\theta}_{P_{k-1},2}\right)}  \hat{u}_{P_{k-1},1,k} \hat{u}_{P_{k-1},2,k} \left(- \sum_{r=1}^k \hat{u}_{P_{k-1},1,r} \hat{u}_{P_{k-1},2,r} \right)}{\sqrt{D_1 D_2}N_1N_2}+ O_p\left(\frac{1}{\theta_1^{1/2} \theta_2^{1/2} m}\right)\\
 &&\hspace{0.5cm}=
\text{sign}\left( \hat{u}_{P_{k-1},1,1} \hat{u}_{P_{k-1},2,2} \hat{u}_{P_{k},1,1} \hat{u}_{P_{k},2,2} \right) 
\left(- \sum_{r=1}^k \hat{u}_{P_{k-1},1,r} \hat{u}_{P_{k-1},2,r} \right)+ O_p\left(\frac{1}{\theta_1^{1/2} \theta_2^{1/2} m}\right)\\
&&\hspace{0.5cm}= \text{sign}\left( \hat{u}_{P_{k-1},1,1} \hat{u}_{P_{k-1},2,2} \hat{u}_{P_{k},1,1} \hat{u}_{P_{k},2,2} \right) \sum_{s=k+1}^m \hat{u}_{P_{k-1},1,s}\hat{u}_{P_{k-1},2,s}+ O_p\left(\frac{1}{\theta_1^{1/2} \theta_2^{1/2} m}\right).
\end{eqnarray*}
\end{minipage}}\\

Using the remark of Theorem \ref{AATheoremcaraceigenstructure}, the sign of the third line is correct with a probability tending to $1$ in $1/m$.
Therefore,  using the convention $\hat{u}_{P_s,i,i}>0$ for $i=1,2,...,s$ and \linebreak $s=1,2,...,k$ leads to \small
\begin{eqnarray*}
\sum_{s=k+1}^m \hat{u}_{P_k,1,s} \hat{u}_{P_k,2,s}&=& \sum_{s=k+1}^m \hat{u}_{P_{k-1},1,s}\hat{u}_{P_{k-1},2,s}+ O_p\left(\frac{1}{\theta_1^{1/2} \theta_2^{1/2} m}\right)\\
&=& \sum_{s=2}^m \hat{u}_{P_{2},1,s}\hat{u}_{P_{2},2,s}+ O_p\left(\frac{1}{\theta_1^{1/2} \theta_2^{1/2} m}\right),
\end{eqnarray*} \normalsize
where we recall that the error $O$ is in probability. 
\end{itemize}
\end{proof}

\subsubsection{Invariant Double Angle Theorem}

\begin{proof}{\textbf{Corollary \ref{AAThInvariantdouble}}}\label{AAproofThInvariantdouble} 
\noindent In order to shorten the equations, we use the following notation 
$$
\begin{array}{lll}
 \hat{\theta}_{P_s,t}=\hat{\theta}_{\hat{\Sigma}_{X,P_s},t}, & \hat{u}_{P_s,t}=\hat{u}_{\hat{\Sigma}_{X,P_s},t}, & \hat{\lambda}_{P_s,t}=\hat{\lambda}_{\hat{\Sigma}_{X,P_s},t}\\
  \hat{\hat{\theta}}_{P_s,t}=\hat{\theta}_{\hat{\Sigma}_{Y,P_s},t}, & \hat{\hat{u}}_{P_s,t}=\hat{u}_{\hat{\Sigma}_{Y,P_s},t}, & \hat{\hat{\lambda}}_{P_s,t}=\hat{\lambda}_{\hat{\Sigma}_{Y,P_s},t}.\\
\end{array}$$
Moreover,   
\begin{eqnarray*}
u^{c_s}&=&\frac{u_{1:s}}{||u_{1:s}||}, \text{ where } u \text{ is a vector of size } m,\\
\hat{\alpha}_{P_s,i}^2 &=& ||\hat{u}_{P_s,i,1:s}||^2,\\
\hat{\hat{\alpha}}_{P_s,i}^2 &=& ||\hat{\hat{u}}_{P_s,i,1:s}||^2.
\end{eqnarray*}
Finally, using the notation \ref{AANot=thetadifferent} and relaxing $\theta_1>\theta_2>...>\theta_k$ allows us to only study $\hat{u}_{P_k,1}$ and $\hat{\hat{u}}_{P_k,1}$ without loss of generality.

\noindent The proof is essentially based on Theorems \ref{AAjointdistribution},\ref{AAThDotproduct}, \ref{AAThInvariantdot} and \ref{AAInvariantth}.
\begin{enumerate}
\item First we investigate $\left\langle \hat{u}_{P_1,1},\hat{\hat{u}}_{P_1,1} \right\rangle^2$:\\ 
\scalebox{0.78}{
\begin{minipage}{1\textwidth}
\begin{eqnarray*}
\left\langle \hat{u}_{P_1,1},\hat{\hat{u}}_{P_1,1} \right\rangle^2&=&
\hat{u}_{P_1,1,1}^2\hat{\hat{u}}_{P_1,1,1}^2+2
\hat{u}_{P_1,1,1}\hat{\hat{u}}_{P_1,1,1} \sum_{i=2}^m \hat{u}_{P_1,1,i}\hat{\hat{u}}_{P_1,1,i}
+\left(\sum_{i=1}^m \hat{u}_{P_1,1,i}\hat{\hat{u}}_{P_1,1,i}\right)^2\\
&=& \underbrace{\hat{u}_{P_1,1,1}^2\hat{\hat{u}}_{P_1,1,1}^2}_{{\rm RV}\left( O_p\left( 1\right),O_p\left(\frac{1}{\theta_1^2 m} \right)\right)}+\underbrace{C_{P_1}}_{{\rm RV}\left( 0,O_p\left(\frac{1}{\theta_1^2 m} \right)\right) }+O_p\left(\frac{1}{\theta_1^2 m}\right).
\end{eqnarray*}
\end{minipage}}
\item Next, we want to prove
\begin{eqnarray*}
&& \left\langle \hat{u}_{P_1,1},\hat{\hat{u}}_{P_1,1} \right\rangle^2= \sum_{i=1}^{k} \left\langle \hat{u}_{P_k,1},\hat{\hat{u}}_{P_k,i} \right\rangle^2+O_p\left( \frac{1}{\theta_1 m} \right).
\end{eqnarray*}
Using Theorem \ref{AAThcomponentdistribution} and \ref{AAjointdistribution},\\ 
\scalebox{0.77}{
\begin{minipage}{1\textwidth}
\begin{eqnarray*}
&&\hspace{-1cm} \left\langle \hat{u}_{P_k,1},\hat{\hat{u}}_{P_k,1} \right\rangle^2  
= \left\langle \hat{u}_{P_k,1,1:k},\hat{\hat{u}}_{P_k,1,1:k}\right\rangle^2 +\underbrace{2
\hat{u}_{P_k,1,1}\hat{\hat{u}}_{P_k,1,1} \sum_{i=k+1}^m \hat{u}_{P_k,1,i}\hat{\hat{u}}_{P_k,1,i}}_{C_{P_k}}
+ O_p\left( \frac{1}{\theta_1 m} \right),\\
&&\hspace{-1cm}\left\langle \hat{u}_{P_k,1},\hat{\hat{u}}_{P_k,s} \right\rangle^2 =
\left\langle \hat{u}_{P_k,1,1:k},\hat{\hat{u}}_{P_k,s,1:k} \right\rangle^2+ O_p\left( \frac{1}{\max(\theta_1, \theta_s) m} \right).
\end{eqnarray*}
\end{minipage}}\\

In this theorem we suppose that Assumption \ref{AAAss=theta} (A4) holds and without loss of generality, we assume that $\theta_1,...,\theta_{k_1}$ are of same order and  $\theta_{k_1+1},...,\theta_{k}$ are also of same order but different from the first group.  Assumption \ref{AAAss=theta} (A4) implies that either all the eigenvalues are proportional or one group has finite eigenvalues. Therefore,  
\begin{eqnarray*}
\sum_{i=1}^k \left\langle \hat{u}_{P_k,1},\hat{\hat{u}}_{P_k,i} \right\rangle^2 &=& \sum_{i=1}^{k_1} \left\langle \hat{u}_{P_k,1},\hat{\hat{u}}_{P_k,i} \right\rangle^2+ O_p\left(\frac{1}{\theta_1 m} \right).
\end{eqnarray*} 
Moreover, we easily see that for $i=1,2,...,k_1$,
\begin{eqnarray*}
\hat{\hat{\alpha}}_{P_k,i}^2 &=& ||\hat{\hat{u}}_{P_s,i,1:k}||^2\\
&=& ||\hat{\hat{u}}_{P_s,i,1:k_1}||^2+ O_p\left( \frac{1}{\theta_1 m} \right).
\end{eqnarray*}

Thus \\ 
\scalebox{0.9}{
\begin{minipage}{1\textwidth}
\begin{eqnarray*}
\sum_{i=1}^{k_1} \left\langle \hat{u}_{P_k,1},\hat{\hat{u}}_{P_k,i} \right\rangle^2 &=& \sum_{i=1}^{k_1} \left\langle \hat{u}_{P_k,1,1:k},\hat{\hat{u}}_{P_k,i,1:k}\right\rangle^2 + C_k + O_p\left( \frac{1}{\theta_1 m} \right)\\
&=& \sum_{i=1}^{k_1} \left\langle \hat{u}_{P_k,1,1:k_1},\hat{\hat{u}}_{P_k,i,1:k_1}\right\rangle^2 + C_k + O_p\left( \frac{1}{\theta_1 m} \right)\\
&=&\sum_{i=1}^{k_1} \hat{\alpha}_{P_k,1}^2 \hat{\hat{\alpha}}_{P_k,i}^2 \left\langle \hat{u}_{P_k,1}^{c_{k_1}},\hat{\hat{u}}_{P_k,i}^{c_{k_1}} \right\rangle^2  + C_k + O_p\left( \frac{1}{\theta_1 m} \right)\\
&=&\hat{\alpha}_{P_k,1}^2 \hat{\hat{\alpha}}_{P_k,1}^2 \sum_{i=1}^{k_1}   \left\langle \hat{u}_{P_k,1}^{c_{k_1}},\hat{\hat{u}}_{P_k,i}^{c_{k_1}} \right\rangle^2 + C_k + O_p\left( \frac{1}{\theta_1 m} \right)\\
&& \hspace{2cm}  + \hat{\alpha}_{P_k,1}^2 \sum_{i=2}^{k_1}  \left(\hat{\hat{\alpha}}_{P_k,i}^2 - \hat{\hat{\alpha}}_{P_k,1}^2 \right) \left\langle \hat{u}_{P_k,1}^{c_{k_1}},\hat{\hat{u}}_{P_k,i}^{c_{k_1}} \right\rangle^2\\
&=&\hat{\alpha}_{P_1,1}^2 \hat{\hat{\alpha}}_{P_1,1}^2 \underbrace{\sum_{i=1}^{k_1}   \left\langle \hat{u}_{P_k,1}^{c_{k_1}},\hat{\hat{u}}_{P_k,i}^{c_{k_1}} \right\rangle^2}_{\text{Part 1}} + \underbrace{C_{P_k}}_{\text{Part 2}} + O_p\left( \frac{1}{\theta_1 m} \right),
\end{eqnarray*}
\end{minipage}}\\

Where the last equality is obtained because for $i=1,2,...,k_1$, $\hat{\hat{\alpha}}_{P_k,i}^2 - \hat{\hat{\alpha}}_{P_k,1}^2=O_p\left( 1/\theta_1\right)$.\\
So,  we just need to show that 
\begin{eqnarray*}
&&\sum_{i=1}^{k_1}   \left\langle \hat{u}_{P_k,1}^{c_{k_1}},\hat{\hat{u}}_{P_k,i}^{c_{k_1}} \right\rangle^2=1+O_p\left( \frac{1}{ \theta_1 m} \right),\\
&&C_{P_k}=C_{P_1}+O_p\left( \frac{1}{\theta_1 m} \right).
\end{eqnarray*}

\paragraph*{Part 1 :}

First we prove that
$$ \sum_{i=1}^{k_1} \left\langle \hat{u}_{P_k,i}^{c_{k_1}},\hat{\hat{u}}_{P_k,i}^{c_{k_1}} \right\rangle^2  = 1  + O_p\left( \frac{1}{\theta_1 m} \right). $$

We apply Gram-Schmidt to $\hat{\hat{u}}_{P_k,1}^{c_{k_1}},\hat{\hat{u}}_{P_k,2}^{c_{k_1}},...,\hat{\hat{u}}_{P_k,k_1}^{c_{k_1}}$, \\ 
\scalebox{0.71}{
\begin{minipage}{1\textwidth}
\begin{eqnarray*}
\hat{\hat{w}}_{P_k,1}&=&\hat{\hat{u}}_{P_k,1}^{c_{k_1}}.\\
\hat{\hat{w}}_{P_k,2}&=&\left(\hat{\hat{u}}_{P_k,2}^{c_{k_1}}-\left\langle \hat{\hat{u}}_{P_k,2}^{c_{k_1}},\hat{\hat{w}}_{P_k,1} \right\rangle\hat{\hat{w}}_{P_k,1}\right) \left(1+O_p \left(\frac{1}{\theta_1^2 m} \right) \right).\\
  && \hspace{-1.5cm}\text{ \normalsize Indeed by Theorems \ref{AAThDotproduct} and \ref{AAThInvariantdot}, } \\
  &&|| \hat{\hat{u}}_{P_k,2}^{c_{k_1}}-\left\langle \hat{\hat{u}}_{P_k,2}^{c_{k_1}},\hat{\hat{w}}_{P_k,1} \right\rangle\hat{\hat{w}}_{P_k,1} || = 1- \underbrace{\left\langle \hat{\hat{u}}_{P_k,1}^{c_{k_1}},\hat{\hat{u}}_{P_k,2}^{c_{k_1}} \right\rangle^2}_{- \alpha_{P_k,1} \alpha_{P_k,2}   \sum_{i=k+1}^m \hat{\hat{u}}_{P_k,1,i}\hat{\hat{u}}_{P_k,2,i} }=1+O_p \left(\frac{1}{\theta_1^2 m} \right).\\
\hat{\hat{w}}_{P_k,p}&=&\left(\hat{\hat{u}}_{P_k,p}^{c_{k_1}}-\sum_{i=1}^{p-1} \left\langle \hat{\hat{u}}_{P_k,p}^{c_{k_1}},\hat{\hat{w}}_{P_k,i} \right\rangle\hat{\hat{w}}_{P_k,i}\right) \left(1+O_p \left(\frac{1}{\theta_1^2 m} \right) \right).
\end{eqnarray*} 
\end{minipage}}\\

However, the norm is more difficult to estimate for $p=3,4,...,k_1$:\\ 
\scalebox{0.71}{
\begin{minipage}{1\textwidth}
\begin{eqnarray*}
  ||\hat{\hat{u}}_{P_k,p}^{c_{k_1}}-\sum_{i=1}^{p-1} \left\langle \hat{\hat{u}}_{P_k,p}^{c_{k_1}},\hat{\hat{w}}_{P_k,i} \right\rangle\hat{\hat{w}}_{P_k,i} || &=& 1- \sum_{i=1}^{p-1}\left\langle \hat{\hat{u}}_{P_k,p}^{c_{k_1}},\hat{\hat{w}}_{P_k,i} \right\rangle^2\\
  &=&1- \sum_{i=1}^{p-1}\left\langle \hat{\hat{u}}_{P_k,p}^{c_{k_1}},\sum_{j=1}^{i} a_j \hat{\hat{u}}_{P_k,j}^{c_{k_1}} \right\rangle^2, \text{ for some } |a_i|<1,\\
  &=&1- \sum_{i=1}^{p-1} \sum_{j_1,j_2=1}^{i} \left\langle \hat{\hat{u}}_{P_k,p}^{c_{k_1}}, a_{j_1} \hat{\hat{u}}_{P_k,j_1}^{c_{k_1}} \right\rangle  \left\langle \hat{\hat{u}}_{P_k,p}^{c_{k_1}}, a_{j_2} \hat{\hat{u}}_{P_k,j_2}^{c_{k_1}} \right\rangle\\
  &=& 1+O_p \left(\frac{1}{\theta_1^2 m} \right) .
\end{eqnarray*}
\end{minipage}}\\

We can express the truncated eigenvectors in a orthonormal basis as,\\ 
\scalebox{0.9}{
\begin{minipage}{1\textwidth}
\begin{eqnarray*}
\Rightarrow && \hat{\hat{u}}_{P_k,1}^{c_k}=\hat{\hat{w}}_{P_k,1},\\
&& \text{For } p=2,...,k_1,\\
&& \hspace{1cm} \hat{\hat{u}}_{P_k,p}^{c_{k_1}}= \left(\hat{\hat{w}}_{P_k,p}+ \sum_{i=1}^{p-1} \left\langle \hat{\hat{u}}_{P_k,p}^{c_{k_1}},\hat{\hat{w}}_{P_k,i} \right\rangle\hat{\hat{w}}_{P_k,i} \right) \left(1+O_p \left(\frac{1}{\theta_1^2 m} \right)  \right).
\end{eqnarray*}
\end{minipage}}\\

Thus\\ 
\scalebox{0.8}{
\begin{minipage}{1\textwidth}
\begin{eqnarray*}
&& \hspace{-1cm} \sum_{p=1}^{k_1} \left\langle \hat{u}_{P_k,1}^{c_{k_1}},\hat{\hat{u}}_{P_k,p}^{c_{k_1}} \right\rangle^2   \\
&=&
\left\langle \hat{u}_{P_k,1}^{c_k},\hat{\hat{w}}_{P_k,1} \right\rangle^2
+ \sum_{p=2}^{k_1} \left\langle \hat{u}_{P_k,1}^{c_{k_1}},\hat{\hat{w}}_{P_k,p}+ \sum_{j=1}^{p-1} \left\langle \hat{\hat{u}}_{P_k,p}^{c_{k_1}},\hat{\hat{w}}_{P_k,j} \right\rangle\hat{\hat{w}}_{P_k,j} \right\rangle^2 +O_p \left(\frac{1}{ \theta_1^2 m} \right)   \\
&=& \sum_{p=1}^{k_1} \left\langle \hat{u}_{P_k,1}^{c_{k_1}},\hat{\hat{w}}_{P_k,p} \right\rangle^2
+ \left( \sum_{p=2}^{k_1} \left\langle \hat{u}_{P_k,1}^{c_{k_1}}, \sum_{j=1}^{p-1} \left\langle \hat{\hat{u}}_{P_k,p}^{c_{k_1}},\hat{\hat{w}}_{P_k,j} \right\rangle\hat{\hat{w}}_{P_k,j} \right\rangle^2  \right.\\
&&\hspace{1cm} \left. +2  \sum_{p=2}^{k_1}  \left\langle \hat{u}_{P_k,1}^{c_{k_1}},\hat{\hat{w}}_{P_k,p} \right\rangle
 \left\langle \hat{u}_{P_k,1}^{c_{k_1}}, \sum_{j=1}^{p-1} \left\langle \hat{\hat{u}}_{P_k,p}^{c_{k_1}},\hat{\hat{w}}_{P_k,j} \right\rangle\hat{\hat{w}}_{P_k,j} \right\rangle  \right)+O_p \left(\frac{1}{ \theta_1^2 m} \right)   \\
 &=& 1+ \sum_{p=2}^{k_1}  \left( \sum_{j=1}^{p-1} \left\langle \hat{\hat{u}}_{P_k,p}^{c_{k_1}},\hat{\hat{w}}_{P_k,j} \right\rangle \left\langle \hat{u}_{P_k,1}^{c_{k_1}}, \hat{\hat{w}}_{P_k,j} \right\rangle \right)^2 \\
&&\hspace{1cm} +2  \sum_{p=2}^{k_1} \sum_{j=1}^{p-1} \left\langle \hat{\hat{u}}_{P_k,p}^{c_{k_1}},\hat{\hat{w}}_{P_k,j} \right\rangle \left\langle \hat{u}_{P_k,1}^{c_{k_1}},\hat{\hat{w}}_{P_k,p} \right\rangle
 \left\langle \hat{u}_{P_k,1}^{c_{k_1}}, \hat{\hat{w}}_{P_k,j} \right\rangle +O_p \left(\frac{1}{ \theta_1^2 m} \right)   \\
  &=&1+A+B+O_p \left(\frac{1}{ \theta_1^2 m} \right). 
\end{eqnarray*}
\end{minipage}}\\

Next we prove separately that $A$ and $B$ are negligible.

\paragraph*{A :} By Theorem \ref{AAThDotproduct}, \ref{AAThInvariantdot}, \\ 
\scalebox{0.71}{
\begin{minipage}{1\textwidth}
\begin{eqnarray*}
&&\hspace{-0.8cm} \left\langle \hat{\hat{u}}_{P_k,p}^{c_{k_1}},\hat{\hat{w}}_{P_k,j} \right\rangle, \  k_1\geqslant p > j:  \\
&&\hspace{-0.3cm}  j=1: \ \left\langle \hat{\hat{u}}_{P_k,p}^{c_{k_1}},\hat{\hat{w}}_{P_k,1} \right\rangle= \left\langle \hat{\hat{u}}_{P_k,p}^{c_{k_1}},\hat{\hat{u}}_{P_k,1}^{c_{k_1}} \right\rangle=O_p\left( \frac{1}{\theta_1 \sqrt{m}} \right),\\
&&\hspace{-0.3cm} j\neq 1: \ \left\langle \hat{\hat{u}}_{P_k,p}^{c_{k_1}},\hat{\hat{w}}_{P_k,j} \right\rangle \left(1+O_p\left(\frac{1}{\theta_1^2 m} \right)\right)= \left\langle \hat{\hat{u}}_{P_k,p}^{c_{k_1}},\hat{\hat{u}}_{P_k,j}^{c_{k_1}} \right\rangle-\sum_{i=1}^{j-1} \underbrace{\left\langle \hat{\hat{u}}_{P_k,j}^{c_{k_1}},\hat{\hat{w}}_{P_k,i} \right\rangle}_{O_p\left( \frac{1}{ \theta_1 \sqrt{m}} \right)} \underbrace{\left\langle \hat{\hat{u}}_{P_k,p}^{c_{k_1}},\hat{\hat{w}}_{P_k,i} \right\rangle}_{O_p\left( \frac{1}{\sqrt{m}} \right)} \\
&&\hspace{6.7cm} =O_p\left( \frac{1}{ \theta_1 \sqrt{m}} \right).
\end{eqnarray*}
\end{minipage}}\\ 
\scalebox{0.71}{
\begin{minipage}{1\textwidth}
\begin{eqnarray*}
&&\hspace{-0.8cm} \left\langle \hat{u}_{P_k,1}^{c_{k_1}},\hat{\hat{w}}_{P_k,j} \right\rangle, k_1 \geqslant j: \\
&&\hspace{-0.3cm} j\neq 1: \ \left\langle \hat{u}_{P_k,1}^{c_{k_1}},\hat{\hat{w}}_{P_k,j} \right\rangle \left(1+O_p\left(\frac{1}{\theta_1^2 m} \right)\right)=
\left\langle \hat{u}_{P_k,1}^{c_{k_1}},\hat{\hat{u}}_{P_k,j}^{c_{k_1}} \right\rangle
-\sum_{i=1}^{j-1} \underbrace{\left\langle \hat{\hat{w}}_{P_k,i},\hat{\hat{u}}_{P_k,j}^{c_{k_1}} \right\rangle}_{O_p\left( \frac{1}{\theta_1 \sqrt{m}} \right)} \underbrace{\left\langle \hat{u}_{P_k,1}^{c_{k_1}},\hat{\hat{w}}_{P_k,i} \right\rangle}_{=O_p\left(1\right)} \\
&&\hspace{6.7cm}=O_p\left( \frac{1}{\sqrt{m}} \right),\\
&&\hspace{-0.3cm} j=1: \ \left\langle \hat{u}_{P_k,1}^{c_{k_1}},\hat{\hat{w}}_{P_k,1} \right\rangle= O_p\left( 1\right).
\end{eqnarray*}
\end{minipage}}\\

Consequently,
\begin{eqnarray*}
\left( \sum_{j=1}^{p-1} \left\langle \hat{\hat{u}}_{P_k,p}^{c_{k_1}},\hat{\hat{w}}_{P_k,j} \right\rangle \left\langle \hat{u}_{P_k,1}^{c_{k_1}}, \hat{\hat{w}}_{P_k,j} \right\rangle \right)^2=O_p\left( \frac{1}{\theta_1^2 m} \right)\\
\end{eqnarray*}
Therefore,  $A=O_p\left( \frac{1}{\theta_1^2 m} \right)$.
\paragraph*{B :} The same estimations as previously lead to\\ 
\scalebox{0.9}{
\begin{minipage}{1\textwidth}
\begin{eqnarray*}
B&=&2  \sum_{p=2}^{k_1} \sum_{j=1}^{p-1} \underbrace{\left\langle \hat{\hat{u}}_{P_k,p}^{c_{k_1}},\hat{\hat{w}}_{P_k,j} \right\rangle}_{O_p\left( \frac{1}{\theta_1 \sqrt{m}} \right)} \underbrace{\left\langle \hat{u}_{P_k,1}^{c_{k_1}},\hat{\hat{w}}_{P_k,p} \right\rangle}_{O_p\left(\frac{1}{\sqrt{m}} \right)}
 \underbrace{\left\langle \hat{u}_{P_k,1}^{c_{k_1}}, \hat{\hat{w}}_{P_k,j} \right\rangle}_{O_p\left( 1\right)}=O_p\left( \frac{1}{\theta_1 m} \right).
\end{eqnarray*}
\end{minipage}}\\

Therefore,  
$$ \sum_{i=1}^{k_1} \left\langle \hat{u}_{P_k,i}^{c_{k_1}},\hat{\hat{u}}_{P_k,i}^{c_{k_1}} \right\rangle^2  = 1  + O_p\left( \frac{1}{\theta_1 m} \right). $$

\paragraph*{Part 2 :} In this part we prove the invariance of $C_{P_1}$. We need to show:
\begin{eqnarray*}
C_{P_k}&=& 2
\hat{u}_{P_k,1,1}\hat{\hat{u}}_{P_k,1,1} \sum_{i=k+1}^m \hat{u}_{P_k,1,i}\hat{\hat{u}}_{P_k,1,i}\\
&=&2
\hat{u}_{P_1,1,1}\hat{\hat{u}}_{P_1,1,1} \sum_{i=2}^m \hat{u}_{P_1,1,i}\hat{\hat{u}}_{P_1,1,i} +O_p\left( \frac{1}{\theta_1 m}\right)\\
&=&C_{P_1}+O_p\left( \frac{1}{\theta_1 m}\right).
\end{eqnarray*}
In order to prove this result we show $C_{P_k}=C_{P_{k-1}}+O_p\left( \frac{1}{\theta_1 m}\right)$ and more precisely,\\ 
\scalebox{0.75}{
\begin{minipage}{1\textwidth}
\begin{eqnarray*}
2
\hat{u}_{P_k,1,1}\hat{\hat{u}}_{P_k,1,1}  \sum_{i=k+1}^m \hat{u}_{P_k,1,i}\hat{\hat{u}}_{P_k,1,i} &=& 2
\hat{u}_{P_{k-1},1,1}\hat{\hat{u}}_{P_{k-1},1,1}  \sum_{i=k}^m \hat{u}_{P_{k-1},1,i}\hat{\hat{u}}_{P_{k-1},1,i}+ O_p\left( \frac{1}{\theta_1 m}\right) .
\end{eqnarray*}
\end{minipage}}\\

The proof is similar to the proofs of invariant eigenvector structure \ref{AAThInvariantdot} and \ref{AAInvariantth}. We use Theorem \ref{AATheoremcaraceigenstructure} in order to estimate each term of the sum. Assuming $P_{k-1}$ satisfies \ref{AAAss=theta}(A4) the last added eigenvalue can be either proportional to $\theta_1$ or to the other group. \\
In this proof we do not use the convention $\hat{u}_{P_k,i,i}>0$ for $i=1,2,...,k$.

\noindent We start by studying $\hat{u}_{P_k,1}.$ As in Theorem \ref{AATheoremcaraceigenstructure}, for $s>k$, \\ 
\scalebox{0.69}{
\begin{minipage}{1\textwidth}
 \begin{eqnarray*}
\hat{u}_{P_k,1,s}&=&\frac{1}{\sqrt{\hat{D}_1} \hat{N}_1} \left(\sum_{i=k}^m \frac{\hat{\lambda}_{P_{k-1},i}}{\hat{\theta}_{P_k,1}-\hat{\lambda}_{P_{k-1},i}} \hat{u}_{P_{k-1},i,s}  \hat{u}_{P_{k-1},i,k} +\frac{\hat{\theta}_{P_{k-1},1}}{\hat{\theta}_{P_k,1}-\hat{\theta}_{P_{k-1},1}} \hat{u}_{P_{k-1},1,s}  \hat{u}_{P_{k-1},1,k} \right. \\
&&\hspace{1.2cm} \left.+ \sum_{i=2}^{k_1} \frac{\hat{\theta}_{P_{k-1},i}}{\hat{\theta}_{P_k,1}-\hat{\theta}_{P_{k-1},i}} \hat{u}_{P_{k-1},i,s}  \hat{u}_{P_{k-1},i,k}+ \sum_{i=k_1+1}^{k-1} \frac{\hat{\theta}_{P_{k-1},i}}{\hat{\theta}_{P_k,1}-\hat{\theta}_{P_{k-1},i}} \hat{u}_{P_{k-1},i,s}  \hat{u}_{P_{k-1},i,k} \right).
\end{eqnarray*} 
\end{minipage}}\\

By a similar proof as part (a), (b) and (c) of Theorem \ref{AATheoremcaraceigenstructure},
\begin{eqnarray*}
&&\hat{A}_s=\sum_{i=k}^m \frac{\hat{\lambda}_{P_{k-1},i}}{\hat{\theta}_{P_k,1}-\hat{\lambda}_{P_{k-1},i}} \hat{u}_{P_{k-1},i,s}  \hat{u}_{P_{k-1},i,k} = O_p\left( \frac{1}{\sqrt{m} \theta_1} \right),\\
&&\hat{B}_s=\frac{\hat{\theta}_{P_{k-1},1}}{\hat{\theta}_{P_k,1}-\hat{\theta}_{P_{k-1},1}} \hat{u}_{P_{k-1},1,s}  \hat{u}_{P_{k-1},1,k} \overset{\scalebox{0.5}{order}}{\sim} \frac{1}{\min(\theta_1,\theta_k)},\\
&& \hat{C}_s=\sum_{i=2}^{k_1} \frac{\hat{\theta}_{P_{k-1},i}}{\hat{\theta}_{P_k,1}-\hat{\theta}_{P_{k-1},i}} \hat{u}_{P_{k-1},i,s}  \hat{u}_{P_{k-1},i,k} = O_p\left(\frac{1}{m \theta_1} \right),\\
&&\hat{C}^G_s=\sum_{i=k_1+1}^{k-1} \frac{\hat{\theta}_{P_{k-1},i}}{\hat{\theta}_{P_k,1}-\hat{\theta}_{P_{k-1},i}} \hat{u}_{P_{k-1},i,s}  \hat{u}_{P_{k-1},i,k} = O_p\left(\frac{1}{m \theta_1} \right),\\
&&\hat{D}_1= \frac{\hat{\theta}_{P_{k-1},1}^2}{ \left(\hat{\theta}_{P_k,1}-\hat{\theta}_{P_{k-1},1}\right)^2}   \hat{u}_{P_{k-1},1,k}^2 + O_p\left( \frac{1}{\theta_1^2} \right)+ O_p\left( \frac{1}{\theta_1^2 m} \right),\\
&&\hat{N}_1=1+O_p\left(\frac{\min(\theta_1,\theta_k)}{\max(\theta_1,\theta_k)m} \right).
\end{eqnarray*}
Thus,
\begin{eqnarray*}
\hat{u}_{P_k,1,s}&=&\frac{1}{\sqrt{\hat{D}_1} \hat{N}_1}  \left(\hat{A}_s+\hat{B}_s+\hat{C}_s+\hat{C}^G_s\right).
\end{eqnarray*}
We now find, \\ 
\scalebox{0.85}{
\begin{minipage}{1\textwidth} 
\begin{eqnarray*}
\sum_{s=k+1}^m \hat{u}_{P_k,1,i}\hat{\hat{u}}_{P_k,1,i} &=& \frac{\sum_{s=k+1}^m \left(\hat{A}_s+\hat{B}_s+\hat{C}_s+\hat{C}^G_s\right) \left(\hat{\hat{A}}_s+\hat{\hat{B}}_s+\hat{\hat{C}}_s+\hat{\hat{C}}^G_s\right)}{\sqrt{\hat{D}_1} \hat{N}_1 \sqrt{\hat{\hat{D}}_1} \hat{\hat{N}}_1} .
\end{eqnarray*} 
\end{minipage}}\\

Many of the terms are negligible,
\begin{eqnarray*}
&&\sum_{s=k+1}^m \hat{A}_s \hat{\hat{A}}_s = O_p\left( \frac{1}{\theta_1^2} \right), \
\sum_{s=k+1}^m \hat{A}_s \hat{\hat{C}}_s = O_p\left( \frac{1}{\sqrt{m} \theta_1^2} \right), \\
&&\sum_{s=k+1}^m \hat{B}_s \hat{\hat{C}}_s = O_p\left( \frac{1}{\theta_1 \min(\theta_1,\theta_k)} \right) , \
\sum_{s=k+1}^m \hat{C}_s \hat{\hat{C}}_s = O_p\left( \frac{1}{m \theta_1^2} \right).
\end{eqnarray*}
Moreover, because $\hat{\hat{u}}_{P_{k-1},1,s}$ is invariant by rotation, we have that \\ 
\scalebox{0.73}{
\begin{minipage}{1\textwidth}
\begin{eqnarray*}
\sum_{s=k+1}^m \hat{A}_s \hat{\hat{B}}_s &=& \frac{\hat{\hat{\theta}}_{P_{k-1},1}}{\hat{\hat{\theta}}_{P_k,1}-\hat{\hat{\theta}}_{P_{k-1},1}}  \hat{\hat{u}}_{P_{k-1},1,k} \underbrace{\sum_{s=k+1}^m  \left( \hat{\hat{u}}_{P_{k-1},1,s} \sum_{i=k}^m \frac{\hat{\lambda}_{P_{k-1},i}}{\hat{\theta}_{P_k,1}-\hat{\lambda}_{P_{k-1},i}} \hat{u}_{P_{k-1},i,s}  \hat{u}_{P_{k-1},i,k}  \right)}_{O_p\left( \frac{1}{m \theta_1^2} \right)}\\
&=&  O_p\left( \frac{1}{\theta_1 \min(\theta_1,\theta_k)} \right).
\end{eqnarray*} 
\end{minipage}}\\

Using the remark of Theorem \ref{AATheoremcaraceigenstructure}, the last term leads to \\ 
\scalebox{0.73}{
\begin{minipage}{1\textwidth}
\begin{eqnarray*}
&&\hspace{-0.5cm}\sum_{s=k+1}^m \hat{u}_{P_k,1,i}\hat{\hat{u}}_{P_k,1,i}\\
&&\hspace{0.5cm} = \frac{1}{\sqrt{\hat{D}_1}\hat{N}_1\sqrt{\hat{\hat{D}}_1} \hat{\hat{N}}_1} \sum_{s=k+1}^m \hat{B}_s \hat{\hat{B}}_s + O_p\left( \frac{1}{\theta_1 m} \right)\\
&&\hspace{0.5cm} =  \frac{\frac{\hat{\theta}_{P_{k-1},1}}{\hat{\theta}_{P_k,1}-\hat{\theta}_{P_{k-1},1}}   \hat{u}_{P_{k-1},1,k} 
 \frac{\hat{\hat{\theta}}_{P_{k-1},1}}{\hat{\hat{\theta}}_{P_k,1}-\hat{\hat{\theta}}_{P_{k-1},1}}   \hat{\hat{u}}_{P_{k-1},1,k}
\sum_{s=k+1}^m \hat{u}_{P_{k-1},1,s} \hat{\hat{u}}_{P_{k-1},1,s}}
{\frac{\hat{\theta}_{P_{k-1},1}}{\left| \hat{\theta}_{P_k,1}-\hat{\theta}_{P_{k-1},1} \right|}   \left|\hat{u}_{P_{k-1},1,k}\right|
 \frac{\hat{\hat{\theta}}_{P_{k-1},1}}{\left|\hat{\hat{\theta}}_{P_k,1}-\hat{\hat{\theta}}_{P_{k-1},1}\right|}   \left| \hat{\hat{u}}_{P_{k-1},1,k} \right|}+ O_p\left( \frac{1}{\theta_1 m} \right)\\
 &&\hspace{0.5cm} = \text{sign} \left( \left( \hat{\theta}_{P_k,1}-\hat{\theta}_{P_{k-1},1} \right)  \hat{u}_{P_{k-1},1,k}  \left( \hat{\hat{\theta}}_{P_k,1}-\hat{\hat{\theta}}_{P_{k-1},1} \right)  \hat{\hat{u}}_{P_{k-1},1,k} \right)  \sum_{s=k+1}^m \hat{u}_{P_{k-1},1,s} \hat{\hat{u}}_{P_{k-1},1,s}\\
 &&\hspace{0.5cm} = \text{ sign}\left( \hat{u}_{P_k,1,1} \right) 
\text{sign}\left( \hat{\hat{u}}_{P_k,1,1} \right) \text{ sign}\left( \hat{\hat{u}}_{P_{k-1},1,1}  \right) \text{ sign}\left( \hat{u}_{P_{k-1},1,1}  \right) \sum_{s=k+1}^m \hat{u}_{P_{k-1},1,s} \hat{\hat{u}}_{P_{k-1},1,s}.
\end{eqnarray*}
\end{minipage}} \\
Finally,\\ 
\scalebox{0.73}{
\begin{minipage}{1\textwidth}  
\begin{eqnarray*}
&&2
\hat{u}_{P_k,1,1}\hat{\hat{u}}_{P_k,1,1}  \sum_{i=k+1}^m \hat{u}_{P_k,1,i}\hat{\hat{u}}_{P_k,1,i} = 2 
\hat{u}_{P_{k-1},1,1}\hat{\hat{u}}_{P_{k-1},1,1}  \sum_{i=k}^m \hat{u}_{P_{k-1},1,i}\hat{\hat{u}}_{P_{k-1},1,i}+ O_p\left( \frac{1}{\theta_1 m}\right)
\end{eqnarray*}
\end{minipage}}\\
and the remark is straightforward assuming the sign convention.
\end{enumerate}
\end{proof}

\bibliographystyle{imsart-nameyear}
\bibliography{biblio4} 
\addcontentsline{toc}{section}{Bibliography}

\end{document}

%% file: macro.tex
\DeclareMathOperator{\cov}{Cov}                   
\DeclareMathOperator{\corr}{Corr}                 
\DeclareMathOperator{\var}{Var}                   
\DeclareMathOperator{\E}{E}                       

\renewcommand{\P}{\mathit{P}}
\renewcommand{\i}{{\,\text{\it \@i}\,}}

\newcommand{\Tr}{{{\rm Trace}}}
\newcommand{\Normal}{\mathbf{N}}

\newcommand{\I}{{\rm I}}

\newcommand{\X}{{\mathbf{X}}}
\newcommand{\Y}{{\mathbf{Y}}}

\newcommand{\1}{{\mathbf{1}}}

\makeatother
